%This is a plain tex paper. 
\magnification=\magstep1
\vbadness=10000
\hbadness=10000
\tolerance=10000

   % complex numbers
\def\Q{{\bf Q}}   % rational numbers
\def\R{{\bf R}}   % real numbers
\def\Z{{\bf Z}}   % integers
\def\Aut{{\rm Aut}}
\def\height{\rm height}

\proclaim The Leech lattice and other lattices. 

This is a corrected (1999) copy of my Ph.D. thesis. I have corrected
several errors, added a few remarks about later work by various people
that improves the results here, and missed out some of the more
complicated diagrams and table 3.

Most of chapter 2 was published as ``The Leech lattice'',
Proc. R. Soc. Lond A398 (365-376) 1985.
The results of  chapter 6 have been published as ``Automorphism
groups of Lorentzian lattices'', J. Alg. Vol 111 No. 1 Nov 1987.
Chapter 7 is out of date, and the paper
``The monster Lie algebra'', Adv. Math Vol 83 No 1 1990 p.30
contains stronger results.

Current address: R. E. Borcherds, Department of Math, Evans Hall \#3840,
University of California at Berkeley, CA 94720-3840, U.S.A.

Current e-mail: reb@math.berkeley.edu

www home page  www.math.berkeley.edu/\hbox{\~{}}reb

\bigskip
\centerline{The Leech lattice and other lattices}
\centerline{Richard Borcherds}
\centerline{Trinity College}
\centerline{June 1984}
\bigskip

\proclaim Preface.

I thank my research supervisor Professor J. H. Conway for his help and
encouragement. I also thank the S.E.R.C. for its financial support
and Trinity College for a research scholarship and a fellowship.

\proclaim Contents.

\item{} Introduction.
\bigskip
\item{} Notation.
\bigskip
\item{0.} Definitions and notation
\item{0.1} Definitions
\item{0.2} Neighbors
\item{0.3} Root systems
\item{0.4} Lorentzian lattices and hyperbolic geometry
\bigskip
\item{1.} Root systems
\item{1.1} Norms of Weyl vectors
\item{1.2} Minimal vectors
\item{1.3} The opposition involution
\item{1.4} Maximal sub root systems
\item{1.5} Automorphisms  of the  fundamental domain
\item{1.7} Norm 4 vectors 
\bigskip
\item{2.} 
The Leech lattice. 
\item{2.1} History
\item{2.2} Some results from modular forms
\item{2.3} Norm 0 vectors in Lorentzian lattices
\item{2.4} Existence of the Leech lattice
\item{2.5} The covering radius of the Leech lattice
\item{2.6} Uniqueness of the Leech lattice
\item{2.7} The deep holes of the Leech lattice
\bigskip
\item{3.} Negative norm vectors in Lorentzian lattices. 
\item{3.1} Negative norm vectors and lattices
\item{3.2} Summary of the classification
\item{3.3} How to find things in $D$
\item{3.4} Vectors of type 0 and 1
\item{3.5} Reducing the norm by 2
\item{3.6} Vectors of $B/nB$
\item{3.7} Type 2 vectors
\item{3.8} Vectors of type at least 3
\item{3.9} Positive norm vectors
\bigskip
\item{4.} 25 dimensional lattices
\item{4.1} The height in $II_{25,1}$
\item{4.2} The space $Z$
\item{4.3} Norm $-2$ vectors
\item{4.4} Norm $-4$ vectors
\item{4.5} 24 dimensional lattices
\item{4.6} 25 dimensional lattices
\item{4.7} Norm $-6$ vectors
\item{4.8} Vectors whose norm is not square-free
\item{4.9} Large norm vectors in $D$
\item{4.10} Norm $-8$ vectors
\item{4.11} Small holes and lattices
\bigskip
\item{5.} Theta functions
\item{5.1} Some standard results about theta functions
\item{5.2} Roots of a 25-dimensional even bimodular lattice
\item{5.3} The $\theta$ function of a 25 dimensional even bimodular lattice
\item{5.4} The $\theta$ function of a 25 dimensional unimodular lattice
\item{5.5} 26 dimensional unimodular lattices
\item{5.6} More on 26 dimensional unimodular lattices
\item{5.6} A 27 dimensional unimodular lattice with no roots
\bigskip
\item{6.} Automorphism  groups of Lorentzian lattices
\item{6.1} General properties of the automorphism  group.
\item{6.2} Notation
\item{6.3} Some automorphisms of the lattice $T$
\item{6.4} Hyperplanes in $T$
\item{6.5} A complex
\item{6.6} Unimodular lattices
\item{6.7} More about $I_{n,1}$
\item{6.8} Other examples
\item{6.9} Higher dimensions
\bigskip
\item{7.} The monster Lie algebra
\item{7.1} Introduction
\item{7.2} Multiplicities of roots of type 0 and 1
\item{7.3} Multiplicities of roots of norm $-2$
\bigskip
\item{}References
\bigskip
\item{} Figure 1. The neighborhood graph for 24 dimensional unimodular 
lattices. 
\item{} Figure 2 The vectors of $II_{25,1}$ of small height and norm.
\item{} Table $-2$ The norm $-2$ vectors of $II_{25,1}$
\item{} Table $-4$ The norm $-4$ vectors of $II_{25,1}$
\item{} Table 3 $B/2B$ for Niemeier lattices $B$

\proclaim Notation.

The following symbols often have the following meanings (but they
sometimes do not).
\item{$A$}
An odd lattice with elements $a$.
\item{$B$}
An even lattice with elements $b$.
\item{$D$}
A fundamental domain, usually of the root system of $II_{25,1}$.
\item{$h$} 
The Coxeter number of a root system or Niemeier lattice. 
\item{$II_{25,1}$} 
The 26-dimensional even unimodular Lorentzian lattice. 
\item{$L$} A lattice
\item{$\Lambda$} The Leech lattice
\item{$N$}
A Niemeier lattice, in other words a 24 dimensional even unimodular 
positive definite lattice.
\item{$\Q$} 
The rational numbers
\item{$\rho$}
The Weyl vector of the norm 2 vectors of some lattice.
\item{$\R$}
The real numbers.
\item{$R$}
A root system.
\item{$R_i(u)$} 
or $U_i$ The simple roots of $R$ that have inner product $-i$ with $u$.
\item{$r,r'$} 
Norm 2 vectors of a lattice.
\item{$S$} 
$=S(R)$ The maximal number of pairwise orthogonal roots of the root
system $R$. See 1.3.
\item{$T$}
A lattice
\item{$t,u,v,z$} 
Negative norm vectors in $II_{25,1}$. Usually $-u^2=2n$, $-v^2=2n-2$, 
$-z^2=0$ or $2n-4$. 
\item{$Z$}
$\Lambda\otimes \Q\cup \infty$ See 4.2.
\item{$\theta,\Delta$}
Modular forms. See 5.1.
\item{$\sigma(R),$}
$\sigma(u^\perp)$ The opposition involutions of $R$ or $u^\perp$. See 1.3. 
\item{$\sigma_u$} 
See 3.7 and 1.5.
\item{$\tau_u$}
See 3.2.
\item{$\cup$}
Union of
\item{$w$}
The Weyl vector of the Weyl chamber $D$ of $II_{25,1}$.
\item{$\cdot0$}
The group of automorphisms of the Leech lattice $\Lambda$.
\item{$\cdot\infty$}
The group of affine automorphisms of the Leech lattice $\Lambda$.
\item{$\cdot 2$,}
$\cdot 3$, $HS$, $M_{12}$, $M_{24}$, $McL$, and so on, stand for some
of the sporadic simple groups.

\proclaim Introduction. 

In this thesis we use the Leech lattice $\Lambda$ and the 26
dimensional even unimodular lattice $II_{25,1}$ to study other
lattices, mainly unimodular ones of about 25 dimensions.

Chapter 0 gives some definitions and quoted results. Chapter 1 gives
several auxiliary results about root systems, most of which are
already known.

In chapter 2 we give new proofs of old results. The main result is a
``conceptual'' proof of Leech's conjecture that $\Lambda$ has covering
radius $\sqrt{2}$. (One proof is said to be more conceptual than
another if it gives less information about the objects being studied.)
This was first proved in [C-P-S] by a rather long calculation. This
result is important because it is used in Conway's proof that
$II_{25,1}$ has a Weyl vector and its Dynkin diagram can be identified
with the points of $\Lambda$. We also give new proofs of the existence
and uniqueness of $\Lambda$, and give a uniform proof that the 23
``holy constructions'' of [C-S b] all give the Leech lattice.

In chapter 3 we give an algorithm for finding all orbits of vectors of
$II_{25,1}$ under $\Aut(II_{25,1})$. Any two primitive vectors of the
same positive norm are conjugate under $\Aut(II_{25,1})$ and orbits of
primitive vectors of norm 0 correspond to the 24 Niemeier lattices, so
the algorithm is mainly concerned with negative norm vectors. There
are 121 orbits of vectors of norm $-2$ (given in table $-2$) and 665
orbits of vectors of norm $-4$ (given in table $-4$). The main reason
for the interest in these vectors is that the norm $-4$ vectors are in
1:1 correspondence with the 665 25-dimensional positive definite
unimodular lattices. Previous enumerations of such lattices were done
in [K] which found the eight 16-dimensional ones, in [C-S a] which
found the 117 23-dimensional ones, and by Niemeier who found the 24 even
24-dimensional ones. All of these used essentially the same method as
that found by Kneser. This method becomes much more difficult to use
in 24 dimensions and seems almost impossible to use in 25 dimensions.

Chapter 4 looks at the vectors of $II_{25,1}$ in more detail and we
prove several curious identities for lattices of dimension at most
25. For example, if $L$ is a unimodular lattice of dimension $n\le 23$
then the norm of the Weyl vector of $L\oplus I^{25-n}$ is a square and
we can use this to give a construction of the Leech lattice for each
such lattice $L$. This is nontrivial even if $L$ is 0-dimensional; in
this case the norm of the Weyl vector is $0^2+1^2+\cdots+24^2=70^2$,
and the construction of the Leech lattice in this case is that given
in [C-S c]. Sections 4.7 to 4.10 examine vectors of norms $\le -6$,
and in 4.11 we use negative norm vectors in $II_{25,1}$ to prove a
rather strange formula (4.11.2) for the smallest multiple of a
``shallow hole'' of $\Lambda$ that is in $\Lambda$.

Chapter 5 gives some applications of theta functions to our
lattices. The height of a primitive vector $u$ of $II_{25,1}$ seems to
depend linearly on the theta function of the lattice $u^\perp$, and in
$5.1$ to $5.4$ we prove this when $u$ has norm $-2$ or $-4$. (The
height of a vector of a fundamental domain of $II_{25,1}$ is $-(u,w)$
where $w$ is the Weyl vector of the fundamental domain.) In 5.5 we
give an algorithm for finding the 26 dimensional unimodular lattices
and use it to show that there is a unique such lattice with no
roots. In 5.6 we show that any 25 dimensional unimodular lattice with
no roots has determinant at least 3, so this means we have found all
unimodular lattices with no roots of dimension at most 26. (Apart from
the Leech lattice and the 26 dimensional one and the trivial
0-dimensional one, there are two others of dimension 23 and 24 which
are both closely related to the Leech lattice.) Finally in 5.7 we
construct a 27-dimensional unimodular lattice with no roots (which is
probably not unique). (Remark added 1999: Bacher and Venkov have shown
that there are exactly 3 27-dimensional unimodular lattices with no
roots.)

In chapter 6 we calculate the automorphism groups of the unimodular
lattices $I_{n,1}$ for $n\le 23$, or more precisely the groups $G_n$
which are the quotients of the automorphism groups by the subgroups
generated by reflections and $-1$. Vinberg showed that $G_n$ was
finite if and only if 
 $n\le 19$, and Conway and Sloane showed that for $n\le 19$
the groups $G_n$ could be identified in a natural way with subgroups
of the Conway group $\cdot 0$ by mapping (the even sublattice of)
$I_{n,1}$ into $II_{25,1}$. Chapter 6 extends Conway and Sloane's work
to $n\le 23$. For $n=20$, 21, and 22, $G_n$ is an amalgamated product
of two subgroups of $\cdot\infty$ and for $n=23$ $G_n$ is a direct
limit of 6 subgroups of $\cdot\infty$. We also find the automorphism
groups of a few other Lorentzian lattices; for example if $L$ is the
even sublattice of $I_{21,1}$ then the reflection group of $L$ has
finite index in the automorphism group of $L$, and has a fundamental
domain with $168+42$ sides. We find the structure of $G_n$ by using
$II_{25,1}$ to construct a complex acted on by $G_n$, whose dimension
is the virtual cohomological dimension of $G_n$.

Vectors of norm $\le 2$ in $II_{25,1}$ can also be thought of as roots
of the ``monster Lie algebra''. In chapter 7 we calculate the
multiplicities of some of these roots. The multiplicity of a root of
norm $-2n$ seems to be closely related to the coefficient $c_n$ of
$q^n$ in
$$
q^{-1}(1-q)^{-24}(1-q^2)^{-24}\cdots=q^{-1}+24+324q+3200q^2+\cdots
$$
We show that for every $n$ there are roots of norm $-2n$ with
multiplicity $c_n$, and show that of the 121 orbits of norm $-2$
vectors, 119 have multiplicity 324. (The other two have multiplicity 0
and 276.) I have recently seen the preprint [F] which proves the
results in 7.1 and 7.2, and proves that $c_n$ is an upper bound for
the multiplicity of nonzero roots by finding a representation of the
monster Lie algebra containing the adjoint representation for which a
nonzero vector of norm $-2n$ has multiplicity $c_n$.

\proclaim Chapter 0  Definitions and assumed results. 

In this chapter we give some standard definitions and results about
lattices and root systems.

\proclaim 0.1 Definitions. 

Most of the results here can be found in [S a].

A {\it lattice} $L$ is a finitely generated free $Z$-module with an
integer valued bilinear form, written $(x,y)$ for $x$ and $y$ in
$L$. The {\it type} of a lattice is even (or II) if the {\it norm}
$x^2=(x,x)$ of every element $x$ of $L$ is even, and odd (or I)
otherwise. If $L$ is odd then the vectors in $L$ of even norm form an
even sublattice of index 2 in $L$. $L$ is called {\it positive
definite, Lorentzian, nonsingular,} and so on if the real vector space
$L\otimes \R$ is. 

A lattice is called {\it indecomposable} if it is not the direct sum of two nonzero sublattices. A result of Eichler states that any definite lattice is the direct sum of its indecomposable sublattices. In particular any positive definite lattice $L$ can be written uniquely as $L=L_1\oplus I^n$ where $L_1$ has no vectors of norm 1 and $I$ is the one dimensional lattice generated by a vector of norm 1. 

If $L$ is a lattice then $L'$ denotes its {\it dual} in $L\otimes \R$;
in other words the vectors of $L\otimes \R$ that have integral inner products
with all elements of $L$. $L'$ contains $L$ and if $L$ is nonsingular
then $L'/L$ is a finite abelian group whose order is called the {\it
determinant} of $L$. (If $L$ is singular we say it has determinant 0.)
The rational valued quadratic form on $L'$ gives a quadratic form on
$L'/L$ which is defined mod 1 if $L$ is odd and mod 2 if $L$ is even.
$L$ is called {\it unimodular}, {\it bimodular}, or {\it trimodular}
if its determinant is 1, 2 or 3.

Even unimodular lattices of a given signature and dimension exist if
and only if there is a real vector space with that signature and
dimension and the signature is divisible by 8. Any two indefinite
unimodular lattices with the same type, dimension, and signature are
isomorphic. $I_{m,n}$ and $II_{m,n}$ $(m\ge 1, n\ge 1)$ are the
unimodular lattices of dimension $m+n$, signature $m-n$, and type $I$
or $II$.

A vector $v$ in a lattice $L$ is called {\it primitive} if $v/n$ is
not in $L$ for any $n>1$. A {\it root} of a lattice $L$ is a primitive
vector $r$ of $L$ such that reflection in the hyperplane $r^\perp$
maps $L$ to itself.  This reflection maps $v$ in $L$ to
$v-2r(v,r)/(r,r)$. Any vector $r$ in $L$ of norm 1 or 2 is a root and
in these cases the reflection in $r^\perp$ fixes all elements of
$L'/L$. We often use ``root'' to mean ``root of norm 2''. A lattice is
said to have {\it many roots} if the roots generate the vector space
$L\otimes \R$.

If $L$ is unimodular then there is a unique element $c$ in $L/2L$ such
that $(c,v)\equiv v^2\bmod L$ for all $v$ in $L$. $c$ or any inverse
image of $c$ in $L$ is called a {\it characteristic vector } of $L$,
and its norm is congruent to the signature of $L\bmod 8$.  (Remark
added 1999: these days the characteristic vector is usually called the
parity vector.)

\proclaim 0.2 Neighbors.

Two lattices $L_1$, $L_2$ in a vector space are called {\it neighbors}
if their intersection has index two in each of them. The main idea in
[K] is to find unimodular lattices by starting with one unimodular
lattice, finding all its neighbors, and repeating this with the new
lattices obtained. It is known that all unimodular lattices of the
same dimension and signature can be obtained up to isomorphism in this
way, starting with one such lattice. In this section we will describe
the neighbors of unimodular lattices.

Let $A$ be an odd unimodular lattice of signature $8n$ and let $A_2$ be
its sublattice of vectors of even norm. If $c$ is a characteristic
vector of $A$ and $a$ is a vector of $A$ of odd norm, then $c/2$ and
$a$ represent different non-zero elements of $A_2'/A_2$ with
$(c/2)^2\equiv 0\bmod 2$, $a^2\equiv 1\bmod 2$, and $(a,c/2)\equiv {1/
2} \bmod 1$, so $A_2'/A_2$ is isomorphic to the group $(Z/2Z)^2$ and
the norms of its non-zero elements are $0,0$, and $1\bmod 2$. If we add a
representative $x$ of some non-zero element of $A_2'/A_2$ to $A_2$ we
get a unimodular lattice which is $A$ if $x$ has odd norm and an even
lattice if $x$ has even norm. Hence $A$ has exactly two even
neighbors.

Conversely if we start with an even unimodular lattice $B$ of
signature $8n$ then its sublattices of index 2 are in 1:1
correspondence with the non-zero elements $b$ of $B/2B$. If $b$ is
such an element we will write $B_b$ for the elements of $B$ that have
even inner product with $b$. Let $a$ be any element of $B$ that has
odd inner product with $b$. $B_b'/B_b$ is isomorphic to $(\Z/2\Z)^2$
and its nonzero elements are represented by $b/2$, $a$, and
$b/2+a$. The norms of these elements are
$$\eqalign{ 0,0,1\bmod 2&\hbox{~~~~if~~~~} b^2\equiv 0\bmod 4,\hbox{ and}\cr
1/2,0,3/2\bmod 2&\hbox{~~~~if~~~~} b^2\equiv 2\bmod 4.\cr }$$ (Note that
$b^2\bmod 4$ only depends on $b\bmod 2B$.) Hence the odd neighbors of
$B$ are in $1:1$ correspondence with the non-zero elements of
$B'_b/B_b$ that have norm $0\bmod 4$. Given such an element $b$ the
corresponding odd neighbor of $B$ is the unimodular odd lattice
containing $B_b$.

If the odd lattice $A$ has a vector $a$ of norm 1 then $2a$ is a
vector of norm 4 in both even neighbors of $A$, and these neighbors
are exchanged by reflection in $a^\perp$. (The two even neighbors of
$A$ can sometimes be exchanged by an automorphism of $A$ even if $A$
has no vectors of norm 1.) If $B$ is an even neighbor of $A$ then the
vector $b=2a$ gives back the lattice $A$ as the odd lattice containing
$B_b$. Conversely if $b$ is any norm 4 vector in the even lattice $B$
then the odd lattice of $B_b$ contains the vector $b/2$ of norm
1. Hence there is a 1:1 correspondence between
\item{1.} 
Pairs $(a,A)$ where $a$ is a vector of norm 1 in the odd unimodular
lattice $A$ of signature $8n$, and
\item{2.} 
Pairs $(b,B)$ where $b$ is a vector of norm 4 in the even lattice $B$
of signature $8n$.

In particular we can find all positive definite unimodular lattices of
dimension less than $8n$ by finding all vectors of norm 4 in $8n$
dimensional even unimodular lattices, because in a positive definite
lattice all vectors of norm 1 are conjugate under its automorphism
group. This was used in [C-S a] to find all 23 dimensional unimodular
lattices from the Niemeier lattices.

Remark. We can use this to show that if $B$ is an even unimodular
lattice that is not positive definite than all vectors of norm 4 in
$B$ are conjugate under $\Aut(B)$. This is generalized later.

\proclaim 0.3 Root systems. 

The results here can be found in [B]. We summarize some definitions
and basic properties of finite root systems. Chapter 1 contains more
detailed information about them.

``Root system'' will mean ``root system all of whose roots have norm
2'' unless otherwise stated, so we only consider components of type
$a_n$, $d_n$, $e_6$, $e_7$, $e_8$. Most of the results about root
systems in this thesis can easily be generalized to root systems
$b_n$, $c_n$, $f_4$ and $g_2$ but this would make everything more
complicated to state. (We use small latter $x_n$ to stand for
spherical Dynkin diagrams, and capital letters $X_n$ to stand for the
corresponding affine Dynkin diagrams.)

The norm 2 vectors in a positive definite lattice $A$ form a root
system which we call the root system of $A$. The hyperplanes
orthogonal to these roots divide $A\otimes \R$ into regions called
{\it Weyl chambers}. The reflections in the roots of $A$ generate a
group called the {\it Weyl group} of $A$, which acts simply
transitively on the Weyl chambers of $A$. Fix one Weyl chamber
$D$. The roots $r_i$ that are orthogonal to the faces of $D$ and
that have inner product at most 0 with the elements of $D$ are called
the simple roots of $D$. (These have opposite sign to what are usually
called the simple roots of $D$. This is caused by the irritating fact
that the usual sign conventions for positive definite lattices are not
compatible with those for Lorentzian lattices. With the convention
used here something is in the Weyl chamber if and only if it has inner
product at most 0 with all simple roots, and a root is simple if and
only if it has inner product at most 0 with all simple roots not equal
to itself.)

The {\it Dynkin diagram} of $D$ is the set of simple roots of $D$. It
is drawn as a graph with one vertex for each simple root of $D$ and
two vertices corresponding to the distinct roots $r,s$ are joined by
$-(r,s)$ lines. (If $A$ is positive definite then two vertices are
always joined by 0 or 1 lines. We will later consider the case that
$A$ is Lorentzian and then its Dynkin diagram may contain multiple
bonds, but these are not the same as the multiple bonds appearing in
$b_n$, $c_n$, $f_4$, and $g_2$.) 
The Dynkin diagram of $A$ is a union of
components of type $a_n$, $d_n$, $e_6$, $e_7$ and $e_8$. The {\it Weyl
vector} $\rho$ of $D$ is the vector in the vector space spanned by
roots of $A$ that has inner product $-1$ with all simple roots of
$D$. It is in the Weyl chamber $D$ and if there are a finite number of roots it
is equal to half the sum of the
positive roots of $D$, where a root is called positive if its inner
product with any element of $D$ is at least 0. The {\it height} of an
element $a$ of $A$ is $-(a,\rho)$, so a root of $A$ is simple if and
only if it has height 1.

The {\it Coxeter number $h$} of $A$ is defined to be (number of
roots of $A$)/(dimension of $A$). The Coxeter number of a
component of $A$ is defined as the Coxeter number of the lattice
generated by the roots of that component. Components $a_n$, $d_n$,
$e_6$, $e_7$ and $e_8$ have Coxeter numbers $n+1$, $2n-2$, 12, 18, and
30 respectively. For each component $R$ of the Dynkin diagram of $D$
there is an orbit of roots of $A$ under the Weyl group, and this orbit
has a unique representative $v$ in $D$, which is called the {\it
highest root} of that component. We can write
$-v=n_1r_1+n_2r_2+\cdots$, where the $r$'s are the simple roots of $R$
and the $n$'s are positive integers called the {\it weights} of the
roots $r_i$. The sum of the $n$'s is $h-1$, where $h$ is the Coxeter
number of $R$, and the height of $v$ is $1-h$. The {\it extended
Dynkin diagram} of $A$ is the simple roots of $A$ together with the
highest roots of $A$, and it is the Dynkin diagram of the positive
semidefinite lattice $A\oplus O$, where $O$ is a one-dimensional
singular lattice. We say that the highest roots of $A$ have weight
one. Any point of weight one of a Dynkin diagram is called a tip of that
Dynkin diagram. We write $X_n$ for the extended Dynkin diagram
corresponding to the Dynkin diagram $x_n$.

%Here are the connected extended Dynkin diagrams with their weights:
%To do: insert picture of extended Dynkin diagrams here. 

The automorphism group of $A$ is a split extension of its Weyl group
by $N$, where $N$ is the group of automorphisms of $A$ fixing $D$. $N$
acts on the Dynkin diagram of $D$ and $\Aut(A)$ is determined by its
Dynkin diagram $R$, the group $N$, and the action of $N$ on $R$.
There is a unique element $i$ of the Weyl group taking $D$ to $-D$,
and $-i$ is called the {\it opposition involution} of $D$ and denoted
by $\sigma$ or $\sigma_D$. $\sigma$ fixes $D$ and has order 1 or
2. (Usually $-\sigma$ is called the opposition involution.)

If $A$ is Lorentzian or positive semidefinite then we can still talk
about its root system and $A$ still has a fundamental domain $D$ for
its Weyl group and a set of simple roots. $A$ may or may not have a
Weyl vector but does not have highest roots or an opposition
involution.
\proclaim 0.4 Lorentzian lattices any hyperbolic geometry. 

We  describe the geometry of Lorentzian lattices and its relation to
hyperbolic space, and give Vinberg's algorithm for finding the
fundamental domains of hyperbolic reflection groups.

Let $L$ be an $(n+1)$-dimensional Lorentzian lattice (so $L$ has
signature $n-1$). Then the vectors of $L$ of zero norm form a double
cone and the vectors of negative norm fall into two components. The
vectors of norm $-1$ in one of these components form a copy of
$n$-dimensional hyperbolic space $H_n$. The group ${\rm Aut}(L)$ is a
product $(\Z/2\Z)\times {\Aut}_+(L)$, where $\Z/2\Z$ is generated by $-1$
and ${\rm Aut}_+(L)$ is the subgroup of ${\rm Aut}(L)$ fixing each
component of negative norm vectors. See [Vi] for more
details. (If we wanted to be more intrinsic we could define $H_n$ to
be the set of nonzero negative definite subspaces of $L\otimes \R$.)

If $r$ is any vector of $L$ of positive norm then $r^\perp$ gives a
hyperplane of $H_n$ and reflection in $r^\perp$ is an isometry of
$H_n$. If $r$ has negative norm then $r$ represents a point of $H_n$
and if $r$ is nonzero but has zero norm then it represents an infinite
point of $H_n$.

The group $G$ generated by reflections in roots of $L$ acts as a
discrete reflection group on $H_n$, so we can find a fundamental
domain $D$ for $G$ which is bounded by reflection hyperplanes. The
group ${\rm Aut}_+(L)$ is a split extension of this reflection group by
a group of automorphisms of $D$.

Vinberg  gave an algorithm for finding a fundamental domain $D$
which runs as follows. Choose a vector $w$ in $L$ of norm at most 0;
we call $w$ the {\it controlling vector}.  The
hyperplanes of $G$ passing through $w$ form a root system which is
finite (or affine if $w^2=0$); choose a Weyl chamber $C$ for this root
system. Then there is a unique fundamental domain $D$ of $G$
containing $w$ and contained in $C$, and its simple roots can be found
as follows.
\item{Step 0:} 
Take all the simple roots of $C$ as some of the simple roots of $D$.
Enumerate the roots $r_1,r_2,\ldots$ of $L$ that have negative inner
product with $w$ in increasing order of the distances of their
hyperplanes from $w$, in other words in increasing order of
$-(r_i,w)/r_i^2$.
\item{Step n:} ($n\ge 1$, $n$ running through some countable ordinals):
Take the root $r_n$ as a simple root for $D$ if and only if it has
inner product at most 0 with all the simple roots of $D$ that we
already have. (It is sufficient to check this for the simple roots
whose hyperplanes are strictly closer to $w$ than the hyperplanes of
$r_n$.)

This algorithm produces a finite or countable set of roots which are
the simple roots of a fundamental domain $D$ of $G$, If at any time
the simple roots we already have bound a region in $H_n$ of finite
volume then no further roots are accepted by the algorithm, so the
Dynkin diagram of $D$ is finite.

If $H_n$ has dimension greater than 1 then a convex region of $H_n$
bounded by hyperplanes has finite volume if and only if it contains a
finite number of infinite points of $H_n$.

\proclaim Chapter 1 Root systems. 

In this chapter we collect several results on root systems and
lattices for which there is no convenient reference. Most of them are
already known.

\proclaim 1.1 Norms of Weyl vectors. 

Here we find a formula for the norm of the Weyl vector of a simple
root system.

{\bf Notation.} Let $R$ be a simple root system with rank $n$, Coxeter
number $h$ and Weyl vector $\rho$. ($\rho$ is half the sum of the
positive roots.) We do not assume that all roots of $R$ have the same
length or norm 2. Let $(,)$ be a bilinear form on $R$ invariant under
the Weyl group (such a form is unique up to scalar multiples) and let
$n^*$ be the sum of the norms of the roots of $R$. ($n^*$ depends on
$(,)$.)  

\proclaim Lemma 1.1.1. 
If $(,)$ is the Killing form on $R$ (``$\hbox{forme bilineaire canonique}$'')
then $n^*=n$.

Proof. By [B, Chapter V no. 6.2 Corollary to theorem 1] we have
$$\sum_\alpha (\alpha,\beta)^2/(\alpha,\alpha) = (\beta,\beta)h$$
for any form $(,)$, where $\sum_\alpha$ means sum over all roots $\alpha$. 
By [B, Chapter VI Section 1 No. 1.12] the Killing form $(,)$ satisfies
$$(x,x)=\sum_\alpha(x,\alpha)^2$$
so
$$\eqalign{
hn^*&=
h\sum_\beta\beta^2\cr
&=\sum_\beta\sum_\alpha(\alpha,\beta)^2/(\alpha,\alpha)\cr
&=\sum_\alpha(\alpha,\alpha)^{-1}\sum_\beta(\alpha,\beta)^2\cr
&= \sum_\alpha1\cr
&= hn\cr
}$$
Q.E.D.

\proclaim Theorem 1.1.2. $\rho^2=(h+1)n^*/24$

Proof. The ration of both sides does not depend on which bilinear form
$(,)$ we choose, so we can assume that $(,)$ is the Killing form. In
this case the strange formula of Freudenthal and de Vries shows that
$\rho^2=\dim(G)/24$ where $G$ is the Lie algebra of $R$. We have
$\dim(G)=(h+1)n=(h+1)n^*$ (by lemma 1.1.1) and this implies the
theorem. Q.E.D.

\proclaim 
Corollary 1.1.3. If all roots of $R$ have norm 2 then
$\rho^2=h(h+1)n/12$. More specifically if $R$ is $a_n$, $d_n$, $e_6$,
$e_7$, or $e_8$ then $\rho^2$ is $n(n+1)(n+2)/12$, $(n-1)n(2n-1)/6$,
78, $399/2$ or 620 respectively.

Proof. $R$ has $hn$ roots so $n^*=2hn$. The result now follows from
theorem 1.1.2. 

Theorem 1.1.2 occurs in [M b formula 8.2]. 

\proclaim 1.2 Minimal vectors. 

If $L$ is a lattice and $R$ its sublattice generated by roots then we
will find a canonical set of representatives for $L/R$, called the
minimal vectors of $L$ (sometimes called minuscule vectors).

{\bf Notation.} Let $R$ be a simple root system all of whose roots
have norm 2. We also write $R$ for the lattice generated by $R$ and
$R'$ for the dual of $R$, so that $R'/R$ is a finite abelian group
whose order is the determinant $d$ of the Cartan matrix of $R$. We will
call a vector of $R'$ minimal if it has inner product 0 or $\pm 1$
with all roots of $R$ and it lies in the Weyl chamber of $R$.

\proclaim Lemma 1.2.1. 
If $r'$ is a nonzero minimal vector of $R$ then there is exactly one
simple root $r$ of $R$ with $(r,r')\ne 0$ and this root has weight
one. Every simple root of weight 1 comes from a unique minimal vector
in this way. The minimal vectors form a complete set of
representatives of the elements of $R'/R$.

Proof. We first show that every element $\hat r$ of $R'/R$ is
represented by a minimal vector. Let $r'$ be an element of $R'$
representing $\hat r$ of smallest possible norm. The Weyl group of $R$
fixes all elements of $R'/R$ so we can assume that $r'$ is in the Weyl
chamber of $R$. For any root $r$, $r'-r$ also represents $\hat r$ so
$(r'-r)^2\ge r'^2$. This implies that $(r,r')\le 1$ for all roots $r$
so $r'$ must be minimal.

Now let $v=-m_1r_1-m_2r_2-\cdots$ be the highest root of $R$, where
the $r_i$'s are the simple roots of $R$ of weights $m_i$. If $r'$ is
nonzero and minimal then $(r',v)=1$ and $(r',r_i)\le 0$ so $(r',r_i)$
must be $-1$ for one simple root $r_i$ and zero for all other simple
roots. Conversely if $r'$ has inner product $-1$ with one simple root
and 0 with all other simple roots then it is in the Weyl chamber and
has inner product 1 with $v$, so it has inner product at most 1 with
any root and hence is minimal.

We have shown that there is a bijection between nonzero minimal
vectors and tips of the Dynkin diagram of $R$. The number of such tips
is always equal to $d-1$ where $d$ is the order of $R'/R$. As every
element of $R'/R$ is represented by a minimal vector, this shows that
it must be represented by a unique minimal element. Q.E.D.

Remark. If $R$ does not have all its roots of norm 2 then minimal
vectors are defined to be those that have inner product 0 or
$\pm r^2/2$ with any root $r$. Minimal vectors correspond to simple
roots of weight 1 in the {\it dual} root system of $R$. $b_n$ and
$c_n$ have one nonzero minimal vector while $f_4$ and $g_2$ have none.

Here are the norms of the minimal vectors of $R$:
\item{$R$} Possible norms
\item{$a_n$} 0, $m(n+1-m)/(n+1)$ ~~~~~$1\le m\le n$
\item{$d_n$} 0, 1, $n/4$, $n/4$
\item{$e_6$} 0, $4/3$, $4/3$
\item{$e_7$} 0, $3/2$
\item{$e_8$} 0

Now let $L$ be a positive definite lattice and $R$ the lattice
generated by norm 2 roots of $L$. We say that a vector of $L$ is
minimal if it is in the Weyl chamber of $L$ and has inner product 0 or
$\pm 1$ with all roots of $R$.

\proclaim Lemma 1.2.2. 
The minimal vectors of $L$ form a complete set of representatives for $L/R$. 

Proof. We can show as in the proof of lemma 1.2.1 that every element
of $L/R$ is represented by a minimal vector of $L$. It also follows
from lemma 1.2.1 that two minimal vectors of $L$ cannot be congruent
mod $R$, so every element of $L/R$ is represented by a unique minimal
vector. Q.E.D.

\proclaim 1.3 The opposition involution. 

Let $R$ be a root system with a Weyl chamber $D$. There is a unique
element $i$ of the Weyl group taking $D$ to $-D$ and $\sigma=-i$ is
called the opposition involution of $D$. $\sigma$ has order 1 or 2 and
fixes $D$ and the Weyl vector $\rho$ of $D$. In this section we will
prove some properties of $\sigma$.

Suppose that $R$ is simple with all roots of norm 2, and write $R$ for
the lattice generated by $R$ and $R'$ for its dual.

\proclaim Lemma 1.3.1. 
$\sigma$ has order 1 if and only if every element of $R'/R$ has order 1 or 2. 

Proof. $\sigma$ is an automorphism of the Dynkin diagram of $R$ and
has order 1 if and only if it fixes all the tips of $R$. By lemma
1.2.1 this is true if and only if $\sigma$ fixes all elements of
$R'/R$. $i$ fixes all elements of $R'/R$ as it is in the Weyl group,
so $\sigma$ acts as $-1$ on $R'/R$ and so it fixes all elements of
$R'/R$ if and only if all elements of $R'/R$ have order 1 or 2.

(This result occurs as one of the exercises in [B].)

\proclaim Theorem 1.3.2. 
Suppose $L$ is a lattice with root system $R$ and Weyl chamber
$D$. Then $\sigma$ acts as 1 on all components of $R$ of type $a_1$,
$d_{2n}$, $e_7$, and $e_8$, as the unique non-trivial automorphism on
components of type $a_n$ $(n\ge 2)$, $d_{2n+1}$ and $e_6$, and as $-1$
on all vectors that are orthogonal to all roots.

Proof. The only simple root systems $X$ for which $X'/X$ has an
element of order greater than 2 are $a_n$ $(n\ge 2)$, $d_{2n+1}$, and
$e_6$. The theorem follows from this and lemma 1.3.1. Q.E.D.

For a positive definite lattice $L$ we define $S(L)$ to be the
dimension of the fixed space of $\sigma$.  \proclaim Lemma
1.3.3. $S(L)$ is equal to the sum of the following number for each
component of the root system of $L$: $n$ for $a_{2n-1}$ and $a_n$,
$2n$ for $d_{2n}$ and $d_{2n+1}$, $4$ for $e_6$, $7$ for $e_7$, and
$8$ for $e_8$. 

Proof. $S(L)$ is the sum of $S(X)$ for each component $X$ of the
Dynkin diagram of $L$. For any component $X$, $S(X)$ is equal to
(number of points of $X$ fixed by $\sigma$)+(1/2)(number of points of
$X$ not fixed by $\sigma$) and by working out this number for each
possible component we get the result stated in the lemma. Q.E.D.

\proclaim Lemma 1.3.4. 
If $r$ is any root of a lattice $L$ then $S(L)=1+S(r^\perp)$.

Proof. We can assume that $r$ is a highest root of $L$. If $\sigma$
and $\sigma_0$ are the opposition involutions of $L$ and $r^\perp$
then $\sigma(r)=r$, $\sigma_0(r)=-r$, and $\sigma=\sigma_0$ on
$r^\perp$. Q.E.D.

\proclaim Lemma 1.3.5. 
All maximal sets of pairwise orthogonal roots in a lattice $L$ are
conjugate under then Weyl group and the number of roots in such a set
is $S(L)$.

Proof. This follows from lemma 1.3.4 by induction on the rank of the
root system of $L$. Q.E.D.

\proclaim Lemma 1.3.6. 
$S(L)+4\rho^2+(\hbox{number of roots of $L$})/2\equiv 0\bmod 4$.

Proof. Choose any root $r$ of $L$. As we go from $r^\perp$ to $L$ we
increase $S$ by 1, $\rho^2$ by $(h-1)^2/2$ and the number of roots by
$4h-6$ where $h$ is the Coxeter number of the component of $r$. Hence
we increase $S+4\rho^2+\hbox{roots}/2$ by $1+2(h-1)^2+2h-3=2h(h-1)$
which is divisible by 4. Lemma 1.3.6 now follows by induction on the
number of roots of $L$. Q.E.D.

\proclaim 1.4 Maximal sub root systems. 

Let $R$ be a simple spherical root system and let $R_0$ be its
extended Dynkin diagram. We find the maximal sub roots systems of $R$
and then prove a few technical lemmas which will be used in 4.4.

\proclaim Theorem 1.4.1.
The Dynkin diagrams $R_1$ of maximal sub root systems of $R$ are
obtained as follows:
\item{(1)} 
Either: Delete one point of prime weight $p$ from $R_0$. In this case the lattice of $R_1$ has index $p$ in the lattice of $R$. 
\item{(2)}
Or: Delete two points of weight 1 from $R_0$. In this case the lattice
of $R_1$ has dimension one less than that of $R$.

Proof. Let $r$ be a root in $R$ but not in $R_1$ that is in the Weyl
chamber of $R_1$ (for some choice of simple roots of $R_1$). Then $r$
has inner product 0 or $-1$ with every simple root of $R_1$, so $r$
together with the simple roots of $R_1$ forms a Dynkin diagram
$R'$. The root system generated by $r$ and $R_1$ is $R$ because $R_1$
is maximal, so $R'$ is either the Dynkin diagram or the extended
Dynkin diagram of $R$. 

First suppose that $R'$ is the extended Dynkin diagram of $R$. $r$
cannot be a point of weight 1 of $R'$ (otherwise $R_1$ would be $R$)
so $r$ has weight $m\ge 2$. Choose a point $r_1$ of weight 1 of $R'$
and take a set of simple roots of $R$ to be $R'-r_1$ (so that $r_1$ is
the highest root). Let $v$ be the vector that has inner product $-1$
with $r$ and 0 with the other simple roots of $R'$ other than
$r_1$. Then the root system of roots in $v^\perp$ has $R_1-r-r_1$ as
its Dynkin diagram (as $v$ is in the Weyl chamber of $R$) and
$(v,r_1)=m$ as $r_1$ is the highest root of $R$, so the root system
$R_1$ is exactly those roots $s$ with $m|(s,v)$. $m$ must be prime as
$R$ contains roots $t$ with $(t,v)=n$ for all $n$ between $1$ and $m$,
so that if $m$ was not prime $R_1$ would not be maximal. Conversely if
$m$ is prime then the lattice of $R_1$ has prime index in the lattice
of $R$, so $R_1$ is a maximal sub root system of $R$.

If $R'$ is the Dynkin diagram of $R$ we can use similar reasoning to
prove case 2 of 1.4.1. Q.E.D.

[D] finds all sub root systems of root systems.

Remark. There is of course a similar theorem for $b_n$, $c_n$, $f_4$,
and $g_2$ except that it is also necessary to consider the extended
Dynkin diagram of the dual root system.

\proclaim Lemma 1.4.2. 
Let $R$ be an extended Dynkin diagram of Coxeter number $h$, and let
$r_1,r_2,\ldots,r_n$ be some subset of its vertices with weights
$w_1,\ldots, w_n$. Then
$$(\sum_iw_ir_i,\rho)=h-\sum_iw_i$$ 
where $\rho$ is the Weyl vector of
the Dynkin diagram $R$ with the points $r_i$ deleted.

Proof. Let $r_i'$ be the other points of $R$ with weights $w_i'$. Then 
$$\eqalign{
\sum_iw_ir_i+\sum_iw_i'r_i'&=0\cr
\sum_iw_i+\sum_iw_i'&=h\cr
(r_i',\rho)&=-1\cr
}$$
The lemma follows from these three equalities. Q.E.D. 

Let $L_2$ be a lattice of index 1 or 2 in the lattice $L$, where the
root system of $L$ is simple. By 1.4.1 the Dynkin diagram of $L_2$ is
obtained from the extended Dynkin diagram of $L$ by deleting one or
two points of weight 1 or a point of weight 2.

\proclaim Lemma 1.4.3. 
Notation is as above. If $\rho$ is the Weyl vector of $L_2$ and $r$ is
one of the (one or two) points deleted from the extended Dynkin
diagram of $L$ then
$$(r,\rho)=h/m-1,$$ 
where $h$ is the Coxeter number of $L$ and $m$ is
the sum of the weights of the deleted points (so $m=1,1+1$, or 2).

Proof. This follows from 1.4.2 and from the fact that any two points
of weight 1 of an extended Dynkin diagram $R$ can be exchanged by an
automorphism of $R$, so all the points deleted from the extended
Dynkin diagram of $L$ have the same inner product with $\rho$. Q.E.D.

Remark. If $m=1$ then 1.4.3 is the usual formula for the height of the
highest root of a root system.

\proclaim 1.5 Automorphisms of the fundamental domain. 

{\bf Notation.} $L$ is a Lorentzian lattice whose reflection group has
a fundamental domain $D$. $u$ is a negative norm vector in $D$.

We construct some automorphisms of $D$ from certain vectors $u$. 

\proclaim Theorem 1.5.1. 
Suppose that reflection in $u^\perp$ is an automorphism  of $L$. Put
$$\sigma_u=r_u\sigma(u^\perp)$$
where $r_u$ is reflection in $u^\perp$
and $\sigma(u^\perp)$ is the opposition involution of the root system
of roots of $L$ in $u^\perp$. Then $\sigma_u$ is an automorphism of $L$
fixing $D$ and $u$ and has order 1 or 2.

Proof. $r_u$ and $\sigma(u^\perp)$ are commuting elements of $\Aut(L)$
of orders 1 or 2, so $\sigma_u$ is an automorphism of $L$ of order 1
or 2. Both $r_u$ and $\sigma(u^\perp)$ multiply $u$ by $-1$ so
$\sigma_u$ fixes $u$. If $D_0$ is the fundamental domain of the finite
reflection group of $u^\perp$ then both $r_u$ and $\sigma(u^\perp)$,
and hence $\sigma_u$, fix $D_0$. By applying Vinberg's algorithm with
$u$ as a controlling vector we see that any automorphism of $L$ fixing
$u$ and $D_0$ must fix $D$, so $\sigma_u$ fixes $D$. Q.E.D.

Remark. In 3.7 we will use a similar but more complicated construction
to get automorphisms of $D$ from certain other vectors $u$.

\proclaim 1.7 Norm 4 vectors of lattices. 

We describe the norm 4 vectors of a positive definite lattice $A$. The
main use of this is that the norm 4 vectors of $8n$ dimensional even
unimodular lattices correspond to the unimodular lattices of dimension
at most $8n-1$. Most results of this section come from [C-S a].

If $r$ is a norm 4 vector of a lattice $A$ then $r$ has inner product
at most $[\sqrt{4}\times \sqrt{2}]=2$ with all roots of $A$ and if it
has inner product 2 with some root $s$ then $r$ is the sum of the
orthogonal roots $s$ and $r-s$. If it has inner product at most 1 with
all roots then it is a minimal vector and hence by 1.2.2 it is not in
the lattice generated by the roots of $A$, so $r$ is a sum of roots of
$A$ if and only if it is the sum of two orthogonal roots. If $r$ is
the sum of two roots $s_1$ and $s_2$ that are in different components
of the root system of $A$ then $s_1$ and $s_2$ are the only roots that
have inner product 2 with $r$. Suppose that $r$ is the sum of two
roots from the same component $R$ of the root system of $A$. If $s$ is
a root of $R$ then the possible orbits of norm 4 vectors $r$ under the
Weyl group of $R$ correspond to components of $R^*$, where $R^*$ is
the roots of $R$ orthogonal to $s$. There are two orbits it $R$ is
$d_n$ $(n\ge 5)$, 3 orbits if $R$ is $d_4$ (permuted by the triality
automorphism of $d_4$), one orbit if $R$ is $a_n$ $(n\ge 3)$, $e_6$,
$e_7$, or $e_8$ and no norm 4 vectors if $R$ is $a_1$ or $a_2$. The
following table lists the possibilities for the norm 4 vectors of $A$.
\halign{
\hfill$#$&~~\hfill$#$&~~\hfill$#$&~~\hfill$#$&~~\hfill$#$&~~\hfill$#$&~~\hfill$#$\cr
R&R_1&R_2&n&t&m&c\cr
\hbox{none}&\hbox{none}&\hbox{none}&\hbox{?}&\hbox{?}&0&\hbox{?}\cr
xy&xy^*&x^*y^*&n_xn_y&h_x+h_y-2&2&2(h_x+h_y-4)\cr
a_n (n\ge 3)&a_{n-2}&a_1a_1a_{n-4}&(n+1)!/4(n-3)!&2(n-1)&4&2(n-3)\cr
d_4 (\times 3)&a_1&a_3&8&6&6&0\cr
d_n (n\ge 5)&a_1&a_{n-1}&2n&2(n-1)&2(n-1)&0\cr
d_n (n\ge 5)&d_{n-2}&d_{n-4}a_3&2\times n!/3(n-4)!&2(2n-5)&6&2(n-4)\cr
e_6&a_5&d_4&270&16&8&2\cr
e_7&d_6&d_5a_1&756&26&10&2\cr
e_8&e_7&d_7&2160&46&14&1\cr
}
\item{}
If $x$ is a Dynkin diagram then $x^*$ is the Dynkin diagram of roots
of $x$ orthogonal to some fixed root of $x$.
\item{$R$} 
is the set of components of the root system of $A$ containing roots
that have inner product 2 with $r$.
\item{$m$}
is the number of roots that have inner product 2 with $r$.
\item{$t$}
is the maximum possible height of $r$. $t=m+2^{m/2-2}c$.
\item{$n$}
is the number of conjugates of $r$ under the Weyl group of $R$. 
\item{$R_1$}
is the component of $R^*$ corresponding to $r$. 
\item{$R_2$}
is the Dynkin diagram of the roots of $R$ orthogonal to $r$. 
\item{$h_x$} 
and $h_y$ are the Coxeter numbers of $x$ and $y$. 
\item{$n_x$}
and $n_y$ are the numbers of roots of $x$ and $y$.
\item{$c$}
has the following meaning. If $B$ is an $8n$-dimensional even positive
definite unimodular lattice and $r$ is a norm 4 vector in $B$ then $r$
gives an $8n $ dimensional unimodular lattice $A$ containing $m+2>0$
vectors of norm 1 and $2^{1+m/2}c$ characteristic vectors of norm
8. If $A_1$ is the lattice obtained from $A$ by discarding the vectors
of norm 1 then $A_1$ has $c$ characteristic vectors of norm
$\dim(A_1)-8(n-1)$. For example, if $A_1$ is $8n+3$ dimensional then
it has 0 or 2 characteristic vectors of norm 3, depending on whether
the corresponding norm 4 vector in the $8n+8$ dimensional even lattice
comes from a $d_5$ or an $e_6$.

This is essentially the same method as that used in [C-S a] to find
all norm 4 vectors in the Niemeier lattices, and hence all 23
dimensional unimodular lattices.

Remark. If $r$ is a norm 4 vector of $A$ that has inner product 2 with
at least 4 roots of $A$ then the projection of the Weyl vector $\rho$
of $R$ into $R_2$ is a Weyl vector $\rho_2$ of $R_2$. This implies
that $\rho^2=\rho_2^2+t^2/4$.

\proclaim Chapter 2 The Leech lattice. 

%To do: references and proof read. 

\proclaim 2.1.~History.

In 1935 Witt found the 8- and 16-dimensional even unimodular lattices
and more than 10 of the 24-dimensional ones. In 1965 Leech found a
24-dimensional one with no roots, called the Leech lattice. Witt's
classification was completed by Niemeier in 1967, who found the twenty
four 24-dimensional even unimodular lattices; these are called the
Niemeier lattices. His proof was simplified by Venkov (1980), who used
modular forms to restrict the possible root systems of such lattices.

When he found his lattice, Leech conjectured that it had covering
radius $\sqrt 2$ because there were several known holes of this
radius. Parker later noticed that the known holes of radius $\sqrt 2$
seemed to correspond to some of the Niemeier lattices, and inspired by
this Conway, Parker, and Sloane [C-P-S] found all the holes of this
radius. There turned out to be 23 classes of holes, which were
observed to correspond in a natural way with the 23 Niemeier lattices
other than the Leech lattice. Conway [C b] later used the fact that
the Leech lattice had covering radius $\sqrt 2$ to prove that the
26-dimensional even Lorentzian lattice $II_{25,1}$ has a Weyl vector,
and that its Dynkin diagram can be identified with the Leech lattice
$\Lambda$.

So far most of the proofs of these results involved rather long
calculations and many case-by-case discussions. The purpose of this
chapter  is to provide conceptual proofs of these results. We give new proofs
of
the existence and uniqueness of the Leech lattice and of the fact that it has
covering radius $\sqrt 2$. Finally we prove that the deep
holes of $\Lambda$ correspond to the Niemeier lattices and give a
uniform proof that the 23 ``holy constructions'' of $\Lambda$ found by
Conway and Sloane [C-S b] all work.

The main thing missing from this treatment of the Niemeier lattices is
a simple proof of the classification of the Niemeier lattices. (The
classification is not used in this chapter.) By Venkov's results it
would be sufficient to find a uniform proof that there exists a unique
Niemeier lattice for any root system of rank 24, all of whose
components have the same Coxeter number.

\proclaim 2.2.2 Some results from modular forms. 

We give some  results about Niemeier lattices, which are
proved using modular forms. From this viewpoint the reason why
24-dimensional lattices are special is that certain spaces of modular
forms vanish.

{\it Notation.} $\Lambda$ is any Niemeier lattice with no roots. A
Niemeier lattice is an even 24-dimensional unimodular lattice.

\proclaim Lemma 2.2.1. 
[C c]. Every element of $\Lambda$ is congruent $\bmod$
$2\Lambda$ to an element of norm at most 8.

We sketch Conway's proof of this. If $C_i$ is the number of elements
of $\Lambda$ of norm $i$ then elementary geometry shows that the
number of elements of $\Lambda/2\Lambda$ represented by vectors of
norm at most 8 is at least

$$C_0+C_4/2+C_6/2+C_8/(2\times\dim(\Lambda)).$$

The $C_i$'s can be worked out using modular forms and this sum
miraculously turns out to be equal to $2^{24}=$ order of
$\Lambda/2\Lambda$. Hence every element of $\Lambda/2\Lambda$ is
represented by an element of norm at most 8. Q.E.D.

(This is the only numerical calculation that we need.)

We now write $N$ for any Niemeier lattice. 

\proclaim Lemma 2.2.2.
(Venkov) If $y$ is any element of the Niemeier lattice $N$ then 
$$\sum(y,r)^2=y^2r^2n/24,$$
where the sum is over the $n$ elements $r$ of some fixed norm.

This is proved using modular forms; the critical fact is that the
space of cusp forms of dimension ${24\over 2}+2=14$ is
zero-dimensional. See [V] for details. 

\proclaim Lemma 2.2.3. 
(Venkov). The root system of $N$ has rank 0 or 24 and all
components of this root system have the same Coxeter number $h$.

This follows easily from lemma 2.2.2 as in [V]. We call $h$
the Coxeter number of the Niemeier lattice $N$. If $N$ has no roots we
put $h=0$.

It is easy to find the 24 root systems satisfying the condition of
lemma 2.2.3. Venkov gave a simplified proof of Niemeier's result that
there exists a unique Niemeier lattice for each root system. The next four
lemmas are not needed for the proof that $\Lambda$ has covering radius
$\sqrt 2$.

\proclaim Lemma 2.2.4.
(Venkov). $N$ has $24h$ roots and the norm of its Weyl vector is
$2h(h+1)$.

{\it Proof.} The root system of $N$ consists of components $a_n$,
$d_n$, and $e_n$, all with the same Coxeter number $h$. For each of
these components the number of roots is $hn$ and the norm of its Weyl
vector is ${1\over 12}nh(h+1)$. The lemma now follows because the rank
of $N$ is 0 or 24. Q.E.D.

\proclaim Lemma 2.2.5.
If $y$ is in $N$ then
$$\sum(y,r)^2=2hy^2,$$
where the sum is over all roots $r$ of $N$.

This follows from lemmas 2.2.2 and 2.2.4. Q.E.D.

\proclaim Lemma 2.2.6. 
If $\rho$ is the Weyl vector of $N$ then $\rho$ lies in $N$. 

{\it Proof.} We show that if $y$ is in $N$ then $(\rho,y)$ is an
integer, and this will prove that $\rho$ is in $N$ because $N$ is
unimodular. We have
$$\eqalign{ (2\rho,y)^2&= (\sum_r r,y)^2\qquad\hbox{(The sums are over
all positive roots $r$.)}\cr &=(\sum_r(r,y))^2\cr
&\equiv\sum_r(r,y)^2\bmod 2\cr &=y^2h\qquad\hbox{by lemma 2.2.5}\cr
&\equiv 0\bmod 2 \qquad \hbox{as $y^2$ is even.}\cr }$$ The term
$(2\rho,y)^2$ is an even integer and $(\rho,y)$ is rational, so
$(\rho,y)$ is an integer. Q.E.D.

\proclaim Lemma 2.2.7.
Suppose $h\ne 0$ and $y$ is in $N$. Then
$$(\rho/h-y)^2\ge 2(1+1/h)$$ and the $y$ for which equality holds form
a complete set of representatives for $N/R$, where $R$ is the
sublattice of $N$ generated by roots.

{\it Proof.} $\rho^2=2h(h+1)$, so 
$$\eqalign{
(\rho/h-y)^2-2(1+1/h)
&= ((\rho-hy)^2-\rho^2)/h^2\cr
&=(hy^2-2(\rho,y))/h\cr
&=(\sum(y,r)^2-\sum(y,r))/r\quad\hbox{by lemma 2.2.5}\cr
&= \sum((y,r)^2-(y,r))/h,\cr
}$$
where the sums are over all positive roots $r$. This sum is greater
than or equal to 0 because $(y,r)$ is integral, and is zero if and
only if $(y,r)$ is 0 or 1 for all positive roots $r$. In any lattice
$N$ whose roots generate a sublattice $R$, the vectors of $N$ that
have inner product 0 or 1 with all positive roots of $R$ (for some
choice of Weyl chamber of $R$) form a complete set of representatives
for $N/R$, and this proves the last part of the lemma. Q.E.D.

{\it Remark.} This lemma shows that if $N$ has roots then its covering
radius is greater than $\sqrt 2$. (In fact is is at least
$\sqrt{(5/2)}$.) In particular the covering radius of the Niemeier
lattice with root system $A_2^{12}$ is at least $\sqrt{(8/3)}$; this
Niemeier lattice has no deep hole that is half a lattice vector.

\proclaim 2.3.~Norm zero vectors in Lorentzian lattices. 

Here we describe the relation between norm 0 vectors in Lorentzian
lattices $L$ and extended Dynkin diagrams in the Dynkin diagram of
$L$.

{\it Notation.} Let $L$ be any even Lorentzian lattice, so it is
$II_{8n+1,1}$ for some $n$. Let $U$ be $II_{1,1}$; it has a basis of
two norm 0 vectors with inner product $-1$.

If $V$ is any $8n$-dimensional positive even unimodular lattice then
$V\oplus U\cong L$ because both sides are unimodular even Lorentzian
lattices of the same dimension, and conversely if $X$ is any
sublattice of $L$ isomorphic to $U$ then $X^\perp$ is an $8n$-dimensional
even unimodular lattice. If $z$ is any primitive norm 0 vector of $L$
then $z$ is contained in a sublattice of $L$ isomorphic to $U$ and all
such sublattices are conjugate under \Aut$(L)$. This gives 1:1
correspondences between the sets:

(1) $8n$-dimensional even unimodular lattices (up to isomorphism);

(2) orbits of sublattices of $L$ isomorphic to $U$ under \Aut$(L)$; and

(3) orbits of primitive norm 0 vectors of $L$ under \Aut$(L)$. 

If $V$ is the $8n$-dimensional unimodular lattice corresponding to the
norm 0 vector $z$ of $L$ then $z^\perp\cong V\oplus O$ where $O$ is the
one-dimensional singular lattice. The Dynkin diagram of $V\oplus O$ is
the Dynkin diagram of $V$ with all components changed to the
corresponding extended Dynkin diagram.

Now choose coordinates $(v,m,n)$ for $L\cong V\oplus U$, where $v$ is
in $V$, $m$ and $n$ are integers, and $(v,m,n)^2=v^2-2mn$ (so $V$ is
the set of vectors $(v,0,0)$ and $U$ is the set of vectors
$(0,m,n)$). We write $z$ for $(0,0,1)$ so that $z$ is a norm 0 vector
corresponding to the unimodular lattice $V$. We choose a set of simple
roots for $z^\perp\cong V\oplus O$ to be the vectors $(r_j^i,0,0)$ and
$(r_j^0,0,1)$ where the $r_j^i$'s are the simple roots of the
components $R_j$ of the root system of $V$ and $r_j^0$ is the highest
root of $R_j$. We let $\bar R_j$ be the vectors $(r_j^i,0,0)$ and
$(r_j^0,0,1)$ so that $\bar R_j$ is the extended Dynkin diagram of
$R_j$. Then
$$\sum m_ir_j^i=z\qquad\hbox{and}\qquad\sum m_i=h_j,$$ where the
$m_i$'s are the weights of the vertices of $\bar R_j$ and $h_j$ is the
Coxeter number of $R_j$. We can apply Vinberg's algorithm to find a
fundamental domain of $L$ using $z$ as a controlling vector and this
shows that there is a unique fundamental domain $D$ of $L$ containing
$z$ such that all the vectors of the $\bar R_j$s are simple roots of
$D$.

Conversely suppose that we choose a fundamental domain $D$ of $L$ and
let $\bar R$ be a connected extended Dynkin diagram contained in the
Dynkin diagram of $D$. If we put $z=\sum m_ir_i$, where the $m_i$'s are
the weights of the simple roots $r_i$ of $\bar R$, then $z$ has norm 0
and inner product 0 with all the $r_i$s (because $\bar R$ is an
extended Dynkin diagram) and has inner product at most 0 with all
simple roots of $D$ not in $\bar R$ (because all the $r_i$ do), so $z$
is in $D$ and must be primitive because if $z'$ was a primitive norm 0
vector dividing $z$ we could apply the last paragraph to $z'$ to find
$z'=\sum m_ir_i=z$. The roots of $\bar R$ together with the simple
roots of $D$ not joined to $\bar R$ are the simple roots of $D$
orthogonal to $z$ and so are a union of extended Dynkin
diagrams. This shows that the following 3 sets are in natural 1:1
correspondence:

(1) equivalence classes of extended Dynkin diagrams in the Dynkin
diagram of $D$, where two extended Dynkin diagrams are equivalent if
they are equal or not joined;

(2) maximal disjoint sets of extended Dynkin diagrams in the Dynkin
diagram of $D$, such that no two elements of the set are joined to each
other;

(3) primitive norm 0 vectors of $D$ that have at least one root
orthogonal to them. 

\proclaim 2.4.~Existence of the Leech lattice. 

In this section we prove the existence of a Niemeier lattice with no
roots (which is of course the Leech lattice). This is done by showing
that given any Niemeier lattice we can construct another Niemeier
lattice with at most half as many roots. It is a rather silly proof
because it says nothing about the Leech lattice apart from the fact
that it exists.

{\it Notation.} Fix a Niemeier lattice $N$ with a Weyl vector $\rho$
and Coxeter number $h$, and take coordinates $(y,m,n)$ for
$II_{25,1}=N\oplus U$ with $y$ in $N$, $m$ and $n$ integers.

By lemma 2.2.6 the vector $w=(\rho,h,h+1)$ is in $II_{25,1}$ and by
lemma 2.2.4 it has norm 0. We will show that the Niemeier lattice
corresponding to $w$ has at most half as many roots as $N$.

\proclaim Lemma 2.4.1.
There are no roots of $II_{25,1}$ that are orthogonal to $w$ and
have inner product $0$ or $\pm 1$ with $z=(0,0,1)$.

{\it Proof.} Suppose that $r=(y,0,n)$ is a root that has inner product
0 with $w$ and $z$. Then $y$ is a root of $N$, and $(y,\rho)=nh$
because $(y,\rho)-nh=((y,0,n),(\rho,h,h+1))=(r,w)=0$. But for any root
$y$ of $N$ we have $1\le |(y,\rho)|\le h-1$, so $(y,\rho)$ cannot be a
multiple of $h$.

If $(y,1,n)$ is any root of $II_{25,1}$ that has inner product $-1$
with $z$ and $0$ with $w$ then
$$(y,\rho)-(h+1)-nh=((\rho,h,h+1),(y,1,n))=(w,r)=0$$
and
$$y^2-2n=(y,1,n)^2=r^2=2$$
so
$$\eqalign{
(y-\rho/h)^2&=y^2-2(y,\rho)/h+\rho^2/h^2\cr
&=2+2n-2(nh+h+1)/h+2h(h+1)/h^2\cr
&=2
}$$
which contradicts lemma 2.2.7. Hence no root in $w^\perp$ can have inner
product 0 or $-1$ with $z$. Q.E.D.

{\it Remark.} In fact there are no roots in $w^\perp$.

\proclaim Lemma 2.4.2. 
The Coxeter number $h'$ of the Niemeier lattice of the vector $w$ is at
most ${1\over 2}h$.

{\it Proof.} We can assume $h'\ne 0$. Let $R$ be any component of the
Dynkin diagram of $w^\perp$ (so $R$ is an extended Dynkin
diagram). The sum $\sum m_ir_i$ is equal to $w$ by section 2.3, where
the $r_i$s are the roots of $R$ with weights $m_i$. Also $\sum
m_i=h'$, $(r_i,z)\le -2$ by lemma 2.4.1, and
$$(\sum m_ir_i,z)=((\rho,h,h+1),z)=-h,$$
so $h'\le {1\over 2} h$. Q.E.D.

\proclaim Theorem 2.4.3. There exists a Niemeier lattice with no roots. 

{\it Proof.} By lemma 2.4.2 we can find a Niemeier lattice with Coxeter
number at most ${1\over 2}h$ whenever we are given a Niemeier lattice
of Coxeter number $h$. By repeating this we eventually get a Niemeier
lattice with Coxeter number 0, which must have no roots. Q.E.D.

Leech was the first to construct a Niemeier lattice with no roots
(called the Leech lattice), and Niemeier later proved that it was
unique as part of his enumeration of the Niemeier lattices. In section 2.6 we
will give another proof that there is only one such lattice.
In section 2.7 we will  show that the Niemeier lattice of the vector
$(\rho,h,h+1)$ never has any roots and so is already the Leech
lattice.

\proclaim 2.5.~The covering radius of the Leech lattice. 

Here we prove that any Niemeier lattice $\Lambda$ with no roots has
covering radius $\sqrt 2$.

{\it Notation.} We write $\Lambda$ for any Niemeier lattice with no
roots, and put $II_{25,1}=\Lambda\oplus U$ with coordinates
$(\lambda,m,n)$ with $\lambda$ in $\Lambda$, $m$ and $n$ integers and
$(\lambda,m,n)^2=\lambda^2-2mn$. We let $w$ be the norm 0 vector
$(0,0,1)$ and let $D$ be a fundamental domain of the reflection group
of $II_{25,1}$ containing $w$. If we apply Vinberg's algorithm, using
$w$ as a controlling vector, then the first batch of simple roots to
be accepted is the set of roots $(\lambda, 1, {1\over 2} \lambda^2-1)$
for all $\lambda$ in $\Lambda$. (Conway proved that no other roots are
accepted. See section 2.6.) We identify the vectors of $\Lambda$ with
these roots of $II_{25,1} $, so $\Lambda$ becomes a Dynkin diagram
with two points of $\Lambda$ joined by a bond of strength 0,1,2,... if
the norm of their difference is 4,6,8,.... We can then talk about
Dynkin diagrams in $\Lambda$; for example an $a_2$ in $\Lambda$ is two
points of $\Lambda$ whose distance apart is $\sqrt 6$.

Here is the main step in the proof that $\Lambda$ has covering radius
$\sqrt 2$.

\proclaim Lemma 2.5.1. 
If $X$ is any connected extended Dynkin diagram in $\Lambda$ then $X$
together with the points of $\Lambda$ not joined to $X$ contains a
union of extended Dynkin diagrams of total rank 24. (The rank of a
connected extended Dynkin diagram is one less than the number of its
points.)

{\it Proof.} $X$ is an extended Dynkin diagram in the set of simple
roots of $II_{25,1}$ of height 1, and hence determines a norm 0 vector
$z=\sum m_ix_i$ in $D$, where the $x_i$s are the simple roots of $X$
of weights $m_i$. The term $z$ corresponds to some Niemeier lattice
with roots, so by lemma 2.2.3 the simple roots of $z^\perp$ form a union
of extended Dynkin diagrams of total rank 24, all of whose components
have the same Coxeter number $h=\sum m_i$. Also
$$h=-(z,w)$$ because $z=\sum m_ix_i$ and $(x_i,w)=-1$. By applying
Vinberg's algorithm with $z$ as a controlling vector we see that the
simple roots of $D$ orthogonal to $z$ form a set of simple roots of
$z^\perp$, so the lemma will be proved if we show that all simple
roots of $D$ in $z^\perp$ have height 1, because then all simple roots
of $z^\perp$ lie in $\Lambda$.

Suppose that $Z$ is any component of the Dynkin diagram of $z^\perp$
with vertices $z_i$ and weights $n_i$ so that $\sum n_iz_i=z$. Then
$$-h=(z,w)=(\sum n_iz_i,w)=\sum n_i(z_i,w).$$ Now $(z_i,w)\le -1$ as
there are no simple roots $z_i$ with $(z_i,w)=0$, and $\sum n_i=h$ as
$Z$ has the same Coxeter number $h$ as $X$, so we must have
$(z_i,w)=-1$ for all $i$. Hence all simple roots of $D$ in $z^\perp$
have height 1. Q.E.D.

It follows easily from this that $\Lambda$ has covering radius $\sqrt
2$. For completeness we sketch the remaining steps of the proof of
this. (This part of the proof is taken from [C-P-S].)

\item{Step 1.}  The distance between two vertices of any hole is at
most $\sqrt 8$. This follows easily from lemma 2.2.1. 
\item{Step 2.} By a property of extended Dynkin diagrams one of the following must hold for the vertices of any hole.
\itemitem{(i)} The vertices form a spherical Dynkin diagram. In this case the
hole has radius $(2-1/\rho^2)^{1\over 2}$, where $\rho^2$ is the norm
of the Weyl vector of this Dynkin diagram.
\itemitem{(ii)} The distance between two points is greater than $\sqrt 8$. This
is impossible by step 1.
\itemitem{(iii)} The vertices contain an extended Dynkin diagram. 
\item{Step 3.} Let $V$ be the set of vertices of any hole of radius
greater than $\sqrt 2$. By step 2, $V$ contains an extended Dynkin
diagram. By lemma 2.5.1 this extended Dynkin diagram is one of a set of
disjoint extended Dynkin diagrams of $\Lambda$ of total rank 24. These
Dynkin diagrams then form the vertices of a hole $V'$ of radius $\sqrt
2$ whose center is the center of any of the components of its
vertices. Any hole whose vertices contain a component of the vertices
of $V'$ must be equal to $V'$, and in particular $V=V'$ has radius
$\sqrt 2$. Hence $\Lambda$ has covering radius $\sqrt 2$. Q.E.D.

\proclaim 2.6.~Uniqueness of the Leech lattice. 

We continue with the notation of the previous section. We show that
$\Lambda$ is unique and give Conway's calculation of $\Aut(II_{25,1})$.

\proclaim Theorem 2.6.1. 
The simple roots of the fundamental domain $D$ of $II_{25,1}$ are just
the simple roots $(\lambda,1,{1\over 2}\lambda^2-1)$ of height one.

{\it Proof} [C b]. If $r'=(v,m,n)$ is any simple root of
height greater than 1 then it has inner product at most 0 with all
roots of height 1. $\Lambda$ has covering radius $\sqrt 2$ so there is
a vector $\lambda$ of $\Lambda$ with $(\lambda-v/m)^2\le 2$. But then
an easy calculation shows that $(r',r)\ge 1/2m$, where $r$ is the
simple root $(\lambda,1,{1\over 2}\lambda^2-1)$ of height 1, and this
is impossible because $1/2m$ is positive. Q.E.D. 

\proclaim Corollary 2.6.2. The Leech lattice is unique, in other words
any two Niemeier
lattices with no roots are isomorphic.

(Niemeier's enumeration of the Niemeier lattices gave the first proof
of this fact. Also see [C c] and [V].)

{\it Proof.} Such lattices correspond to orbits of primitive norm 0
vectors of $II_{25,1}$ that have no roots orthogonal to them, so
it is sufficient to show that any two such vectors are conjugate under
$\Aut(II_{25,1})$, and this will follow if we show that $D$ contains
only one such vector. But if $w$ is such a vector in $D$ then theorem
2.6.1 shows that $w$ has inner product $-1$ with all simple roots of
$D$, so $w$ is unique because these roots generate $II_{25,1}$. Q.E.D.

\proclaim Corollary 2.6.3. 
$\Aut(II_{25,1})$ is a split extension $R.(\cdot\infty)$, where $R$ is
the subgroup generated by reflections and $\pm 1$, and $\cdot \infty$
is the group of automorphisms of the affine Leech lattice.

{\it Proof} [C b]. Theorem 2.6.1 shows that $R$ is simply
transitive on the primitive norm 0 vectors of $II_{25,1}$ that are not
orthogonal to any roots, and the subgroup of $\Aut(II_{25,1})$
fixing one of these vectors is isomorphic to $\cdot\infty$. Q.E.D.

\proclaim 2.7.~The deep holes of the Leech lattice. 

Conway, Parker, and Sloane [C-P-S] found the 23 orbits of ``deep holes'' in
$\Lambda$ and observed that they corresponded to the 23 Niemeier
lattices other than $\Lambda$. Conway  Sloane [C-S b,c] later gave
a ``holy construction'' of the Leech lattice for each deep hole, and
asked for a uniform proof that their construction worked. In this
section we give uniform proofs of these two facts.

We continue with the notation of the previous section, so
$II_{25,1}=\Lambda\oplus U$. We write $Z$ for the set of nonzero
isotropic subspaces of $II_{25,1}$, which can be identified with the
set of primitive norm 0 vectors in the positive cone of
$II_{25,1}$. $Z$ can also be thought of as the rational points at
infinity of the hyperbolic space of $II_{25,1}$. We introduce
coordinates for $Z$ by identifying it with the space $\Lambda\otimes
Q\cup \infty$ as follows. Let $z$ be a norm 0 vector of $II_{25,1}$
representing some point of $Z$. If $z=w=(0,0,1)$ we identify it with
$\infty$ in $\Lambda\otimes \Q\cup \infty$. If $z$ is not a multiple of
$w$ then $z=(\lambda,m,n)$ with $m\ne 0$ and we identify $z$ with
$\lambda/m$ in $\Lambda\otimes \Q\cup \infty$. It is easy to check that
this gives a bijection between $Z$ and $\Lambda\otimes \Q\cup \infty$.

The group $\cdot\infty$ acts on $Z$. This action can be described
either as the usual action of $\cdot\infty$ on $\Lambda\otimes \Q\cup
\infty$ (with $\cdot\infty$ fixing $\infty$) or as the action of
$\cdot\infty=\Aut(D)$ on the isotropic subspaces of $II_{25,1}$. We now
describe the action of the whole of $\Aut(II_{25,1})$ on $Z$.

\proclaim Lemma 2.7.1.
Reflection in the simple root $r=(\lambda,1,{1\over 2} \lambda^2-1)$
of $D$ acts on $\Lambda\otimes \Q\cup \infty$ as inversion in a sphere
of radius $\sqrt 2$ around $\lambda$ (exchanging $\lambda$ and
$\infty$).

{\it Proof.} Let $v$ be a point of $\Lambda\otimes \Q$ corresponding to
the isotropic subspace of $II_{25,1}$ generated by $z=(v,1,{1\over
2}v^2-1)$. The reflection in $r$ maps $z$ to $z-2(z,r)r/r^2$,
$$(z,r)={1\over 2} (z^2+r^2-(z-r)^2)=1-{1\over 2}(v-\lambda)^2,$$
so reflection in $r^\perp$ maps $z$ to 
$$\eqalign{ z-(1-{1\over 2}(v-\lambda)^2)r 
&=(v,1,{1\over 2}v^2-1)-([1-{1\over 2}(v-\lambda)^2]\lambda,1-{1\over
2}(v-\lambda)^2,?) \cr
&\propto (\lambda+2(v-\lambda)/(v-\lambda)^2,1,?),\cr
}$$
which is a norm 0 vector corresponding to the point
$\lambda+2(v-\lambda)/(v-\lambda)^2$ of $\Lambda\otimes \Q\cup \infty$,
and this is the inversion of $v$ in the sphere of radius $\sqrt 2$
about $\lambda$. Q.E.D. 

\proclaim Corollary 2.7.2.
The Niemeier lattices with roots are in natural bijection with the
orbits of deep holes of $\Lambda$ under $\cdot\infty$. The vertices of
a deep hole form the extended Dynkin diagram of the corresponding
Niemeier lattice.

(This was first proved by Conway, Parker, and Sloane [C-P-S], who explicitly
calculated all the deep holes of $\Lambda$ and observed that they
corresponded to the Niemeier lattices.)

{\it Proof.} By lemma 2.7.1 the point $v$ of $\Lambda\otimes \Q$
corresponds to a norm 0 vector in $D$ if and only if it has distance
at least $\sqrt 2$ from every point of $\Lambda$, i.e. if and only if
$v$ is the center of some deep hole of $\Lambda$. In this case the
vertices of the deep hole are the points of $\Lambda$ at distance
$\sqrt 2$ from $v$ and these correspond to the simple roots $(\lambda,
1, {1\over 2} \lambda^2-1)$ of $D$ in $z^\perp$.

Niemeier lattices $N$ with roots correspond to the orbits of primitive
norm 0 vectors in $D$ other than $w$ and hence to deep holes of
$II_{25,1}$. The simple roots in $z^\perp$ form the Dynkin diagram of
$z^\perp\cong N\oplus O$, which is the extended Dynkin diagram of
$N$. Q.E.D.

Conway and Sloane [C-S b] gave a construction for the Leech lattice
from each Niemeier lattice. They remarked ``The fact that this
construction always gives the Leech lattice still quite astonishes us,
and we have only been able to give a case by case verification as
follows. $\ldots$ We would like to see a more uniform proof''. Here is
such a proof. Let $N$ be a Niemeier lattice with a given set of simple
roots whose Weyl vector is $\rho$. We define the vectors $f_i$ to be
the simple roots of $N$ together with the highest roots of $N$ (so
that the $f_i$s form the extended Dynkin diagram of $N$), and define
the glue vectors $g_i$ to be the vectors $v_i-\rho/h$, where $v_i$ is
any vector of $N$ such that $g_i$ has norm $2(1+1/h)$. By lemma 2.2.7
the $v_i$s are the vectors of $N$ closest to $\rho/h$ and form a
complete set of coset representatives for $N/R$, where $R$ is the
sublattice of $N$ generated by roots. In particular the number of
$g_i$s is $\sqrt {\det(R)}$.

The ``holy construction'' is as follows. 

(i) The vectors $\sum m_if_i+\sum n_ig_i$ with $\sum n_i=0$ form the
lattice $N$.

(ii) The vectors $\sum m_if_i+\sum n_ig_i$ with $\sum m_i+\sum n_i=0$
form a copy of $\Lambda$.

Part (i) follows because the vectors $f_i$ generate $R$ and the
vectors $v_i$ form a complete set of coset representatives for
$N/R$. We now prove (ii).

We will say that two sets are isometric if they are isomorphic as
metric spaces after identifying pairs of points whose distance apart
is 0. For example $N$ is isometric to $N\oplus O$. Let $z$ be a
primitive norm 0 vector in $D$ corresponding to the Niemeier lattice
$N$. The sets of vectors $f_i$ and $g_i$ are isometric to the sets of
simple roots $f_i'$ and $g_i'$ of $D$, which have inner product 0 or
$-1$ with $z$. (This follows by applying Vinberg's algorithm with $z$
as a controlling vector.)

\proclaim Lemma 2.7.3.
The vectors $f_i'$ and $g_i'$ generate $II_{25,1}$. 

{\it Proof.} The vectors $\sum m_if_i'+n_ig_i'$ with $\sum n_i=0$ are
just the vectors in the lattice generated by $f_i'$ and $g_i'$, which
are in $z^\perp$, and they are isometric to $N$ because the vectors
$\sum m_i f_i+\sum n_ig_i$ with $\sum n_i=0$ form a copy of $N$. The
vector $z$ is in the lattice generated by the $f_i'$s, so the whole of
$z^\perp$ is contained in the lattice generated by the vectors $f_i$
and $g_i$. The vectors $g_i'$ all have inner product 1 with $z$, so
the vectors $f_i'$ and $g_i'$ generate $II_{25,1}$. Q.E.D.

The lattice of vectors $\sum m_if_i'+\sum n_ig_i'$ with $\sum m_i+\sum
n_i=0$ is $w^\perp$ by lemma 2.2.7 because
$(f_i',w)=(g_i',w)=1$. $w^\perp$ is $\Lambda\oplus O$, which is
isometric to $\Lambda$ and therefore isomorphic to $\Lambda$
because it is contained in the positive definite space $N\otimes \Q$.
This proves (ii). 

The holy construction for $\Lambda$ is equivalent to the
$(\rho,h,h+1)^\perp$ construction of section 2.4, so that $(\rho,h,h+1)$
is always a norm 0 vector corresponding to $\Lambda$.

\proclaim Chapter 3 Negative norm vectors in Lorentzian lattices. 

To classify lattices we use the following idea. If $A$ is an
$n$-dimensional positive definite unimodular lattice then $A\oplus
(-I)$ is isomorphic to $I_{n,1}$ and the lattice $-I$ determines a
pair of norm $-1$ vectors in $I_{n,1}$. Conversely the orthogonal
complement of a norm $-1$ vector in $I_{n,1}$ is an $n$-dimensional
unimodular lattice, and in this way we get a bijection between such
lattices and orbits of norm $-1$ vectors in $I_{n,1}$.

In practice we use the even lattices $II_{8n+1,1}$, and in particular
$II_{25,1}$, instead of $I_{n,1}$ because their automorphism groups
behave better. Unimodular lattices correspond to norm $-4$ vectors in
these lattices which is the initial motivation for trying to find
negative norm vectors in them.

Sections 3.2 to 3.5 of this chapter contain an algorithm for finding
some of the norm $-2m$ vectors of $II_{8m+1,1}$ from those of norm
$-(2m-2)$, which for $8n\le 24$ finds all the norm $-2m$ vectors.
Hence to find the norm $-4$ vectors we start with the norm 0 vectors
(which are known because they correspond to the Niemeier lattices) and
then find the norm $-2$ vectors and then the norm $-4$ vectors. The
norm $-2$ and $-4$ vectors of $II_{25,1}$ are given in tables $-2$ and
$-4$; there are 121 orbits of norm $-2$ vectors and 665 orbits of norm
$-4$ vectors. In chapters 4 and 5 we will use some of the vectors of
norms $-6$, $-8$, $-10$, and $-14$. The number of orbits of norm $-2m$ vectors
of $II_{8n+1,1}$ increases fairly rapidly with $m$, and extremely
rapidly with $n$ for $n\ge 4$.

Sections 3.6 to 3.8 contain more details about the algorithm and show
how it simplifies for vectors of ``type $\ge 3$''. Finally in 3.9 we
classify the positive norm vectors of $II_{8n+1,1}$; this turns out to
be much easier than classifying the negative norm ones because there is
a unique orbit of primitive vectors of any given even positive norm
(for $n\ge 1$).

\proclaim 3.1 Negative norm vectors and lattices. 

We now investigate the lattices that occur as orthogonal complements
of negative norm vectors in $II_{8n+1,1}$.  Norm $-4$ vectors of
$II_{8n+1,1}$ give unimodular lattices, and norm $-2$ vectors given
even bimodular lattices.

\proclaim Lemma 3.1.1. There is a bijection between 
\item{(1)} Orbits of primitive vectors $u$ of norm $-2m<0$ in 
$II_{8n+1,1}$ under $\Aut(II_{8n+1,1})$.
\item{(2)} 
Pairs $(B,b)$ where $B$ is an $8n+1$-dimensional even positive
definite lattice with determinant $2m$ and $b$ is a generator of
$B'/B$ of norm $1/2m\bmod 2$.
\item{} The group $\Aut(B,b)$ of automorphisms of $B$ fixing $b$ in $B'/B$ is isomorphic to the group $\Aut(II_{8n+1,1},u)$ of automorphisms of
$II_{8n+1,1}$ fixing $u$. 

Proof. To get from (1) to (2) we take $B$ to be $u^\perp$. We get $b$
by taking any vector $x$ of $II_{8n+1,1}$ with $(x,u)=-1$ and projecting
it into $B\otimes \Q$. This gives a well defined element $b$ of $B'/B$
which is independent of the choice of $x$.

Conversely if we are given $(B,b)$ we form the lattice $B\oplus U$
where $U$ is a one dimensional lattice generated by $u$ of norm
$-2m$. We add $b\oplus u/2m$ to this lattice and this gives a copy of
$II_{8n+1,1}$. ($b\oplus u/2m$ is well defined mod $B\oplus U$ and has
even norm.)

The two maps given above are inverses of each other so this gives a
bijection between (1) and (2). The statement about automorphism groups
follows easily. Q.E.D. 

\proclaim Lemma 3.1.2. There is a bijection between 
\item{(1)} Orbits of norm $-2$ vectors $u$ of $II_{8n+1,1}$, and
\item{(2)} Even $8n+1$-dimensional positive definite bimodular lattices $B$. 
\item{} We have $\Aut(B)\cong \Aut(II_{8n+1,1},u)$. 

Proof. This follows from 3.1.1 provided we show that the nonzero
element $b$ of $B'/B$ has norm $1/2$ (rather than $3/2$) $\bmod 2$. If
it had norm $3/2\bmod 2$ we could form the lattice $B\oplus a_1$ where
$a_1$ is a one dimensional lattice generated by $d$ of norm 2. Adding
the vector $b\oplus d/2$ to this would give an $8n+2$-dimensional even
positive definite unimodular lattice, which is impossible. (Note that
$b\oplus d/2$ is well defined $\bmod B\oplus D$ and has even norm.)
Q.E.D.

\proclaim Lemma 3.1.3. There is a bijection between
\item{(1)}
Orbits of norm $-4$ vectors $u$ of $II_{8n+1,1}$, and
\item{(2)} $8n+1$-dimensional positive definite unimodular lattices $A$.
\item{}$\Aut(A)\cong (\Z/2\Z)\times \Aut(II_{8n+1,1},u)$ where
$\Z/2\Z$ is generated by $-1$. 

Proof. Let $c$ be a characteristic vector of $A$, so $c$ has norm
$1\bmod 8$.  If $B$ is the sublattice of even elements of $A$ then the
element $b=c/2$ has norm $1/4\bmod 2$ and generates the group $B'/B$
which is cyclic of order 4.  The result follows from this and
3.1.1. To prove the result about $\Aut(A)$, note that the automorphism
$-1$ of $A$ exchanges the two generators of $B'/B$, so $\Aut(A)\cong
(\Z/2\Z)\times \Aut(B,b)$ where $\Z/2\Z$ is generated by $-1$. Now use
3.1.1 again. Q.E.D.

{\bf Remark.} In 4.9 we will show that given any orbit of norm
$2m$ vectors in $e_8$ we can find an interpretation of the norm $-2m$
vectors of $II_{8n+1,1}$. $e_8$ has unique orbits of vectors of norms
2 and 4 which give 3.1.2 and 3.1.3.

\proclaim 3.2 Summary of the classification of negative norm vectors. 

The next few section give an algorithm for finding the negative norm
vectors in $II_{25,1}$. Here we summarize this method and indicate how
far it can be generalized to work for other Lorentzian lattices.

First we take a fundamental domain $D$ of the reflection group $R$ of
$II_{25,1}$. Every negative norm vector of $II_{25,1}$ is conjugate
under $R$ to a unique vector in $D$, so it is enough to classify
orbits of vectors in $D$ under $\Aut(D)$ (the subgroup of
$\Aut(II_{25,1})$ fixing $D$, which is isomorphic to $\cdot\infty$).

The properties of $II_{25,1}$ that we use are the following (which
characterize it):
\item{(1)}
$II_{25,1}$ is even.
\item{(2)}
$II_{25,1}$ is unimodular.
\item{(3)}
$II_{25,1}$ contains a norm 0 vector $w$ which has inner product 
$-(r,r)/2$ with all simple roots $r$ of $D$. 
\item{(4)}
The primitive norm 0 vectors of $II_{25,1}$ are known. 

Of these properties (1) and (2) are not essential. The main effect of
not having them would be that we would have to work with arbitrary
root systems instead of those that have all their roots of norm
2. This would not create any major new difficulties but would make
most results messier to state. Property (3) is not used in large parts
of the machinery, but some version of it is needed for getting hold of
the vectors of $D$ that are not orthogonal to any roots. For most
purposes it would be sufficient to assume that $D$ contains a nonzero
vector $w$ (preferably of norm 0) that had inner product 0 or
$-(r,r)/2$ with all simple roots $r$ of $D$. The lattices $II_{9,1}$
and $II_{17,1}$ have such a vector $w$ (in fact $II_{17,1}$ has two
orbits of such vectors) and we can use the machinery of this chapter
with little change to classify the negative norm vectors in these
lattices. (However, classifying vectors in these lattices is fairly
easy anyway because the fundamental domains for their reflection
groups have finite volumes and automorphism groups of order 1 or 2.)

We now define several important invariants of the vector $u$ of $D$.
\item{}
The {\bf norm} of $u$ is $(u,u)=-2n\le 0$.
\item{}
The {\bf type} of $u$ is the smallest possible inner product of $u$
with a norm 0 vector of $II_{25,1}$. $u$ has type 0 if and only if it
has norm 0 because $II_{25,1}$ is Lorentzian and $u^2\le 0$.
\item{}
The {\bf height} of $u$ is $-(u,\rho)$; it is positive if $u$ is not a
multiple of $w$.
\item{}
The {\bf root system} of $u$ is the root system of $u^\perp$; its
Dynkin diagram is a union of $a$'s, $d$'s, and $e$'s if $u^2<0$ and is
the extended Dynkin diagram of a Niemeier lattice if $u^2=0$.

The algorithm of the next few sections for classifying negative norm
vectors runs roughly as follows:

Given $u$ of norm $-2n$ in $D$ we try to reduce it to a vector of norm
$-(2n-2)$ in $D$ by adding a root $r$ of $u^\perp$ to $u$. One of three
things can happen:
\item{(1)}
There are no roots in $u^\perp$. In this case the vector $v=u-w$ is
still in $D$ because $u$ has inner product at most $-1$ with all
simple roots of $D$ and $w$ has inner product (0 or) $-1$ with these
roots. If $u$ is not a multiple of $w$ then $v$ has smaller norm than
$u$.
\item{(2)}
There are roots $r$ of $u^\perp$, but $r+u$ is never in $D$. In this
case we can describe $u$ precisely: $u$ has norm 0, or there is a norm
0 vector $z$ with $(u,z)=-1$ (in other words $u$ has type 1) and the
orbit of $u$ depends only on $z$ and the norm of $u$.
\item{(3)}
There is a root $r$ of $u^\perp$ such that $v=u+r$ is in $D$. In this
case we have two problems to solve:
\itemitem{(a)}
Given $v$, what vectors $u$ can we obtain from it like this?
\itemitem{(b)}
When are two $u$'s, obtained from different $v$'s, conjugate under
$\Aut(D)$?

For $u$ in $D$ we let $R_i(u)$ be the simple roots of $D$ that  have
inner product $-i$ with $u$, so $R_i(u)$ is empty for $i<0$ and
$R_0(u)$ is the Dynkin diagram of $u^\perp$. We write $S(u)$ for
$R_0(u)\cup R_1(u)\cup R_2(u)$. Then given $S(v)$ we can find all
vectors $u$ of $D$ that come from $v$ as in (3) above, and $S(u)$ is
contained in $S(v)$. By keeping track of the action of $\Aut(D,u)$ on
$S(u)$ for vectors $u$ of $D$ we can solve problems (a) and (b).

In 3.6 to 3.8 we will look at properties of vectors that depend on
their type. A vector of type $m$ can be ``reduced'' to a vector of
$-$norm at most $2m$, and vectors of type $m$ are closely related to
elements of $B/mB$ as $B$ runs through the Niemeier lattices. More
specifically: \item{} Type 0 vectors are just norm 0 vectors.  \item{}
Type 1 vectors are closely related to type 0 vectors.  \item{} Type 2
vectors can be reduced to vectors of norms $-2$ and $-4$. Knowing the
type 2 vectors of $D$ is essentially the same as knowing $B/2B$ for
all Niemeier lattices $B$, or knowing all 24 dimensional unimodular
lattices and all 25 dimensional even bimodular lattices. For any type
2 vector $u$ we construct a special automorphism $\sigma_u$ of $D$
fixing $u$ and use these automorphisms to prove several identities for
the lattices $u^\perp$. (For example, twice the norm of the Weyl
vector of a 25 dimensional even bimodular lattice is a square.) 

We also find some properties of vectors whose type is at least 3; in
particular the algorithm simplifies a bit for these vectors.

For each vector $u$ of $D$ we define the involution $\tau_u$ by
$\tau_u(v)=u+\sigma(v)$, where $\sigma$ is the opposition involution
of the root system of $u^\perp$. $\tau_u$ is not in $\Aut(L)$ because
it does not fix 0.

\proclaim Lemma 3.2.1. 
\item{(a)} $(u,v)+(u,\tau_u(v))=u^2$ and $\tau_u^2=1$. 
\item{(b)} $v$ has negative inner product with all simple roots of $u^\perp$
if and only if $\tau_u(v)$ does.

Proof. (a) follows easily from the fact that $\sigma(u)=-u$. (b)
follows from the fact that $\sigma$ fixes the Weyl chamber of
$u^\perp$ while $u$ has inner product 0 with all simple roots of
$u^\perp$. Q.E.D.

$\tau_u$ is a sort of duality map on vectors of $II_{25,1}$. We will
use the fact that it identifies the following pairs of sets:
\item{(1)}
Vectors $z$ with $z^2=0$ and $(z,u)=-1$; and vectors $v$ with
$-v^2=2n-2$, $-(v,u)=2n-1$.
\item{(2)}
Roots in $r^\perp$; and vectors $v$ with $-v^2=2(n-1)$ and $-(v,u)=2n$. 
\item{(3)}
Roots that have inner product $-1$ with $u$; and vectors $v$ with
$-v^2=2n-4$ and $-(v,u)=2n-1$.

$\tau_u(v)$ can be thought of as a sort of twisted version of $u-v$;
the twisting is so that (b) in 3.2.1 holds.

\proclaim 3.3 How to find things in $D$. 

Given a vector $u$ in the fundamental domain $D$ of $II_{25,1}$ we
often construct other vectors $v$ from it and we would like to be able
to prove that these vectors $v$ are also in $D$, preferably without
checking that their inner products with all the (infinite number of)
simple roots of $D$ are $\le 0$. This section given conditions for
$v$ to be in $D$ which only depend on checking the inner products of
$v$ with a finite number of simple roots. 

Let $H$ be hyperbolic space acted on by a discrete reflection group
$G$, and choose a fundamental domain $D$ for $G$ and a point $d$ in
$D$. We take an orbit $S$ of points of $H$ and find the set $S_d$ of
points of $S$ closest to $D$. (As $G$ is discrete there are only a
finite number of points of $S$ within any given distance from $d$, so
there exist points of $S$ closest to $d$.) If $s$ is in $S$ and there
is a reflection hyperplane $P$ of $G$ passing strictly between $s$ and
$d$ then the reflection of $s$ in $P$ is closer to $d$ than $s$ is, so
the points of $S_d$ are in the conjugates of $D$ containing $d$. If
$G_d$ is the subgroup of $G$ fixing $d$ then it acts transitively on
the conjugates of $D$ containing $d$, so each such conjugate contains
exactly one point of $S_d$ and the points of $S_d$ form a single orbit
under $G_d$.

We translate this into the language of Lorentzian lattices. $d$ is a
vector in the fundamental domain $D$ of the group $G$ generated by
reflections of the Lorentzian lattice $L$, so $d$ has norm $\le 0$.
We let $S$ be an orbit of vectors of norm $\le 0$ whose inner product
with $d$ is $\le 0$.  We regard all these points as points of the
hyperbolic space $H$ of $L$. The distance of $s$ and $d$ regarded as
points of $H$ increases with $-(s,d)$ and the group $G_d$ is generated
by the reflections of $L$ fixing $d$. The results of the paragraph
above show that the vectors of $S$ whose inner product with $d$ is
minimal (in absolute value) form a single orbit under $G_d$, exactly
one element of which is in $D$. We rephrase this as follows:

\proclaim Theorem 3.3.3. 
If $L$, $G$, $D$, and $d$ are as above and $S_d$ is the set of vectors
of some fixed norm $\le 0$ in the same cone as $d$ that have minimal
$|$inner product$|$ with $d$ then every orbit of $S_d$ under $G_d$
contains a unique point in $D$. A point of $S_d$ is in $D$ if and only
if it has inner product $\le 0$ with all simple roots of $D$ in
$d^\perp$. (Such simple roots form the Dynkin diagram of $d^\perp$.)

\proclaim 3.4 Vectors of types 0 and 1. 

The main part of the algorithm breaks down for vectors $u$ of type 0
and 1, because in these cases certain vectors constructed from $u$ are
not in $D$. Fortunately these vectors are easy to classify in
$II_{25,1}$: the primitive type 0 vectors correspond to the Niemeier lattices,
and the type 1 vectors can easily be reduced to type 0 vectors.

{\bf Notation.} $L$ is the lattice $II_{8n+1,1}$, and $D$ is a
fundamental domain for the reflection group $G$ of $L$.

\proclaim Lemma 3.4.1.
Suppose that $Z_1$, $Z_2$ are norm 0 vectors of $L$ that have inner
product $-1$. Then $Z_1$ and $Z_2$ lie in adjacent Weyl chambers of
$G$ and the reflection of the root $r=z_1-z_2$ maps $z_1$ to
$z_2$. $z_1$ and $z_2$ cannot both be in $D$.

Proof. It is easy to check that the reflection in $r^\perp$ exchanges
$z_1$ and $z_2$, so $z_1$ and $z_2$ cannot lie in the same Weyl
chamber, and in particular they cannot both be in $D$. To show that
$z_1$ and $z_2$ lie in adjacent Weyl chambers it is enough to show
that there is no root other then $\pm r$ whose hyperplane passes
strictly between $z_1$ and $z_2$.

Write $L=B\oplus U$ where $U$ is the lattice generated by $z_1$ and
$z_2$. If $r'=b+n_1z_1+n_2z_2$ (with $b$ in $B$) is any root that has
positive inner product with $z_1$ and negative inner product with
$z_2$ then $n_1>$, $n_2<0$, and $b^2-2n_1n_2=r'^2=2$, so $b=0$ as $B$
is positive definite which implies $n_1=1$, $n_2=-1$, and
$r=r'$. Q.E.D.

Vectors $z$ of type 0 in $D$ are just the vectors of norm 0, so they
correspond to $8n$-dimensional even unimodular lattices $B$ in such a
way that $z^\perp\cong B\oplus O$. We will say that $z$ is of type
$B$.

\proclaim Lemma 3.4.2. 
Let $u$ be any primitive vector of $D$ with negative norm. $u$ has
type 1 if and only if $u^\perp$ is the sum of a unimodular lattice $B$
and a one dimensional lattice $N$.  If $u^2=-2n$ then $N$ has determinant
$2n$, and $u$ can be written uniquely as $u=nz_1+z_2$ where $z_1$ and
$z_2$ are norm 0 vectors with $(z_1,z_2)=-1$ and $z_1$  in
$D$. $z_1$ and $z_2$ are both of type $B$ and $z_2$ is not in $D$.

Proof. Let $z_1$ be a norm 0 vector with $(z_1,u)=-1$ and put
$z_2=u-nz_1$, so that $z_2$ is a norm 0 vector with $(z_1,z_2)=-1$. It
is easily checked that $z_1$ is the only norm 0 vector that has inner
product $-1$ with $u$, unless $u^2=-2$ in which case $z_1$ and $z_2$
are the only two such vectors. In either case 3.3.1 implies that there
is a unique such vector in $D$, which we can assume is $z_1$. By 3.4.1
$z_2$ is then not in $D$. 

The lattice  $U=\langle z_1,z_2\rangle$ is unimodular so $L=B\oplus U$ for
some even unimodular lattice $B$. $u^\perp$ is therefore the sum of
$B$ and a one dimensional lattice. Conversely if $u^\perp$ is of this
form it is easy to check that there is a norm 0 vector that has inner
product $-1$ with $u$. Q.E.D.

\proclaim 3.5 Reducing the norm by 2. 

This section gives the main part of the algorithm for finding norm
$-2n$ vectors $u$ from norm $-(2n-2)$ vectors $v$. We consider pairs
of such vectors that have inner product $-2n$, and show how given one
element of such a pair, we can find all possibilities for the other
element. The main theorem is 3.5.4 which shows how to find the
possible $u$'s from the set of simple roots of $D$ that have small
inner product with $v$.

{\bf Notation.} $D$ is a fundamental domain of the reflection group
$G$ of an even unimodular Lorentzian lattice $L$. $u$ is a vector in
$D$ of negative norm $-2n$ and $R_i$ is the set of simple roots of $D$
that have inner product $-i$ with $u$, so $R_0$ is the Dynkin diagram
of $u^\perp$. $v$ is a vector of norm $-2(n-1)$ that is in the same
component of negative norm vectors as $u$. In this section we assume
that $R_0$ is not empty, and choose one of its components $C$.

\proclaim Lemma 3.5.1. 
$-(u,v)\ge 2n-1$ and if $u$ or $v$ has type at least 2 then $-(u,v)\ge 2n$. 

Proof. $(u,v)^2\ge u^2v^2=2n\times 2(n-1)=(2n-1)^2-1$ so $(u,v)\ge
2n-1$ unless $n=1$ and $(u,v)=0$, but this is impossible as it would
give a norm 0 vector $v$ orthogonal to the negative norm vector $u$.

If $-(u,v)=2n-1$ then $u-v$ is a norm 0 vector that has inner product
$-1$ with $u$ and 1 with $v$, so $u$ and $v$ both have types 0 or
1. Q.E.D.

\proclaim Lemma 3.5.2. There is a bijection between 
\item{(1)}
Roots $r$ of $u^\perp$, and
\item{(2)}
Norm $-2(n-1)$ vectors $v$ with $-(v,u)=2n$.
\item{}
given by $v=\tau_u(r)$, $r=\tau_u(v)$. If $u$ has type at least 2 then
$v$ is in $D$ if and only if $r$ is a highest root of some component
of $R_0$ and in this case $v=r+u$.

Proof. The bijection follows from 3.2.1. If $u$ has type at least 2
then by 3.5.1 we have $-(u,v)\ge 2n$, so by 3.3.1 the vector $v$ is in
$D$ if and only if it has inner product $\le 0$ with all roots of
$R_0$, and by 3.2.1 this is true if and only if $\tau_u(v)=r$ has
inner product $\le 0$ with all roots of $R_0$, in other words $r$ is a
highest root of some component of $R_0$. If $r$ is a highest root
then it is fixed by the opposition involution of $u^\perp$ so
$v=\tau_u(r)=u+r$. Q.E.D.

\proclaim Lemma 3.5.3. 
Suppose $v$ is in $D$, $v^2=-2(n-1)$ and $(v,u)=-2n$. Then 
$$\eqalign{
R_0(u)&\subseteq R_0(v)\cup R_1(v)\cup R_2(v)=S(v)\cr
R_i(u)&\subseteq R_0(v)\cup R_1(v)\cup \cdots\cup R_i(v)\hbox{ for } i\ge 1.\cr
}$$

Proof. $v$ is in $D$, so $v=u+r$ for some highest root $r$ of
$u^\perp$. $r$ has inner product 0, $-1$, or $-2$ with all simple
roots of $u^\perp$, and $-r$ is a sum of roots of $R_0(u)$ with
positive coefficients, so $r$ has inner product $\ge 0$ with all
simple roots of $D$ not in $R_0(u)$. The lemma follows from this and
the fact that $(v,s)=(u,s)+(r,s)$ for any simple root $s$ of
$D$. Q.E.D.

{\bf Remark.} If $s$ is in $R_0(u)$ and $R_2(v)$ then we must have
$s=-r$ as $s^2=r^2=-(r,s)=2$, so $s$ must be a component of type
$a_1$ of $R_0(u)$.

We now start with a vector $v$ of norm $-2(n-1)$ and try to
reconstruct $u$ from it. $u-v$ is a highest root of some component of
$R_0(u)$, and $R_0(u)$ is contained in $S(v)$, so we should be able to
find $u$ from $S(v)$. By 3.5.3 $S(u)$ is contained in $S(v)$, so we
can repeat this process with $u$ instead of $v$.

\proclaim Theorem 3.5.4.
Suppose that $v$ has norm $-2(n-1)$ and is in $D$ (so $n\ge 1$). Then
there are bijections between
\item{(1)}
Norm $-2n$ vectors $u$ of $D$ with $(u,v)=-2n$. 
\item{(2)}
Simple spherical Dynkin diagrams $C$ contained in the Dynkin diagram
$\Lambda$ of $D$ such that if $r$ is the highest root of $C$ and $c$
in $C$ satisfies $(c,r)=-i$, then $c$ is in $R_i(v)$.
\item{(3)}
Dynkin diagrams $C$ satisfying one of the following three conditions: 
\itemitem{Either}
$C$ is an $a_1$ and is contained in $R_2(v)$, 
\itemitem{or}
$C$ is an $a_n$ ($n\ge 2$) and the two endpoints of $C$ are in $R_1(v)$ while the other points of $C$ are in $R_0(v)$, 
\itemitem{or}
$C$ is $d_n$ ($n\ge 4$), $e_6$, $e_7$, or $e_8$ and the unique point of
$C$ that has inner product $-1$ with the highest root of $C$ is in
$R_1(v)$ while the other points of $C$ are in $R_0(v)$.

Proof. Let $u$ be as in (1) and put $r=v-u$. $r$ is orthogonal to $u$
and has inner product $\le 0$ with all roots of $R_0$ (because $v$
does) so it is a highest root of some component $C$ of $R_0(u)$. $r$
therefore determines some simple spherical Dynkin diagram $C$
contained in $\Lambda$. Any root $c$ of $C$ has $(c,v-r)=(c,u)=0$, so
$c$ is in $R_i(v)$ where $i=-(c,r)$. This gives a map from (1) to (2).

Conversely if we start with a Dynkin diagram $C$ satisfying (2) and
put $u=v-r$ (where $r$ is the highest root of $C$) then $(c,u)=0$ for
all $c$ in $C$, so $(r,u)=0$ as $r$ is a sum of the $c$'s. This
implies that $u^2=-2n$ and $(u,v)=-2n$. We now have to show that $u$
is in $D$. Let $s$ be any simple root of $D$. If $s$ is in $C$ then
$(s,r)=(s,v)$ and if $s$ is not in $C$ then $(s,r)\ge 0$, so in any
case $(s,u)=(s,v-r)\le 0$ and hence $u$ is in $D$. This gives a map
from (2) to (1) and shows that (1) and (2) are equivalent.

Condition (3) is just the condition (2) written out explicitly for
each possible $C$, so (2) and (3) are also equivalent. Q.E.D.

\proclaim 3.6 Vectors of $B/nB$. 

Here we show that the type $n$ vectors of $II_{25,1}$ are closely
related to vectors of $B/nB$ for Niemeier lattices $B$. The type $n$
vectors of $II_{25,1}$ can be found from those of small norm; in
particular the type 2 vectors can be found from those of norm $-2$ or
$-4$ so all type 2 vectors of $II_{25,1}$ can be read off from tables
$-2$ and $-4$.

{\bf Notation.} $L$ is an even unimodular Lorentzian lattice. 

\proclaim Theorem 3.6.1. 
There is a 1:1 correspondence between isomorphism classes of
\item{(1)}
Triples $(\hat b, B, k)$ where $\hat b$ is an element of $B/nB$ with
$\hat b^2\equiv k\bmod 2n$ and $B$ is an even unimodular positive
definite lattice whose dimension is 2 less than that of $L$, and
\item{(2)}
Pairs $(z,u)$ of vectors of $L$ with $z^2=0$, $(z,u)=-n$, and $z$ primitive.
\item{} We have $u^2=k$ and $B$ is the Niemeier lattice corresponding to $z$. 

Proof. Suppose we are given $z$ and $u$. We choose a norm 0 vector
$z'$ with $(z,z')=-1$ so that $L=B\oplus\langle z,z'\rangle$ for some
lattice $B$. Then $u=b+mz+nz'$ for some $b$ in $B$ and some integer
$m$. All other norm 0 vectors $z''$ that have inner product $-1$ with
$z$ are conjugate to $z'$ under the ``translation'' automorphism of
$L$ that maps $b+mz+nz'$ to $(b+nb_1)+(m+(b_1,b)+nb_1^2/2)z+nz'$ for
some $b_1$ in $B$, so changing $z'$ does not change $b\bmod nB$. Hence
given $z$ and $u$, $B$ and $\hat b$ are well defined. The norm of
$\hat b$ is well defined mod $2n$ and is congruent to $u^2$, because
$u^2=b^2-2mn$, so we can put $k=u^2$.

Conversely if we are given $(\hat b, B, k)$ as in (1) then we choose a
representative $b$ of $\hat b$ in $B$. We put $L=B\oplus\langle
z,z'\rangle$ where $z$ and $z'$ are norm 0 vectors with inner product
$-1$, and take $u$ to be $b+mz+nz'$, where $m$ is chosen so that $u$
has norm $k$, in other words $m=(b^2-k)/2n$. If we choose a different
representative of $\hat b$ in $B$ then it is easy to check that the
new $u$ we get is conjugate to the old one under a ``translation''
automorphism of $L$ fixing $z$, so the isomorphism class of $(z,u)$ is
well defined. The maps we have constructed from (1) to (2) and from
(2) to (1) are inverses of each other. Q.E.D.

We can restrict $k=u^2$ to run through a set of representatives of the
even integers mod $2n$. For example vectors of $B/2B$ for Niemeier
lattices $B$ correspond to pairs $(z,u)$ of $II_{25,1}$ as in (2) with
$u^2=-2$ or $-4$.

\proclaim 3.7 Type 2 vectors. 

If $u$ is a norm $-2$ vector of $II_{25,1}$ in $D$ then reflection in
$u^\perp$ is an automorphism of $II_{25,1}$ and from this we can
construct an automorphism of $D$ as in 1.5. If $u$ has larger norm
then reflection in $u^\perp$ is no longer an automorphism of
$II_{25,1}$, but if $u$ has type 2 then we can still construct an
automorphism of $D$ associated to $u$ by a more complicated
construction. The existence of this automorphism will be used in
chapter 4 to prove some properties of the root system of $u^\perp$.

We will also describe the relation between type 2 vectors of $II_{25,1}$ and vectors of $B/2B$ in more detail. (See table 3.)

{\bf notation.} $L$ is an even unimodular Lorentzian lattice whose
reflection group $G$ has fundamental domain $D$.

\proclaim Theorem 3.7.1. 
Suppose that $u$ is a norm $-2n\le -2$ vector of $D$ such that there
is a norm 0 vector $z$ with $(z,u)=-2$. Then there is an element
$\sigma_u$ of $\Aut(L)$ not depending on $z$ with the following
properties:
\item{(1)}
$\sigma_u$ fixes $u$ and $D$.
\item{(2)}
$\sigma_u^2=1$.
\item{(3)}
Any element of $L$ fixed by $\sigma_u$ is in the vector space
generated by $u$, $z$, and the norm 2 roots of $u^\perp$.

Proof. If $u$ has norm $-2$, so that reflection in $u^\perp$ is an
automorphism of $L$, then the existence of $\sigma_u$ follows from
1.5.1, so we will assume that $-u^2=2n\ge 4$. We construct $\sigma_u$
as the product of 3 commuting involutions $r(u)$, $\sigma(u^\perp)$,
and $r_z$, where $r(u)$ is reflection in $u^\perp$, $\sigma(u^\perp)$
is the opposition involution of $u^\perp$, and $r_z$ is reflection in
a certain vector of $u^\perp$. (Note that $r(u)$ and $r_z$ are not
automorphisms of $L$.) $\sigma_u$ fixes $u$ because $r_z$ does and
$r(u)$ and $\sigma(u^\perp)$ multiply $u$ by $-1$. $r(u)$ obviously
commutes with $r_z$ and $\sigma(u^\perp)$.

The vector $h=u-nz$ has norm $2n$ and is in $u^\perp$. We let $H$ be
the hyperplane $h^\perp$. Let $G_0$ be the group generated by
reflections of the roots of $u^\perp$, and let $G_1$ be the group
generated by $G_0$ and $r(H)$, where $r(H)$ is reflection in $H$.

We now check that the orbit of $H$ under $G_0$ does not depend on
$z$. If $u^2\le -6$ then there is at most one norm 0 vector $z$ that
has inner product $-2$ with $u$, while if $u^2=-4$ then the projection
of any such vector $z$ into $u^\perp$ is a norm 1 vector of the
unimodular lattice containing $u^\perp$. Any two such norm 1 vectors
are conjugate under the group generated by $G_0$ and $-1$, so the
hyperplane orthogonal to these norm 1 vectors are all conjugate under
$G_0$.

Let $x$ be any vector of $L$. Then
$$\eqalign{
r(u)r(H)x&= x-2(x,u)u/(u,u)-2(x,h)h/(h,h)\cr
&= x+2(x,u)u/2n-2(x,u-nz)(u-nz)/2n\cr
&= x+(x,z)(u-nz)+(x,u)z\cr
}$$
so $r(u)r(H)x$ is in $L$ and hence $r(u)r(H)$ fixes $L$.

$r(H)$ normalizes $G_0$ and is not in $G_0$ because it is not an
automorphism of $L$, so $G_0$ is a reflection group of index 2 in the
reflection group $G_1$. In particular the fundamental domain $D_0$ of
$G_0$ containing $D$ contains exactly 2 fundamental domains for $G_1$,
so there is a unique conjugate $H'$ of the hyperplane $H$ passing
through $D_0$. We define $r_z$ to be reflection in $H'$. $r_z$ is the
unique element of $G_1$ not in $G_0$ that fixes $D_0$, so $r_z$
commutes with $\sigma(u^\perp)$.

$\sigma_u=r(u)r_z\sigma(u^\perp)$ fixes $D$ because $r_z$ and
$r(u)\sigma(u^\perp)$ do, so (1) holds. $\sigma_u^2=1$ because
$\sigma_u$ is the product of three commuting involutions, so (2)
holds. On the space orthogonal to $u$, $z$, and all roots of
$u^\perp$, $\sigma_u$ acts as $-1$, and this implies that any element
of $L$ fixed by $\sigma_u$ is in the vector space generated by $u$,
$z$, and roots of $u^\perp$, so (3) holds. $\sigma_u$ does not depend
on $z$ because $G_1$ and hence $r_z$ are independent of $z$. Q.E.D.

\proclaim Lemma 3.7.2.
Let $u$ and $v$ be vectors of $D$ with $-u^2=4$, $-v^2=2$,
$-(u,v)=4$. If there is a norm 0 vector $z$ that has inner product
$-2$ with $u$ and $v$ then $\sigma_u=\sigma_v$.

Proof. The existence of $z$ implies that $\sigma_u$ is well
defined. $\sigma_u$ and $\sigma_v$ both have period 2 and both fix the
fundamental domain $D$, so to prove they are equal it is sufficient to
prove that their product is in the reflection group of $L$. We have
$\sigma_u=-r(u)r(z_u)g_1$ and $\sigma_v=-r(v)g_2$ where $g_1,g_2$ are
in the reflection group and $z_u$ is the projection of $z$ into
$u^\perp$, so we just have to prove that $r(z_u)r(u) r(v)$ is in the
reflection group of $L$. (It is a product of reflections, but these so
not come from norm 2 vectors.)

Put $r_1=z-v$, $r_2=u-z-v$, $r_3=u-v$. Then $r_1$, $r_2$, and $r_3$
all have norm 2, and $r(z_u)r(u)r(v)=r(r_3)r(r_2)r(r_1)$ because both
sides map $z_u$ to $-z_u$, $v$ to $2u-v$, $u$ to $u-4v$, and fix the
subspace of $L$ orthogonal to $z_u$, $u$, and $v$. Q.E.D.

\proclaim Corollary 3.2.2${1\over 2}$.
If $u$ is a norm $-4$ vector in $D$ of type $\le 2$ and $v$ is a norm
$-2$ vector of $D$ with $-(u,v)=4$ then
\item{Either:} $\sigma_u=\sigma_v$
\item{Or:} 
$v$ has type 1 and the $8n+1$-dimensional unimodular lattice
corresponding to $u$ has exactly 4 vectors of norm 1.

Proof. Suppose $\sigma_u\ne \sigma_v$ and let $z_1$, $z_2$ be norm 0
vectors with $-(z_i,u)=2$, $z_1+z_2=u$. We have
$(z_1,v)+(z_2,v)=(u,v)=-4$, $(z_i,v)<0$, and by 3.7.2 $(z_i,v)\ne -2$,
so one of $(z_i,v)$ is $-1$ and hence $v$ has type 1. The projection
$-r$ of $v$ into $u^\perp$ is a root of the unimodular lattice $A$
containing $u^\perp$ and $(r,i)=\pm 1$ for all norm 1 vectors $i$ of
$A$ because these are the projections of norm 0 vectors $z$ with
$-(z,u)=-2$ and $v$ has inner product $-1$ or $-3$ with all such
vectors. This implies that $A$ has exactly 4 vectors of norm 1 (as it
certainly has some vectors of norm 1). Q.E.D.

{\bf Remarks.} The second case of 3.2.2${1\over 2}$ does sometimes occur. In
fact the dimensions of the fixed spaces of $\sigma_u$ and $\sigma_v$
can be different so $\sigma_u$ and $\sigma_v$ are not even
conjugate. 3.2.2${1\over 2}$ implies that if $v_1$ and $v_2$ are type 2 norm
$-2$ vectors of $D$ that both have inner product $-4$ with a norm $-4$
vector $u$ of $D$, then $\sigma_{v_1}=\sigma_{v_2}$ when $u$ has type
2. This is not always true if $u$ has type 3. If $u$ is in one of the
three orbits of norm $-4$ vectors of $D$ of height $\le 2$ (so that
$u$ corresponds to a unimodular lattice with no roots) then $\sigma_u$
is not equal to $\sigma_v$ for any norm $-2$ vector $v$ of $D$, so
3.7.1 does give a few new elements of $\Aut(D)$.

We now describe the relation between type 2 vectors and orbits of
$B/2B$, where $B$ runs through the even unimodular lattices of the
right dimension. By 3.6.1 there is a bijection between
\item{(1)}
Pairs $(\hat b,B)$ where $B$ is as above and $\hat b$ is in $B/2B$, and
\item{(2)}
Pairs $(z,u)$ with $z^2=0$, $-u^2=4$, $-(z,u)=2$ (where $z$ and $u$
are in $L$).

If $\hat b$ is in $B/2B$ we will write $\hat b^2$ for the minimum
possible norm $b^2$ of a representative $b$ of $\hat b$ in $B$. For
each lattice $B$ we draw the following graph. (See table 3 for all
examples with $\dim(B)\le 24$.) The vertices of the graph are the
orbits of $B/2B$ under $\Aut(B)$ and they will lie on rows according
to the value of $\hat b^2$, with one row for each possible value of
$\hat b^2$. Two vertices are joined if and only if there are
representatives of these two orbits whose difference is a root of
$B$. If $\hat b$ is in $B/2B$ we write $R_b$ for the roots of $B$ that
have even inner product with $\hat b$. To the vertex corresponding to
$\hat b$ we attach the Dynkin diagram of $R_b$ arranged into orbits
under the subgroup $G$ of $\Aut(B)$ fixing $\hat b$ and the Dynkin
diagram of $R_b$. We write the order of $G$ next to this Dynkin
diagram.

\proclaim Lemma 3.7.3.
\item{(1)} 
If $b$ is a minimal norm representative of $\hat b$ then all roots of
$B$ have inner product 0, $\pm 1$, $\pm 2$ with $b$.
\item{(2)}
If two vertices of the graph are joined then they lie on different
adjacent rows.

Proof.
\item{(1)} 
If $r$ is any root of $B$ then $(2r+b)^2\ge b^2$ because $b$ is a
minimal norm representative of $\hat b$, so $-(b,r)\le 2$; this proves
(1).
\item{(2)}
Let $b_1$ and $b_2$ be minimal norm representatives of the vertices
with $b_1^2\le b_2^2$. By assumption there is a root $r$ such that
$r+b_1$ is conjugate to a representative of $\hat b_2$. By (1),
$(r,b_1)$ is 0, $\pm 1$, or $\pm 2$, so we can assume that $(r,b_1)$ is
0, $-1$ or $-2$. If it was $-1$ then $r+b_1$ would be conjugate to
$b_1$ by a reflection which is impossible as the vertices $b_1$ and
$b_2$ are different, while if it was $-2$ then $r+b_1$ would have
smaller norm than $b_1$ which is impossible as we assumed that
$b_1^2\le b_2^2$. Hence $(r,b_1)=0$ so $(r+b_1)^2=b_1^2+2$ and $r+b_1$
is a minimal norm representative for $\hat b_2$, so $\hat b_1^2+2=\hat
b_2^2$. Q.E.D.

By the remark at the end of 3.6 the vectors of $B/2B$ correspond to
pairs $(z,u)$ where $(z,u)=-2$ and $z^2=0$ and $u^2=-2$ or $-4$, and
hence either to (some) roots $r$ of 25 dimensional bimodular lattices
$C$ or to 24 dimensional unimodular lattices $A$ that have $B$ as a
neighbor. We now describe the lattice corresponding to the vector
$\hat b$ of $B/2B$ in terms of the norm of $\hat b$ and the Dynkin
diagram $R$ of the roots of $B$ that have even inner product with
$\hat b$. If $\hat b^2=0\bmod 4$ then $\hat b$ corresponds to some 24
dimensional unimodular lattice $A$, and if $\hat b^2=2\bmod 4$ then it
corresponds to a root $r$ of a 25 dimensional bimodular lattice
$C$. Let $b$ be a minimal norm representative for $\hat b$ in $B$. See
table 3 for examples. Proofs are omitted as they are easy and
uninteresting.

\item{$\hat b^2=0$}
$A$ is isomorphic to $B$.
\item{$\hat b^2=2$}
$C$ is isomorphic to $B\oplus a_1$ and $r$ corresponds to the root $b$ of $B$.
\item{$\hat b^2=4$}
$A$ is an odd unimodular lattice containing vectors of norm 1. If it
has $2n$ vectors of norm 1 then the components of $R$ joined to
vertices of norm 2 in table 3 form a $d_n$ Dynkin diagram (where $d_1$
is empty, $d_2=a_1^2$, $d_3=a_3$).
\item{} 
Example. If $B=A_1^{24}$ then there are two classes of norm 4 vectors
in $B/2B$, which have root systems $a_1^2a_1^{22}$ and $a_1^{16}$. The
corresponding lattices $A$ have 4 and 2 vectors of norm 1, so they
give lattices of dimension 22 or 23 with no vectors of norm 1 and root
systems $a_1^{22}$ and $a_1^{16}$.
\item{$\hat b^2=6$}
$r$ is a root in $C$ such that the component of the Dynkin diagram $S$
of $C$ containing $r$ is not an $a_1$. If the Dynkin diagram of $C$ is
$SX$ then the Dynkin diagram $R$ is $S'X$ where $S'$ is the Dynkin
diagram of the roots of $S$ orthogonal to $r$. A component of $R$ is
joined to a vertex of norm 4 in table 3 if it is in $S'$ and to a
vertex of norm 8 otherwise. This almost determines the Dynkin diagram
of $C$ in terms of $R$ and the set $S'$ of components of $R$ joined to
vectors of norm 4; we change $S'$ as follows:
\halign{
\hfill$#$&#\hfill&\qquad$#$\hfill\cr
S'&&S\cr
\hbox{empty}&&a_2\cr
a_n&($n$ not 5)&a_{n+2}\cr
a_5&&a_7\hbox{ or } e_6\cr
a_1d_n&($n$ at least 2)&d_{n+2}\qquad(d_2=a_1^2, d_3=a_3)\cr
d_6&&e_7\cr
e_7&&e_8\cr
}
\item{} 
Example. Let $B$ be $A_4^6$. There is a $\hat b$ of norm 6 in $B/2B$
such that $R$ is $a_6^3a_5$ and the $a_5$ is joined to a vertex of
norm 4. Hence the Dynkin diagram of $C$ is either $a_6^3a_7$ or
$a_6^3e_6$. By 4.3.4 the norm of the Weyl vector of $C$ is half a
square, so the Dynkin diagram must be $a_6^3e_6$.
\item{$\hat b^2=8$}
$A$ is an odd unimodular lattice with no vectors of norm 1 with two
even neighbors $B$ and $D$ such that $D$ has roots not in $A$. A
component of $R$ is joined to a vertex of norm 10 in table 3 if and
only if it is a component of the Dynkin diagram of $D$.
\item{} Example. 
For $B=A_9^2D_6$ there are two vectors of norm 8 with $R$ equal to
$a_5^2d_6a_3^2$. The other even neighbors of the lattices $A$
corresponding to them are $A_9D_6^2$ and $D_6^4$.
\item{$\hat b^2=10$}
$r$ is a root of $C$ such that the component of the Dynkin diagram $S$
of $C$ is an $a_1$. $S$ is $a_1R$.
\item{$\hat b^2=12$}
$A$ is an odd 24 dimensional unimodular lattice such that the other
neighbor of $A$ has the same roots as $A$. See 4.5.3.

For Niemeier lattices $\hat b^2$ is at most 12. In higher dimensions
vectors $\hat b$ with norm greater than 12 behave like those of norm
10 or 12.

\proclaim 3.8 Vectors of type at least 3. 

Vectors of type $\ge 3$ form a sort of general case. By 3.5,  3.6
and 3.7 all vectors of type $\le 2$ in $II_{25,1}$ are known, so we
show how the algorithm for finding negative norm vectors simplifies
for vectors of type $\ge 3$. (3.8.6.) We also show that for norm $-2n$
vectors $u$ of type $\ge 3$ the correspondence between roots that have
inner product $-1$ with $u$ and norm $-(2n-4)$ vectors $v$ with
$-(z,u)=2n-1$ behaves well. For example if $u$ in $D$ has norm $-4$
then there is a bijection between simple roots having inner product
$-1$ with $u$, and norm 0 vectors of $D$ having inner product 3 with
$u$. See 4.5 for some examples.

{\bf Notation.} $L$ is an even unimodular Lorentzian lattice. $u$ is a
vector of norm $-2n\le -4$ in the fundamental domain $D$ of the
reflection group of $L$. $R_i(u)$ is the set of simple roots of $D$
that have inner product $-i$ with $u$.

\proclaim Lemma 3.8.1. There is a bijection between 
\item{(1)}
Roots $r$ with $(r,u)=-i$, and
\item{(2)} Vectors $z$ with $z^2=-2(n-2)$, $-(z,u)=2n-1$
\item{}
given by $z=\tau_u(r)$, $r=\tau_u(z)$. If $z$ is in $D$ then $r$ is a
simple roots of $D$. If $u$ has type at least 3 then $z$ is in $D$ if
and only if $r$ is a simple root of $D$.

Proof. By 3.2.1 $\tau_u$ exchanges the sets in (1) and (2) and $r$ has
inner product $\le 0$ with all roots of $R_0(u)$ if and only if $z$
does. By Vinberg's algorithm using $u$ as a controlling vector a root
$r$ with $(r,u)=-1$ is simple if and only if it has inner product $\le
0$ with all roots of $R_0(u)$, so if $z$ is in $D$ then $r$ is simple.

Now suppose that $u$ has type t least 3. If $z$ is any vector of norm
$-2(n-2)$ that has negative inner product with $u$ then $(z,u)^2\ge
z^2u^2=2n\times 2(n-2)=(2n-2)^2-4$, so either $-(z,u)\ge 2n-2$ or
$n=2$. If $n=2$ then $z$ has norm 0 and $u$ has type at least 3, so
$-(u,z)\ge 3$. If $-(u,z)\ge 2n-2$ then $-(u,z)$ cannot be equal to
$2n-2$ because then $u-z$ would be a norm 0 vector having inner product
$-2$ with $u$, which is impossible because $u$ has type at least 3, so
$-(u,z)\ge 2n-1$. In both cases we have $-(u,z)\ge 2n-1$, so by 3.3.1
a vector $z$ with $z^2=-2(n-2)$ and $-(u,z)=2n-1$ is in $D$ if and
only if it has inner product $\le 0$ with all simple roots of
$R_0(u)$. Hence $z$ is in $D$ if and only if $r$ is simple. Q.E.D.

From now on $u$ and $v$ are vectors of $D$ with $u^2=-2n$,
$v^2=-2(n-1)$, and $(u,v)=-2n$ (as in 3.5). We also assume that $u$
has type at least 3.

\proclaim Lemma 3.8.2. $v$ has type at least 2. 

Proof. If $v$ had type 0 or 1 then by 3.4.2 we would have
$v=(n-1)z+z'$ for norm 0 vectors $z$ and $z'$. Then
$(n-1)(z,u)+(z',u)=(v,u)=-2n$ and $z$ and $z'$ have negative inner
products with $u$, so $z$ and $z'$ cannot both have inner products
less than $-2$ with $u$. This is impossible as $u$ has type at least
3. Q.E.D.

If $r$ is a simple root in $R_1(u)$ and $R_0(v)$ then there are two
ways to construct a norm $-2(n-2)$ vector from it:
\item{(1)}
$r$ is in $R_1(u)$ and so determines a norm $-2(n-2)$ vector $z$ as in
3.8.1.
\item{(2)}
$r$ is in some component of $R_0(u)$ and this component of $R_0(u)$
determines a vector $z'$ of norm $-2(n-2)$ with $-(v,z')=2n-2$ as in
lemma 3.5.2.

\proclaim Lemma 3.8.3 $z=z'$. 

Proof. By 3.8.2 $v$ has type at least 2, so by 3.5.2 $z'$ is in
$D$. By 3.8.1 $z$ is in $D$, so to show that $z=z'$ it is sufficient
to show that they are conjugate under the group generated by $-1$ and
the reflections of $L$, because they are both in the same fundamental
domain $D$ of this group.

By 3.8.1 the opposition involution $\sigma(u^\perp)$ maps $z$ to
$r-u$. The vector $r'=u-v$ satisfies $r'^2=2$, $(r',u)=0$,
$(r',r)=-1$, so reflection in $r'$ maps $r-u$ to 
$$\eqalign{ &(r-u)-(r',r-u)r'\cr =& (r-u)-(-1)(u-v)\cr =&r-v\cr }$$
The vector $z'$ is equal to $v+r''$, where $r''$ is the highest root
of the component of $R_0(v)$ containing $r$, and hence $r''$ is
conjugate to $-r$ by reflections leaving $v$ fixed. This implies that
$r-v$ is conjugate to $-r''-v=-z'$ under the reflection group of $L$,
so we have shown that $z$ is conjugate to $r-u$, $r-v$, $-z'$, and
$z'$ under the group generated by reflections of $L$ and $-1$. Q.E.D.

{\bf Remark.} If $u$ has type less than 3 then we can still define $z$
and $z'$ and they are conjugate under the reflection group of $L$, but
they can be different.

We will now show how the method in section 3.5 for constructing norm
$-2n$ vectors $u$ from norm $-2(n-1)$ vectors $v$ can be described
more explicitly for vectors $u$ of type at least 3; in fact, there are
only 13 different ways to construct vectors $u$. This is useful
because all vectors of type 0, 1, and 2 are known once those of norms
0, $-2$, or $-4$ are known (which they are for $II_{25,1}$).

Let $v$ have norm $-2(n-1)$ and be in $D$. Recall from 3.5.4 that norm
$-2n$ vectors $u$ of $D$ with $-(u,v)=-2n$ correspond to certain
simple spherical Dynkin diagrams $C$ whose points are simple roots of
$D$. Let $c$ be the highest root of $C$, so $v=c+u$.

\proclaim Lemma 3.8.5. 
$u$ has type at least 3 if and only if $v$ has type at least 2 and $c$
has inner product 0 or $-1$ with all highest roots of the root system
of $v^\perp$.

Proof. Step 1. If $c_v$ is the projection of $c$ into $v^\perp$ then
$$ c_v^2=c^2-(c,v)^2/(v,v)=2+4/-v^2\le 4$$ so if $r$ is any root of
$v^\perp$ then $(c,r)^2=(c_v,r)^2\le c_v^2r^2\le 4\times 2=8$, so
$(c,r)$ is $0$, $\pm 1$ or $\pm 2$. If $r$ is a highest root of
$v^\perp$ then $(r,c)=(r,v-u)=-(r,u)\le 0$ as $-r$ is a sum of simple
roots of $D$, so $(r,c)=0$,$- 1$, or $-2$.

Step 2. Now assume that $u$ has type at least 3. By 3.8.2 $v$ has type
$\ge 2$. If $c$ does not has inner product 0 or $-1$ with all highest
roots of $v^\perp$ then by step 1 there is a highest root $r$ with
$(r,c)=-2$. Then $r+c$ is a norm 0 vector that has inner product 2
with $u$, which is impossible. This proves that if $u$ has type $\ge
3$ then both conditions of 3.8.5 must hold.

Step 3. We now assume that $v$ has type at least 2 and $c$ has inner
product 0 or $-1$ with all highest roots of $v^\perp$ and prove that
$u$ has type at least 3. If not, there is a norm 0 vector $z$ with
$(z,u)=-1$ or $-2$; by changing $z$ to $2z$ if necessary we can assume
that $(z,u)=-2$. If $(z,c)\le 0$ we reflect $z$ in $c$, so we can also
assume that $(z,c)\ge 0$. ($(c,u)=0$ so this does not affect $(z,u)$.)
We have $v=c+u$, $(z,u)=-2$, $(z,c)\ge 0$ and $(z,v)\le -2$ (because
$v$ has type at least 2) so $(z,v)=-2$ and $(z,c)=0$. $r=z+c$ is then
a norm 2 vector in $v^\perp$ with $(r,c)=2$, so by step 1 the highest
root of the component of $v^\perp$ containing $r$ has inner product
$-2$ with $c$, contradicting our assumptions. Hence $u$ has type at
least 3. Q.E.D.

\proclaim Theorem 3.8.6. 
Let $v$ be a vector in $D$ of type $\ge 2$ and norm $-2(n-1)$. Write
$V_0$, $V_1$ and $V_2$ for the simple roots of $D$ that have inner
products $0$, $-1$, or $-2$ with $v$ (and likewise for $U_i$). Let $X$
be one of $a_i$, $d_i$ or $e_i$. Then the vectors $u$ of $D$ of type
at least 3 with $-u^2=2n=-(u,v)$ such that the component $C$ of
$u^\perp$ corresponding to $v$ is an $X$ can be obtained as follows:
\item{$X$ is $a_1$:}
$C$ is a point of $V_2$ joined to at most one point in each component
of $V_0$, and any point to which it is joined is a tip of its
component. The Dynkin diagram $U_0$ of $u^\perp$ is $C$ together with
the points of $V_0$ not joined to $C$.
\item{}
From now on the points of $V_2$ are not used, and $C$ is always the
component of $U_0$ containing a certain point $x$. Black dots
$\bullet$
represent points of $V_0$ and white dots $\circ$ represent points of $V_1$.
\item{$X$ is $a_2$:}
$C$ consists of two points $x$, $y$ of $V_1$, which are joined and
such that there is no component of $V_0$ joined to both $x$ and
$y$. $U_0$ is $C$ together with the points of $V_0$ not joined to $C$.
\item{$X$ is $a_i$:}
$(i\ge 3)$ Case 1. There is a component $a_j$ $(j\ge i-2)$ of $V_0$
and two points $x$, $y$ of $V_1$ such that $x$ is joined to one end of
this $a_j$, $y$ is joined to some point of the $a_j$, and there is no
other component of $V_0$ joined to both $x$ and $y$.
$U_0$ is $x$, $y$, and $V_0$ with the following points deleted:
\itemitem{(1)} 
Any points of $V_0$ not in the $a_j$ that are joined to $x$ or $y$.
\itemitem{(2)} A point of $a_j$ if $j>i-2$. 
\item{} $C$ is the component of $U_0$ containing $x$ and $y$.
\item{$X$ is $a_i$:}
$(i\ge 5)$ Case 2. There is a component $d_{i-1}$ of $V_0$ and two points $x$, $y$ of $V_1$ joined to 2 tips of $d_{i-1}$ not both joined to the same point of
$d_{i-1}$. $U_0$ is $x$, $y$, and $V_0$ with the following points deleted:
\itemitem{(1)} 
Any points of $V_0$ not in the $d_{i-1}$ that are joined to $x$ or $y$.
\itemitem{(2)}
The third tip of the $d_{i-1}$. 
\item{$X$ is $d_4$:}
There is a point $x$ of $V_1$ joined to an end of each of three
components $a_i$, $a_j$, $a_k$ of $V_0$ and joined to no other points
of $V_0$.  $U_0$ is $x$ and $V_0$ with a point deleted from any of the
$a$'s that have more than one point.
\item{$X$ is $d_5$:}
There is a point $x$ of $V_1$ joined to one end of an $a_i$ $(i\ge 1)$
of $V_0$ and to the point next to one end of an $a_j$ $(j\ge 3)$ of
$V_0$ and to no other points of $V_0$.
$U_0$ is $x$ and $V_0$ with a point of the $a_i$ deleted if $i\ge 2$ 
and a point of the $a_j$ deleted if $j\ge 4$. 
\item{$X$ is $d_i$:}
$(i\ge 6)$ (There is a second way to get $d_6$'s and $d_7$'s.) Find a
point $x$ of $V_1$ joined to one end of an $a_j$ ($j\ge 1$) of $V_0$,
to a tip of a 
$d_{i-2}$ of $V_0$, and to no other points of $V_0$. $U_0$ is $x$ and
$V_0$, with a point of the $a_j$ deleted if $j\ge 2$.
\item{$X$ is $d_6$:}
Case 2. Find a point $x$ of $V_1$ joined to one end of an $a_i$ $(i\ge
1$), to a tip of a $d_5$, and to no other points of $V_0$. $U_0$ is
$V_0$ and $x$ with a point of the $d_5$ deleted and a point of the
$a_i$ deleted if $i\ge 2$.
\item{$X$ is $d_7$:}
Case 2. Find a point $x$ of $V_1$ joined to a tip of an $e_6$, to one
end of an $a_i$ of $V_0$ ($i\ge 1$) and to no other points of
$V_0$. $U_0$ is $x$ and $V_0$ with a point deleted from the $e_6$ and
a point deleted from the $a_i$ if $i\ge 2$.
\item{$X$ is $e_6$:}
Find a point $x$ of $V_1$ joined to a point third from the 
end of an $a_i$ of $V_0$ ($i\ge 5$) and to no other points of
$V_0$. $U_0$ is $x$ and $V_0$ with a point deleted from the $a_i$ if $i\ge 6$.
\item{$X$ is $e_7$:}
Case 1. Find a point $x$ of $V_1$ joined to a tip of a $d_6$ of
$V_0$  and to no other points of
$V_0$. $U_0$ is $x$ and $V_0$.
\item{$X$ is $e_7$:}
Case 2. Find a point $x$ of $V_1$ joined to a tip of a $d_7$ of
$V_0$  and to no other points of
$V_0$. $U_0$ is $x$ and $V_0$ with a point deleted from the $d_7$.
\item{$X$ is $e_8$:}
Find a point $x$ of $V_1$ joined to an $e_7$ of $V_0$ and to no other
points of $V_0$. $U_0$ is $x$ and $V_0$. (This case only happens once
for norm $-4$ vectors in $II_{25,1}$.)

Proof. Apply lemma 3.8.5 to every possible $X$. Q.E.D.

\proclaim Corollary 3.8.7. 
Notation as in 3.8.6. If the component $C$ is not of type $a_n$ then
$S(u^\perp)=S(v^\perp)+1$.

Proof. Check each case of 3.8.6 for $C$ a $d_n$ ($n\ge 4$), $e_6$,
$e_7$, or $e_8$. Q.E.D.

\proclaim 3.9 Positive norm vectors in Lorentzian lattices. 

Classifying vectors of zero or negative norm in Lorentzian lattices is
usually difficult. Classifying positive norm vectors turns out to be
much easier. Here we will find the positive norm vectors of
$II_{8n+1,1}$; the same method can be used for vectors in other
lattices whose orthogonal complement is indefinite although the result
will usually be more complicated.

\proclaim Theorem 3.9.1. 
If $m\ge 1$ and $n\ge 1$ then $II_{8n+1,1}$ has exactly one orbit of
primitive vectors of norm $2m$.

Proof. By an obvious generalization of 3.1.1 such orbits are in
bijection with pairs $(a,A)$ such that
\item{(1)}
$A$ is an $8n+1$-dimensional even Lorentzian lattice, and
\item{(2)}
$A'/A$ is cyclic, generated by $a$ in $A'/A$ of norm $-1/2m\bmod 2$. 

We will first show that $A$ is unique by finding its genus and then
show that $\Aut(A)$ is transitive on the elements $a$ satisfying
(2). $A$ is indefinite and has dimension at least 3, so by a theorem
of Eichler it is determined by its spinor genus. $A'/A$ is cyclic so
the spinor genus of $A$ is determined by its genus. Condition (2)
determines $A\otimes \Z_p$ if $p$ is odd and $A\otimes \R$ is
determined by the fact that $A$ is Lorentzian, so $A$ is determined by
$A\otimes \Z_2$. ($\Z_p$ is the ring of $p$-adic integers and $\R$ is
the reals.)

We now work out $A\otimes \Z_2$. Let $x$ be an element of $A\otimes
\Z_2$ such that $x/2m$ represents $a$, so that $x^2\equiv -2m\bmod
8m^2$. We will show that $A\otimes \Z_2$ is determined up to
isomorphism by $x^2$. The $\Z_2$ lattice $x^\perp$ is even and self
dual and its determinant is determined by $x^2$, so it is
unique. $A\otimes \Z_2$ is the sum of this lattice and the lattice
generated by $x$, so it is determined by $x^2$. Let $\hat A$ be the
lattice $e_8^n\oplus X$ where $X$ is generated by an element of norm
$-2m$. Then for any 2-adic integer congruent to $-2m\bmod 8m^2$ it is
easy to check that $\hat A\otimes \Z_2$ contains an element $\hat x$ of
this norm such that $(\hat x, \hat a)\equiv 0\bmod 2m$ for all $\hat
a$ in $\hat A$ (because $e_8$ contains elements of any positive even
norm). This implies that $A\otimes \Z_2$ is isomorphic to $\hat
A\otimes \Z_2$, and so the genus of $A$, and hence $A$ are uniquely
determined by the conditions (1) and (2). The lattice $\hat A$
satisfies these conditions, so $A\cong\hat A$.

We can now find the orbits of $\Aut(A)$ on $A'/A$. We can put
$A=(e_8)^n\oplus X$ where $X$ is generated by $x$ of norm $-2m$. The
elements $a$ of $A'/A$ with $a^2\equiv -1/2m\bmod 2$ are represented
by the elements $kx/2m$ where $k$ runs over a set of representatives
of the elements of $\Z/2m\Z$ with $k^2\equiv 1\bmod 4m$. For any such
$k$ we construct an automorphism of $A$ mapping $x/2m$ in $A'/A$ to
$kx/2m$. To do this we can and will assume that $A\cong e_8\oplus X$,
in other words $n=1$. We can also assume that $k\ge 1$ so we can find
an element $e$ in $e_8$ with $e^2=(k^2-1)/2m$ as this number is
even. If $x'=2me+kx$ then $x'^2=-2m$ and $2m$ divides $(x',y)$ for all
$y$ in $A$, so $x'^\perp$ is an even positive definite 8 dimensional
unimodular lattice and hence is isomorphic to $e_8$ as $e_8$ is the
only such lattice. This implies that $A\cong \langle x'\rangle \oplus
x'^\perp
\cong
\langle x'\rangle \oplus e_8$, so there is an automorphism of
$A$ taking $x$ to $x'$, and hence $x/2m$ to $kx/2m$ in $A'/A$. This
implies that all the elements $a$ satisfying (2) are conjugate under
$\Aut(A)$, so $II_{8n+1,1}$ has only one orbit of primitive vectors of
norm $2m$. Q.E.D.

{\bf Remark.}  $I_{n,1}$ ($n\ge 4$) has (at least) two orbits of
primitive vectors of norm $m$ whenever $m\equiv n-1\bmod 8$.

\proclaim Chapter 4 25 dimensional lattices. 

In chapter 3 we found properties of negative norm vectors most of
which held for vectors in all Lorentzian lattices $II_{8n+1,1}$. In
this chapter we look in more detail at the vectors of $II_{25,1}$. In
sections 4.1 to 4.4 we look at vectors of norm $-2$ or $-4$ and find
several strange identities for the root systems of the lattices
corresponding to them, and use these to give several hundred
constructions of the Leech lattice. The section from 4.5 onwards give
a collection of isolated (and rather uninteresting) results on vectors
of larger negative norm. Finally in 4.11 we give an application of
some of these results to shallow holes of the Leech lattice.

\proclaim 4.1 The height in $II_{25,1}$. 

Unlike the lattices $II_{8n+1,1}$ for $8n\ge 32$, $II_{25,1}$ contains
a Weyl vector $w$, so we can define the height of a vector $u$ in
$II_{25,1}$ to be $-(u,w)$. In this section we show how to calculate
the heights of vectors of $II_{25,1}$ that have been found with the
algorithm in chapter 3. (In chapter 5 we will show that for vectors
$u$ in $D$ of fixed norm $-2$ or $-4$, and probably for all other
norms, the height depends linearly on the theta function of $u^\perp$
and it might be possible to find some substitute for the height like
this in $II_{8n+1,1}$ for $8n\ge 32$.)

{\bf Notation.} $D$ is a fundamental domain of the reflection group of
$II_{25,1}$ and $w$ is the norm 0 vector that has inner product $-1$
with all simple roots of $D$. The height of a vector $u$ is
$-(u,w)$. There are no roots of height 0, and a root is simple if and
only if it has height 1.

\proclaim Lemma 4.1.1. 
Suppose $u$, $v$ are vectors in $D$ of norms $-2n$, $-2(n-1)$ with
$-(u,v)=-2n$ (as in 3.5) and suppose that $v$ corresponds to the
component $C$ of $R_0(u)$ as in 3.5.2. Then
$$\height(u)=\height(v)+h-1$$
where $h$ is the Coxeter number of the component $C$. 

Proof. By 3.5.2 $v=r+u$ where $r$ is the highest root of $C$, so
$\height(u)=\height(v)+(r,w)$. $r=-\sum_im_ic_i$ where the $c_i$ are
the simple roots of $C$ with weights $m_i$ and $\sum_im_i=h-1$. All
the $c_i$ have inner product $-1$ with $w$, so $(r,w)=h-1$. Q.E.D.

\proclaim Lemma 4.1.2. 
Let $u$ be a primitive vector of $D$ of height 0 or 1, so that there
is a norm 0 vector $z$ with $(z,u)=0$ or $-1$, and suppose that $z$
corresponds to a Niemeier lattice $B$ with Coxeter number $h$.
\item{(a)} 
If $u$ has type 0 then its height is $h$. The Dynkin diagram of
$u^\perp$ is the extended Dynkin diagram of $B$.
\item{(b)}
If $u$ has type 1 then $\height(u)=1+(1-u^2/2)h$. The Dynkin diagram
of $u^\perp$ is the Dynkin diagram of $B$ if $u^2<-2$ and the Dynkin
diagram of $B$ plus an $a_1$ if $u^2=-2$.

Proof.
\item{(a)}
The Dynkin diagram of $u^\perp$ is a union of extended Dynkin
diagrams. If this union is empty then $u$ must be $w$ and therefore
has height $0=h$. If not then let $C$ be one of the
components. $u=\sum_im_ic_i$ where the $c_i$'s are the simple roots of
$C$ with weights $m_i$. $\sum_im_i=h$ because $C$ is an extended
Dynkin diagram and all the $c_i$'s have height 1, so $u$ has height
$h$.
\item{(b)}
By 3.4.2 $u=nz+z'$ with $u^2=-2n$, where $z$ is a norm 0 vector in
$D$. By part (a) $z$ has height $h$. By 3.4.1 $z'=z+r$ where $r$ is a
simple root of $D$, so $\height(z')=\height(z)+\height(r)=h+1$. Hence
$\height(u)=nh+h+1=1+(1-u^2/2)h$. The lattice $u^\perp$ is $B\oplus N$
where $N$ is a one dimensional lattice of determinant $2n$, so the
Dynkin diagram is that of $B$ plus that of $N$, and the Dynkin diagram
of (norm 2 roots of) $N$ is empty unless $2n=2$ in which case it is
$a_1$. This proves that the Dynkin diagram of $u^\perp$ is what it is
stated to be. Q.E.D.

\proclaim 4.2 The space $Z$. 

{\bf Notation.} $L$ is $II_{25,1}$, and $w$ is the Weyl vector of a
fundamental domain $D$ of $L$.

There are several correspondences between some vectors of $II_{25,1}$
and some points of $\Lambda\otimes \Q$. For example simple roots of
$II_{25,1}$ correspond to points of $\Lambda$, and primitive norm 0
vectors in $D$ other than $w$ correspond to deep holes of $\Lambda$
(which are points of $\Lambda\otimes Q$). In this section we
generalize these two correspondences to give a map from a larger class
of vectors of $II_{25,1}$ to points of $\Lambda\otimes \Q$, or more
precisely to the set $\Lambda\otimes \Q\cup \infty$ which can be
identified with the non-zero isotropic subspaces of $II_{25,1}$.

We define a space $Z$ as follows:

Points of $Z$ are maximal isotropic subspaces of $L$, in other words
pairs $(z,-z)$ of primitive norm 0 vectors. For the sake of confusion
we will use the same letter $z$ for a norm 0 vector of $L$ and for the
isotropic sublattice corresponding to it. We define the distance
between two points by
$$\eqalign{
d(w,w)&=0,\cr
d(z,w)=d(w,z)&=\infty\hbox{ if $z$ is not $w$},\cr
d(z_1,z_2)&=(z_1/\height (z_1)-z_2/\height(z_2))^2\cr
&= -2(z_1,z_2)/\height(z_1)\height(z_2)\cr
&\hbox{ if neither $z_1$ nor $z_2$ is $w$}.\cr
}$$

\proclaim Lemma 4.2.1. 
$Z$ is isometric to the affine rational Leech lattice together with a
point at infinity. This isometry is compatible with the action of
$\cdot\infty$ on $Z$ and the affine Leech lattice. ($\cdot\infty$ acts
on $Z$, since it is the subgroup of $\Aut(L)$ fixing $w$.)

Proof. The space $Z$ can be identified with $\Lambda\otimes \Q\cup
\infty$ as in section 1.8. If we choose coordinates $(\lambda,m,n)$
for $II_{25,1}=\Lambda\oplus U$ so that $w=(0,0,1)$ then the norm 0
vector $(\lambda, m,n)$ is identified with $\lambda/m$ in
$\Lambda\otimes \Q\cup \infty$. The formula for $d(z_1,z_2)$ follows
from this by an easy calculation. Q.E.D.

If $u$ is any point of $L$ with $u^2\le 0$ or $(u,w)\ne 0$ then there is
a unique non-zero isotropic subspace of $L\otimes \Q$ containing a
non-zero point of the form $u+nw$ for some rational $n$. ($n$ is
uniquely determined if $u$ is not a multiple of $w$.) This gives a map
$Z$ from such vectors $u$ to the space $Z$ such that $Z(u)=\infty$ if
and only if $u$ is a multiple of $w$. (In the coordinates of 4.2.1
$Z((\lambda,m,n))=\lambda/m$ even when $(\lambda,m,n)$ does not have
norm 0.) The images $Z(r)$ of the simple roots $r$ of $II_{25,1}$ form
a copy of the Leech lattice in $Z$.

\proclaim Theorem 4.2.2. 
There are natural 1:1 correspondences between the following three
sets:
\item{(1)}
Points of the space $Z$, in other words non-zero isotropic subspaces
of $II_{25,1} $.
\item{(2)}
The vector $w$ together with all primitive vectors $u$ of $D$ such
that $u^\perp$ contains a root of $II_{25,1}$.
\item{(3)}
Rational (finite or infinite) points on the boundary of $D$ when $D$
is considered as a subset of hyperbolic space.
\item{}
(The map between (1) and (3) is not smooth, because $Z$ looks like 24
dimensional space with a point at infinity, while the boundary of $D$
has corners.)

Proof. We first check that (2) and (3) are the same. Any point $d$ on
the boundary of $D$ defines a vector $u$ of $II_{25,1}$ where $u$ is
the primitive vector of $D$ representing $d$. ($u$ has norm 0 if and
only if $d$ is a point at infinity.) $d$ is on the boundary of $D$ if
and only if it has norm 0 or lies on some hyperplane of the boundary
of $D$, in other words if and only if $u$ has norm 0 or $u^\perp$
contains a root of $II_{25,1}$. But if $u$ has norm 0 then either
$u^\perp$ contains a root of $u=w$, so (2) and (3) are the same.

To show that (1) and (3) are the same we ``project from $w$'', in other
words a non-zero isotropic subspace of $II_{25,1}$ represented by an
infinite point $z$ of hyperbolic space corresponds to a point $d$ on
the boundary of $D$ if and only if $z=w=d$ or if $z\ne w\ne d$ and
$z$, $w$, and $d$ are co-linear. It is easy to check that this gives a
1:1 correspondence between (1) and (3) using the fact that $D$ is
convex. Q.E.D.

\proclaim 4.3 Norm $-2$ vectors. 

{\bf Notation. } $u$ is a norm $-2$ vector in the fundamental domain
$D$ of $II_{25,1}$. $D$ has Weyl vector $w$ Recall that orbits of such
vectors $u$ correspond to even 25 dimensional bimodular lattices $B$,
where $B\cong u^\perp$. Roots $r$ of $u^\perp$ correspond to norm 0 vectors
$z$ with $(z,u)=-2$ by $r=\tau_u(v)$. $r$ is simple if and only if $z$
has inner product $\le 0$ with all simple roots of $u^\perp$, and if
$u$ has type $\ge 2$ this is equivalent to saying that $z$ is in $D$.

In this section we find some properties of norm $-2$ vectors $u$ of
$II_{25,1}$ and the lattices $B$ corresponding to them. Such vectors
have the special property that reflection in $u^\perp$ is an
automorphism of $II_{25,1}$, and using this we can construct an
automorphism of $D$ for each vector $u$. This implies the important
result that the projection of $w$ into $u^\perp$ is the Weyl vector of
$u^\perp$ (4.3.3) and several strange properties of the root system of
$u^\perp$ follow from this. (4.3.3 is almost the only result of this
chapter that is used later.)

\proclaim Lemma 4.3.1. 
$u^\perp$ contains roots, or in other words every 25 dimensional even
bimodular lattice has a root.

Proof. If $u^\perp$ contains no roots then $u=w+u_1$ for some $u_1$ in
$D$. We have $u_1^2=u^2+2\height (u)$, so $u_1^2=0$ and $u$ has height
1 because $u_1^2\le 0$, $u^2=-2$ and the height of $u$ is
positive. Then $\height(u_1)=\height(u)=1$, so $u$ is a norm 0 vector
in $D$ that has inner product $-1$ with the norm 0 vector $w$ of $D$,
but this is impossible by 3.4.1. Q.E.D.

\proclaim Corollary 4.3.2. $u$ has type 1 or 2. 

This follows from 4.3.1 and the correspondence between roots of
$u^\perp$ and norm 0 vectors that have inner product $-2$ with
$u$. Q.E.D.

\proclaim Theorem 4.3.3. 
The projection $\rho$ of $w$ into $u^\perp$ is the Weyl vector of $u^\perp$. 

Proof. Recall from 1.5.1 that there is an automorphism $\sigma_u$ of
$D$ fixing $u$ and acting as $\sigma(u^\perp)$ on
$u^\perp$. $\sigma_u$ fixes $w$ because it fixes $D$, so it also fixes
$\rho$. On the subspace of $u^\perp$ orthogonal to all roots of
$u^\perp$, $\sigma_u$ acts as $-1$, so any element of $u^\perp$ fixed
by $\sigma_u$, and in particular $\rho$, must be in the space
generated by roots of $u^\perp$. However $\rho$ also has inner product
$-1$ with all simple roots of $u^\perp$ because $w$ does, so it is the
Weyl vector of $u^\perp$. Q.E.D.

\proclaim Corollary 4.3.4. If $\rho$ is the Weyl vector of $u^\perp$
then $2\rho^2=\height(u)^2$. 

Proof. By 4.3.3 $\rho$ is the projection of $w$ into $u^\perp$, so
$$\eqalign{
\rho^2&=(w-(w,u)u/(u,u))^2\cr
&=w^2+(w,u)^2+(w,u)^2u^2/4 \hbox{ because } (u,u)=-2\cr
&= (w,u)^2/2\hbox{ because } w^2=0\cr
&=\height(u)^2/2.\cr
}$$
Q.E.D.

We find all the norm $-2$ vectors (which are listed in table $-2$) by
taking all the norm 0 vectors of $D$ and applying the algorithm of
chapter 3, and in particular 3.5.4, to them.

{\bf Example.}  Let $z$ be a primitive norm 0 vector of $D$ with
Dynkin diagram $D_{24}$. The simple roots of $D$ that have inner
product 0, $-1$, or $-2$ with $z$ are
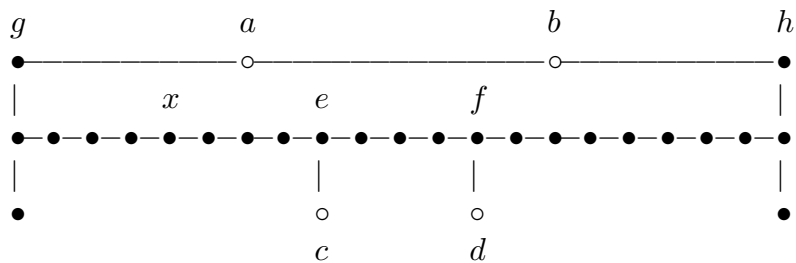
\halign
{
$#$&&$#$\cr
\cr
g&&&&&&&&&&&&a&&&&&&&&&&&&&&&&b&&&&&&&&&&&&h\cr
\bullet&\!-&\!-&\!-&\!-&\!-&\!-&\!-&\!-&\!-&\!-&\!-&\circ
&\!-&\!-&\!-&\!-&\!-&\!-&\!-&\!-&\!-&\!-&\!-&\!-&\!-&\!-&\!-&\circ
&\!-&\!-&\!-&\!-&\!-&\!-&\!-&\!-&\!-&\!-&\!-&\bullet\cr
|&&&&&&&&x&&&&&&&&e&&&&&&&&f&&&&&&&&&&&&&&&&|\cr
\bullet&\!-&\bullet&\!-&\bullet&\!-&\bullet&\!-&\bullet&\!-&\bullet&\!-&
\bullet&\!-&\bullet&\!-&\bullet&\!-&\bullet&\!-&\bullet&\!-&\bullet&\!-&
\bullet&\!-&\bullet&\!-&\bullet&\!-&\bullet&\!-&\bullet&\!-&\bullet&\!-&
\bullet&\!-&\bullet&\!-&\bullet\cr
|&&&&&&&&&&&&&&&&|&&&&&&&&|&&&&&&&&&&&&&&&&|\cr
\bullet&&&&&&&&&&&&&&&&\circ&&&&&&&&\circ&&&&&&&&&&&&&&&&\bullet\cr
&&&&&&&&&&&&&&&&c&&&&&&&&d&&&&&&&&&&&&&&&&
\cr}

We obtain 3 vectors $u$ in $D$ of norm $-2$ form it as follows:
\item{(1)}
Add $c$ and delete $e$. The Dynkin diagram of $u$ is $a_1d_{10}d_{14}$
and the height of $u$ is $\height(z)+$ Coxeter number of $a_1 -1=46+2-1=47$.
\item{(2)} 
Add $a$ and $b$ and delete $g$ and $h$. The Dynkin diagram of $u$ is
$a_2a_{23}$ and its height is 48.
\item{(3)}
Add $a$ and delete $c$ and $x$. The Dynkin diagram of $u$ is
$e_7d_{18}$ and its height is 63.  

{\bf Remark.} The Leech lattice can be constructed from any 25
dimensional even bimodular lattice $B$. We form the lattice generated
by $B\oplus U$ and $b'\oplus u/2$, where $u$ is a one dimensional
lattice generated by $u$ of norm $-2$ and $b$ is in $B'-B$. This
lattice is isomorphic to $II_{25,1}$ and the vector
$w=\rho+\sqrt{\rho^2/2}u$ is in the lattice and is a primitive norm 0
vector corresponding to $\Lambda$, so $\Lambda$ is isomorphic to
$w^\perp/w$.

Given a norm $-2$ vector $u$ we may wish to find the norm $-4$ vectors
that can be obtained from it as in 3.5, and to do this we need to know
the simple roots of $D$ that have inner product $-1$ with $u$. These
can be found from the last column of table $-2$ as follows.  \proclaim
Lemma 4.3.5. Suppose that $u$ has type 2 and does not have height
2. If $r$ is a simple root of $D$ with $-(r,u)=1$ then $r$ is in
$z^\perp$ for one of the norm 0 vectors $z$ of $D$ with $-(z,u)=2$.

Proof. $(r,w)=(r,u)=-1$ and $\rho=w-\height(u)u/2$ by 4.3.3, so
$(r,\rho)$ is not 0 because $u$ does not have height 2. Therefore
there is some highest root $r'$ of $u^\perp$ with $(r,r')\ne 0$. It is
easy to show that for any such $r'$ we have $0\le (r,r')\le \sqrt 5$,
and $(r,r')$ cannot be 2 because then $r-r'$ would be a norm 0 vector
having inner product $-1$ with $u$, so $(r,r')=1$. Hence $r$ is in
$z'$ where $z$ is the norm 0 vector $r'+u$ of $D$. Q.E.D.

{\bf Remark.} There is a curious ``duality'' for some norm $-2$ vectors. If the norm $-2$ vector $u$ has Dynkin diagram $X$ then there is often another norm $-2$ vector $u'$ with Dynkin diagram $X'$ where $X'$ is obtained from $X$ by making some of the following changes to it:
$$a_{2n-1}\leftrightarrow a_{2n},
\quad a_1^4\leftrightarrow d_4\leftrightarrow a_2^4,
\quad a_8\leftrightarrow e_6\leftrightarrow d_5\leftrightarrow a_4^2.$$
We have $S(u^\perp)=S(u'^\perp)$ and the subgroups of $\Aut(D)$ fixing
$u$ and $u'$ are isomorphic. For example there are norm $-2$ vectors
with Dynkin diagrams $a_3a_2^5a_1^6$ and $a_4a_2^6a_1^5$, or
$a_8a_5^3$ and $e_6a_6^3$. This duality seems to be connected with the
existence of norm 0 vectors $z$ having inner product $-2$ with $u$,
such that the Cartan matrix of the corresponding Niemeier lattice has odd
determinant.

{\bf Remark.} The automorphism $\sigma_u$ of $II_{25,1}$ fixes $w$ and
so can be considered as an element of $\cdot\infty$. In particular it
determines a conjugacy class of elements of period 2 in $\cdot
0$. $\cdot 0$ has 5 conjugacy classes of elements of order 1 or 2;
they are characterized by the dimension of their fixed subspace which
can be 0, 8, 12, 16, or 24 (the possible weights of words in the
binary Golay code of the Steiner system $S(5,8,24)$). If the elements
of $II_{25,1}$ fixed by $\sigma(u) $ form an $n$-dimensional subspace
then the subspace of elements fixed by the corresponding element of
$\cdot0$ is $n-2$-dimensional. $\sigma_u=\sigma(u^\perp)r_u$, so the
dimension of the space of elements fixed by $\sigma_u$ is 1 plus the
dimension of the subspace of $u^\perp$ fixed by $\sigma(u^\perp)$,
which is $1+S(u^\perp)$. Hence if $B$ is any 25 dimensional even
bimodular lattice then $S(B)$ must be 1, 9, 13, 17, or 25, and in
particular if $N$ is a Niemeier lattice then by looking at $B=N\oplus
a_1$ we see that $S(N)$ must be 0, 8, 12, 16, or 24. (The case
$S(N)=8$ does not occur because $S(N)\ge rank(N)/2$ and $N$ has rank
0 or 24.)

\proclaim 4.4 Norm $-4$ vectors. 

{\bf Notation.} In the next few sections $u$ is a norm $-4$ vector in
the fundamental domain $D$ of $II_{25,1}$. Such vectors correspond to
25 dimensional unimodular lattices $A$, where $u^\perp$ is the lattice
of even elements of $A$. The odd vectors of $A$ can be taken as the
projections of the vectors $y$ with $-(y,u)=2$ into $u^\perp$.

Such a vector $u$ can behave in 4 different ways, depending on whether
the unimodular lattice $A_1$ with no norm 1 vectors corresponding to
$u$ is at most 23 dimensional, or 24 dimensional and odd, or 24
dimensional and even, or 25 dimensional. (4.4.1.) Reflection in
$u^\perp$ is not an automorphism of $L$, but if $u$ has type at most 2
(in other words $A_1$ has dimension at most 24) then we can still
construct an automorphism of $D$ fixing $u$ as in 3.7, and this allows
us to prove identities for the root systems of such lattices similar
to those for the lattices of norm $-2$ vectors. We will use these
identities to give a construction of the Leech lattice for each
unimodular lattice of dimension at most 23; the case of the 0
dimensional lattice was found by Curtis, Conway and Sloane. The list
of the 665 orbits of norm $-4$ vectors appears in table $-4$ and the
unimodular lattices of dimension at most 25 can be read off from
this. In this section we examine the norm $-4$ vectors that correspond
to lattices $A_1$ of dimension at most 23.

\proclaim Theorem 4.4.1. 
Norm 1 vectors of $A$ correspond to norm 0 vectors $z$ of $II_{25,1}$
with $-(z,u)=-2$. Write $A=A_1\oplus I^n$ where $A_1$ has no vectors
of norm 1. Then $u$ is in exactly one of the following four classes:
\item{(A)} $u$ has type 1. $A_1$ is a Niemeier lattice. 
\item{(B)}
$u$ has type 2 and $A$ has at least 4 vectors of norm 1, so that $A_1$
is at most 23 dimensional (but may be even). There is a unique norm 0
vector $z$ of $D$ with $(z,u)=-2$ and this vector $z$ is of the same
type as either of the two even neighbors of $A_1\oplus I^{n-1}$.
\item{(C)}
$U$ has type 2 and $A$ has exactly two vectors of norm 1, so $A_1$ is
24 dimensional and odd. There are exactly two norm 0 vectors that have
inner product $-2$ with $u$, and they are both in $D$. They have the
types of the two even neighbors of $A_1$.
\item{(D)}
$u$ has type at least 3 and $A=A_1$ has no vectors
of norm 1. (We will later show that $u$ has type at most 3.)

Proof. $z$ is a norm 0 vector with $(z,u)=-2$ if and only if $u/2-z$
is a norm 1 vector of $A$. Most of 4.4.1 follows from this, and part
(A) follows from 3.4. The only things to check are the statements
about norm 0 vectors  that are in $D$.

If $u$ has type $\ge 2$ then by 3.3.1 a norm 0 vector $z$ with
$-(z,u)=-2$ is in $D$ if and only if it has inner product $\le 0$ with
all simple roots of $u^\perp$, so there is one such vector in $D$ for
each orbit of such norm 0 vectors under the reflection group of
$u^\perp$. If $A$ has at least 4 vectors of norm 1 then they form a
single orbit under the Weyl group of (the norm 2 vectors of)
$u^\perp$, which proves (B), while if $A$ has only two vectors of norm
1 then they are both orthogonal to all norm 2 vectors of $A$ and so
form two orbits under they Weyl group of $u^\perp$. Q.E.D.

In the rest of this section we consider the norm $-4$ vectors $u$
corresponding to lattices $A_1$ of dimension at most 23. We let $A_1$
have dimension $25-n$, so $A=A_1\oplus I^n$ and $n\ge 2$. The (norm 2)
root system of $u^\perp$ is $Xd_n$ where $X$ is the Dynkin diagram of
$A_1$. We let $z$ be the unique norm 0 vector of $D$ with $-(z,u)=2$,
write $i$ for the projection of $z$ into $u^\perp$ so $i$ is the norm
1 vector of $A$ in the Weyl chamber of $u^\perp$. $z$ corresponds to
the even neighbors of $A_1\oplus I^{n-1}$; we let $h$ be their Coxeter
number which by 4.1.2 is equal to $-(z,w)$. Write $\rho$ for the Weyl
vector of $u^\perp$.

\proclaim Theorem 4.4.2. 
\item{(1)} 
The projection of $w$ into $u^\perp$ is $\rho$, and
$\rho^2=(h+n-1)^2$.
\item{(2)}
$\height(u)=-(w,u)=2(h+n-1)$. 
\item{(3)}
$A_1$ has  $16(h-1)-2(n-1)(n-10) $ roots. In particular the
number of roots of a unimodular lattice of dimension at most 23 is
determined mod 16 by its dimension.

Proof. $A$ has at least 4 vectors of norm 1, so any vector of norm 1
and in particular $i$ is in the vector space generated by vectors of
norm 2. Hence by 3.7.1 and the same argument as in 4.3.3 the
projection $w-\height(u)u/4$ of $w$ into $u^\perp$ is $\rho$. The norm
4 vector $-2i$ of $A$ is the sum of $2(n-1)$ simple roots of the $d_n$
component of the Dynkin diagram of $u^\perp$, so
$(-2i,w)=(-2i,\rho)=-2(n-1)$.

$i$ is the projection of $z$ into $u^\perp$, so $i=z-u/2$, and hence
$$\eqalign{
\height(u)&=-(w,u)\cr
&= 2(w,i-z)\cr
&= 2(\height(z)+n-1)\cr
&=2(h+n-1).\cr
}$$
If we calculate the norms of both sides of $\rho=w-\height(u)u/4$ we
find that $\rho^2=(h+n-1)^2$, which proves (1) and (2). We will prove
later (5.4.1) that a 25 dimensional unimodular lattice with $2n$ norm 1
vectors corresponding to the norm $-4$ vector $u$ has
$8\height(u)-20+4n$ norm 2 vectors, and (c) follows from this and the
fact that $d_n$ has $2n(n-1)$ norm 2 vectors. Q.E.D.

This theorem is not trivial even if $A_1$ is the 0 dimensional
lattice. $n$ is then 25 and $\rho$ can be taken as
$(0,1,2,\ldots,24)$. The even neighbors of $I^{24}$ are both $D_{24}$
with $h=46$, so we find that
$\rho^2=0^2+1^2+2^2+\cdots+24^2=(46+25-1)^2$. Watson showed that the
only solution of $0^2+1^2+\cdots+k^2=m^2$ with $k\ge 2$ is $k=24$. If
$A$ is a 25 dimensional lattice with 0 or 2 vectors of norm 1 then its
Weyl vector does not usually have square norm.

For every unimodular lattice $A_1$ of dimension at most 23 there is a
construction of the Leech lattice. The lattice $A=A_1\oplus I^n$
(where $n=25-\dim(A_1)$) is 25 dimensional and its Weyl vector $\rho$
has square norm $(h+n-1)^2$. We form the lattice $L=A\oplus(-I)$. The
vector $w=(\rho,h+n-1)$ is a norm 0 vector of $L$ and in either of the
two even neighbors of $L$ (which are isomorphic to $II_{25,1}$) $w$ is
a norm 0 vector corresponding to the Leech lattice.

For example, if $A_1$ is the 0 dimensional lattice then $w$ is
$(0,1,\ldots, 24|70)$ and this is the construction in [C-S c]. They
chose the vector $w$ because it was lexicographically the first one
which did not obviously have a root orthogonal to it, and showed that
$w$ corresponded to the Leech lattice by writing down an explicit
basis for $\Lambda$ in $w^\perp$.

{\bf Remark.} If $u$ is any norm $-4$ vector of type $\le 2$ then we
can find restrictions on $S(u^\perp)$ as in the last paragraph of
4.3. $S(u^\perp)$ must be 0, 8, 12, 16, or 24 if the unimodular
lattice $A_1$ corresponding to $u$ has even dimension, and must be 2,
10, 14, or 18 otherwise. (10 does not occur. If $u$ has type 3 then
$S(u^\perp)$ is 2, 6, 8, 10, 12, 14, or 18 but I do not know a uniform
proof of this.)

\proclaim 4.5 24 dimensional unimodular lattices. 

{\bf Notation.} We deal with the third case of 4.4.1, so $u$ is a norm
$-4$ vector of $D$ with exactly two norm 0 vectors $z_1$, $z_2$ that
have inner product $-2$ with $u$. $z_1$ and $z_2$ are both in $D$ and
have Coxeter numbers $h_1$, $h_2$ where $h_i=-(z_i,w)$. The Niemeier
lattices corresponding to $z_1$ and $z_2$ are the two even neighbors
of the 24 dimensional unimodular lattice corresponding to $u$ and
their Coxeter numbers are also $h_1$ and $h_2$.  
We write $\rho$ for the Weyl vector of $u^\perp$. 

We give a formula for the Weyl vector of the 24 dimensional unimodular
lattice $A$ corresponding to $u$, then show that a Niemeier lattice
$B$ has $S(B)=24$ if and only if it has a hole of radius $\sqrt{3}$
that is halfway between two lattice points, and finally describe how
the simple roots that have small inner product with $u$ are arranged.

\proclaim Theorem 4.5.1. 
If $w_u$ is the projection of $w$ into $u^\perp$ then
$$w_u=\rho-(h_1-h_2)(z_1-z_2)/4.$$ 
$u=z_1+z_2$, $\height(u)=h_1+h_2$,
and $\rho^2=h_1h_2$. $u^\perp$ has $8(h_1+h_2-2)$ roots.

Proof. $u-z_1$ is a norm 0 vector which has inner product $-2$ with
$u$ and so must be $z_2$. Hence $u=z_1+z_2$ and
$\height(u)=\height(z_1)+\height(z_2)=h_1+h_2$.

By 3.7.1 $w_u-\rho$ is a multiple of the projection of $z_1$ into
$u^\perp$, so $w_u=\rho+k(z_1-z_2)$ for some $k$. Also
$w_u=w-(w,u)u/(u,u)=w-(h_1+h_2)(z_1+z_2)/4$ as $w_u$ is the projection
of $w$ into $u^\perp$. Taking the inner product of both sides of
$$\rho+k(z_1-z_2)=w-(h_1+h_2)(z_1+z_2)/4$$ with $z_1$ shows that $k$
is $(h_2-h_1)/4$ which proves the formula for $w_u$. If we calculate
the norms of both sides of this equality and simplify we find that
$\rho^2=h_1h_2$. The formula for the number of roots of $u^\perp$ is
proved as in 4.4.2. Q.E.D.

\proclaim Corollary 4.5.2. 
If $A$ is an odd 24 dimensional unimodular lattice with no vectors of
norm 1 and whose even neighbors have Coxeter numbers $h_1$ and $h_2$,
then $\rho^2=h_1h_2$ where $\rho$ is the Weyl vector of $A$, and $A$
has $8(h_1+h_2-2)$ roots.

Proof. See 4.5.1. Q.E.D. 

{\bf Remark.} The condition $\rho^2=h_1h_2$ nearly characterizes the
root systems of these lattices $A$. In fact if $B_1$ and $B_2$ are two
{\sl non-isomorphic} Niemeier lattices with Coxeter numbers $h_1$ and
$h_2$ then there is a 1:1 correspondence between 24 dimensional
lattices with neighbors $B_1$ and $B_2$, and root systems satisfying
$\rho^2=h_1h_2$ and that are isomorphic to root systems of index 1 or
2 in the root systems of $B_1$ and $B_2$.

\proclaim Proposition 4.5.3.  
Notation is as in 4.5.2. Let $B_1$, $B_2$ be the two even neighbors
of $A$. Then $h_2\le 2h_1+2$ and the following conditions are
equivalent:
\item{(1)} $h_2=2h_1+2$
\item{(2)}
The element $x$ of $B_2/2B_2$ corresponding to $A$ has $n(x)\ge
12$. ($n(x)$ is the minimum norm of an element of $B_2$ representing
$x$. In fat $n(x)$ is always $\le 12$ for Niemeier lattices.)
\item{(3)}
$A$ has the same root system as $B_1$. 
\item{}
Any of these conditions imply that $S(B_2)=24$. Conversely any
Niemeier lattice $B_2$ with $S(B_2)=24$ has a unique orbit of vectors
$x$ in $B_2/2B_2$ with $n(x)=12$.

Proof. (2) and (3) are equivalent because roots of $B_1$ not in $A$
correspond to points of $B_2$ whose distance from $x/2$ is
$\sqrt{2}$. The root system of any odd neighbor $A$ of $B_1$ is a
subset of the root system of $B_1$, so $\rho^2(A)\le \rho^2(B_1)$ with
equality if and only if $A$ and $B_1$ have the same root system (where
$\rho^2(L)$ is the norm of the Weyl vector of the norm 2 vectors of
$L$). If $h_2\ge 2h_1+2$ then
$$ \rho^2(A)=h_1h_2\ge 2h_1(h_1+1)=\rho^2(B_1)\ge \rho^2(A)$$ so
$h_2=2h_1+2$ and $A$ and $B_1$ have the same root system. Conversely
if $A$ and $B_1$ have the same root system then it follows easily that
$h_2=2h_1+2$ so (1) and (3) are equivalent.

If a simple root system $R$ of Coxeter number $2h+2$ and rank $n$ has
a sub-root system all of whose components have Coxeter number $h$ and
whose rank is still $n$ then $R$ must be $a_1$, $d_{2n}$, $e_7$ or
$e_8$ and all these have $S(R)=n$. This shows that (1) implies that
$S(B_2)=24$. Conversely observation shows that any of the Niemeier
lattices $N$ with $S(N)=24$ has a unique orbit of vectors $x$ in
$N/2N$ with $n(x)=12$. (It is probably easy to prove this without
using the classification of Niemeier lattices.) Q.E.D.

{\bf Remark.} The Niemeier lattices $B_2$ with $S(B_2)=24$ are
$A_1^{24}(=D_2^{12})$, $D_4^6$, $D_6^4$, $D_8^3$, $D_{12}^2$, $D_{24}$,
$D_{10}E_7^2$, $E_8^3$, and $D_{16}E_8$. The corresponding lattices
$A$ have Dynkin diagrams $\emptyset(=d_1^{24})$, $a_1^{24}(=d_2^{12})$,
$a_3^8(=d_3^8)$, $d_4^6$, $d_6^4$, $d_{12}^2$, $a_7^2d_5^2$, $d_8^3$,
and $d_8^3$. These are the only odd 24 dimensional unimodular lattices
that have the same root system as a Niemeier lattice.

We now describe the simple roots of $D$ that have inner product 0,
$-1$, or $-2$ with $u$; write $R_0$, $R_1$, $R_2$ for these sets. Let
$Z_i$ be the set of simple roots in $z_i^\perp$, let $r$ be a root of
$Z_1$ nor in $R_0$, and let $C$ be the component of $Z_1$ containing
$r$. There are three possibilities:
\item{(1)}
$r$ is  a point of weight 1 of $C$ and all other points of $C$ are in $R_0$. 
\item{(2)}
$r$ is a point of weight 2 of $C$ and all other points of $C$ are in $R_0$. 
\item{(3)}
There are two points of $C$ not in $R_0$, both of which have weight 1
and one of which is $r$.

We will say that $r$ has weight 1, 2, or 11 and total weight $m=1,2,2$
in these three cases.

\proclaim Lemma 4.5.4. $(r,u)=(r,z_2)=-2/m$. 

Proof. By 1.4.3 $(r,\rho)=h_1/m-1$ where $\rho$ is the Weyl vector of
$u^\perp$. We also have $u=z_1+z_2$, $(r,w)=-1$, $(r,z_1)=0$, and by
4.5.1 $w-\rho=(h_2z_1-h_1z_2)/2$, so taking inner products of this
with $r$ we find $-h_1/m=h_1(r,z_2)/2$. Q.E.D.

In particular if a simple root $r$ of $D$ is in $Z_1$ or $Z_2$ then it
has inner product 0, $-1$, or $-2$ with $u$. Conversely if a simple
root has inner product 0 or $-1$ with $u$ then it has inner product 0
with either $z_1$ or $z_2$.

Now let $r_i$ be simple roots in $Z_i$. If $r_1$ and $r_2$ are both
joined to some component of $R_0$ then they are not connected in the
Dynkin diagram of $D$, otherwise the number of bonds between them is
given in the following table.
\halign{#\hfill&&~~~\hfill#\hfill\cr
&$r_2$ of weight 1&$r_2$ of weight 11&$r_2$ of weight 2\cr
$r_1$ of weight 1&2&1&1\cr
$r_1$ of weight 11&1&0 or 1&0\cr
$r_1$ of weight 2&1&0&0\cr
}
The ambiguity in the middle is because if $r_1$ and $r_2$ are both of
weight 11 then there are other points $r_1'$ and $r_2'$ of weight 11
in the same components of $Z_1$ and $Z_2$ as $r_1$ and $r_2$, and each
of $r_1$ and $r_1'$ is joined to exactly one of $r_2$ and $r_2'$.

{\bf Remark.} If $u$ is a norm $-4$ vector of $D$ of type 3 then by
3.8.1 the simple roots that have inner product $-1$ with $u$ are in
1:1 correspondence with the norm 0 vectors of $D$ that have inner
product $-3$ with $u$. This is not true if $u$ has type 2; in fact
there are no norm 0 vectors $z$ of $D$ that have inner product $-3$
with $u$ because $u=z_1+z_2$ for some norm 0 vectors $z_i$ and $z$
must have inner product $-1$ with some $z_i$ and is hence conjugate to
it.

\proclaim 4.6 25 dimensional lattices. 

{\bf Notation.} $u$ is a norm $-4$ vector in $D$ of type $\ge
3$. $u^\perp$ is the sublattice of even vectors of a strictly 25
dimensional unimodular lattice $A$. We write $t$ for the height of $u$
and $\rho$ for the Weyl vector of $A$.

We continue sections 4.4 and 4.5 by examining norm $-4$ vectors of
type $\ge 3$. There do not seem to be many nice properties of these
vectors, so instead we give several examples of calculations with them
to illustrate some of the lemmas in 3.8.

\proclaim Theorem 4.6.1.
\item{(1)}
$A$ has $8t-20$ roots and there are $8t+80$ norm 0 vectors that have
inner product 3 with $u$. (In particular $A$ has roots and $u$ has
type exactly $3$.)
\item{(2)}
$4\rho^2\le t^2$. (Equality holds for 17 vectors $u$, with heights 4, 6,
8, 10, 10, 10, 10, 14, 14, 18, 22, 22, 22, 34, 34, 58.)
\item{(3)}
$S(A)$ is even, and divisible by 4 if and only if $\rho^2$ is not an
integer. If $A$ has a component $d_n$ ($n\ge 4$) $e_6$, $e_7$, or
$e_8$ then $S(A)$ is 10, 14, or 18.

Proof. (1) follows from 5.4.1 and 5.4.3. The projection of $w$ into
$u^\perp$ is equal to $\rho$ plus something orthogonal to all roots of
$u^\perp$, so its norm $(w,u)/u^2$ is at least $\rho^2$. This implies
(2). By 1.3.6, $S(A)+4\rho^2 +({\rm number~ of ~ roots~ of~ A})/2
\equiv 0\bmod 4$ and the number of roots of $A$ is $4\bmod 8$ by (1),
so $S(A)$ is even, and divisible by 4 if and only if $\rho^2$ is not
an integer. If $A$ has a component that is not an $a_n$, then by 3.8.7
$S(A)=1+S(v^\perp)$ for some norm $-2$ vector $v$, and $S(v^\perp)$ is
1, 9, 13, 17, or 25. $S(A)$ cannot be 2 or 26, so it is 10, 14, or
18. Q.E.D.

{\bf Observation 4.6.2.} $S(A)$ is always 2, 6, 8, 10, 12, 14, or
18. These are the exponents of the $e_7$ root system. If $S(A)=18$
then the Dynkin diagram of $A$ contains no $a_{2n}$'s.

We find all the 368 norm $-4$ vectors of type 3 by taking each of the
97 norm $-2$ vectors of type 2 in turn and applying 3.8.6 to it.

{\bf Example 4.6.2.} Let $v$ be the norm $-2$ vector in $D$ of height
36 such that $v^\perp$ has Dynkin diagram $a_{18}e_6$. The roots that
have inner product 0, $-1$, or $-2$ with $v$ are
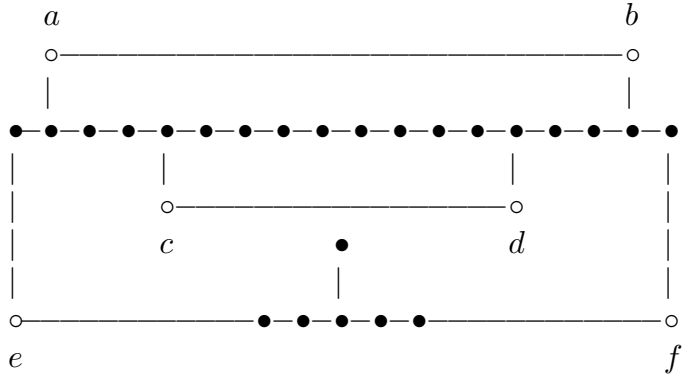
\halign
{
$#$&&$#$\cr
&&a&&&&&&&&&&&&&&&&&&&&&&&&&&&&&&b&&\cr
&&\circ&\!-&\!-&\!-&\!-&\!-&\!-&\!-&\!-&\!-&\!-&\!-&\!-
&\!-&\!-&\!-&\!-&\!-&\!-&\!-&\!-&\!-&\!-&\!-&\!-&\!-&\!-&\!-&\!-&\!-&\circ
\cr
&&|&&&&&&&&&&&&&&&&&&&&&&&&&&&&&&|&&\cr
\bullet&\!-&\bullet&\!-&\bullet&\!-&\bullet&\!-&\bullet&\!-&\bullet&\!-&
\bullet&\!-&\bullet&\!-&\bullet&\!-&\bullet&\!-&\bullet&\!-&\bullet&\!-&
\bullet&\!-&\bullet&\!-&\bullet&\!-&\bullet&\!-&\bullet&\!-&\bullet\cr
|&&&&&&&&|&&&&&&&&&&&&&&&&&&|&&&&&&&&|\cr
|&&&&&&&&\circ&\!-&\!-&\!-&\!-&\!-&\!-&\!-&\!-&\!-&\!-&\!-&\!-&\!-&\!-&\!-&\!-&\!-&\circ&&&&&&&&|\cr
|&&&&&&&&c&&&&&&&&&\bullet&&&&&&&&&d&&&&&&&&|\cr
|&&&&&&&&&&&&&&&&&|      &&&&&&&&&&&&&&&&&|\cr
\circ&\!-&\!-&\!-&\!-&\!-&\!-&\!-&\!-&\!-&\!-&\!-&\!-&\bullet&\!-&\bullet&\!-&\bullet&\!-&\bullet&\!-&\bullet&\!-&\!-&\!-&\!-&\!-&\!-&\!-&\!-&\!-&\!-&\!-&\!-&\circ\cr
e&&&&&&&&&&&&&&&&&&&&&&&&&&&&&&&&&&f\cr
\cr}
In addition there should be lines from $c$ to $f$ and from $d$ to $e$, 
which I cannot be bothered to make \TeX  draw. 
$c$ and $d$ have inner product $-2$ with $v$, while $a$, $b$, $c$, $d$
have inner product $-1$ with $v$. The deep hole of the $a_{18}$
component is $A_{17}E_7$ and contains $a$, $b$, $e$, and $f$, while
the deep hole of the $e_6$ component is $A_{24}$ and contains $e$ and
$f$. The automorphism group has order 2 and acts in the obvious way.

There are 4 strictly 25 dimensional unimodular lattices that can be
obtained from this as in 3.8.6:
\item{(1)}
$d_7a_{16}$ Add the point $e$ and delete a point from the $a_{18}$ and $e_6$ components. The height is $36+($Coxeter number of $d_7) -1=47$.
\item{(2)} $a_{19}d_5$. Add $b$ and $e$ and delete two other points. 
The height is 55.
\item{(3)} $a_4a_{15}d_5$. Add $a$ and $e$ and delete two other points. 
The height is 40.
\item{(4)} $a_1a_4e_6a_{13}$. Add $c$ and delete a point from the $a_{18}$.
The height is 37.

{\bf Example 4.6.3.} Let $u$ be the norm $-4$ vector of height 37 from
example 4.6.2. The simple roots of $D$ that have inner product 0, $-1$
or $-2$ with $u$ are as follows. (There are no simple roots that have
inner product $-2$.)
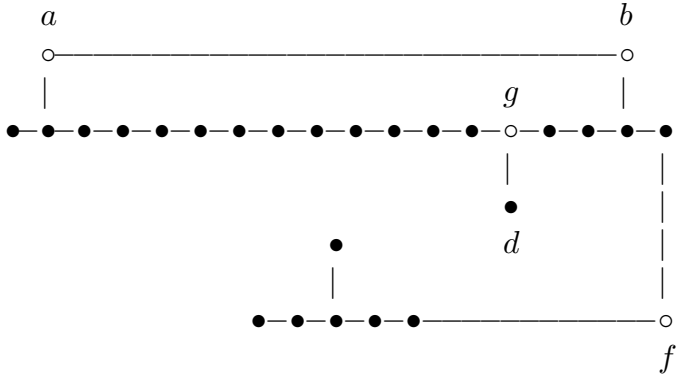
\halign
{
$#$&&$#$\cr
&&a&&&&&&&&&&&&&&&&&&&&&&&&&&&&&&b&&\cr
&&\circ&\!-&\!-&\!-&\!-&\!-&\!-&\!-&\!-&\!-&\!-&\!-&\!-
&\!-&\!-&\!-&\!-&\!-&\!-&\!-&\!-&\!-&\!-&\!-&\!-&\!-&\!-&\!-&\!-&\!-&\circ
\cr
&&|&&&&&&&&&&&&&&&&&&&&&&&&g&&&&&&|&&\cr
\bullet&\!-&\bullet&\!-&\bullet&\!-&\bullet&\!-&\bullet&\!-&\bullet&\!-&
\bullet&\!-&\bullet&\!-&\bullet&\!-&\bullet&\!-&\bullet&\!-&\bullet&\!-&
\bullet&\!-&\circ  &\!-&\bullet&\!-&\bullet&\!-&\bullet&\!-&\bullet\cr
&&&&&&&&&&&&&&&&&&&&&&&&&&|&&&&&&&&|\cr
&&&&&&&&&&&&&&&&&&&&&&&&&&\bullet&&&&&&&&|\cr
&&&&&&&&&&&&&&&&&\bullet&&&&&&&&&d&&&&&&&&|\cr
&&&&&&&&&&&&&&&&&|      &&&&&&&&&&&&&&&&&|\cr
&&&&&&&&&&&&&\bullet&\!-&\bullet&\!-&\bullet&\!-&\bullet&\!-&\bullet&\!-&\!-&\!-&\!-&\!-&\!-&\!-&\!-&\!-&\!-&\!-&\!-&\circ\cr
&&&&&&&&&&&&&&&&&&&&&&&&&&&&&&&&&&f\cr
\cr}
Using 3.5.4 we can find the norm $-2$ vectors $v$ corresponding to the
components of $u^\perp$:
\item{$a_1$:}
$v^\perp=a_{18}e_6$. $\height(v)=37-{\rm Coxeter~ number~ of~} a_1+1=36$.
\item{$a_4$:}
$v^\perp=a_{15}e_7a_3$. $\height(v)=33$. $v^\perp$ contains $g$, $b$,
and $f$ (and other points).
\item{$e_6$:}
$v^\perp=a_{13}a_{10}a_1$. $\height(v)=26$. $v^\perp$ contains $f$.
\item{$a_{13}$:}
$v^\perp=a_{12}e_6a_6$. $\height(v)=24$. $v^\perp$ contains $a$ and  $g$.

By 3.8.1 the roots $a$, $b$, $f$, $g$ that have inner product $-1$
with $u$ correspond to the norm 0 vectors of $D$ that have inner
product $-3$ with $u$. The norm 0 vectors are:
\item{$a$:} 
$A_{11}D_7E_6$. The total number of vectors conjugate to this under
the Weyl group of $u^\perp$ is 91.
\item{$b$:} $E_8D_{16}$. Number of conjugates is 10. 
\item{$f$:} $A_{15}D_9$. Number of conjugates is 135. 
\item{$g$:} $A_{17}E_7$. Number of conjugates is 140. 

The norm 0 vector corresponding to $f$ is the same as the norm 0 vector
corresponding to the $a_{10}$ component of the norm $-2$ vector
$a_{13}a_{10}a_1$ in accordance with 3.8.3, and similarly for the
other norm 0 vectors. The total number of norm 0 vectors that have
inner product $-3$ with $u$ is $91+10+135+140=376=8\height(u)+80$.

Notice that (for example) a simple root $r$ of $u^\perp$ is in
$z^\perp$ (where $z$ is the norm 0 vector corresponding to $g$) if and
only if $\sigma(r)$ is not joined to $g$, where $\sigma$ is the
opposition involution of $u^\perp$. This shows that the opposition
involution in 3.8.1 really is necessary and is not caused by a bad
choice of sign conventions.

We can find norm $-6$ vectors from $u$ using 3.8.6. The ones of type
$\ge 3$ are:
$$\eqalign{
a_2a_1a_{11}a_2a_1a_1e_6 & \qquad\height= 36\cr
a_{14}a_3e_6&\qquad\height=51\cr
a_5a_{12}e_6&\qquad\height=42\cr
a_4a_{13}a_1a_1d_5&\qquad\height=41\cr
a_6a_{12}d_5&\qquad\height=43\cr
d_4a_{11}a_2e_6&\qquad\height=42\cr
d_7a_{13}a_2a_1&\qquad\height=48\cr
}$$

{\bf Example 4.6.4} Let $u$ be the norm $-4$ vector of height 3 with
$u^\perp=a_1^2$. There are 101 simple roots $r$ with $-(r,u)=1$, 100
giving norm 0 vectors of type $A_1^{24}$ and one giving a norm 0
vector of Leech type. There are $104=8\height(u)+80$ norm 0 vectors
that have inner product $-3$ with $u$, as the norm 0 vector of Leech
type has 4 conjugates under the Weyl group of $u^\perp$. This is the
only type 3 norm $-4$ vector $u$ such that one of the norm 0 vectors
$z$ with $-(u,z)=3$ is of Leech type. The automorphism group of the
lattice $u^\perp$ is $\Z/2\Z\times 2^2.HS.2$ where $HS$ is the
Higman-Sims simple group.

{\bf Example 4.6.5.} Using 3.8.6 we can find all pairs $(r,A)$ where
$r$ is a root in a strictly 25 dimensional unimodular lattice $A$. We
have to match up pairs with the same lattice $A$ to get a list of such
lattices. If $A$ is determined by its root system this is easy, but
sometimes there are several lattices $A$ with the same root system and
then some care is needed in matching up the pairs $(r,A)$. For example
there are 3 lattices $A$ with root system $a_6a_5a_4a_3a_3$. Here we
give the simple roots of $D$ that have inner products $0$ or $-1$ with
three norm $-4$ vectors $u$ corresponding to these three lattices.

We give these three diagrams as subsets of the following diagram,
which gives all the roots that have inner product 0 or $-1$ with a
norm $-2$ vector in $D$ with Dynkin diagram $a_7a_6a_5a_4a_1$.

(4 rather complicated diagrams omitted here as they are too hard to do
in \TeX.)

\proclaim 4.7 Norm $-6$ vectors. 

Here we look briefly at norm $-6$ vectors of $II_{25,1}$ and show how
they can be used to find the even lattices of determinant 3 and
dimension 26. We omit most proofs because they are similar to previous
ones and the results do not seem important.

\proclaim Lemma 4.7.0. 
There are natural bijections between the following three set:
\item{(1)}
Orbits of norm $-6$ vectors $u$ in $D$. 
\item{(2)}
Equivalence classes of roots $r$ in 26 dimensional trimodular lattices
$T$. (There is a unique such lattice with no roots. See 5.7.)
\item{(3)}
Maps $e_6a_1\mapsto B$ where $B$ is a 32 dimensional even unimodular lattice. 

Proof. The correspondence between (2) and (3) follows by taking $T$ to
be the orthogonal complement of $e_6$ in $B$ and $r$ to be the root
corresponding to the $a_1$. The correspondence between (1) and (3)
follows from 4.9.1. Q.E.D.

Remark added 1999: Each of the three sets above has 2825 elements.

\proclaim Lemma 4.7.1. The following assertions are equivalent.
\item{(1)} $u$ has type 1.
\item{(2)}
The 26 dimensional lattice $T$ is the sum of a Niemeier lattice and an
$a_2$ and the root $r$ is in the $a_2$ component.
\item{(3)}
The map $e_6a_1\mapsto B$ factors as $e_6a_1\mapsto e_8\mapsto B$.

Proof omitted. 

\proclaim Lemma 4.7.2. The following conditions are equivalent:
\item{(1)} $u$ has type 1 or 2.
\item{(2)}
$T$ contains a root of norm 6, or equivalently $T'$ contains a vector
of norm $2/3$.
\item{(3)}
The map $e_6a_1\mapsto B$ factors as $e_6a_1\mapsto e_7\mapsto B$.

Proof omitted. 

The norm $-6$ vectors $u$ of type $\le 2$ can be found from table $-2$
using 3.7. They correspond to roots of 25 dimensional even bimodular
lattices. There are 24 orbits of such $u$'s of type 1 and 297+32 of
type 2. From now on we will assume that $u$ is a norm $-6$ vector of
$D$ of type $\ge 3$. By 4.7.2 this is equivalent to assuming that all
roots of $T$ have norm 2, or that the $e_6$ is a component of the
Dynkin diagram of $T$.

We now look at the set of norm 0 vectors that have inner product $-3$
with $u$. Such vectors come in pairs whose sum is $u$. Write $C$ for
the component of the Dynkin diagram of $T$ containing $r$.

\proclaim Lemma 4.7.3. 
The orbits of norm 0 vectors $z$ with $-(z,u)=3$ under the Weyl group
of $u^\perp$ depend on $C$ as follows:
\halign{
$#$\hfill&~~~~\hfill#\hfill&#\cr
C&Size of orbits\cr
a_1&none\cr
a_n(n\ge 2)&$(n-1)+(n-1)$&(Two orbits)\cr
d_n(n\ge 4)&$4n-8$\cr
e_6&$20$\cr
e_7&$32$\cr
e_8&$56$\cr
}
The number of such norm 0 vectors is $2h-4$  where $h$ is the Coxeter number 
of $C$. 

Proof. (Sketch) Such norm 0 vectors correspond to roots of $T$ that have inner product 1 with $r$. Q.E.D. 

{\bf Remark 4.7.4.} The number of roots of $r^\perp$ appears to be
$$6\height(u)-20+z_3$$ where $z_3$ is the number of norm 0 vectors
that have inner product $-3$ with $u$. I have a proof of this when the
Dynkin diagram of $u^\perp$ has a component other than $a_1$ or $a_2$,
but it is too long and disgusting to be worth reproducing. I have
check a few cases when the root system of $u^\perp$ consists of just
$a_1$'s and $a_2$'s.

(Remark added 1999: the formula for the number of roots of $u^\perp$ is
correct: see Adv. Math. Vol 83 no. 1 1990, p.30.)

Assume that the number of roots of $u^\perp$ is always given by the
expression above. Let $T$ be any 26 dimensional even trimodular
lattice with no roots of norm 6 and define $t$ to be $h-5+\height(u)$
where $u$ is a norm $-6$ vector of $D$ corresponding to some root $r$
of $T$, and where $h$ is the Coxeter number of the component of $r$ in
$T$. (If $T$ has no roots then define $t$ to be 0.) Then the number of
roots of $T$ is $6t$, and in particular $t$ is well defined, This
follows from 4.7.3 and 4.7.4.

{\bf Example 4.7.5.} Here is a short table of some 26 dimensional even
trimodular lattices $T$ with no roots of norm 6. It is complete for
$t=0$ or 1.

\halign{
#\hfill&~~~~~~$#$\hfill\cr $t$&\hbox{Dynkin diagram of }T\cr 0&\hbox{no roots}
\cr
1&a_2\cr 1&a_1^3\cr 2&a_3\cr 2&a_1^3a_2\cr 4&d_4\cr
48&a_{13}d_7a_4a_1\cr 92&d_{15}a_{11}\cr } The $d_{15}$ and $a_{11}$
components of the one with $t=92$ give norm $-6$ vectors of $D$ of
heights 69 and 85 and Dynkin diagrams $d_{13}a_1a_{11}$ and
$d_{15}a_9$.

Remark added 1999: the norm $-6$ vectors have been classified and
there are 2825 orbits. The 26 dimensional even lattices of determinant
3 have almost been classified: there are 4 or 5 ambiguities left to
resolve, and the number of them is about 678.

\proclaim 4.8 Vectors whose norm is not squarefree.

Given a vector of norm $-2m^2n$ in $II_{25,1}$ we show how to find an
orbit of vectors of norm $-2n$. This is used in 4.10 to construct norm
$-2$ vectors from norm $-8$ vectors, and in 4.11 to solve a problem
about shallow holes in the Leech lattice. 4.8.1 shows that for some
norm $-2m^2n$ vectors we can find a canonical norm $-2n$ vector in the
orbit associated to it.

{\bf Notation.} Let $u$ be a primitive vector of norm $-2m^2n$ in
$D$. $u^\perp$ is a 25 dimensional even lattice $B$ of determinant
$2m^2n$ such that $B'/B$ is generated by an element $b$ of norm
$1/2m^2n\bmod 2$. If $C$ is the lattice generated by $B$ and $2mnb$
then $C$ is a 25 dimensional even lattice of determinant $2n$ such
that $C/C'$ is generated by $c=mb$ of norm $1/2n\bmod 2$. Hence by
3.1.1 $C$ determines an orbit of primitive vectors of norm $-2n$ in
$D$. We will say that any vector of this orbit is associated to
$u$. There is an orbit of vectors associated to $u$ for every positive
integer $m$ such that $2m^2$ divides $u^2$.

\proclaim Theorem 4.8.1. 
Suppose that the primitive vector $u$ has norm $-2m^2n$ with $m\ge
2$. If there is a norm $-2n$ vector $v$ of $II_{25,1}$ with
$-(v,u)=2mn+1$ then there is such a vector in $D$, which is unique if
$m\ge 3$ or if $u$ has type at least 3 and norm $-8$. Any such vector
$v$ is associated to $u$.

Proof. For such vectors $v$ we have $(v,u)^2\ge v^2u^2=4m^2n^2$ and
equality cannot hold as $u$ is primitive, so $-(v,u)\ge 2mn+1$. Hence
by 3.3.1 if such a vector $v$ exists then there is one in $D$.

There is at most one such vector $v$ in $D$ if for any two such
vectors $v_1$, $v_2$ we have $(v_1-v_2)^2<4$, because then either
$v_1=v_2$ or $v_1$ and $v_2$ are conjugate by reflection in the root
$(v_1-v_2)$ of $u^\perp$. The projection of $v_i$ into $u^\perp$ has
norm $2/m+1/2m^2n$ so if $m\ge 3$ then $(v_1-v_2)^2<4$. If $m=2$ then
$(v_1-v_2)^2\le 4$ and if equality holds and $u^2=-8$ then
$(v_1+v_2-u)$ is a norm 0 vector that has inner product 2 with
$u$. Hence if $m\ge 3$ or $u$ has type $\ge 3$ and norm $-8$ there is
at most one such vector $v$ in $D$.

To show that $v$ is associated to $u$ we need to show that $u^\perp$
is isomorphic to a sublattice of $v^\perp$. Let $r$ be the vector
$u-mv$. Then $r^2=2m$, $(r,u)=m$, $(r,v)=-1$. Reflection in $r^\perp$
maps $u$ to $mv$ and maps vectors of $II_{25,1}$ in $u^\perp$ to
$II_{25,1}$ and so maps $u^\perp$ into $v^\perp$. (Note that
reflection in $r^\perp$ does not map all vectors of $II_{25,1}$ to
$II_{25,1}$!) Q.E.D.

\proclaim 4.9 Large norm vectors in $II_{25,1}$. 

We have given several correspondences between negative norm vectors in
$II_{25,1}$ and certain classes of lattices; for example, norm $-6$
vectors correspond to 26 dimensional even trimodular lattices. In this
section we show how to assign a class of lattices to vectors of norm
$-2n$ in $II_{25,1}$. In fact there are several ways to do this, one
for each orbit of norm $2n$ vectors in $e_8$.

\proclaim Theorem 4.9.1. 
Let $x$ be a primitive vector in $e_8$ of norm $2n$. Then there is a
bijection between isomorphism classes of
\item{(1)}
Norm $-2n$ vectors of $D$ (or $II_{25,1}$).
\item{(2)}
Maps of $x^\perp$ to $B$, where $B$ is a 32 dimensional even
unimodular lattice.

Proof. From a norm $-2n$ vector $u$ of $II_{25,1}$ we get a 25
dimensional even lattice $U$ with $U'/U$ cyclic of order $2n$
generated by some element $c$ of norm $1/2n\bmod 2$. $B$ is generated
by $U\oplus x^\perp$ and the vector $c\oplus d$ where $d$ is a vector
of ${x^\perp}'/x^\perp$ which is the projection in ${x^\perp}'$ of a
vector in $e_8$ having inner product 1 with $x$. It is routine to
check that this gives a bijection between (1) and (2). Q.E.D.

{\bf Remark.} A norm 0 vector in $II_{33,1}$ corresponding to $B$ can
be constructed as $x\oplus u$ in $e_8\oplus II_{25,1}\cong II_{33,1}$.

There are unique orbits of vectors of norms 0, 2, 4, 6 in $e_8$ and
their orthogonal complements have root systems $e_8$, $e_7$, $d_7$,
and $e_6a_1$. This gives the correspondences of 3.1 between vectors of
norms 0, $-2$, $-4$, and $-6$ in $II_{25,1}$ and certain
lattices. Vectors of norms 8, 10, and 14 in $e_8$ will be used later.

\proclaim 
Appendix to 4.9: An algorithm for writing down the orbits of vectors in $e_8$. 

Any orbit of vectors on $e_8$ has a unique representative in the Weyl
chamber of $e_8$. As $e_8$ is unimodular, these representatives are in
1:1 correspondence the sequences of 8 non negative integers giving
their inner products with the 8 simple roots of $e_8$. Let $C$ be the
set of all points of $e_8\otimes \Q$ that are at least as close to 0
as to any other point of $e_8$, and for any $x$ in $e_8$ write $n(x)$
for the smallest integer $n$ such that $x/n$ is in $C$. As $e_8$ is
generated by its norm 2 vectors, $x/n$ is in $C$ if and only if it has
inner product $\le 1$ with all roots of $C$, so if $x$ is in the Weyl
chamber of $e_8$ and has inner products $x_i$ with the simple roots of
$e_8$, then $n(x)=\sum_iw_ix_i$ (where the $w_i$'s are the weights of
the simple roots) and $x/n(x)$ lies on the boundary of $C$. We have
$$n(x)^2/2\le x^2\le n(x)^2.$$

We list the orbits of $e_8$ in order of $n(x)$. We can find the orbits
with $n(x)=n$ from those with $n(x)\le n-1$ as follows:

\item{(1)} If $n$ is even take the vector 
$x=$
\halign{
$#$&&$#$\cr
&&&&\bullet\cr
{n\over 2}
&&&&|\cr
\bullet&\!-&\bullet&\!-&\bullet&\!-&\bullet&\!-&\bullet&\!-&\bullet&\!-&\bullet\cr
}
of norm $n^2$. There are 16 lattice points of $e_8$ nearest to $x$. 
\item{} We number the points of $e_8$ as follows:
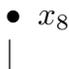
\halign{
$#$&&$#$\cr
&&&&\bullet&x_8\cr
&&&&|\cr
\bullet&\!-&\bullet&\!-&\bullet&\!-&\bullet&\!-&\bullet&\!-&\bullet&\!-&\bullet\cr
x_7&&x_6&&x_5&&x_4&&x_3&&x_2&&x_1\cr
}
\item{(2)}
For each $y$ corresponding to $(y_1,\ldots,y_8)$ giving an orbit with
$n(y)=n-2$ take the vector corresponding to
$(y_1+1,y_2,\ldots,y_8)$. This gives a vector $x$ of norm
$y^2+2(n-1)$. There are 2 lattice points nearest to $x/n$.
\item{(3)} For each vector $y$ with $n(y)=n-1$ and $y_1=y_2=y_3=y_4=y_5=0$, 
$y_8=1$ take a vector $x$ corresponding to
$(0,0,0,0,0,y_6+1,y_7,0)$. $x$ has norm $y^2+2(n-1)$ and there are 8
lattice vectors of $e_8$ closest to $x/n$.
\item{(4)}
For each vector $y$ with $n(y)=n-1$ and $y_1=y_2=\ldots y_{i-1}=0$,
$y_i=1$ ($1\le i\le 7$) take a vector $x$ with $x_8=y_8$ if $i\le 4$
and $x_8=y_8+1$ if $i\ge 5$, $x_{i+1}=1+y_{i+1}$ if $i\le 6$,
$x_1=x_2=\ldots=x_i=0$, $x_j=y_j$ for $i+2\le j\le 7$. $x$ has norm
$y^2+2(n-1)$ and there are $i+2$ points of $e_8$ closest to $x/n$.

By applying (1), (2), (3), and (4) we get all vectors $x$ in $D$ with
$n(x)=n$. After a bit of practice the vectors of $e_8$ whose norm is
at most (say) 100 can be worked out as fast as they can be written
down.

The orbits of vectors in $e_8/ne_8$ are represented by the vectors $x$
of the Weyl chamber with $n(x)\le n$. A vector of $e_8/ne_8$
represented by a vector $x$ with $n(x)\le x$ is represented by a
unique such vector if $n(x)\le n-1$, and by (the number of vectors of
$e_8$ nearest to $x/n$) such vectors is $n(x)=n$.  For example
$e_8/3e_8$ contains 5 orbits, of sizes 1, 240, $240\times 9$
$240\times 28/3$, and $240\times 72/9$, represented by vectors of norm
0, 2, 4, 6, 8. (See the table below.)

Here is a table of the orbits of vectors of $e_8$ with $n(x)$ at most 6. 
It includes all orbits of norm at most 24. 
\halign{
\hfill$#$~~&$#$\hfill~&\hfill#~&\hfill$#$&$#$\hfill\cr
n(x)&x_i(i=765(8)4321)&norm&\hbox{ number of }&\hbox{vectors}\cr
\cr
0&000(0)0000&0&1&/1\cr
\cr
2&000(0)0001&2&240&/2\cr
2&100(0)0000&4&9.240&/16\cr
\cr
3&000(0)0010&6&28.240&/3\cr
3&000(1)0000&8&72.240&/9\cr
\cr
4&000(0)0002&8&240&/2\cr
4&100(0)0001&10&126.240&/2\cr
4&000(0)0100&12&252.240&/4\cr
4&010(0)0000&14&288.240&/8\cr
4&200(0)0000&16&9.240&/16\cr
\cr
5&000(0)0011&14&56.240&/2\cr
5&000(1)0001&16&576.240&/2\cr
5&100(0)0010&18&756.240&/3\cr
5&000(0)1000&20&1008.240&/5\cr
5&100(1)0000&22&576.240&/8\cr
\cr
6&000(0)0003&18&240&/2\cr
6&100(0)0002&20&126.240&/2\cr
6&000(0)0101&22&756.240&/2\cr
6&010(0)0001&24&2016.240&/2\cr
6&200(0)0001&26&126.240&/2\cr
6&000(0)0020&24&28.240&/3\cr
6&000(1)0010&26&2016.240&/3\cr
6&100(0)0100&28&2520.240&/4\cr
6&001(0)0000&30&2016.240&/6\cr
6&110(0)0000&32&576.240&/8\cr
6&000(2)0000&32&72.240&/9\cr
6&300(0)0000&36&9.240&/16\cr
}

The number of vectors in the orbit of $x$ is written in the form
$m/n$, where $m$ is the number of vectors of $e_8$ in that orbit and
$n$ is the number of vectors of $e_8$ nearest to $x/n(x)$.

{\bf Example.} From the orbit of norm 30 in the table above we see that
we could find all 32 dimensional even unimodular lattices with $a_4
a_2a_1$ in their Dynkin diagrams by finding the norm $-30$ vectors in
$II_{25,1}$. More generally we could find all 32 dimensional even
unimodular lattices whose root system contains $e_8$, $e_7$, $d_7$, $e_6a_1$,
$a_7$, $d_5a_2$, $a_6a_1$, $a_4a_3$, or $a_4a_2a_1$ from certain
negative norm vectors in $II_{25,1}$. (These are the cases when
$x^\perp$ is generated by norm 2 roots for some vector $x$ of $e_8$.)

\proclaim 4.10 Norm $-8$ vectors. 

We describe a few properties of norm $-8$ vectors of $II_{25,1}$ and
give an example. Norm $-8$ vectors could be used to find the 32
dimensional even unimodular lattices containing an $a_n$ for $n\ge 7$.

$e_8$ has 2 orbits of norm 8 vectors which look like
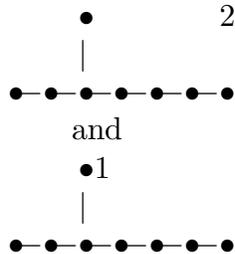
\halign{
$#$&&$#$\cr
&&&&\bullet&&&&&&&&2\cr
&&&&|\cr
\bullet&\!-&\bullet&\!-&\bullet&\!-&\bullet&\!-&\bullet&\!-&\bullet&\!-&\bullet\cr
}
and
\halign{
$#$&&$#$\cr
&&&&\bullet&1\cr
&&&&|\cr
\bullet&\!-&\bullet&\!-&\bullet&\!-&\bullet&\!-&\bullet&\!-&\bullet&\!-&\bullet\cr
}
Using 4.9.1 we see that there are natural correspondences between 
\item{(1)} Norm $-8$ vectors of $II_{25,1}$,
\item{(2)} 
Certain index 2 sublattices of 25 dimensional even bimodular lattices $B$, and
\item{(3)}
$a_n$ ($n\ge 7$) components in the Dynkin diagrams of 32 dimensional
even unimodular lattices $C$.

The index two sublattice of $B$ corresponding to $u$ is isomorphic to
$u^\perp$, so the norm $-2$ vector of $D$ corresponding to $B$ is
associated to $u$ (see 4.8). $u$ has type 1 or 2 if and only if this
norm $-2$ vector has type 1. Suppose that $u$ has type exactly 3. Norm
$-2$ vectors that have inner product $-5$ with $u$ correspond under
$\tau_u$ to norm 0 vectors that have inner product $-3$ with $u$, so
by 4.8.1 there is a unique norm 0 vector $z$ in $D$ with $(z,u)=-3$ and
$\tau_u(z)$ is a norm $-2$ vector associated to $u$. The number of
roots of $\tau_u(z)$ is equal to the number of roots of $u^\perp$ plus
$z_4$, where $z_i$ is the number of norm 0 vectors that have inner
product $-i$ with $u$.

If $u$ corresponds to a component $a_{7+k}$ in a 32 dimensional
unimodular lattice $C$ then there are $k$ norm 0 vectors that have
inner product $-3$ with $u$. If $C$ has Dynkin diagram $a_{7+k}X$ then
$u^\perp$ has Dynkin diagram $a_kX$ (where $a_0$ and $a_{-1}$ are
empty).

If $u$ is primitive then it seems that
$$5n=24t+z_4+24z_3+96z_2-32z_1-114$$ where $n$ is the number of roots
of $u^\perp$ and $t$ is the height of $u$. We have $n+z_4=12t'-18$ if
the norm $-2$ vector associated to $u$ has type 2 and height $t'$.

{\bf Example. } $u$ has height 48, norm $-8$, and the Dynkin diagram
of $u^\perp$ is $a_6a_{12}e_6$. Some of the simple roots that have
small inner product with $u$ are

(Diagram omitted for the moment.)

The numbers are minus the inner products of the simple roots with $u$.
$z_3$ is 13 and $z_4$ is 0. There is a unique norm 0 vector $z$ of $D$
with $-(z,u)=3$ and its Dynkin diagram is $A_{11}D_7E_6$ and contains
the points $a$, $c$, $d$, and $g$. (All the simple roots shown have
inner products 0 or $-1$ with $z$.) By 4.8.1 the norm $-2$ vector
$v=\tau_u(z)$ is associated to $u$. $v^\perp$ has Dynkin diagram
$a_{12}a_6e_6$ which is obtained form the Dynkin diagram of $u$ by
deleting $x$ and adding $b$. (The Dynkin diagram of $v^\perp$ is
isomorphic to that of $u^\perp$ because $z_4=0$. This is not usually
the case.)

$z_3=13$ so $u$ corresponds to a 32 dimensional even unimodular
lattice with Dynkin diagram $a_{20}e_6a_6$. This lattice can also be
obtained via its $e_6$ component from norm $-6$ vectors of $D$ with
heights $61$, $75$ and Dynkin diagrams $a_{18}a_6$ and $a_{20}a_4$ as
in 4.7.

\proclaim 4.11 Shallow holes in the Leech lattice. 

A hole of $\Lambda$ is a point of $\Lambda\otimes \Q$ that is locally
at a maximal distance from points of $\Lambda$. By 2.5 a hole has
radius at most $\sqrt{2}$ where the radius of a hole is the distance
of $x$ from the points of $\Lambda$ nearest $x$; these points are
called the vertices of the hole $x$. The holes of radius $\sqrt 2$ are
the 23 orbits of deep holes in the Leech lattice and their vertices
form an extended Dynkin diagram. The holes of radius less than $\sqrt
2$ are called the shallow holes and their vertices form a spherical
Dynkin diagram of rank 25. It is useful to think of a hole as the
convex hull of its vertices. In [B-C-Q] the shallow holes are found by
searching for sets of 25 points in $\Lambda$ that form a spherical
Dynkin diagram, and it is proved that the list of such holes is
complete by adding up the volumes of all the holes found in some
fundamental domain of $\Lambda$ and showing that this is 1. There are
284 orbits of shallow holes. Any Dynkin diagram in $\Lambda$ all of
whose components are affine or spherical is contained in the vertices
of a deep or shallow hole.

The smallest multiple of a deep hole that lies in $\Lambda$ is the
Coxeter number of the corresponding Niemeier lattice. In this section
we calculate the smallest multiple of any shallow hole that is in
$\Lambda$. We do this by identifying shallow holes with certain
vectors of $D$ whose height is the smallest multiple of the hole in
$\Lambda$, and calculating the height using 4.8. This gives a supply
of vectors of $II_{25,1}$ of large negative norm.

Let $v/k$ be a shallow hole of $\Lambda$, where $v$ is in $\Lambda$ and $k$ is the smallest multiple of $v/k$ in $\Lambda$. We assume that 0 is one of the vertices of the hole, so that it has radius $v^2/k^2$. We write $\rho^2$ for the norm of the Weyl vector of the Dynkin diagram of the vertices of the hole and $d$ for the determinant of the Cartan matrix of this Dynkin diagram. Then the hole has radius $2-1/\rho^2$ and the volume of the convex hull of its vertices is $\sqrt{\rho^2d}/25!$. We have $v^2=k^2(2-1/\rho^2)$ so $k^2/\rho^2$ is an even integer. Define integers $m$ and $n$ by
$$k^2=2m^2n\rho^2$$ where $n$ is squarefree. $n$ is determined by
$\rho^2$, so to find $k$ we only have to find $m$.

\proclaim Lemma 4.11.1. There are 1:1 correspondences between
\item{(1)} Orbits of shallow holes in $\Lambda$, with
$m$, $n$ as above. 
\item{(2)}
Orbits of primitive vectors $u$ of $D$ of norm $-2m^2n$ and height $k$
such that the root system of $u^\perp$ has rank 25.  
\item{(3)}
Pairs $(b,B)$ where $B$ is a 25 dimensional even unimodular lattice
whose root system has rank 25 and such that $B'/B$ is cyclic of order
$2m^2n$ and generated by an element $b$ of norm $1/2m^2n\bmod 2$. (In
fact the automorphism group of any such lattice acts transitively on
the set of $b$'s.)

Proof. The correspondence between (1) and (2) is that between
$\Lambda\otimes \Q\cup\infty$ and points on the boundary of $D$, where
we identify points on the boundary of $D$ with some primitive vectors
of $II_{25,1}$. The correspondence between (2) and (3) is just
3.1.1. Q.E.D.

\proclaim Theorem 4.11.2. 
Let $X$ be the Dynkin diagram of the hole $v/k$. Then $m$ is 1, 2, 3,
or 6. $m$ is divisible by 2 if and only if the number of $a_7$'s in
$X$ is not divisible by 2. $m$ is divisible by 3 if and only if the
number of $a_8$'s in $X$ is not divisible by 3, with the following two
exceptions: $m=3$ for one of the two holes with $X=a_{17}d_7a_1$ (the
one with trivial automorphism group) and $m=1$ for one of the two
holes with $X=a_{17}a_8$. Hence $k$ is determined by $X$, except for a
small ambiguity when $X$ contains an $a_{17}$.

(Most holes have $m=1$. There is only one hole with $m=6$; it has
Dynkin diagram $a_8a_7d_5^2$.)  

Proof. We give some methods for finding upper and lower bounds on
$m$. The proof consists of using these to work out $m$ for each hole
in turn. Let $u$ be the vector of $D$ corresponding to $v/k$ as in
4.11.1, (3).

To find lower bounds for $m$ we find simple roots $r$ of $D$ whose
inner product with $u$ cannot be integral unless $m$ is divisible by
something. The projection of $w$ into $u^\perp$ is the Weyl vector
$\rho$ of $u^\perp$ because the root system of $u^\perp$ has rank 25,
so
$$\eqalign{w=&\rho+(w,u)u/(u,u)\cr
&= \rho+ku/2mn^2,\cr
}$$
so for any simple root $r$ of $D$ we have
$$\eqalign{(u,r)&=(2m^2n/k)((r,w)-(r,\rho))\cr &=(m^2n/k)(-2-(r,2\rho))\cr }$$
$(r,2\rho)$ is an integer and can be calculated if $(r,r')$ is known
for all simple roots $r'$ of $X$ as $\rho$ is a sum of these roots
with known coefficients. The fact that $(u,r)$ is an integer sometimes
implies that $m$ must be divisible by 2, 3, or 6.

{\bf Example.}  Let $v/k$ be the shallow hole with Dynkin diagram
$a_8a_7d_5^2$. Then $\rho^2=162$, so $n=1$ and $k=18m$. Some of the
simple roots of $D$ are shown below, with the $a_8a_7d_5^2$ in
black. (Note the $A_7^2D_5^2$ deep hole.)

(Diagram omitted)

The points $a$, $b$, $c$, $d$, $e$ have inner products 22, 22, 22, 22,
31 with $2\rho$, so their inner products with $u$ are $m/18$ times 24,
24, 24, 24, 33. As these inner products are integers $m$ must be
divisible by 6.

Now we give some ways of finding upper bounds for $m$. 

\item{(1)}
The lattice $L=u^\perp$ has determinant $2m^2n$ and $L'/L$ is cyclic,
so if $M$ is the lattice generated by norm $2$ vectors of $L$ then
$M'/M$ contains elements of order $2m^2$. This implies that if $q=p^i$
($i\ge 2$) is a prime power dividing $2m^2$ then the Dynkin diagram of
$L$ contains a component $a_{qj-1}$ for some integer $j$. $qj-1\le
25$, so $m$ is squarefree and the only primes that can divide $m$ are
2, 3, and 5. If $2$ divides $m$ then the Dynkin diagram $X$ contains
an $a_7$, $a_{15}$, or $a_{23}$; if $3$ divides $m$ then $X$ contains
an $a_8$ or $a_{17}$, and if $5$ divides $m$ then $X$ contains an
$a_{24}$. (In fact $m$ is never divisible by 5.) For example, the hole
with Dynkin diagram $a_{10}e_8e_7$ must have $m=1$, and the hole with
Dynkin diagram $d_{14}a_7a_2a_1a_1$ has $m=1$ or 2.  \item{(2)} Now
suppose that $d$ divides $m$ for some integer $d$. We apply 4.8 to $u$
to get a vector $v$ associated to $u$ of norm $-2m^2n/d^2$. The
lattice $v^\perp$ contains a copy of $u^\perp$ so its Dynkin diagram
has rank 25, and hence by 4.11.1 $v$ corresponds to another shallow
hole such that the sublattice of $v$ generated by roots has a lattice
of index dividing $d$ whose Dynkin diagram is $X$. (Warning: this new
shallow hole is not necessarily conjugate to $dv/k$.) Here is a list
of the sub root systems of simple root systems that have the same rank
and prime index. 
\halign{$#$&~$#$&\hfill$#$&\hfill~$#$&\hfill~$#$\cr
&\hbox{Sub root systems of index = }&2&3&5\cr
a_n&&\hbox{none}&\hbox{none}&\hbox{none}\cr
d_n&(n\ge 4)\hfill&d_m d_{n-m}(2\le m\le n-m)&\hbox{none}&\hbox{none}\cr
&&(d_2=a_1^2,d_3=a_3)\cr
e_6&&a_1a_5&a_2^3&\hbox{none}\cr
e_7&&a_7,a_1d_6&a_2a_5&\hbox{none}\cr
e_8&&d_8,e_7a_1&e_6a_2,a_8&a_4^2\cr
}

{\bf Example.} By (1), the shallow hole with $a_{15}e_7a_3$ has $m=1$
or 2. $m$ cannot be 2, because the table above shows that the only
root system containing $a_{15}e_7a_3$ with index 1 or 2 is
$a_{15}e_7a_3$ and there is only one shallow hole with Dynkin diagram
$a_{15}e_7a_3$.

{\bf Another example.} The hole $a_8a_7d_5^2$ has $m=6$. The shallow
holes corresponding to the divisors 1, 2, 3, or 6 of $m$ are
$a_8a_7d_5^2$, $d_{10}a_8e_7$, $e_8a_7d_5^2$, and $d_{10}e_8e_7$.

\item{(3)} Finally we can sometimes work out $m$ by using the list of vectors of $D$ of norms $-2$ or $-4$. Any such vector $u$ such that the Dynkin diagram of $u^\perp$ has rank 25 gives a shallow hole with $m=1$. 

{\bf Example.} There are two orbits of shallow holes with Dynkin
diagram $a_{17}d_7a_1$. Such a hole has $m=1$ or 3, and $m=1$ if and
only if the vector $u$ corresponding to it has norm $-4$. There is
only one orbit of norm $-4$ vectors with Dynkin diagram
$a_{17}d_7a_1$, so one of the holes has $m=1$ and the other has
$m=3$. (The construction in (2) applied to the hole with $m=3$ gives
the hole with $m=1$.)

{\bf Remark.} In particular if $X$ is the Dynkin diagram of a
unimodular lattice of dimension $i\le 23$ whose vector space is
spanned by its roots then $Xd_{25-i}$ is the Dynkin diagram of a
shallow hole.

By applying these we can work out $m$ for all holes and we find the
result stated in 4.11.2. Q.E.D.

\proclaim Chapter 5 Theta functions. 

This chapter gives some result that are obtained by combining the
theory in chapters 3 and 4 with the theory of theta functions. In
sections 5.1 to 5.4 we calculate the theta functions of the lattices
associated with norm $-2$ and $-4$ vectors of $II_{25,1}$. The height
of such vectors $u$ in the fundamental domain $D$ turns out to be a
linear function of the theta function of $u^\perp$, and there is some
evidence that this is true for all primitive negative norm vectors of
$D$. In sections 5.5 and 5.6 we give an algorithm for finding the 26
dimensional unimodular lattices and use it to show that there is a
unique such lattice with no roots, and then show that the number of
norm 2 roots of a 26 dimensional unimodular lattice is divisible by
4. In section 5.7 we construct a 27 dimensional unimodular lattice
with no roots and show that it is the unique such lattice with a
characteristic vector of norm 3.

{\bf Remark.} Conway and Thompson showed that there are unimodular
lattices with no roots in all dimensions $\ge 37$ by a non
constructive argument. For dimensions at most 27 the only unimodular
lattices with no roots are the 0 dimensional one, one 23 dimensional
one, two 24-dimensional ones, one 26 dimensional one, and at least one
27 dimensional one. There is some evidence that there is more than one
in each of the dimensions 27 and 28, and it is easy to construct odd
31 and 32 dimensional ones from the known 32 dimensional even
unimodular lattices with no roots. Notice the gap at 25 dimensions; in
fact any lattice of determinant 1 or 2 in 25 dimensions has a vector
of norm 1 or 2. (Remark added in 1999: Conway and Sloane have shown
that there are unimodular lattices with no roots in all dimensions at
least 26. Bacher and Venkov have classified all unimodular lattices
with no roots in dimensions 27 and 28.)

\proclaim 5.1 Notation and quoted results. 

We quote some results about theta functions of unimodular
lattices. The theta function $\theta_L$ of an $n$ dimensional
unimodular lattice $L$ can be written in the form
$$\theta_L(q)=\sum_{r=0}^{[n/8]} a_r\theta(q)^{n-8r}\Delta_8(q)^r$$
where $\theta(q)=1+2q+2q^4+2q^9+\cdots$ and $\Delta_8(q)=q+\cdots$ is
a modular form of weight 4. In this chapter $[n/8]$ will usually be 3,
so the theta function is determined by $a_0$, $a_1$, $a_2$, and
$a_3$. The coefficient of $q^0$ is 1, so $a_0=1$. If $L_0$ is the sub
lattice of even elements of $L$, then the theta function of $L_0'$ is
equal to
$$\theta_{L_0'}(q)=\sum_{r=0}^{[n/8]} (-1)^r2^{n-12r}a_rf_r(q)$$
where $f_r(q)$ is a power series with integral coefficients and leading 
term $q^{n/4-2r}$. 

\proclaim Lemma 5.1.1. 
Suppose $[n/8]=3$. Then the number of characteristic vectors of $L$ of
norm $n-24$ is equal to the number of vectors of $L_0'$ of norm
$(n-24)/4$, which is $-2^{n-36}a_3$, and if there are no such vectors
then the number of characteristic vectors of norm $n-16$ is
$2^{n-24}a_2$.

Proof. This follows from the formula for $\theta_{L_0'}(q)$. Q.E.D. 

\proclaim 5.2 Roots of a 25 dimensional even bimodular lattice. 

We will calculate the number of roots of $u^\perp$, where $u$ is a
norm $-2$ vector of $D$, in terms of the height of $u$.

{\bf Notation.} Let $u$ be a norm $-2$ vector of $D$ of height $t$ and
let $r$ be a highest root of the root system of $u^\perp$. (Such a
root always exists because the root system of $u^\perp$ is nonempty.)
Let $z_1$ be the number of norm 0 vectors that have inner product $-1$
with $u$, so $z_1$ is 2 if $u$ has type 1 and 0 otherwise.

\proclaim Lemma 5.2.1. The number of roots of $u^\perp$ is $12t-18+4z_1$. 

Proof. If $z_1\ne 0$ then $t=2h+1$ where $h$ is the Coxeter number of
the Niemeier lattice $N$ such that $u^\perp\cong N\oplus a_1$. The
number of roots of $u^\perp$ is $24h+2$ so the lemma is correct in this
case. From now on we assume that $z_1=0$.

The vector $z=r+u$ is a primitive norm 0 vector of $D$ because $u$ has
type $\ge 2$; let $h$ be its Coxeter number (so $h=-(z,w)$). Write $X$
for the lattice $\langle u,r\rangle^\perp$. We will find the number of
roots of $u^\perp$ by comparing it with the number of roots of
$X$. The Coxeter number of the component of $u^\perp$ containing $r$
is $t-h+1$ by 4.1.1, so
$$(\hbox{number of roots of~}u^\perp)=(4(t-h+1)-6)
+(\hbox{number of roots of~}X)$$

If $N$ is the Niemeier lattice corresponding to $z$ choose coordinates
$(x,m,n)$ for $II_{25,1}\cong N\oplus U$ with $x$ in $N$ and $m$ and
$n$ integers, so that $z=(0,0,1)$. In these coordinates we have
$w=(\rho,h,h+1)$ and $u=(v,2,n)$ for some $v$ in $N$ and some integer
$n$. We have $u^2=-2$ and $(u,w)=-t$, so $v^2-4n=-2$ and
$(v,\rho)-nh-2(h+1)=-t$.

The roots in $X$ are just the vectors of the form $(x,0,m)$ that have
inner product 0 with $u$ and norm 2, so they are in 1:1 correspondence
with the roots $x$ of $N$ that have even inner product with $v$. Let
$r_i$ be the number of roots of $N$ that have inner product $i$ with
$v$. If we choose $v$ to have minimum possible norm then $v$ has inner
product at most 2 with all roots $r$ of $N$, otherwise we could change
$u=(v,2,n)$ to $(v-2r,2,\hbox{?})$ and $v-2r$ would have smaller norm
than $v$. Hence the number of roots of $N$ is $r_0+2r_1+2r_2$ and the
number of roots of $X$ is $r_0+2r_2$.

By 2.2.5 we have $\sum_r(v,r)^2=2v^2h$ where the sum is over all roots
$r$ of $N$, and $\sum_{r>0}(v,r)=-2(v,\rho)$ where the sum is over all
positive roots $r$ of $N$, so
$$\eqalign{
2r_1+8r_2&=\sum_r(v,r)^2=2v^2h=2h(4n-2)\cr
\hbox{and~}r_1+2r_2&=-\sum_{r>0}(v,r)=2(v,\rho)=2(-t+nh+2(h+1)).\cr
}$$
This implies that 
$$\eqalign{
(\hbox{number of roots of~} X)&= r_0+2r_2\cr
&= (r_0+2r_1+2r_2)+(2r_1+8r_2)-4(r_1+2r_2)\cr
&= 24h+2h(4n-2)+8(t-nh-2(h+1))\cr
&=4h+8t-16,\cr
}$$
so the number of roots of $u^\perp$ is equal to this plus $4(t-h+1)-6$ 
which is $12t-18$. Q.E.D. 

\proclaim 5.3 The theta function of a 25 dimensional bimodular lattice. 

{\bf Notation.} $u$ is a norm $-2$ vector in $D$ and there are $z_1=0$
or 2 norm 0 vectors that have inner product $-1$ with $u$.

We will use the result of the previous section to calculate the theta
function of $u^\perp$.

\proclaim Theorem 5.3.1. 
The theta function of the 25 dimensional bimodular lattice $u^\perp$
is equal to
$$\theta_a+z_1\theta_b+12t\theta_c$$ where $\theta_a$, $\theta_b$,
$\theta_c$ are theta functions with integral coefficients that are
independent of $u$. We have
$$\eqalign{
\theta_c&=(\sum_n\tau(n)q^{2n})(1+2q^2+2q^8+2q^{18}+\cdots)\cr
&= q^2-22q^4\cdots\cr
\cr
\theta_a+2\theta_b&= \theta_\Lambda\times (1+2q^2+2q^8+\cdots)\cr
&\hbox{where } \theta_\Lambda=1+196560q^4+\cdots \cr
&\hbox{is the theta function of the Leech lattice.}\cr
\cr
\theta_a&=1-18q^2+143496q^4+\cdots\cr
}$$

Proof. We first show that the space generated by the theta functions
of the 25 dimensional bimodular lattices is three dimensional. Let $L$
be the lattice $u^\perp$. Then by 5.2.1 the theta function of $L$ is
$1+(12t-18+4z_1)q^2+\cdots$ and the theta function of $L'$ is
$1+z_1q^{1/2}+(12t-18+4z_1)q^2+\cdots$. Let $M$ be the 26 dimensional
unimodular lattice containing $L\oplus a_1$. Then $M$ has
$12t-18+4z_1+2$ vectors of norm 2, $2z_1$ of norm 1, and $2+z_1$
characteristic vectors of norm 2, so its theta function depends
linearly on $z_1$ and $t$. The theta function of its even sublattice
is therefore also a linear combination of $z_1$ and $t$, and is the
product of the theta functions of $L$ and $a_1$, so the theta function
of $L$ is equal to $\theta_a+z_1\theta_b+12t\theta_c$ for some
$\theta_a$, $\theta_b$, $\theta_c$ not depending on $u$.

If we let $L_1$ and $L_2$ be any two even bimodular lattices of the
form $L_i=N_i\oplus a_i$ for some Niemeier lattices $N_1$, $N_2$ with
different Coxeter numbers $h_i$, then the norm $-2$ vectors of $L_i$
have heights $t_i=2h_i+1$. Therefore
$$\eqalign{ 24(h_1-h_2)\theta_c&=12(t_1-t_2)\theta_c\cr &=
\theta_{L_1}-\theta_{L_2}\cr &=
24(h_1-h_2)(\sum_n\tau(n)q^{2n})(1+2q^2+2q^8+\cdots)\cr }$$ so $\theta$
has its stated value. The first few terms of $\theta_a$ can be found by
calculating the theta function of some lattice $u^\perp$ where $u$ has
type 2. (The one with root system $a_{17}e_8$ seems easiest to do
calculations with.)  Q.E.D.

{\bf Remark.} In particular if $u$ has type 2 then the number of norm
4 vectors of $u^\perp$ is $143496-264t$. The number of such norm 4
vectors can also be calculated from the norm $-4$ vectors of $II_{25,1}$
that have inner product $-4$ with $u$, and this gives a useful check
when enumerating these norm $-4$ vectors.

\proclaim 5.4 The theta function of a 25 dimensional unimodular lattice. 

{\bf Notation.} $u$ is a norm $-4$ vector in $D$ with height
$t$. $z_i$ is the number of norm 0 vectors that have inner product $i$
with $u$.

We will express the theta function of $u^\perp$ in terms of the height
of $u$.

\proclaim Lemma 5.4.1. The number of norm 2 roots of $u^\perp$ is 
$$8t-20+2z_2+8z_1.$$

Proof. If $u$ is one of the norm $-4$ vectors with no roots on
$u^\perp$ then either $t=1$, $z_1=1$, $z_2=2$ or $t=2$, $z_1=0$,
$z_2=2$, so the lemma is correct in these cases. If $u$ has type 1
then $t=3h+1$ where $h$ is the Coxeter number of the Niemeier lattice
$N$ such that $u^\perp\cong N$ plus a one dimensional lattice of
determinant 4, so $z_1=1$, $z_2=2$, and the number of roots of
$u^\perp$ is $24h$, so the lemma is also correct in these cases.

From now on we can assume that $u$ has type $\ge 2$ and that $u^\perp$
contains some highest root $r$. Then $v=r+u$ is a norm $-2$ vector of
$D$. Let $h$ be the Coxeter number of the component of $u^\perp$
containing $r$ so that $v$ has height $t-(h-1)$, and let $X$ be the
lattice $\langle u,v\rangle^\perp$. The number of roots of $u^\perp$
is equal to the number of roots of $X$ plus $4h-6$.

The projection $u_v$ of $u$ into $v^\perp$ has norm 4, so $u$ has
inner product 0, $\pm1$, or $\pm 2$ with all roots of $v^\perp$. Let
$r_i$ be the number of roots of $v^\perp$ that have inner product $i$
with $u$, so that the number of roots of $v^\perp$ is
$r_0+2r_1+2r_2$. Let $\rho$ be the projection of $w$ into $v^\perp$,
which by 4.3.3 is equal to half the sum of the positive roots of
$v^\perp$. Then
$$\eqalign{
(r_1+2r_2)/2=(u,\rho)&= (u_v,\rho)\cr
&=(u_v,w)\cr
&=(u-2v,w)\cr
&= 2(t-(h-1))-t\cr
&= t-2h+2\cr
}$$

$r_2$ is the number of roots $r'$ of $II_{25,1}$ that have inner
product 0 with $v$ and $2$ with $u$, and $r'$ is such a root if and
only if $r'+v$ is a norm 0 vector that has inner product 2 with $u$
and 0 with $r=v-u$, so $r=z_2$ is $r$ is not the sum of two norm 1
vectors in the unimodular lattice containing $u^\perp$ and $r=z_2-4$
if $r$ is the sum of two norm 1 vectors. The vector $v$ has type 1 if
and only if $r$ is the sum of two norm 1 vectors, so
$$\eqalign{
\hbox{number of roots of $X$}=r_0&=(r_0+2r_1+2r_2)-2(r_1+2r_2)+2r_2\cr
&=(\hbox{roots of }u^\perp)-4(t-2h+2)+2r_2\cr
&= 12(t-(h-1))-18(+8 \hbox{ if $v$ has type 1})\cr
&\qquad -4(t-2h+2)+2z_2(-8 \hbox{ if $r$ is the sum of}\cr
&\qquad\qquad\hbox{ two norm 
$-1$ vectors})\cr
&= 8t-4h-14+2z_2.\cr
}$$
Therefore the number of roots of $u^\perp$ is equal to this plus
$4h-6$ which gives the result of the lemma (as $z_1=0$). Q.E.D.

\proclaim Theorem 5.4.2. 
The theta function of the unimodular lattice of $u^\perp$ is given by 
$$\theta_a+2z_1\theta_b+z_2\theta_c+8t(\sum_n\tau_nq^{2n})(\sum_nq^{4n^2})$$
where the $\theta$'s are fixed theta functions with integer coefficients. 

Proof. It follows from section 5.1 that this theta function must be a
linear function of $z_1$, $z_2$, and the number of roots of $u^\perp$,
so by 5.4.1 it is a linear function of $z_1$, $z_2$, and $t$. The
explicit expression for the coefficient of $t$ is proved in exactly
the same way as in 5.3.1. Q.E.D.

{\bf Example 5.4.3.} Suppose that the unimodular lattice of $u^\perp$
has no vectors of norm 1, so that $z_1=z_2=0$. Then this lattice has
$2a_2=2(8t-20-26a_1-1200a_0)=16t+160$ characteristic vectors of norm 9
(where the $a_i$'s are as in 5.1). This is twice the number of norm 0
vectors that have inner product 3 with $u$, so in particular any norm
$-4$ vector $u$ has type at most 3.

{\bf Remark.} 5.4.2 and 5.3.1 suggest that the heights of primitive
vectors $u$ in $D$ of fixed norm depend linearly on the theta function
of $u^\perp$. The method used to prove 5.4.1 and 5.2.1 is too messy to
generalize easily to vectors of larger norm, though I have used it to
give a rather long proof that the height depends linearly on the theta
function for norm $-6$ $u$'s such that the Dynkin diagram of $u^\perp$
contains a component with at least 4 points.

(Remark added in 1999: the height of a primitive vector $u$ of $D$
depends linearly on the theta functions of all the cosets of $L$ in
$L'$ where $L=u^\perp$. For vectors of norms $-2$, $-4$, or $-6$ the
theta function of these cosets depend linearly on the theta function
of $u^\perp$, so the height of vectors of these norms is a linear
function of the theta function. See ``Automorphic forms with
singularities on Grassmannians''.)

\proclaim 5.5 26 dimensional unimodular lattices. 

In this section we combine some results about theta functions with the
algorithm of chapter 3 to give a practical algorithm for finding all
26 dimensional unimodular lattices. We use this to find the least
uninteresting such lattice: the unique one with no roots.

\proclaim Lemma 5.5.1. 
A 26-dimensional unimodular lattice $L$ with no vectors of norm 1 has
a characteristic vector of norm 10.

Proof. If $L$ has a characteristic vector $x$ of norm 2 then $x^\perp$ is a 25
dimensional even bimodular lattice and therefore has a root $r$ by
4.3.1; $2r+x$ is a characteristic vector of norm 10. If the lemma is
not true we can therefore assume that $L$ has no vectors of norm 1 and
no characteristic vectors of norm 2 or 10. Its theta function is
determined by section 5.1 and these conditions and turns out to be
$1-156q^2+\cdots$ which is impossible as the coefficient of $q^2$ is
negative. Q.E.D.

\proclaim Lemma 5.5.2. There is a bijection between isomorphism classes of
\item{(1)} 
Norm 10 characteristic vectors $c$ in 26-dimensional unimodular
lattices $L$, and
\item{(2)}
Norm $-10$ vectors $u$ in $II_{25,1}$
\item{}
given by $c^\perp\cong u^\perp$. 
\item{}We have $\Aut(L,c)=\Aut(II_{25,1}, u)$. 

Proof. Routine. Note that $-1$ is a square mod 10. Q.E.D. 

5.5.1 and 5.5.2 give an algorithm for finding 26 dimensional
unimodular lattices $L$, in case anyone finds a use for such
things. $L$ has a vector of norm 1 if and only if there is a norm 0
vector that has inner product $-1$ or $-3$ with $u$, and $L$ has a
characteristic vector of norm 2 if and only if there is a norm 0
vector that has inner product $-2$ with $u$, so it is sufficient to
find all norm $-10$ vectors $u$ of type at least 4.

\proclaim Lemma 5.5.3. 
Take notation as in 5.5.2. $L$ has no roots if and only if $u^\perp$
has no roots and $u$ has type at least $5$.

Proof. If $u^\perp$ has roots then obviously $L$ has too. If there is
a norm 0 vector $z$ that has inner product 1, 2, 3, or 4 with $u$ then
the projection $z_u$ of $z$ into $u^\perp$ has norm $1/10$, $4/10$,
$9/10$, or $16/10$. The lattice $L$ contains $u^\perp+c$, and the
vector $z_u\pm 3c/10$, $z_u\pm 4c/10$, $z_u\pm c/10$, or $z_u\pm
2c/10$ is in $L$ for some choice of sign and has norm 1, 2, 1, or
2. Hence if $u$ has type at most 4 then $L$ has roots. Conversely if
$L$ has a root $r$ then either $r$ has norm 2 and inner product 0,
$\pm 2$, $\pm 4$ with $c$ or it has norm 1 and inner product $\pm 1$,
$\pm 3$ with $c$, and each of these cases implies that $u^\perp$ has
roots or that $u$ has type at most 4. Q.E.D.

Now let $L$ be a 26 dimensional unimodular lattice with no roots
containing a characteristic vector $c$ of norm 10, and let $u$ be a
norm $-10$ vector of $D$ corresponding to it as in 5.5.2.

\proclaim Lemma 5.5.4. 
$u=z+w$, where $z$ is a norm 0 vector of $D$ corresponding to a
Niemeier lattice with root system $A_4^6$, and $w$ is the Weyl vector
of $D$. In particular $u$ is determined up to conjugacy under
$\Aut(D)$.

Proof. $u^\perp$ has no roots so $u=w+z$ for some vector $z$ of
$D$. By 5.5.3 $u$ has type at least 5, so $-(z,w)=-(u,w)\ge 5$. Hence
$$-10=u^2=z^2+2(z,w)\le 2(z,w)\le -10$$ so $(z,w)=-5$ and $z^2=0$. The
only norm 0 vectors $z$ in $D$ with $-(z,w)=5$ are the primitive ones
corresponding to $A_4^6$ Niemeier lattices. Q.E.D.

\proclaim Lemma 5.5.5. 
If $u=z+w$ is as in 5.5.4 then the 26 dimensional unimodular lattice
corresponding to $u$ has no roots.

Proof. $u^\perp$ obviously has no roots so by 5.5.3 we have to check
that there are no norm 0 vectors that have inner product $-1$, $-2$,
$-3$, or $-4$ with $u$. Let $x$ be any norm 0 vector in the positive
cone. If $x$ has type $A_4^6$ then $(x,u)\le (x,w)\le -5$; if $x$ has
Leech type then $(x,u)\le -(x,z)\le -5$; if $x$ has type $A_1^{24}$
then $(x,u)=(x,w)+(x,z)\le -2+-3$ ($(x,z)$ cannot be $-2$ as there are
no pairs of norm 0 vectors of types $A_1^{24}$ and $A_4^6$ that have
inner product $-2$ by the classification of 24 dimensional unimodular
lattices); and if $x$ has any other type then $(x,u)=(x,w)+(x,z)\le
-3+-2$. Q.E.D.

\proclaim Theorem 5.5.6. 
There is a unique 26 dimensional unimodular lattice $L$ with no
roots. Its automorphism group acts transitively on the 624
characteristic norm 10 vectors of $L$ and the stabilizer of such a
vector has order $5^3.2.120$, so the automorphism group has order
$2^8.3^2.5^4.13$.

Proof. By 5.5.1 $L$ has a characteristic vector of norm 10, so by
5.5.3 and 5.5.4 $L$ is unique and its automorphism group acts
transitively on the characteristic vectors of norm 10. By 5.5.5 $L$
exists. The theta function is determined by the conditions that $L$
has no vectors of norm 1 or 2 and no characteristic vectors of norm 2,
and it turns out that the number of characteristic vectors of norm 10
is 624. The stabilizer of such a vector is isomorphic to
$\Aut(II_{25,1}, u)$, which is a group of the form $5^3.2.S_5$ where
$S_5$ is the symmetric group on 5 letters. Q.E.D.

In the rest of this section we show that $G=\Aut(L)$ is isomorphic to
$O_5(5)$.

\proclaim Lemma 5.5.7. 
If $c_1$ and $c_2$ are two characteristic vectors of $L$ then write
$c_1\cong c_2$ if and only if $(c_1,c_2)\equiv 2\bmod 4$. This is an
equivalence relation on such vectors.

Proof. $\cong$ is obviously symmetric and reflexive. Suppose that
$c_1\cong c_2$ and $c_2\cong c_3$. Any two characteristic vectors are
congruent mod $2L$, so $c_1-c_3=2x$ for some $x$ in
$L$. $(x,c_2)\equiv 0\bmod 2$ so $x^2$ is even because $c_2$ is
characteristic. Hence $(c_1-c_3)^2\equiv 0\bmod 8$ and this implies
that $(c_1,c_3)\equiv 2\bmod 4$ because $c_1^2\equiv c_3^2\equiv 2\bmod
8$. Hence $\cong$ is transitive and is therefore an equivalence
relation. Q.E.D.

\proclaim Lemma 5.5.8. 
Let $c$ be a fixed norm 10 characteristic vector of $L$. The subgroup
of $\Aut(L)$ fixing $c$ has 12 orbits when acting on the
characteristic vectors of norm 10, as follows:
\halign{\hfill#&#\hfill&\hfill#&#\hfill\cr
Inner product& with $c$.&~~~~Number of vectors&~in orbits.\cr
10&&1\cr
4&&30\cr
2&&30, 125\cr
0&&1, 1, 125, 125\cr
$-2$&&30, 125\cr
$-4$&&30\cr
$-10$&&1\cr
}

Proof. Calculate the number of things in various orbits of 
$\Aut(II_{25,1},u)$. 
Q.E.D. 

By 5.5.7 and 5.5.8 there are exactly two equivalence classes of
characteristic norm 10 vectors under $\cong$. Let $G_1$ be the
subgroup of $G=\Aut(L)$ of index 2 fixing one of these classes. $G_1$
contains $-1$, so write $G_2$ for $G_1/\langle -1\rangle$.

\proclaim Lemma 5.5.9. $G_2$ is the simple group $P\Omega_5(5)$. 

Proof. $G_2$ has order $2^6.3^2.5^4.13$, and by 5.5.8 it has a rank 3
permutation representation on $1+30+125$ points. Let $N$ be a minimal
nontrivial normal subgroup of $G_2$. The rank 3 permutation
representation is primitive and faithful, so $N$ acts transitively on
the points of it and hence the order of $N$ is divisible by 13. By the
Fratteni argument $G_2$ is the product of $N$ and the normalizer $H$
of a Sylow 13 subgroup of $N$. $H$ cannot have order divisible by 5
because then $G_2$ would have a subgroup of order $13.5$ (which must
be cyclic) and this is impossible as $G$ has a faithful 26 dimensional
rational representation. Therefore $N$ has order divisible by
$5^4$. As $N$ is characteristically simple and its order is divisible
by 13 but not by $13^2$, $N$ must be simple. The order of $N$ is
divisible by $5^4.156=2^2.3.5^4.13$ and divides $2^6.3^2.5^4.13$, and
the only simple group with one of these orders is $P\Omega_5(5)$ of
order $2^6.3^2.5^4.13$ which must therefore be isomorphic to
$G_2$. Q.E.D.

If the group $G_1=2.G_2$ were perfect it would have to be $Sp_4(5)$,
but this group has no 26 dimensional representation with the center
acting as $-1$, so $G_1$ is not perfect and must therefore be
$(\Z/2\Z)\times G_2$, where the $\Z/2\Z$ is generated by $-1$. $G_1$
has a unique faithful rational 26 dimensional representation $L$ which
is the sum of two complex conjugate irreducible representations. $L$
can only be extended to a faithful representation of a group $G$
containing $G_1$ with index 2 if $G$ is isomorphic to $O_5(5)\cong
2\times G_2.2$ and $L$ is an irreducible representation of $G$.

\proclaim 5.6 More on 26 dimensional unimodular lattices. 

Unimodular lattices with dimension at most 23 all behave similarly and
have many nice properties. As the dimension is increased to 24, 25, or
26 the lattices behave more and more badly. For 26-dimensional ones we
will prove the only nice property I know of: the number of roots is
divisible by 4. We also show that for any 25-dimensional bimodular
lattice (not necessarily even) the number of norm 2 vectors is $2\bmod
4$, so that any 25-dimensional lattice of minimum norm at least 3 must
have determinant at least 3.

\proclaim Lemma 5.6.1. 
If $L$ is a 25-dimensional bimodular lattice then the number of norm 2
roots of $L$ is $2\bmod 4$.

Proof. The even vectors of $L$ form a lattice isomorphic to the
vectors that have even inner product with some vector $b$ in an even
25-dimensional bimodular lattice $B$. (Note that $b$ is not in
$B'-B$.) The number of roots of $B$ is $12t-10$ or $12t-18$ where $t$
is the height of the norm $-2$ vector of $D$ corresponding to $B$, so
it is sufficient to prove that the number of norm 2 vectors of $B$
that have odd inner product with $b$ is divisible by 4.

The projection of $w$ into $v^\perp$ is the Weyl vector $\rho$ of $B$
by 4.3.3,
so $b$ has integral inner product with $\rho$. Hence $b$ has even
inner product with the sum of the positive roots of $B$, so it has odd
inner product with an even number of positive roots. This implies that
the number of roots of $B$ that have odd inner product with $b$ is
divisible by 4. Q.E.D.

\proclaim Corollary 5.6.2. 
If $L$ is a 26 dimensional unimodular lattice then the number of norm
2 vectors of $L$ is divisible by 4.

Proof. The result is obvious if $L$ has no norm 2 roots, so let $r$ be
a norm 2 vector of $L$. $r^\perp$ is a 25 dimensional bimodular
lattice so by 5.6.1 the number of roots of $r^\perp$ is $2\bmod 4$. The
number of roots of $L$ not in $r^\perp$ is $4h-6$ where $h$ is the
Coxeter number of the component of $L$ containing $r$, so the number
of norm 2 vectors of $L$ is divisible by 4. Q.E.D.

{\bf Remark.} There are strictly 26 dimensional unimodular lattices
with no roots or with 4 roots, and there are strictly 27 dimensional unimodular
lattices with 0, 6, and 200 roots, so 5.6.2 is the best possible
congruence for the number of roots. For unimodular lattices of
dimension less than 26 there are congruences modulo higher powers of 2
for the number of roots.

\proclaim 5.7 A 27-dimensional unimodular lattice with no roots. 

Recall that the isomorphism classes of the
following sets are isomorphic: 
\item{(1)} 32 dimensional even unimodular lattices with root system $e_6$.
\item{(2)} 26 dimensional even trimodular lattices with no roots.
\item{(3)} 27 dimensional unimodular lattices with no roots 
but with characteristic vectors of norm 3. 

We will show that each of these sets has a unique member. (Note that
26 dimensional even trimodular lattices with roots can be classified
by finding all norm $-6$ vectors in $II_{25,1}$.) The automorphism
group of the unique lattice in (2) or (3) is a twisted Chevalley group
$2\times {}^3D_4(2).3$.

Let $L$ be a 27 dimensional unimodular lattice with no roots and a
characteristic vector $c$ of norm 3. $L$ has no vectors of norm 1 or 2
and exactly two characteristic vectors of norm 3, so its theta
function is determined and is $1+1640q^3+119574q^4+1497600q^5+\cdots$
and in particular $L$ has vectors of norm 5; let $v$ be such a
vector. $(v,c)$ is odd and cannot be $\pm3$ as then $v\pm c$ would
have norm 2, so $(v,c)=\pm 1$ and we can assume that
$(v,c)=1$. $\langle v,c\rangle$ is a lattice $V$ of determinant 14
such that $V'/V$ is generated by an element of norm $3/14$, so
$V^\perp$ is a 25 dimensional even lattice $X$ of determinant 14 such
that $X'/X$ is generated by an element of norm $1/14\bmod 2$. Such
$X$'s correspond to norm $-14$ vectors $x$ in the fundamental domain
$D$. The condition that $L$ has no vectors of norm 1 or 2 implies that
there are exactly two possibilities for $x$: $x$ is either the sum of
$w$ and a norm 0 vector of height 7 corresponding to $A_6^4$, or $x$
is the sum of $w$ and a norm $-2$ vector of height 6 corresponding to
the 25 dimensional bimodular lattice with root system $a_2^9$. Both of
these $x$'s turn out to give the same lattice $L$, which therefore has
two orbits of norm 5 vectors and is the unique lattice satisfying the
condition (3) above. The number of vectors in the two orbits of norm 5
vectors are 419328 and 1078272 and the groups fixing a vector from
these orbits have orders 3024 and $7^2.24$, so $\Aut(L)$ has order
$2^{12}.3^3.5^2.7^2.13$ which is the order of $2\times
{}^3D_4(2).3$. $\Aut(L)$ does not act irreducibly on $L$ because it
fixes the one dimensional subspace generated by the characteristic
vectors of norm 3.

{\bf Remark.}  In the Atlas of finite groups there is given a 27
dimensional representation $V$ of $2\times {}^3D_4(2).3$ which
contains a vector $c$ of norm 3 fixed by ${}^3D_4(2).3$ and a set of
819 vectors of norm 1 acted on by ${}^3D_4(2).3$ that have inner
product 1 with $c$. The lattice generated by $c$ and twice these
vectors of norm 1 is isomorphic to $L$. $V$ can be given the structure
of an exceptional Jordan algebra such that $L$ is closed under
multiplication (although with the Jordan algebra product given in the
Atlas $L$ is not closed).

Atlas coordinates for $L$: $V$ is the vector space of $3\times 3$
Hermitian matrices $y$
$$\pmatrix{a&C&\hat B\cr
\hat C&b&A\cr
B&\hat A&c\cr
}$$
over the real Cayley algebra with units $i_\infty, i_0,\ldots, i_6$ in
which $i_\infty\mapsto 1$, $i_{n+1}\mapsto i$, $i_{n+2}\mapsto j$,
$i_{n+4}\mapsto k$ generate a quaternion subalgebra. The inner product is
given by $Norm(y)=\sum Norm(y_{ij})$. $c$ is the identity matrix, and
$L$ is generated by $c$ and the 819 images of the norm 3 vectors
$$
\pmatrix{1&0&0\cr
0&-1&0\cr
0&0&-1\cr}
\qquad
\pmatrix{-1&0&0\cr
0&0&1\cr
0&1&0\cr}
\qquad
\pmatrix{0&s&s\cr
\hat s&-{1\over 2}&{1\over 2}\cr
\hat s&{1\over 2}&-{1\over 2}\cr}
$$
under the group generated by the maps taking $(a,b,c,A,B,C)$ to
$(a,b,c,iAi,iB,iC)$, $(b,c,a,B,C,A)$, and $(a,c,b,\hat A,-\hat C,
-\hat B)$ where $i=\pm i_n$ and $s=\sum_ni_n/4$.

{\bf Remark.} Let $L$ be the 27 dimensional lattice above. The number
of elements in $L/2L$ of norm $3\bmod 4$ is $2^{25}-2^{12}$. I have
calculated that the number of elements of norm $3\bmod 4$ in $L/2L$
represented by elements of norm 3 or 7 is less than this, so that $L$
contains an element $x$ of norm $3\bmod 4$ that is not congruent to
any element of norm 3 or 7 mod $2L$. Unfortunately this calculation
is rather messy (partly because elements of different orbits of norm 7
vectors can be congruent $\bmod 2L$ to different numbers of norm 7
vectors) and I have not yet thought of an independent way to check
it. Supposing this result, we can construct a lattice $M$ generated by
${1\over 2}(x\oplus i)$ and the elements of $L\oplus I$ that have even
inner product with $x\oplus i$ (where $I$ is generated by $i$ of norm
1). $M$ is a 28 dimensional unimodular lattice with no roots that
contains a characteristic norm 4 vector $c\oplus i$. By 1.7 this
implies that there is a 32 dimensional even unimodular lattice with
root system $d_5$, and again by 1.7 this implies the existence of a 27
dimensional unimodular lattice with no roots and no characteristic
vectors of norm 3.  (Remark added 1999: Bacher and Venkov have shown
that there are exactly two 27 dimensional unimodular lattices with no
roots and no characteristic vectors of norm 3, and there are 38
28-dimensional unimodular lattices with no roots.)

\proclaim Chapter 6 Automorphism  groups of Lorentzian lattices. 

(Remark added in  1999: Most of this chapter was published in 
Journal of Algebra, Vol. 111, No. 1, Nov 1987, 133--153.)

The study of automorphism groups of unimodular Lorentzian lattices
$I_{n,1}$ was started by Vinberg. These lattices have an infinite
reflection group (if $n\ge 2$) and Vinberg showed that the quotient of
the automorphism group by the reflection groups was finite if and only
if $n\le 19$. Conway and Sloane rewrote Vinberg's result in terms of
the Leech lattice $\Lambda$, showing that this quotient (for $n\le
19$) was a subgroup of $\cdot 0 = \Aut(\Lambda)$. In this paper we
continue Conway and Sloane's work and describe $\Aut(I_{n,1})$ for
$n\le 23$. In these cases there is a natural complex $U$ associated to
$I_{n,1}$, whose dimension is the virtual cohomological dimension of
the ``non-reflection part'' $G_n$ of $\Aut(I_{n,1})$, and which is a
point if and only if $n\le 19$. 

\proclaim 6.1 General properties of the automorphism  group. 

In this section we give some properties of the automorphism groups of
Lorentzian lattices that do not depend much on their dimension. The
results given here are not used later except in the examples.

$L$ is a Lorentzian lattice whose reflection group $R$ has fundamental
domain $D$. $G$ is the subgroup of index 2 of $\Aut(R)$ of
automorphisms fixing the two cones of norm 0 vectors. $G$ is a split
extension of $R$ by $\Aut(D)$. We write $H$ for the hyperbolic space
of $L$, including the points at infinity. We will say that a group has
the {\bf fixed point property} if any action of that group on a
hyperbolic space has a (finite or infinite) fixed point. (This is
similar to the property that any action of the group on a tree has a
fixed point.)

\proclaim Lemma 6.1.1.
\item{(1)} 
Any finite group has the fixed point property. 
\item{(2)} 
$\Z$ has the fixed point property.

Proof. 
\item{(1)} 
We embed the hyperbolic space in a Lorentzian space. If $x$ is any
point of hyperbolic space then the sum of all conjugates of $x$ under
the finite group is a point of Lorentzian space representing a point
of hyperbolic space fixed by the group.
\item{(2)} 
The finite and infinite points of a hyperbolic space form a compact
contractible manifold (with boundary), so any action of $\Z$ has a
fixed point.

(Remark added 1999. The original version of the thesis claimed that an
extension of groups with the fixed point property also has the fixed
point property, but this is false for the dihedral group acing on one
dimensional hyperbolic space.)

We will say that a group is {\bf virtually free abelian} if it has a
free abelian subgroup of finite index.

\proclaim Corollary 6.1.2. 
\item{(1)} 
A subgroup of $G$ fixes a finite point of $H$ if and only if it is finite.
\item{(2)}
A subgroup of $G$ fixes a point of $H$ only if it is virtually free
abelian. (Warning: This point of $H$ is not necessarily rational.)

Proof. By 6.1.1 any finite subgroup of $G$ fixes a finite point of
$H$.  Conversely the subgroup of $G$ fixing a finite point of $H$ is
finite and the subgroup of $G$ fixing an infinite point is a split
extension of $\Z^n$ by a finite group. Q.E.D.

\proclaim Corollary 6.1.3. 
Any infinite subgroup of $G$ or $\Aut(D)$ fixing an infinite point is
contained in a unique maximal virtually free abelian subgroup of $G$
or $\Aut(D)$ fixing a finite or infinite point.

Proof. This follows easily from 6.1.2. 

Example. If $L$ is $I_{24,1}$ then Aut (D) contains 76 classes of
maximal virtually free abelian subgroups corresponding to rational
norm 0 vectors in $L$, of ranks 1, 2, 3, 4, 5, 7, and 23, (and an
infinite number of other classes of maximal virtually free abelian
subgroups corresponding to irrational norm 0 vectors of $L$).

\proclaim 6.2.~Notation. 

We define notation for the rest of chapter 6. $L$ is $II_{25,1}$, with
fundamental domain $D$ and Weyl vector $w$. We identify the simple
roots of $D$ with the affine Leech lattice $\Lambda$.

$R$ and $S$ are two sublattices of $II_{25,1}$ such that $R^\perp=S$,
$S^\perp=R$, and $R$ is positive definite and generated by a nonempty
set of simple roots of $D$. The Dynkin diagram of $R$ is the set of
simple roots of $D$ in $R$ and is a union of $a$'s, $d$'s, and
$e$'s. The finite group $R'/R$ is naturally isomorphic to $S'/S$ (in
more than one way) because $II_{25,1}$ is unimodular, and subgroups of
$R'/R$ correspond naturally to subgroups of $S'$ containing $S$ (in
just one way). We fix a subgroup $G$ of $R'/R$ and write $T$ for the
subgroup of $S'$ corresponding to it, so that an element $s$ of $S'$
is in $T$ if and only if it has integral inner product with all
elements of $G$. It is $T$ that we will be finding the automorphism
group of in the rest of this paper. For each component $R_i$ of the
Dynkin diagram of $R$ the nonzero elements of the group $\langle
R_i\rangle'/\langle R_i \rangle$ can be identified with the tips of
$R_i$, and $R'/R$ is a product of these groups. Any automorphism of
$D$ fixing $R$ acts on the Dynkin diagram of $R$ and on $R'/R$, and
these actions are compatible with the map from tips of $R$ to
$R'/R$. In particular we can talk of the automorphisms of $D$ fixing
$R$ and $G$. We will also write $R$ for the Dynkin diagram of $R$. We
write $x'$ for the projection of any vector $x$ of $II_{25,1}$ into
$S$.

{\sl Example.} $\Lambda$ contains a unique orbit of $d_{25}$'s; let
$R$ be generated by a $d_n$ $(n\ge 2)$ contained in one of these
$d_{25}$'s. (Note that  $\Lambda$ contains two classes of $d_{16}$'s and
$d_{24}$'s.) Then $R^\perp=S$ is the even sublattice of $I_{25-n,1}$
and $R'/R$ has order 4. We can choose $G$ in $R'/R$ to have order 2 in
such a way that $T$ is $I_{25-n,1}$. We will use this to find ${\rm
Aut} (I_{m,1})$ for $m\le 23$.

We write $\Aut(T)$ for the group of automorphisms of $T$ induced by
automorphisms of $II_{25,1}$. This has finite index in the group of
all automorphisms of $T$ and is equal to this group in all the
examples of $T$ we give. This follows from the fact that if $g$ is an
automorphism of $T$ fixing $S$ and such that the automorphism of
$T/S$ it induces is induced by an automorphism of $R$ (under the
identification of $G$ with $S'/S$) then $g$ can be extended to an
automorphism of $II_{25,1}$.

\proclaim 6.3.~Some automorphisms of $T$.

The group of automorphisms of $D$ fixing $(R,G)$ obviously acts on
$T$. In this section we construct enough other automorphisms of $T$ to
generate $\Aut(T)$ and in the next few sections we find how these
automorphisms fit together. Recall that $\Lambda$ is the Dynkin diagram
of $D$ and $R$ is a spherical Dynkin diagram of $\Lambda$.

Let $r$ be any point of $\Lambda$ such that $r\cup R$ is a spherical
Dynkin diagram and write $R'$ for $R\cup r$. If $g'$ is any element of
$\cdot\infty$ such that $\sigma(R')g'$ fixes $r$ and $R$ then we
define an automorphism $g=g(r,g')$ of $II_{25,1}$ by
$g=\sigma(R)\sigma(R')g'$. (Recall that $\sigma(X)$ is the opposition
involution of $X$, which acts on the Dynkin diagram $X$ and acts as
$-1$ on $X^\perp$.) These automorphisms will turn out to be a sort of
generalized reflection in the sides of a domain of $T$.

\proclaim Lemma 6.3.1. 
\item{(1)} If $g=g(r,g')$ fixes the group $G$ then $g$ restricted to 
$T$ is an automorphism  of $T$. $g$ fixes $R$ and $G$ if and only if 
$\sigma(R')g'$ does.
\item{(2)} $g$ fixes the space generated by $r'$ and $w'$ and acts 
on this space as reflection in $r'^\perp$. 
\item{(3)} If $g'=1$ then $g$ acts on $T$ as reflection in $r'^\perp$. 
($g'$ can only be 1 if $\sigma(R')$ fixes $r$.)

{\sl Proof.} Both $\sigma(R)$ and $\sigma(R')g'$ exchange the two cones
of norm 0 vectors in $L$ and fix $R$, so $g$ fixes both the cones of
norm 0 vectors and $R$ and hence fixes $R^\perp=S$. If $g$ also fixes
$G$ then it fixes $T$ as $T$ is determined by $G$ and $S$. $\sigma(R)$
acts as $-1$ on $G$ and fixes $R$, so $g$ fixes $R$ and $G$ if and
only if $\sigma(R')g'$ does. This proves (1).

$\sigma(R')g'$ fixes $r$ and $R$ and so fixes $r'$. $\sigma(R)$ acts
as $-1$ on anything perpendicular to $R$, and in particular on $r'$,
so $g(r')=[\sigma(R)\sigma(R')g'](r')=-r'$. If $v$ is any vector of
$L$ fixed by $g'$ and $v'$ is its projection into $T$, then
$g(v)=\sigma(R)\sigma(R')v$ so that $g(v)-v$ is in the space generated
by $R'$ and hence $g(v')-v'$ is in the space generated by $r'$. As
$g(r')=-r'$, $g(v')$ is the reflection of $v'$ in $r'^\perp$. In
particular if $v=w$ or $g'=1$ then $g'$ fixes $v$ so $g$ acts on $v'$
as reflection in $r'^\perp$. This proves (2) and (3). Q.E.D.

\proclaim Lemma 6.3.2. 
\item{(0)} The subgroup of $\cdot\infty=\Aut(D)$ 
fixing $(R,G)$ maps onto the subgroup of $\Aut_+(T)$ fixing $w'$. 
\item{(1)} If $\sigma(R')$ fixes $(R,G)$ then $g(r,1)$ acts on $T$ 
as reflection in $r'^\perp$ and this is an automorphism  of $T$. 
\item{(2)} If $\sigma(R')$ does not fix $(R,G)$ then we define a map $f$ from 
a subset of $\Aut(D)$ to $\Aut(T)$ as follows:
\itemitem{}If $h$ in $\Aut(D)$ fixes all points of $R'$ then we put $f(h)=h$ 
restricted to $T$.
\itemitem{} If $h'$ in $\Aut(D)$ acts as $\sigma(R')$ on $R'$ then we put 
$f(h')=g(r,h')$ restricted to $T$. (This is a sort of twisted reflection 
in $r'^\perp$.)
\item{} Then the elements on which we have defined $f$ form a subgroup of 
$\Aut(D)$ and $f$ is an isomorphism from this subgroup to its image in 
$\Aut(T)$. $f(h)$ fixes $r'^\perp$ and $w'$ while $f(h')$ fixes $r'^\perp$ 
and acts on $w'$ as reflection in $r'^\perp$. (So $f(h)$ fixes the two 
half-spaces of $r'^\perp$ while $f(h')$ exchanges them.)

{\sl Proof.} Parts (0) and (1) follow from 6.3.1. 

If $h$ in $\Aut(D)$ fixes all points of $R'$ then it certainly fixes
all points of $R$ and $G$ and so acts on $T$. $h$ also fixes $w$ and
$w'$.

If $h'$ acts as $\sigma(R')$ on $R'$ then $\sigma(R')g'$ fixes all
points of $R$, so by 6.3.1(1), $g(r,h')$ is an automorphism of $T$, and
by 6.3.1(2), $g(r,h')$ maps $w'$ to the reflection of $w'$ in
$r'^\perp$. (There may be no such automorphisms $g'$, in which case the
lemma is trivial.) It is obvious that $f$ is defined on a subgroup of
$\Aut(D)$, so it remains to check that it is a homomorphism.

We write $h$, $i$ for elements of $\Aut(D)$ fixing all points of $R'$
and $h'$, $i'$ for elements acting as $\sigma(R')$. All four of these
elements fix $R'$ and so commute with $\sigma(R')$. $\quad h$, $i$,
$\sigma(R')h'$, and $\sigma(R')i'$ fix $R$ and so commute with
$\sigma(R)$. $\sigma(R)^2=\sigma(R')^2=1$. Using these facts it follows
that
$$\eqalign{
f(hi)&= hi=f(h)f(i)\cr
f(hi')&=\sigma(R)\sigma(R')hi' = h\sigma(R)\sigma(R')i'=f(h)f(i')\cr
f(h'i)&= \sigma(R)\sigma(R')h'i=f(h')f(i)\cr
f(h'i')&= h'i' \cr
&=\sigma(R)^2\sigma(R')^2h'i'\cr
&=\sigma(R)^2\sigma(R')h'\sigma(R')i'\cr
&=\sigma(R)\sigma(R')h'\sigma(R)\sigma(R')i'\cr
&=f(h')f(i')\cr
}$$
so $f$ is a homomorphism. Q.E.D.

\proclaim 6.4.~Hyperplanes of $T$. 

We now consider the set of hyperplanes of $T$ of the form $r'^\perp$,
where $r$ is a root of $II_{25,1}$ such that $r'$ has positive
norm. These hyperplanes divide the hyperbolic space of $T$ into
chambers and each chamber is the intersection of $T$ with some chamber
of $II_{25,1}$. We write $D'$ for the intersection of $D$ with $T$.

In the section we will show that $D'$ is often a sort of fundamental
domain with finite volume. It is rather like the fundamental domain of
a reflection group, except that it has a nontrivial group acting on
it, and the automorphisms of $T$ fixing sides of $D'$ are more
complicated than reflections.

\proclaim Lemma 6.4.1. $D'$ contains $w'$ in its interior and, in particular, is 
nonempty. 

{\sl Proof.} To show that $w'$ is in the interior of $D'$ we have to
check that no hyperplane $r^\perp$ of the boundary of $D$ separates
$w$ and $w'$, unless $r$ is in $R$. $r$ is a simple root of $D$ with
$-(r,w)=1$ so it is enough to prove that $(r,\rho)\ge 0$, where
$\rho=w-w'$ is the Weyl vector of the lattice generated by
$R$. $-\rho$ is a sum of simple roots of $R$, so $(r,\rho)\ge 0$
whenever $r$ is a simple root of $D$ not in $R$ because all such
simple roots have inner product $\le 0$ with the roots of $R$. This
proves that $w'$ is in the interior of $D'$. Q.E.D.

\proclaim Lemma 6.4.2. 
The faces of $D'$ are the hyperplanes $r'^\perp$, where $r$ runs
through the simple roots of $D$ such that $r\cup R$ is a spherical
Dynkin diagram and $r$ is not in $R$. In particular $D'$ has only a
finite number of faces because $R$ is not empty.

{\sl Proof.} The faces of $D'$ are the hyperplanes $r'^\perp$ for the
simple roots $r$ of $D$ such that $r'^\perp$ has positive norm, and
these are just the simple roots of $D$ with the property in 6.4.2. $D'$
has only a finite number of faces because the Leech lattice
(identified with the Dynkin diagram of $II_{25,1}$) has only a finite
number of points at distance at most $\sqrt 6$ from any given point in
$R$. Q.E.D.

\proclaim Lemma 6.4.3. 
If $R$ does not have rank 24 (i.e., $T$ does not have dimension 2)
then $D'$ has finite volume.

{\sl Proof.} $D'$ is a convex subset of hyperbolic space bounded by a
finite number of hyperplanes, and this hyperbolic space is not one
dimensional as $T$ is not two dimensional, so $D'$ has finite volume
if and only if it contains only a finite number of infinite
points. $R$ is nonempty so it contains a simple root $r$ of $D$. The
points of $D'$ at infinity correspond to some of the isotropic
subspaces of $II_{25,1}$ in $r^\perp$ and $D$. The hyperplane
$r^\perp$ does not contain $w$ as $(r,w)=-1$, so the fact that $D'$
contains only a finite number of infinite points follows from 6.4.4
below (with $V=r^\perp$). Q.E.D.

\proclaim Lemma 6.4.4. 
If $V$ is any subspace of $II_{25,1}$ not containing $w$ then $V$
contains only a finite number of isotropic subspaces that lie in $D$.

{\sl Proof.} Let $II_{25,1}$ be the set of vectors $(\lambda,m,n)$
with $\lambda$ in $\Lambda$, $m$ and $n$ integers, with the norm given
by $(\lambda,m,n)^2=\lambda^2-2mn$. We let $w$ be $(0,0,1)$ so that
the simple roots are $(\lambda,1,\lambda^2/2-1)$. As $V$ does not
contain $w$ there is some vector $r=(v,m,n)$ in $V^\perp$ with
$(r,w)\ne 0$, i.e., $m\ne 0$. We let each norm 0 vector $z=(u,a,b)$
which is not a multiple of $w$ correspond to the point $u/a$ of
$\Lambda\otimes \Q$. If $z$ lies in $V$ then $(z,r)=0$, so
$(u/a-v/m)^2=(r/m)^2$, so $u/a$ lies on some sphere in $\Lambda\otimes
Q$. If $z$ is in $D$ then $u/a$ has distance at least $\sqrt 2$ from
all points of $\Lambda$ (i.e., it is a ``deep hole''), but as
$\Lambda$ has covering radius $\sqrt 2$ these points form a discrete
set so there are only a finite number of them on any sphere. Hence
there are only a finite number of isotropic subspaces lying in $V$ and
$D$. Q.E.D.

{\sl Remark.} If $w$ is in $V$ then the isotropic subspaces of
$II_{25,1}$ in $V$ and $D$ correspond to deep holes of $\Lambda$ lying
on some affine subspace of $\Lambda\otimes \Q$. There is a universal
constant $n_0$ such that in this case $V$ either contains at most
$n_0$ isotropic subspaces in $D$ or contains an infinite number of
them. If $w$ is not in $V$ then $V$ can contain an arbitrarily large
number of isotropic subspaces in $D$.

\proclaim 6.5.~A complex. 

We have constructed enough automorphisms to generate $\Aut(T)$, and
the problem is to fit them together to give a presentation of
$\Aut(T)$. We will do this by constructing a contractible complex
acted on by $\Aut(T)$. For example, if this complex is one dimensional
it is a tree, and groups acting on trees can often be written as
amalgamated products.

{\sl Notation.} We write $\Aut_+(T)$ for the group of automorphisms of
$T$ induced by $\Aut_+(II_{25,1})$. $D_r$ is a fundamental domain of
the reflection subgroup of $\Aut_+(T)$ containing $D'$. The hyperbolic
space of $T$ is divided into chambers by the conjugates of all
hyperplanes of the form $r'^\perp$ for simple roots $r$ of $D$.

\proclaim Lemma 6.5.1. 
Suppose that for any spherical Dynkin diagram $R'$ containing $R$ and
one extra point of $\Lambda$ there is an element of $\cdot\infty$
acting as $\sigma(R')$ on $R'$. Then $\Aut_+(T)$ acts transitively on
the chambers of $T$ and $\Aut(D_r)$ acts transitively on the chambers
of $T$ in $D_r$.

{\sl Proof. } By lemma 6.3.2 there is an element of $\Aut_+(T)$ fixing
any face of $D'$ corresponding to a root $r$ of $D$ and mapping $D'$
to the other side of this face. Hence all chambers of $T$ are
conjugates of $D'$. Any automorphism of $T$ mapping $D'$ to another
chamber in $D_r$ must fix $D_r$, so $\Aut(D_r)$ acts transitively on
the chambers in $D_r$. Q.E.D.

$D_r$ is decomposed into chambers of $T$ by the hyperplanes
$r'^\perp$ for $r$ a root of $II_{25,1}$. We will write $U$ for the
dual complex of this decomposition and $U'$ for the subdivision of
$U$. $U$ has a vertex for each chamber of $D_r$, a line for each pair
of chambers with a face in common, and so on. $U'$ is a simplicial
complex with the same dimension as $U$ with an $n$-simplex for each
increasing sequence of $n+1$ cells of $U$. $U$ is not necessarily a
simplicial complex and need not have the same dimension as the
hyperbolic space of $T$; in fact it will usually have dimension 0, 1,
or 2. For example, if $D_r=D'$ then $U$ and $U'$ are both just points.

\proclaim Lemma 6.5.2. $U$ and $U'$ are contractible. 

{\sl Proof. } $U$ is contractible because it is the dual complex of
the contractible space $D_r$. ($D_r$ is even convex.) $U'$ is
contractible because it is the subdivision of $U$. Q.E.D.

\proclaim Theorem 6.5.3. 
Suppose that $\Aut(D_r)$ acts transitively on the maximal simplexes of
$U'$, and let $C$ be one such maximal simplex. Then $\Aut(D_r)$ is the
sum of the subgroups of $\Aut(D_r)$ fixing the vertices of $C$
amalgamated over their intersections.

{\sl Proof.} By 6.5.2, $U'$ is connected and simply connected. $C$ is
connected and by assumption is a fundamental domain for $\Aut(D_r)$
acting on $U'$. By a theorem of Macbeth ([S b, p. 31]) the group
$\Aut(D_r)$ is given by the following generators and relations:

{\sl Generators:} An element $\hat g$ for every $g$ in $\Aut(D_r)$
such that $C$ and $g(C)$ have a point in common.

{\sl Relations:} For every pair of elements $(s,t)$ of $\Aut(D_r)$
such that $C$, $s(C)$, and $u(C)$ have a point in common (where
$u=st$) there is a relation $\hat s\hat t=\hat u$.

Any element of $\Aut(D_r)$ fixing $C$ must fix $C$ pointwise. This
implies that $C$ and $g(C)$ have a point in common if and only if $g$
fixes some vertex of $C$, i.e., $g$ is in one of the groups $C_0$,
$C_1,\ldots$ which are stabilizers of the vertices of $C$, so we have
a generator $\hat g$ for each $g$ that lies in (at least) one of these
groups. There is a point in all of $C$, $s(C)$, and $u(C)$ if and only
if some point of $C$, and hence some vertex of $C$, is fixed by $s$
and $t$. This means that we have a relation $\hat s\hat t=\hat u$
exactly when $s$ and $t$ both lie in some group $C_i$. This is the
same as saying that $\Aut(D_r)$ is the sum of the groups $C_i$
amalgamated over their intersections. Q.E.D.

{\sl Example.} If $C$ is one dimensional then $\Aut(D_r)$ is the free
product of $C_0$ and $C_1$ amalgamated over their intersection. If the
dimension of $C$ is not 1 then $\Aut(D_r)$ cannot usually be written
as an amalgamated product of two nontrivial groups.

\proclaim 6.6.~Unimodular lattices. 

In this section we apply the results of the previous section to find
the automorphism group of $I_{m,1}$ for $m\le 23$.

\proclaim Lemma 6.6.1. 
Let $X$ be the Dynkin diagram $a_n$ ($1\le n\le 11$), $d_n$ ($2\le n\le
11$), $e_n$ ($3\le n\le 8$), or $a_2^2$ which is contained in
$\Lambda$. Then any automorphism of $X$ is induced by an element of
$\cdot\infty$. If $X$ is an $a_n$ ($1\le n\le 10$), $d_n$ ($2\le n\le
25$, $n\ne 16$ or $24$), $e_n$ ($3\le n\le 8$), or $a_2^2$ then
$\cdot\infty$ acts transitively on Dynkin diagrams of type $X$ in
$\Lambda$.

{\sl Proof.} A long, unenlightening calculation. See section 6.9. Q.E.D.

{\sl Remark.} $\cdot\infty$ acts simply transitively on ordered
$a_{10}$'s in $\Lambda$. There are two orbits of $d_{16}$'s and
$d_{24}$'s (see section 6.9 and example 2 of section 6.8) and many orbits
of $a_n$'s for $n\ge 11$.

{\sl Notation.} We take $R$ to be a $d_n$ contained in a $d_{25}$ for
some $n$ with $2\le n\le 23$. If $R$ is $d_4$ we label the tips of $R$
as $x$, $y$, $z$ in some order, and if $R$ is $d_n$ for $n\ne 4$ we
label the two tips that can be exchanged by an automorphism of $R$ as
$x$ and $y$. If $n=3$ or $n\ge 5$ we label the third tip of $R$ as
$z$. We let $G$ be the subgroup of $\langle R'\rangle/\langle R
\rangle $ of order 2 which corresponds to the tip $z$ if $n\ge 3$ and
to the sum of the elements $x$ and $y$ if $n=2$. An automorphism of
$R$ fixes $G$ if $n\ne 4$ or if $n=4$ and it fixes $z$. The lattice
$S=R^\perp$ is isomorphic to the even sublattice of $I_{25-n,1}$. We
let $T$ be the lattice corresponding to $G$ that contains $S$, so that
$T$ is isomorphic to $I_{25-n,1}$.

\proclaim Lemma 6.6.2. 
Any root $r$ of $V$ such that $r\cup R$ is a spherical Dynkin diagram
is one of the following types:
\itemitem{Type $a$:} $r\cup R$ is $d_na_1$ (i.e., $r$ is not joined to any 
point of $R$.) $r'$ is then a norm 2 vector of $T$. 
\itemitem{Type $d$:} 
$r$ is joined to $z$ if $n\ge 3$ or to $x$ and $y$ if $n=2$, so that
$r\cup R$ is $d_{n+1}$. $r'$ has norm 1.
\itemitem{Type $e$:} 
$r$ is joined to just one of $x$ or $y$, so that $r\cup R$ is $e_{n+1}$. 
$2r'$ is then a characteristic vector of norm $8-n$ in $T$. 

{\sl Proof.} Check all possible cases. Q.E.D.

In particular if $r$ is of type $a$ or $d$, or of type $e$ with $n=6$
or 7, then $r'$ (or $2r'$) has norm 1 or 2 and so is a root of
$T$. Note that in these cases $\sigma(R\cup r)$ fixes $R$ and $z$ and
therefore $G$, so by 6.3.2(2), $r'^\perp$ is a reflection of $T$. In the
remaining four cases ($r$ of type $e$ with $2\le n\le 5$)
$\sigma(r\cup R)$ does not fix both $R$ and $G$.

\proclaim Corollary 6.6.3. 
If $n\ge 6$ then then reflection group of $T=I_{25-n,1}$ has finite
index in $\Aut(T)$. Its fundamental domain has finite volume and a
face for each root $r$ of $\Lambda$ such that $r\cup R$ is a spherical
Dynkin diagram.

{\sl Proof.} This follows from the fact that all walls of $D'$ give
reflections of $T$ so $D'$ is a fundamental domain for the reflection
group. By 6.4.3, $D'$ has finite volume. Q.E.D.

{\sl Remarks.} The fact that a fundamental domain for the reflection
group has finite volume was first proved by Vinberg [Vi] for $n\ge 8$
and by Vinberg and Kaplinskaja [V-K] for $n=6$ and $7$ (i.e., for
$I_{18,1}$ and $I_{19,1}$). Conway and Sloane [C-S d] show implicitly that
for $n\ge 6$ the non-reflection part of $\Aut(T)$ is the subgroup of
$\cdot\infty$ fixing $R$ and their description of the fundamental
domain of the reflection group in these cases is easily seen to be
equivalent to that in 6.6.3. $I_{18,1}$ and $I_{19,1}$ have a ``second
batch'' of simple roots of norm 1 or 2; from 6.6.2 we see that this
second batch consists of the roots which are characteristic vectors
and they exist because of the existence of $e_8$ and $e_7$ Dynkin
diagrams. The non-reflection group of $I_{20,1}$ is infinite because of
the existence of $e_6$ Dynkin diagrams and the fact that the
opposition involution of $e_6$ acts non-trivially on the
$e_6$. ($\dim(I_{20,1})=1+\dim(II_{25,1})+\dim(e_6)$.)

From now on we assume that $n$ is 2, 3, 4, or 5. Recall that $D_r$ is
the fundamental domain of the reflection group of $T$ containing $D$
and $U$ is the complex which is the dual of the complex of conjugates
of $D$ in $D_r$.

\proclaim Lemma 6.6.4. $\Aut(D_r)$ acts transitively on the vertices of $U$. 

{\sl Proof. } This follows from 6.6.1 and 6.5.1. Q.E.D.

\proclaim Lemma 6.6.5. 
If $n=3$, $4$, or $5$ then $U$ is one dimensional and if $n=2$, then $U$
is two dimensional. (By the remarks after 6.6.3, $U$ is zero
dimensional for $n\ge 6$.)

{\sl Proof.} Let $s$ and $t$ be simple roots of $\Lambda$ such that
$s\cup R$ and $t\cup R$ are $e_{n+1}$'s, so that $s'^\perp$ and
$t'^\perp$ give two faces of $D'$ inside $D_r$. Suppose that these
faces intersect inside $D_r$. We have $s'^2=t'^2=2-n/4$, $(s',t')\le
0$, and $r=s'+t'$ lies in $T$ (as $2s'$ and $2t'$ are both
characteristic vectors of $T$ and so are congruent mod $2T$). $r'^2$
cannot be 1 or 2 as then the intersection of $s'^\perp$ and $t'^\perp$
would lie on the reflection hyperplane $r^\perp$, which is impossible
as we assumed that $s'^\perp$ and $t'^\perp$ intersected somewhere in
the interior of $D_r$. Hence
$$ 3\le r^2=(s'+t')^2\le 2(2-n/4)\le 2(2-2/4)=3$$
so $r^2=3$, $n=2$, and $(s',t')=0$. 

If $n=3$, 4, or 5 this shows that no two faces of $D'$ intersect in
the interior of $D_r$, so the graph whose vertices are the conjugates
of $D'$ in $D_r$ such that two vertices are joined if and only if the
conjugates of $D'$ they correspond to have a face in common is a
tree. As it is the 1-skeleton of $U$, $U$ must be one dimensional.

If $n=2$ then it is possible for $s'^\perp$ and $t'^\perp$ to
intersect inside $D_r$. In this case they must intersect at right
angles, so $U$ contains squares. However, in this case $s$ and $t$
cannot be joined to the same vertex of $R=a_1^2$, and in particular it
is not possible for three faces of $D'$ to intersect inside $D_r$, so
$U$ is two dimensional. (Its two-dimensional cells are squares.)
Q.E.D.

If $X$ is a Dynkin diagram in $\Lambda$ we will write $G(X)$ for the
subgroup of $\cdot\infty$ fixing $X$. If $X_1\subset X_2\subset
X_3\cdots$ is a sequences of Dynkin diagrams of $\Lambda$ we write
$G(X_1\subset X_2\subset X_3\cdots)$ for the sum of the groups
$G(X_i)$ amalgamated over their intersections.

\proclaim Theorem 6.6.6. 
The non-reflection part of $\Aut_+(T)=\Aut_+(I_{25-n,1})$ is given by
$$
\eqalign{
G(d_5\subset e_6) & \hbox{ if } n=5,\cr
G(d_4^*\subset d_5) & \hbox{ if } n=4 
\hbox{ (where $d_4^*$ means that one of the tips of the $d_4$ is labeled)},\cr
G(a_3\subset a_4) & \hbox{ if } n=3,\cr
G(a_1^2\subset a_1a_2\subset a_2^2) & \hbox{ if } n=2.\cr
}
$$

{\sl Proof.} By 6.6.4, $\Aut(D_r)$ acts transitively on the vertices of
$U$.  Using this and 6.6.1 it is easy to check that it acts transitively
on the maximal flags of $U$, or equivalently on the maximal simplexes
of the subdivision $U'$ of $U$. For example, if $n=5$ this amounts to
checking that the group $G(d_5)$ acts transitively on the $e_6$'s
containing a $d_5$.

By 6.5.3 the group $\Aut(D_r)$ is the sum of the groups fixing each
vertex of a maximal simplex $C$ of $U'$ amalgamated over their
intersections. By 6.3.2 the subgroup of $\Aut(D_r)$ fixing the vertex
$D'$ of $U$ can be identified with $G(R)=G(d_n)$ and the subgroup
fixing a face of $D'$ can be identified with $G(e_{n+1})$ for the
$e_{n+1}$ corresponding to this face. These groups are the groups
fixing two of the vertices of $C$, and if $n=2$ it is easy to check
that the group fixing the third vertex can be identified with
$G(a_2^2)$. Hence $\Aut(D_r)$ is $G(d_n\subset e_{n+1})$ if $n=3$, 4,
or 5 (where if $n=4$ we have to use a subgroup of index 3 in
$G(d_4)$), and $G(d_2\subset e_3\subset a_2^2)$ if $n=2$. Q.E.D.

{\sl Examples.} $n=5$: $\Aut(I_{20,1})$ The domain $D'$ has 30 faces of
type $a$, 12 of type $d$, and 40 of type $e$. $\Aut(D')$ is
$\Aut(A_6)$ of order 1440 (where $A_6$ is the alternating group of
order 360). $\Aut(D_r)$ is
$$\Aut(A_6)*_H(S_3\wr Z_2),$$
where $S_3\wr Z_2=G(e_6)$ is a wreathed product and has order
72. $H$ is its unique subgroup of order 36 containing an element of
order 4. (Warning: $\Aut(A_6)$ contains two orbits of subgroups
isomorphic to $H$. The image of $H$ in $\Aut(A_6)$ is not contained in
a subgroup of $\Aut(A_6)$ isomorphic to $S_3\wr Z_2$.) $\Aut(D_r)$ has Euler
characteristic $1/1440+1/72-1/36=-19/1440$.

$n=4$: $\Aut(I_{21,1})$. $D'$ has 42 faces of type $a$, 56 of type
$d$, and 112 of type $e$. $\Aut(D')$ is $L_3(4)\cdot 2^2$. $\Aut(D_r)$
is
$$L_3(4).2^2 *_{M_{10}}\Aut(A_6).$$
$\Aut(A_6)$ has 3 subgroups of index 2, which are $S_6$, $PSL_2(9)$, and 
$M_{10}$. $\Aut(D_r)$ has Euler characteristic $-11/2^8.3^2.7.$

$n=3$: $\Aut(I_{22,1})$. $D'$ has 100 faces of type $a$, 1100 of type
$d$, and 704 of type $e$. $\Aut(D')$ is $HS.2$, where $HS$ is the
Higman-Sims simple group, and $\Aut(D_r)$ is
$$HS.2*_HH.2$$
where $H$ is $PSU_3(5)$ of order $2^4.3^2.5^3.7$. Its Euler
characteristic is $3.13/2^{10}.5^5.7.11$.

$n=2$: $\Aut(I_{23,1})$. $D'$ has 4600 faces of type $a$, 953856 of
type $d$, and 94208 of type $e$. $\Aut(D')$ is $\cdot 2\times2$, where
$\cdot 2$ is one of Conway's simple groups.  $\Aut(D_r)$ is the direct
limit of the groups
$$
\matrix
{
&&&\cdot2\times 2&=&G(a_1^2)\cr
&&\nearrow&&\nwarrow\cr
&McL&\leftarrow&PSU_4(3)&\rightarrow&PSU_4(3).2\cr
&\downarrow&&\downarrow&&\downarrow\cr
G(a_1a_2)=&McL.2&\leftarrow&PSU_4(3).2
  &\rightarrow&PSU_4(3)\cdot D_8&=G(a_2^2)\cr
}
$$
$McL$ is the McLaughlin simple group and $D_8$ is the dihedral group
of order 8. The direct limit is generated by $\cdot 2\times 2$ and an
outer automorphism of $McL$; the $PSU_4(3)$'s are there to supply one
additional relation. The Euler characteristic of $\Aut(D_r)$ is the
sum of the reciprocals of the orders of the groups in the center and
the vertices of the diagram above minus the sum of the reciprocals of
the groups on the edges (see [S c]), which is
$3191297/2^{20}.3^6.5^3.7.11.23$.

{\sl Remark.} For $n=2$, 3, 4, 5, or 7 the number of faces of $D'$ of
type $e$ is $(24/(n-1)-1)2^{12/(n-1)}$, which is $23\cdot2^{12}$,
$11\cdot2^6$, $7\cdot2^4$, $5\cdot 2^3$, or $3\cdot 2^2$. For $n=6$
this expression is $20.06\ldots$ and there are 20 faces of type
$e$. Table 3 of [C-S d] gives the number of faces of
type $a$ and $d$ of $I_{m,1}$ for $m\le 23$.

\proclaim 6.7.~More about $I_{n,1}$. 

Here we give more information about $I_{n,1}$ for $20\le n\le 23$. In
the tables in [C-S d] the heights of the simple roots
they calculate appear to lie on certain arithmetic progressions; we
prove that they always do. We then prove that the dimension of the
complex $U$ is the virtual cohomological dimension of the
non-reflection part of $\Aut(I_{n,1})$.

{\sl Notation.} $D$ is a $d_n$ of $\Lambda$ contained in a
$d_{25}$. Let $C$ be the sublattice of all elements of $T=I_{25-n,1}$
which have even inner product with all elements of even norm. $C$
contains $2T$ with index 2 and the elements of $C$ not in $2T$ are the
characteristic vectors of $T$. $w'$ is the projection of $w$ into $T$
and $D_r$ is the fundamental domain of the reflection group of $T$
containing $w'$.

\proclaim Lemma 6.7.1. 
Suppose $2\le n\le 5$ Then all conjugates of $w'$ in $D_r$ are
congruent $\bmod 2^{n-2}C$. (If $n\ge 6$ then $D_r=D'$ so there are no
other conjugates of $w'$ in $D_r$.)

{\sl Proof.} It is sufficient to prove that any two conjugates of $w'$
which are joined as vertices of the graph of $D_r$ are congruent
$\bmod 2^{n-2}C$ because this graph is connected, and we can also
assume that one of these vertices is $w'$ because $\Aut(D_r)$ acts
transitively on its vertices. Let $w''$ be a conjugate of $w'$ joined
to $w'$.

$w''$ is the reflection of $w'$ in some hyperplane $e^\perp$, where
$e$ is a characteristic vector of $T$ of norm $8-n$, so $e$ is in $C$.
Therefore
$$w''=w'-2(w',e)e/(e,e).$$
$e=2r'$ for a vector $r$ of $\Lambda$ such that $r\cup R$ is an
$e_{n+1}$ diagram. The projections of $r$ and $w$ into the lattice
$I^n$ containing $R$ are $({1\over 2},{1\over 2},\ldots,{1\over 2})$
and $(0,1,\ldots,n-1)$, which have inner product $(n-1)n/4$. Hence
$(w',r')=(w,r)-(\hbox{inner product of projections of $w$ and $r$ into
$\langle R\rangle$})=-1-(n-1)n/4$, so
$$-2(w',e)/(e,e)=(4+(n-1)n)/(8-n).$$
For $2\le n\le 5 $ the expression on the right is equal to $2^{n-2}$, so 
$w'$ and $w''$ are congruent ${} \bmod 2^{n-2}C$. Q.E.D.

\proclaim Theorem 6.7.2. 
Suppose $n=2$, $3$, $4$, or $5$ and let $r$ be a simple root of the
fundamental domain $D_r$ of $T=I_{25-n,1}$.
\item{} If $r^2=1$ then $(r,w')\equiv -n\bmod 2^{n-2}$. 
\item{} If $r^2=2$ then $(r,w')\equiv -1\bmod 2^{n-1}$. 

{\sl Proof.} $(r,w')=(s',w'')$ for some conjugate $w''$ of $w$ in
$D_r$ and some simple root $s'$ of $I_{25-n,1}$ that is the projection
of a simple root $s$ of $D$ into $I_{25-n,1}$. By 6.7.1, $w'$ is
congruent to $w'' \bmod 2^{n-2}C$ so $(r,w')$ is congruent to
$(s',w')\bmod 2^{n-2}$ if $r^2=1$ and ${} \bmod 2^{n-1}$ if $r^2=2$
(because elements of $C$ have even inner product with all elements of
norm 2). $(s',w')$ is equal to $(s,w)-$(inner product of projections
of $s$ and $w$ into $R$), which is $-1$ if $s'$ has norm 2 and
$-1-(n-1)$ if $s'$ has norm 1. Q.E.D.

This explains why the heights of the simple roots for $I_{m,1}$ with
$20\le n\le 23$ given in table 3 of [C-S d] seem to lie
on certain arithmetic progressions.

Now we show that the dimension of the complex $U$ is the virtual
cohomological dimension of $\Aut(D_r)$.

\proclaim Lemma 6.7.3. 
The cohomological dimension  of any torsion free subgroup of
$\Aut(D_r)$ is at most $\dim(U)$. 

{\sl Proof.} Any such subgroup acts freely on the contractible complex $U$. 
Q.E.D.

\proclaim Corollary 6.7.4. 
The virtual cohomological dimension  of $\Aut(D_r)$ is equal to $\dim(U)$.

{\sl Proof.} $\Aut(D_r)$ contains torsion-free subgroups of finite
index, so by 6.7.3 the virtual cohomological dimension of $\Aut(D_r)$
is at most $\dim(U)$. For $n\le 5$, $\Aut(D_r)$ is infinite and so has
virtual cohomological dimension at least 1, while for
$n=2$. $\Aut(D_r)$ contains subgroups isomorphic to $Z^2$ (because
there are 22-dimensional unimodular lattices whose root systems
generate a vector space of codimension 2) so $\Aut(D_r)$ has virtual
cohomological dimension at least 2. Q.E.D.

Lemma 6.7.3 implies that $\Aut(D_r)$ contains no subgroups of the form
$Z^i$ with $i>\dim(U)$, and this implies that if $L$ is a
$(24-n)$-dimensional unimodular lattice then the space generated by
roots of $L$ has codimension at most $\dim(U)$. This can of course
also be proved by looking at the list of such lattices. (Vinberg used
this in reverse: he showed that the non-reflection part of
$\Aut(I_{m,1})$ was infinite for $m\ge 20$ from the existence of
19-dimensional unimodular lattices with root systems of rank 18. There
are two such lattices, with root systems $a_{11}d_7$ and $e_6^3$; they
are closely related to the two Niemeier lattices $a_{11}d_7e_6$ and
$e_6^4$ containing an $e_6$ component.)

\proclaim 6.8.~Other examples. 

We list some more examples of (not necessarily unimodular) Lorentzian
lattices with their automorphism groups.

{\sl Example 1.} $R$ is an $e_8$ in $\Lambda$ so that $T$ is
$II_{17,1}$. The fundamental domain $D'$ of the reflection group has
finite volume and its Dynkin diagram is the set of points of $\Lambda$
not connected to $e_8$, which is a line of 17 points with 2 more
points joined onto the 3rd and 15th points. This diagram was found by
Vinberg [Vi].

{\sl Example 2.} Similarly if $R$ is one of the $d_{16}$'s of
$\Lambda$ not contained in a $d_{17}$ then $T$ is $II_{9,1}$ and the
points of $\Lambda$ not joined to $R$ form an $e_{10}$ which is the
Dynkin diagram of $II_{9,1}$.

{\sl Example 3.} All $e_7$'s of $\Lambda$ are conjugate; if $R$ is one
of them then $T$ is the 19-dimensional even Lorentzian lattice of
determinant 2. There are 3+21 roots $r$ for which $r\cup R$ is a
spherical Dynkin diagram and these 24 points are arranged as a ring of
18 points with an extra point joined on to every third point. The
three roots joined to the $e_7$ correspond to norm 2 roots $r$ of $T$
with $r^\perp$ unimodular, while the other 21 roots correspond to norm
2 roots of $T$ such that $r^\perp$ is not unimodular. The
non-reflection part of $\Aut_+(T)$ is $S_3$ of order 6 acting in the
obvious way on the Dynkin diagram. (Remark added 1999: this example
was first found by Nikulin.)

{\sl Example 4.} Let $R$ be the unique orbit of $e_6$'s in
$\Lambda$. Then $T$ is the 20-dimensional even Lorentzian lattice of
determinant 3. $D'$ is a fundamental domain for the reflection group of
$T$ and has $12+24$ faces, coming from 12 roots of norm 6 and 24 of
norm 2. The non-reflection group of $\Aut_+(T)$ is a wreath product
$S_3\wr Z_2$ of order 72. (Remark added in 1998: this example was first
found by Vinberg in ``The two most algebraic $K3$ surfaces'',
Math. Ann. 265 (1983), no. 1, 1--21.)

{\sl Example 5.} $R$ is $d_4$ and $T$ is the even sublattice of
$I_{21,1}$ so that $T$ has determinant 4. The domain $D'$ is a
fundamental domain for the reflection group of $T$ and has 168 walls
corresponding to roots of norm 4 and 42 walls corresponding to roots
of norm 2. $\Aut(D)$ is isomorphic to $L_3(4).D_{12}$ of order
$2^8.3^3.5.7$. $T$ is a 22-dimensional Lorentzian lattice whose
reflection group has finite index in its automorphism group; I do not
know of any other such lattices of dimension $\ge 21$. (Remark added
in 1998: Esselmann recently proved in ``\"Uber die maximale Dimension
von Lorentz-Gittern mit coendlicher Spiegelungsgruppe'', Number Theory
61 (1996), no. 1, 103--144, that the lattice $T$ is essentially the
only example of a Lorentzian lattice of dimension at least 21 whose
reflection group has finite index in its automorphism group.) $T$ is
contained in three lattices isomorphic to $I_{21,1}$ each of whose
automorphism groups has index 3 in $\Aut(T)$. However the reflection
groups of these lattices do not have finite index in their
automorphism groups.

{\sl Example 6.} $R$ is $a_1$ and $T$ is the 25-dimensional even
Lorentzian lattice of determinant 2. This time $D'$ is not a
fundamental domain for the reflection group. It has 196560 faces
corresponding to norm 2 roots and 16773120 faces perpendicular to norm
6 vectors (which are not roots). However, the simplicial complex of $T$
is a tree so $\Aut(D')$ is $(\cdot 0)*_{(\cdot 3)}(2\times \cdot 3)$,
i.e., it is generated by $\cdot 0$ and an element $t$ of order 2 with
the relations that $t$ commutes with some $\cdot 3$ of $\cdot 0$. $D'$
has finite volume but if any of its 16969680 faces are removed the
resulting polyhedron does not!

\proclaim 6.9.~The automorphism  groups of high-dimensional Lorentzian lattices. 

{\sl Notation.} $L$ is $II_{8n+1,1}$ $(n\ge 1)$ and $D$ is a
fundamental domain of the reflection group of $L$.

$X$ is the Dynkin diagram of $D$. In this section we will show that if
$8n\ge 24$ then $\Aut(L)$ acts transitively on many subsets of $X$,
and use this to generalize some of the results of the previous
sections to higher dimensional lattices.

\proclaim Lemma 6.9.1. 
Let $R$ be a spherical Dynkin diagram. Suppose that whenever $R'$ is a
spherical Dynkin diagram in $X$ which is isomorphic to $R$ plus one
point $r$ there is an element $g$ of $\Aut(D)$ such that $g\sigma(R')$
fixes $R$ (resp. fixes all points of $R$).
\item{} For any map $f:R\mapsto X$ we construct $(M,f',C)$, where
\item{$M$}
is the lattice $f(R)^\perp$,
\item{$f'$} 
is the map from $R'/R$ to $M'/M$ such that $f(r)\equiv-f'(r)\bmod L$
for $r$ in $R'/R$,
\item{$C$}
is the cone of $M$ contained in the cone of $L$ containing $D$.
\item{} If $f_1$, $f_2$ are two such maps then the images $f_1(R)$,
$f_2(R)$ (resp. $f_1$ and $f_2$) are conjugate under $\Aut(L,D)$ if
the two pairs $(M_1,f_1',C_1)$, $(M_2,f_2',C_2)$ are isomorphic.

{\sl Proof.} It is sufficient to show that a triple $(M,f',C)$
determines $f(R)$ (resp. $f$) up to conjugacy under
$\Aut(L,D)$. Given $M$ and $f'$ we can recover $L$ as the lattice
generated by $R\oplus M$ and the elements $r\oplus f'(r)$ for $r$ in
$R'$. We have a canonical map from $R$ to this $L$, so we have to show
that the Weyl chamber $D$ of $L$ is determined up to conjugacy by
elements of the group fixing $R$ and $M$ (resp. fixing $M$ and fixing
all points of $R$.) This Weyl chamber is determined by its
intersection with $R$ and $M$, and its intersection with $R$ is just
the canonical Weyl chamber of $R$. Its intersection with $M$ is in the
cone $C$ and is in some Weyl chamber of the norm 2 roots of $M$. All
such Weyl chambers of $M$ in $C$ are conjugate under automorphisms of
$L$ fixing $M$ and all points of $R$, so we can assume that the
intersection with $M$ is contained in some fixed Weyl chamber $W$ of
$M$.

By 6.3.1 and the assumption on $R$ all the Weyl chambers of $L$ whose
intersection with $M$ is in $W$ are conjugate under the group of
automorphisms of $L$ fixing $R$ (resp. fixing all points of $R$) and
hence $(M,f',C)$ determines $f(R)$ (resp. $f$). Q.E.D.

\proclaim Corollary 6.9.2. 
If $R$ is $e_6$, $e_7$, $e_8$, $d_4$, or $d_m$ $(m\ge 6)$ then two copies of
$R$ in $X$ are conjugate under $\Aut(D)$ if and only if their
orthogonal complements are isomorphic lattices.

{\sl Proof. } If $R'$ is any Dynkin diagram containing $R$ and one
extra point then $\sigma(R')$ fixes $R$. The result now follows from
6.9.1. Q.E.D.

{\sl Remark.} If $L$ is $II_{9,1}$ or $II_{17,1}$ then $\Aut(D)$ is
not transitive on $d_5$'s. $\sigma(e_6)$ does not fix the $d_5$'s in
$e_6$.

\proclaim Lemma 6.9.3. 
$\Aut(D)$ is transitive on $e_8$'s. The simple roots of $D$
perpendicular to an $e_8$ form the Dynkin diagram of $II_{8n-7,1}$ and
the subgroup of $\Aut(D)$ fixing the $e_8$ is isomorphic to the
subgroup of $Aut(II_{8n-7,1})$ fixing a Weyl chamber.

{\sl Proof.} The transitivity on $e_8$'s is in 6.9.2. The rest of 6.9.3
follows easily. Q.E.D. 

\proclaim Lemma 6.9.4. 
If $8n\ge 24$ then for any $e_6$ in $X$ there is an element of
$\Aut(D)$ inducing $\sigma(e_6)$ on it.

{\sl Proof. } By 6.9.2, $\Aut(D)$ is transitive on $e_6$'s so it is
sufficient to prove it for one $e_6$. It is true form $8n=24$ by
calculation and using 6.9.3 it follows by induction for $8n>24$. Q.E.D.

\proclaim Theorem 6.9.5. Classification of $d_m$'s in $X$. 
\item{(1)} 
$\Aut(D)$ acts transitively on $e_6$'s $e_7$'s, and $e_8$'s in $X$.
if $8n\ge 24$ then for any $e_6$ there is an element of $\Aut(D)$
inducing the nontrivial automorphism of this $e_6$.
\item{(2)} 
For any $m$ with $4\le m\le 8n+1$ there is a unique orbit of $d_m$'s in
$X$ such that $d_m^\perp$ is not unimodular, unless $m=5$ and $8n=8$
or $16$.  Any automorphism of such a $d_m$ is induced by an element of
$\Aut(D)$ if and only if $m\le 8n-13$.
\item{(3)} 
For any $m$ with $16\le 8m\le 8n$ there is a unique orbit of
$d_{8m}$'s such that $d_{8m}^\perp$ is unimodular, and these are the
only $d$'s whose orthogonal complement is unimodular. There is no
element of $\Aut(D)$ inducing the nontrivial automorphism of $d_{8m}$.

{\sl Proof. } Part (1) follows from 6.9.2 and 6.9.4 because there is only
one isomorphism class of lattices of the form $e_i^\perp$ for
$i=6,7,8$. From 6.9.4, 6.9.2, and 6.9.1 it follows that two $d_m$'s of $X$
are conjugate under $\Aut(D)$ if and only if their orthogonal
complements are isomorphic, unless $m=5$ and $n\le 2$. $d_m^\perp$ is
either the even sublattice of $I_{8n+1-m,1}$, or $m$ is divisible by 8
and $d_m^\perp$ is $II_{8n+1-m,1}$. In the second case we must have
$m\ge 16$ because if $m$ was 8 the Dynkin diagram of
$(II_{8n+1-m,1})^\perp$ would be $e_8$ and not $d_8$. This shows that
there is one orbit of $d_m$'s unless $8|m$, $m\ge 16$ or $m=5$, $n\le
2$ in which case there are two orbits.

If $d_m^\perp$ is unimodular then $d_m$ is contained in an even
unimodular sublattice of $L$, and there are no automorphisms of this
lattice acting non-trivially on $d_m$, so there are no elements of
$\Aut(D)$ inducing a nontrivial automorphism of $d_m$.

If $d_m^\perp$ is not unimodular then there is an element of $\Aut(D)$
inducing a nontrivial automorphism of $d_m$ if and only if there is an
automorphism of the Dynkin diagram of $I_{8n+1-m,1}$ acting
non-trivially on $M'/M$, where $M$ is the sublattice of even elements
of $I_{8n+1-m,1}$. There is no such automorphism of $I_k$ for $k\le
13$ and there is such an automorphism for $k=14$, so there is an
element of $\Aut(D)$ inducing a nontrivial automorphism of $d_m$ for
$m=8n-13$ and there is no such element if $m\ge 8n-12$. If $m<8n-13$
then our $d_m$ is contained in a $d_{8n-13}$ so there is still a
nontrivial automorphism of $d_m$ induced by $\Aut(D)$. Finally, if
$m=4$ and $8n\ge 24$ then as in the proof of 6.9.4 we see that there is
some $d_4$ such that $\Aut(D)$ induces all automorphisms of $d_4$. As
$\Aut(D)$ is transitive on $d_4$'s, this is true for any $d_4$. Q.E.D.

\proclaim Lemma 6.9.6. 
If $2\le m\le 11$ and $8n\ge 24$ then for any $a_m$ in $L$ there is an
element of $\Aut(D)$ inducing the nontrivial automorphism of $a_m$.

{\sl Proof. } This is true for $8n=24$ by calculation. There is an
element of $\Aut(D)$ acting as $\sigma(a_m)$ on $a_m$ if and only if
there is an automorphism of the Weyl chamber of the lattice
$M=a_m^\perp$ which acts as $-1$ on $M'/M$. The lattice $M$ is isomorphic to
$N\oplus e_8^{n-3}$, where $N$ is $a_m^\perp$ for some $a_m$ in
$II_{25,1}$ and $N$ has an automorphism of its Weyl chamber, so $M$
has one too. Hence for $m\le 11$ there is an element of $\Aut(D)$
inducing $\sigma(a_m)$ on $a_m$. Q.E.D.

\proclaim Corollary 6.9.7. 
If $8n\ge 24$ and $m\le 10$ then $\Aut(D)$ is transitive on $a_m$'s in $D$. 

{\sl Proof.} It follows from 6.9.6 and 6.9.5 that if $R'$ is any spherical
Dynkin diagram in $X$ containing $a_m$ and one extra point (so $R'$ is
$a_{m+1}$, $d_{m+1}$, $e_{m+1}$, or $a_ma_1$) then there is an element
$g$ of $\Aut(D)$ such that $g\sigma(R')$ fixes $a_m$. Hence by 6.9.1,
$\Aut(D)$ is transitive on $a_m$'s. Q.E.D.

\proclaim Corollary 6.9.8. 
If $n\ge 20 $ and $n\equiv 4,5$ or $6\bmod 8$ then the non-reflection
part of $\Aut(I_{n,1})$ can be written as a nontrivial amalgamated
product.

{\sl Proof. } The results of this section show that the analogue of
6.1 is true for $II_{8i+1,1}$ for $8i\ge 24$. This is all that is
needed to prove the analogue of 6.6.6. Q.E.D.

(If $n\ge 23$ then the group cannot be written as an amalgamated
product of finite groups.)

{\sl Remark.} If $n\ge 10$ and $n\equiv 2$ or $3\bmod 8$ and $G$ is
the subgroup of $\Aut(I_{n,1})$ generated by the reflections of
non-characteristic roots, then $\Aut(I_{n,1})/G$ is a nontrivial
amalgamated product.

\proclaim Chapter 7 The monster Lie algebra. 

\proclaim 7.1 Introduction. 

The monster Lie algebra  $\hat M$ was defined in [B-C-Q-S] as the Kac-Moody
algebra whose Dynkin diagram is the Dynkin diagram of $II_{25,1}$, so
$\hat M$ has a simple root for each point of $\Lambda$. See [M] for a
summary of Kac-Moody algebras. (He only deals with ones of finite
rank, but the extension of many results to algebras of infinite rank
like $\hat M$ is trivial.) The center of $\hat M$ is infinite
dimensional and if we quotient out $\hat M$ by its center the
resulting algebra $M$ has a 26 dimensional Cartan subalgebra which can
be identified with $II_{25,1}\otimes \R$. The roots of $M$ can be
identified with some points of $II_{25,1}$ such that the simple roots
of $M$ are the simple roots of a fundamental domain $D$ of
$II_{25,1}$. From tables $-2$ and $-4$ we know all the roots of norm
2, 0, $-2$, and $-4$ of $M$. In this chapter we will work out the
multiplicities of the roots of norms 2, 0, and $-2$ and the roots of
type 1. Note that if $u$ is a root of $M$ then in general it will
split into several roots of $\hat M$.

Remark added 1999: the results of this chapter are generalized to all roots in
``The monster Lie algebra'', Adv. Math Vol 83 No. 1 (1990). 

\proclaim Lemma 7.1.1. Any norm 2 vector of $II_{25,1}$ has multiplicity  1. 

Proof. Any simple root of $II_{25,1}$ has multiplicity  1 and any norm 2 vector
of $II_{25,1}$ is conjugate under the Weyl group to a simple
root. Q.E.D.

The zero vector of $II_{25,1}$ has mult 26. Any root that is not 0 and
does not have norm 2 has norm $\le 0$ and so is conjugate under the
Weyl group to a vector in the fundamental domain $D$ of
$II_{25,1}$. From now on we assume that $u$ is a vector of $D$.

\proclaim Lemma 7.1.2. If $u$ has height 0 or 1 then it has multiplicity 0. 

Proof. $u$ cannot be written as a sum of simple roots, as all simple
roots have height 1. Q.E.D.

\proclaim 7.2 Vector of types 0 and 1. 

We will calculate the multiplicities of vectors of types 0 and 1
considered as roots of $M$. We do this by using the fundamental
representation of the affine subalgebras of $M$.

The Kac-Weyl formula form $M$ states that 
$$
\prod_r(1-e^r) = \sum_{\sigma\in G}\det(\sigma)e^{\sigma(w)-w}$$
where the product is over the roots $r$ of $M$ with $-(r,w)>0$ taken
with their multiplicities, and the sum is over the elements $\sigma$
of the Weyl group $G$ of $II_{25,1}$.

Let $z$ be a primitive norm 0 vector in $D$ not equal to $w$. Then the
Dynkin diagram of $z^\perp$ is the Dynkin diagram of an affine Lie
algebra. By quotienting out part of the center of this Lie algebra and
adding an outer derivation we get a Lie algebra $Z$ with Cartan
subalgebra $II_{25,1}\otimes \R$ whose simple roots are the simple
roots of $D$ in $z^\perp$. $Z$ acts on $M$ by the adjoint
representation of $M$ restricted to $Z$ and $Z$ preserves each of the
subspaces $M_i$ where $M_i$ is the sum of the roots spaces of $M$ for
the roots that have inner product $-i$ with $z$. $M_0$ is the adjoint
representation of $Z$ (7.2.1), while $M_i$ for $i$ positive is a sum
of a finite number of highest weight representations of $Z$ and for
$i$ negative is a sum of a finite number of lowest weight
representations of $Z$. $M_1$ is a sum of fundamental representations
of $Z$ (7.2.5), and as the characters of fundamental representations
are known this allows us to compute the multiplicity of roots of $M$
that have inner product $-1$ with $z$.

\proclaim Lemma 7.2.1. 
$M_0$ is the adjoint representation of $Z$, and any norm 0 vector of
$D$ not conjugate to a multiple of $w$ has multiplicity 24.

Proof. If we take the terms of the Kac-Weyl formula of the form $e^r$
for roots $r$ with $(r,z)=0$ we find
$$\prod_{(r,z)=0}(1-e^r)=\sum_{\sigma\in
H}\det(\sigma)e^{\sigma(w)-w}$$ where $H$ is the Weyl group of
$z^\perp$. The right hand side of this is the right hand side of the
Kac-Weyl formula for $Z$, so the roots $r$ on the left hand side must
have the same multiplicity as roots of $Z$ or as roots of $M$. This
implies that $M_0$ is the adjoint representation of $M$.

The non-zero norm 0 roots of the adjoint representation of $Z$ have
multiplicity 24 which implies the second part of 7.2.1. Q.E.D.

\proclaim Lemma 7.2.2. 
$$\prod_{(r,z)<0}(1-e^r)=
{\sum_{\sigma\in G}\det(\sigma)e^{\sigma(w)-w}
\over
\sum_{\sigma\in H}\det(\sigma)e^{\sigma(w)-w}
}$$

Proof. This follows from the Kac-Weyl formula and 7.2.1. Q.E.D. 

\proclaim Lemma 7.2.3. 
$$
\sum_{-(z,r)=1}-e^r
=
{\sum_{\sigma\in G\atop (\sigma(w)-w,z)=-1}\det(\sigma)e^{\sigma(w)-w}
\over
\sum_{\sigma\in H}\det(\sigma)e^{\sigma(w)-w}
}$$

Proof. This is obtained by taking just the terms $e^r$ of 7.2.2 with
$(r,z)=-1$. Q.E.D.

Let $d$ be the square root of the Cartan matrix of the Dynkin diagram
of $z^\perp$. Then by 2.7 there are $d$ simple roots of $II_{25,1}$
with $-(r,z)=-1$ and any vector $x$ with $(x,z)=-1$ can be written
uniquely as $x=r+y$ for one of these $d$ roots $r$ and some $y$ that
is the sum of roots of $z^\perp$.

\proclaim Lemma 7.2.4. 
If $(\sigma(w)-w,z)=-1$ then $\sigma$ can be written uniquely as
$\sigma=\sigma_1\sigma_2$ where $\sigma $ is in $H$ and $\sigma_2$ is
reflection in one of the $d$ simple roots $r$ of $D$ with
$-(r,z)=1$. Conversely for any such $\sigma_1$ and $\sigma_2$ we have
$(\sigma_1\sigma_2(w)-w,z)=-1$.

Proof. Suppose that $(\sigma(w)-w,z)=-1$, so that
$-(\sigma^{-1}(z),w)=1+-(z,w)$. As $z$ is in the fundamental domain
$D$ of $w$ this implies that there is some root $r$ with $(r,w)\le 0$
such that the reflection of $\sigma^{-1}(z)$ in $r^\perp$ has inner
product $-(z,w)$ with $w$, and is therefore equal to $z$ as it is
conjugate to $z$ under the Weyl group and in the same fundamental
domain as $z$. Hence $z-(z,r)r=\sigma^{-1}(z)$. Taking inner products
with $w$ shows that $(w,z)-(z,r)(w,r)=(w,\sigma^{-1}(z))$, so
$(z,r)(w,r)=1$ and hence $(w,r)=(z,r)=-1$. This implies that $z$ is
one of the $d$ simple roots that have inner product $-1$ with $z$, and
if we put $\sigma_2=$ reflection in $r^\perp$,
$\sigma_1=\sigma\sigma_2$ then we find that
$$\eqalign{ (\sigma_1(w),z)&= (\sigma_2(w),\sigma^{-1}(z))\cr &=
(w+r,z-(z,r)r)\cr &= (w,r)\cr }$$ so $\sigma_1$ is in $H$. Hence
$\sigma$ has a decomposition as $\sigma=\sigma_1\sigma_2$.

If some element $\sigma $ of $G$ is equal to $\sigma_1\sigma_2$ for
some $\sigma_1$, $\sigma_2$ as above then
$(\sigma(w)-w,z)=(\sigma_2(w),\sigma_1^{-1}(z))-(w,z)=(w+r,z)-(w,z)=-1$. Also
$\sigma(w)-w=\sigma_1(w+r)-w=r+$ a sum of roots of $z^\perp$, so $r$
and hence $\sigma_1$ and $\sigma_2$ are determined by $\sigma$. Q.E.D.

\proclaim Theorem 7.2.5. 
$M_1$ considered as a representation of $Z$ is a sum of $d$
irreducible representations whose highest weights are the simple roots
$r$ with $-(r,z)=1$.

Proof. The character of $M_1$ is $\sum e^r$ where the sum is taken
over all roots $r$ of $M$ with $-(r,z)=1$, so by 7.2.3 it is equal to
$${
-\sum_{(\sigma(w)-w,z)=-1}\det(\sigma)e^{\sigma(w)-w}
\over
\sum_{\sigma\in H} \det(\sigma)e^{\sigma(w)-w}
}$$
By 7.2.4 this is equal to 
$${
-\sum_{\sigma_2}\sum_{\sigma_1\in H}\det(\sigma_1\sigma_2)
e^{\sigma_1\sigma_2(w)-w}
\over
\sum_{\sigma\in H} \det(\sigma)e^{\sigma(w)-w}
}$$
(where the first sums are over the $\sigma_1, \sigma_2$ as in 7.2.4)
$${
=\sum_r e^r\sum_{\sigma_1\in H}\det(\sigma_1)
e^{\sigma_1(w+r)-w-r}
\over
\sum_{\sigma\in H} \det(\sigma)e^{\sigma(w)-w}
}$$
where the first sum is over the $d$ roots $r$. Each term of the sum
over $r$ is the character of the irreducible representation of $Z$
with highest weight $r$, so 7.2.5 follows. Q.E.D.

\proclaim Corollary 7.2.6. 
If $u$ has inner product $-1$ with a norm 0 vector $z$ such that $z$
is not a conjugate of $w$ then the multiplicity of $u$ is the
coefficient of $q^{-u^2/2} $ in
$$\eqalign{
\Delta(q)^{-1}&= q^{-1}(1-q)^{-24}(1-q^2)^{-24}\cdots\cr
&= q^{-1}+24+324q+3200q^2+\cdots.\cr
}$$

Proof. We can assume that there is a norm 0 vector $z$ in $D$ not
equal to $w$ such that $-(u,z)=1$. The multiplicity of $u$ as a root
of $M$ is its multiplicity as a weight of $M_1$. If $M_r$ is the
representation with highest weight $r$ then $M_r$ has level $-(r,z)=1$
so its character is known. In fact, by formula 6.1 of [M], the
multiplicity of $u$ in $M_r$ is the coefficient of $q^{-u^2/2}$ in
$\Delta(q)^{-1}$ if $u-r$ is a sum of roots of $z^\perp$ and 0
otherwise. Any vector $u$ with $-(u,z)=1$ can be written uniquely as
$u=y+r$ for some $r$ and some $y$ which is a sum of roots of
$z^\perp$. This implies 7.2.6. Q.E.D.

\proclaim 7.3 Multiplicities of norm $-2$ vectors. 

$u$ is a vector of $D$ of norm $-2$. In this section we will show that
the multiplicity of $u$ as a root of $M$ is 0 if $u$ has height 0, 276
if $u$ has height 2, and 324 otherwise. This is proved by expressing
the multiplicity as a sum in terms of the geometry of $II_{25,1}$ near
$u$ and then evaluating this sum for each of the 121 orbits of norm
$-2$ vectors in $D$.

\proclaim Lemma 7.3.1. 
Let $X$ be the Dynkin diagram of a finite dimensional semisimple Lie
algebra $A$ all of whose roots have norm 2, and let $R$ be a finite
dimensional irreducible representation of $A$. $R$ is given by
assigning a non-negative integer to each point of $X$. Recall that
$R$ is said to be real, complex or quaternionic depending on whether
$(1,S^2R)=1$, $(1,R^2)=0$, or $(1,\Lambda^2R)=1$. Let $r$ be the
highest weight vector of $R$. Then
\item{(1)} 
The dual representation $R^*$ of $R$ has highest weight $\sigma(r)$,
so $R$ is complex if and only if $\sigma(r)\ne r$. (For $\sigma(r)$
see 1.3.)
\item{(2)}
If $\sigma(r)=r$ then $R$ is real or quaternionic depending on whether
$2(r,\rho)$ is even or odd, where $\rho$ is the Weyl vector of $A$.

Proof. These are well known facts about representations of Lie
algebras. Q.E.D.

\proclaim Lemma 7.3.2. 
If $A$ is a finite dimensional simple Lie algebra of rank $n$ all of
whose roots have norm 2 then there is a representation $R$ of $A$ not
containing 1 whose weights of non-zero multiplicity are as follows:
\item{(1)}
Every vector of the root lattice of $A$ that has norm 4 and is the sum
of two roots of $A$ has multiplicity 1 (as a weight of $R$).
\item{(2)}
Every root of $A$ has multiplicity $n-1$. 
\item{(3)}
0 has multiplicity $n(n+1)/2-1$. 

Computational proof. Check for each case. For $e_6$, $e_7$, $e_8$ $R$
is irreducible; for $d_n$ ($n\ge 5$) it is the sum of two irreducible
representations; for $d_4$ it is the sum of 3 irreducible
representations; for $a_n$ ($n\ge 3$) it is the sum of the adjoint
representation and another irreducible representation; for $a_2$ it is
the adjoint representation and for $a_1$ $R$ is 0.

Conceptual proof. From the theory of affine Lie algebras $A$ has a
series of representations $A_i$ such that
$\sum x_i\dim(A_i)=\theta_A(x)/\prod(1-x^i)^{rank(A)}$ where $\theta_A$
is the theta function of the root lattice of $A$. $A_0$ is 1, $A_1$ is
the adjoint representation of $A$, and $A_2$ is isomorphic to $R\oplus
A_1\oplus A_0$. Q.E.D.  

\proclaim Lemma 7.3.3. 
If $A$ is a finite dimensional semisimple Lie algebra of rank $n$ all
of whose roots have norm 2 then $A$ has a representation $R$ not
containing 1 whose weights are as follows:
\item{(1)}
Every vector of the root lattice of $A$ that has norm 4 and is the sum
of two roots of $A$ has multiplicity 1 (as a weight of $R$).
\item{(2)}
Every root of $A$ has multiplicity $n-1$. 
\item{(3)}
0 has multiplicity $n(n+1)/2-$ number of components of the Dynkin
diagram of $A$.

Proof. Take $R$ to be a sum of the representations of 7.3.2 for each
simple component of $A$, added to the sum of the product of the
adjoint representations of all pairs of simple components of $A$. It
is easy to check that $R$ has the required properties. Q.E.D.

\proclaim Lemma 7.3.4. 
Let $A$ be a finite dimensional reductive Lie algebra of rank $n$
whose semisimple part has rank $n'$. Then $A$ has a representation $R$
not containing 1 whose weights are as follows:

\item{(1)}
Every vector of the root lattice of $A$ that has norm 4 and is the sum
of two roots of $A$ has multiplicity 1 (as a weight of $R$).
\item{(2)}
Every root of $A$ has multiplicity $n-1$. 
\item{(3)}
0 has multiplicity $n'(n'+1)/2+(n-n')n'-$ (number of components of the Dynkin
diagram of $A$).

Proof. $R$ is the representation of 7.3.3 corresponding to the
semisimple part of $A$ added to $n-n'$ copies of the adjoint
representation of $A$ on its semisimple part. Q.E.D.

\proclaim Theorem 7.3.5. The multiplicity of $u$ as a root of $M$ is
$$\eqalign{
&325\cr
&-24\times(\hbox{rank of the roots $r$ of $u^\perp$ 
such that $u+r$ is a norm 0}\cr
&\qquad \hbox{vector corresponding to $\Lambda$})\cr
&-(d(d+1)/2+\hbox{number of components of the Dynkin diagram of $u^\perp$}\cr
&\qquad - \hbox{the number of orbits of $R_1$ under $\sigma$.})\cr
}$$
Here $d$ is $25-$ (the rank of the Dynkin diagram of $u^\perp$), $R_1$
is the set of simple roots of $D$ that have inner product $-1$ with
$u$, and $\sigma$ is the opposition involution of $u^\perp$ (which
acts on $R_1$).

Proof. We write the Kac-Weyl formula in the form
$$
\prod_{(r,u)<0}(1-e^r)
= {
\sum_{\sigma\in G}\det(w)e^{\sigma(w)-w}
\over
\prod_{(r,u)=0}(1-e^r)
}
$$
where the products are over positive roots $r$ of $II_{25,1}$ (with
multiplicities). Let $X_i$ be the representation of the Lie algebra of
$u^\perp$ whose weights are the projections into $u^\perp$ of the
roots $r$ of $M$ with $-(r,u)=i$, and let $U_i$ be the (virtual)
representation of $u^\perp$ whose weights are the projections of the
vectors $r$ of $II_{25,1}$ with $-(r,u)=i$, with multiplicity equal to
the number of times that $e^r$ occurs in the expansion of the Kac-Weyl
character formula above. Looking at the right hand side of the
Kac-Weyl character formula we see that any irreducible representation
occurring in $U_2$ has a highest weight vector with norm at least
$-u^2-(\sigma(w)-w)^2\ge 4$ (as $u^2=-2$ and $\sigma(w)-w$ is in $D$),
and a glance at the left hand side shows that $U_2$ is equal to
$\Lambda^2X_1-X_2$. The multiplicity of $u$ as a root of $M$ is equal
to the multiplicity of 0 as a weight of $X_2$.

Therefore we have
$$\eqalign{ X_2&= \Lambda^2(X_1)-U_2\cr &= \Lambda^2(X_1)^{(2)}+Y\cr
}$$ where $\Lambda^2(X_1)^{(2)}$ is the sum of the irreducible sub
representations of $\Lambda^2(X_1)$ whose highest weight vectors have
norm at most 2, and $Y$ is a sum of irreducible representations of
$u^\perp$ whose highest weight vectors have norm at least 4 because
all components of $U_2$ have this property.

If $R$ is a representation of $u^\perp$ then write $R'$ for the
representation of $u^\perp$ that is the sum of the irreducible sub
representations of $R$ whose weights are sums of roots of
$u^\perp$. Then
$$X_2'=\Lambda^2(X_1)^{(2)'}\oplus Y'.$$ $X_2$ and hence $X_2'$ has no
highest weight vectors of norm $\ge 6$, and all norm 4 vectors of
$u^\perp$ that are sums of two roots have multiplicity 1 as weights of
$X_2$. Hence $Y'$ is the sum of the irreducible representations of
$u^\perp$ whose highest weights are the norm 4 vectors in the Weyl
chamber of $u^\perp$ that are the sum of two roots of $u^\perp$. As
weights of $X_2'$ the norm 2 vectors $r$ of $u^\perp$ have
multiplicity equal to the multiplicity of $u+r$ as a root of $M$,
which by 7.2.1 and 7.1.2 is 24 if $z+r$ does not correspond to
$\Lambda$ and 0 otherwise.

Let $R$ be the representation of 7.3.4 for $A=u^\perp$, and let $S$ be
$R$ minus the adjoint representation of every component of $u^\perp$
such that if $r$ is a root of that component then $r+u$ is a norm 0
vector of $II_{25,1}$ corresponding to the Leech lattice. $u^\perp$
has rank 25 and the semisimple part of $u^\perp$ has rank $25-d$, so
the weights of $S$ are as follows:
\item{(1)} 
All norm 4 vectors of $u^\perp$ that are the sum of two roots have
multiplicity 1.
\item{(2)}
All roots $r$ of $u^\perp$ have multiplicity 24 or 0 equal to the
multiplicity of $r+u$ as a root of $M$.
\item{(3)}
0 has multiplicity 
$$\eqalign{
&(25-d)(26-d)/2+d(25-d)\cr
-&\hbox{number of components of }u^\perp\cr
-&24\times (\hbox{rank of roots of $u^\perp$ corresponding to $\Lambda$})\cr
\cr
=&325-d(d+1)/2 -\hbox{number of components of }u^\perp\cr &- 24\times
(\hbox{rank of roots of $u^\perp$ corresponding to $\Lambda$}).\cr }$$
$X_2'$ and therefore $S$ have the same multiplicities for all nonzero
roots, and the representation 1 does not occur in $Y'$ or $S$. Hence
the multiplicity of 0 in $X_2$, which is equal to the multiplicity of
$u$ as a root of $M$, is equal to the multiplicity of 0 as a weight of
$S$ plus the number of times that 1 occurs as an irreducible sub
representation of $\Lambda^2(X_1)^{(2)'}$. The number of times that 1
occurs in $\Lambda^2(X_1)^{(2)'}$ is equal to the number of times that
1 occurs in $\Lambda^2(X_1)$, so to complete the proof of 7.2.5 we
have to show that this number is equal to the number of orbits of
$\sigma$ on $R_1$.

Let $t_1, t_1,\ldots$ be the points of $R_1$, so that they are the
simple roots of $D$ that have inner product $-1$ with $u$. $X_1$ is a
sum of irreducible representations of $u^\perp$ whose highest weights
are the projections of the $t_i$ into $u^\perp$; call these
representations $T_1, T_2,\ldots$. Then
$$\Lambda^2(X_1)=
\bigoplus_i\Lambda^2(T_i)\oplus \bigoplus_{i<j} T_i\otimes T_j$$
so the number of $1$'s in $\Lambda^2(X_1)$ is equal to the number of
$T_i$'s that are quaternionic plus the number of pairs $(i<j)$ such
that $T_i$ and $T_j$ are dual.

By 7.3.1 $T_i$ is dual to $T_j$ for $i\ne j$ if and only if $T_i$ and
$T_j$ are exchanged by $\sigma$, so the number of pairs $(i<j)$ such
that $T_i$ is dual to $T_j$ is equal to the number of orbits of
$\sigma$ on $R_1$ of size 2.

If $T_i$ is quaternionic then $t_i$ is fixed by $\sigma$. Now suppose
that $t_i$ is fixed by $\sigma$. Let $x$ be the projection of $t_i$
into $u^\perp$ and let $\rho$ be the Weyl vector of $u^\perp$. $x$ is
a minimal vector of $u^\perp$ and is in the subspace of $u^\perp$
generated by roots as it is fixed by $\sigma$, so $x^2\equiv
(x,\rho)\bmod 1$. However $x^2=5/2$, so by 7.3.1 the representation
$T_i$ is quaternionic. Hence the number of $1$'s in $\Lambda^2(X_1)$
is equal to the number of orbits of $\sigma$ on $R_1$. Q.E.D.

{\bf Remark.} In the last paragraph of this proof, $\rho$ is the
projection of $w$ into $u^\perp$ by 4.3.3. Hence
$(x,\rho)=(t_i,\rho)=(t_i,w-(w,u)u/(u,u))=-1+\height(u)/2$, so if
there exists a $t_i$ in $R_1$ fixed by $\sigma$ then the height of $u$
must be odd.

\proclaim Corollary 7.3.6. 
The multiplicity of the norm $-2$ vector $u$ of $D$ considered as a
root of $M$ is 0 if $u$ has height 1, 276 if $u$ has height $2$, and
324 otherwise.

Proof. The rank of the roots of $u^\perp$ corresponding to $\Lambda$
is 1 or 2 if $u$ has height 1 or 2, and 0 otherwise. (See table $-2$.)
We can evaluate the term in parentheses in 7.3.5 and we find that it
is 301 if $u$ has height 1 and 1 for the other 120 orbits of roots
$u$. (The sets $R_1$ had to be found for the enumeration of the norm
$-4$ vectors in $II_{25,1}$ so calculating the term in parentheses is
very easy.) Hence 7.3.6 follows from 7.3.5. Q.E.D.

{\bf Example.} Let $u$ be the norm $-2$ vector of height 31 with
Dynkin diagram $a_{15}d_8a_1$. 

(Diagram missed out.)

$R_1$ has 5 point which form 3 orbits under $\sigma$, indicated by
$x$, $y$, and $z$. $d=25-(15+18+1)=1$ and there are 3 components of
$u^\perp$, so by 7.3.5 the multiplicity of $u$ is $325-((1\times
2)/2+3-3)=324$.

\proclaim References.

\item{[B-C-Q-S]} Borcherds, R. E.; Conway, J.
   H.; Queen, L.; Sloane, N. J. A. A monster Lie algebra?
   Adv. in Math. 53 (1984), no. 1, 75--79.
\item{[B-C-Q]} Borcherds, R. E.; Conway, J.
   H.; Queen, L.; The cellular structure  of the Leech lattice, 
Chapter 25 of Conway, J. H.;
   Sloane, N. J. A. Sphere packings, lattices and groups. Third
   edition. Grundlehren der Mathematischen Wissenschaften, 290.
   Springer-Verlag, New York, 1999.
\item{[B]}
Bourbaki, Nicolas \'El\'ements de
   math\'ematique.  Groupes et alg\'ebres
 de Lie. Chapitres 4, 5 et 6. 
 Masson, Paris, 1981.  ISBN: 2-225-76076-4
\item{[C a]} Conway, J. H. A group of  
   order $8,315,553,613,086,720,000$. Bull. London Math. Soc. 1 
   1969 79--88.
\item{[C b]}Conway, J. H. The automorphism
   group of the $26$-dimensional even unimodular Lorentzian lattice.
   J. Algebra 80 (1983), no. 1, 159--163.
\item{[C c]}Conway, J. H. A characterisation of Leech's lattice. Invent. Math. 7 1969
   137--142.
\item{[C-P-S]} Conway, J. H.; Parker, 
R. A.; Sloane, N. J. A. The covering radius of the Leech lattice.
   Proc. Roy. Soc. London Ser. A 380 (1982), no. 1779, 261--290.
\item{[C-S a]}Conway, J. H.; Sloane, N.
J. A. The unimodular lattices of dimension up to $23$ and the
   Minkowski-Siegel mass constants. European J. Combin. 3
   (1982), no. 3, 219--231. 
\item{[C-S b]}Conway, J. H.; Sloane, N. J. A.
 Twenty-three constructions for the Leech lattice. Proc.
   Roy. Soc. London Ser. A 381 (1982), no. 1781, 275--283.
\item{[C-S c]}Conway, J. H.; Sloane, N.
J. A. Lorentzian forms for the Leech lattice. Bull. Amer. Math. 
   Soc. (N.S.) 6 (1982), no. 2, 215--217
\item{[C-S d]}Conway, J. H.; Sloane, N.
J. A. Leech roots and Vinberg groups. Proc. Roy. Soc. London Ser.
   A 384 (1982), no. 1787, 233--258.
\item{[D]}Dynkin, E. B. Semisimple
 subalgebras of semisimple Lie algebras.  Mat. Sbornik
   N.S. 30(72), (1952). 349--462 

\item{[F]}Frenkel, I. B. Representations 
   of Kac- Moody algebras and dual resonance models. Applications of 
   group theory in physics and mathematical physics (Chicago, 1982),
   325--353, Lectures in Appl. Math., 21, Amer. Math. Soc.,
   Providence, R.I., 1985.
\item{[K]}Kneser, Martin Klassenzahlen
   definiter quadratischer Formen.  Arch. Math. 8 (1957),
   241--250. 
\item{[M]}Macdonald, I. G. Affine Lie
   algebras and modular forms. Bourbaki Seminar, Vol. 1980/81, pp.
   258--276, Lecture Notes in Math., 901, Springer, Berlin-New York,
   1981.
\item{[M a]}Macdonald, I. G. Affine root
   systems and Dedekind's $\eta $-function. Invent. Math. 15 (1972),
   91--143. 
\item{[S a]}Serre, J.-P. A course in
   arithmetic.  Graduate Texts in Mathematics,
   No. 7. Springer-Verlag, New York-Heidelberg, 1973.
\item{[S b]}Serre, Jean-Pierre Trees.
   Springer-Verlag,
   Berlin-New York, 1980. ix+142 pp. ISBN: 3-540-10103-9
\item{[S c]}Serre, Jean-Pierre Cohomologie
   des groupes discrets.  Prospects in mathematics (Proc.
   Sympos., Princeton Univ., Princeton, N.J., 1970), pp. 77--169. Ann. of
   Math. Studies, No. 70, Princeton Univ. Press, Princeton, N.J., 1971.
\item{[V]}Venkov, B. B. Classifications
   of integral even unimodular $24$-dimensional quadratic forms.
    Algebra, number theory and their applications. Trudy
   Mat. Inst. Steklov. 148 (1978), 65--76, 273.
\item{[Vi]} Vinberg, E. B. Some arithmetical discrete groups in Loba\v cevski\u\i spaces. Discrete
   subgroups of Lie groups and applications to moduli (Internat. Colloq.,
   Bombay, 1973), pp. 323--348. Oxford Univ. Press, Bombay, 1975.
\item{[V-K]}Vinberg, E. B.;
Kaplinskaja, I. M. The groups $O_{18,1}(Z)$ and
   $O_{19,1}(Z)$. Dokl. Akad. Nauk SSSR 238 (1978), no.
   6, 1273--1275. Translation in Sov. Math. 19 no. 1 (1978) p. 194--197.

\proclaim Figure 1 The neighborhood graph for 8, 16, and 24 dimensions. 

Each circle of this graph represents an $8n$ dimensional even
unimodular lattice ($8n=8, 16, 24$), and each line represents an
$8n$-dimensional odd unimodular lattice with no vectors of norm 1. The
two circles which a line is joined to are the two even neighbors of
the odd lattice. ($8n$-dimensional lattices with vectors of norm 1
would correspond to some more lines both of whose ends joined the same
circle.) The 8 and 16 dimensional cases were worked out in [K], where
it was also proved that the graph is connected in any given
dimension. In 32 dimensions the graph has at least $10^8$ circles and
$10^{17}$ lines.

The thick lines represent the odd lattices with ``$h_2=2h_1+2$'' (see
4.5.3.)

(This figure has not yet been convertex to \TeX. See Conway and Sloane, 
``Sphere packings, lattices, and groups'', chapter 17 to see a copy of it.
Alternatively the reader can draw their own copy using table $-4$, 
which lists the 24 dimensional unimodular lattices and their even neighbors.)

\proclaim Figure 2 Some vectors of $II_{25,1}$. 

The diagram gives the orbits of vectors of $II_{25,1}$ of height at
most 4 and $-$norm at most 40. For larger norms of height $i$ just
repeat the last $i$ rows of the diagram. Each orbit of vectors of a
give height and norm is represented by a Dynkin diagram (preceded by a
number $n$ if the orbit is $n$ times a primitive orbit for $n>1$) or
by $\circ$ when the Dynkin diagram is empty. Let $u$ be a vector in
$D$. If $r$ is a highest root of $u^\perp$ and $u$ does not have type
0 or 1 then $r+u$ is also in $D$ and also appears in the table one
column to the left of $u$. If there are no roots in $u^\perp$ then
$u=v+w$ for some $v$ in $D$ in the same column as $u$ and $i$ rows
higher. The diagram can be continued in the obvious way for norm 0
vectors, and tables $-2$ and $-4$ contain a complete description of it
for norm $-2$ and $-4$ vectors.

The orbits of height $i$ with non-empty Dynkin diagram correspond
naturally to the orbits $x$ of $\Lambda\bmod i\Lambda$, and the Dynkin
diagram describes the points of $\Lambda$ nearest to $x/i$. For
example, from the column with $i=2$ we see that there are 4 orbits of
elements of $\Lambda\bmod 2\Lambda$ represented by elements of norms
0, 4, 6, and 8, such that the number of vectors nearest to half these
vectors is 1, 2, 2, and 48 and these vectors form a Dynkin diagram
$a_1$, $a_1^2$, $a_2$, or $A_1^{24}$.

\halign{
\hfill$#$&&~~~\hfill$#$\cr
&\hbox{Height=}0&1&2&3&4\cr
\hbox{Norm}\cr
0&nw(n\ge3)\cr
&w,2w&&A_1^{24}&A_2^{12}&A_3^8,2A_1^{24}\cr
-2&&a_1&a_2&a_1^9&a_2a_1^{12}\cr
-4&&\circ&a_1^2~\circ&a_1^2&a_2^2,a_1^8,a_1^6\cr
-6&&\circ&\circ&a_1^3,a_1~\circ&a_1^2,a_1^3\cr
-8&&\circ&2a_1~\circ~\circ&a_1~\circ&2a_2,a_1^4,a_1^2,a_1,a_1~\circ~\circ\cr
-10&&\circ&\circ&a_1~\circ&a_1^2,a_1~\circ\cr
-12&&\circ&\circ~\circ~\circ&a_1~\circ~\circ~\circ
&a_1^2,a_1,a_1~\circ~\circ~\circ\cr
-14&&\circ&~\circ&a_1~\circ~\circ
&a_1,a_1~\circ~\circ\cr
-16&&\circ&\circ~\circ~\circ&~\circ~\circ
&2a_1^2,a_1,a_1~\circ~\circ~\circ~\circ~\circ~\circ~\circ\cr
-18&&\circ&~\circ&3a_1~\circ~\circ~\circ~\circ
&a_1~\circ~\circ~\circ\cr
-20&&\circ&\circ~\circ~\circ&~\circ~\circ~\circ
&a_1,a_1~\circ~\circ~\circ~\circ~\circ~\circ\cr
-22&&\circ&~\circ&\circ~\circ
&a_1~\circ~\circ~\circ~\circ\cr
-24&&\circ&\circ~\circ~\circ&\circ~\circ~\circ~\circ~\circ
&a_1~\circ~\circ~\circ~\circ~\circ~\circ~\circ~\circ~\circ~\circ\cr
-26&&\circ&~\circ&\circ~\circ~\circ
&a_1~\circ~\circ~\circ~\circ\cr
-28&&\circ&\circ~\circ~\circ&~\circ~\circ
&a_1~\circ~\circ~\circ~\circ~\circ~\circ~\circ~\circ\cr
-30&&\circ&~\circ&\circ~\circ~\circ~\circ~\circ
&~\circ~\circ~\circ~\circ~\circ\cr
-32&&\circ&\circ~\circ~\circ&~\circ~\circ~\circ
&3a_1~\circ~\circ~\circ~\circ~\circ~\circ~\circ~\circ~\circ~\circ~\circ\cr
-34&&\circ&~\circ&\circ~\circ
&~\circ~\circ~\circ~\circ~\circ\cr
-36&&\circ&\circ~\circ~\circ&\circ~\circ~\circ~\circ~\circ
&~\circ~\circ~\circ~\circ~\circ~\circ~\circ~\circ~\circ\cr
-38&&\circ&~\circ&\circ~\circ~\circ
&~\circ~\circ~\circ~\circ~\circ\cr
-40&&\circ&\circ~\circ~\circ&~\circ~\circ
&\circ~\circ~\circ~\circ~\circ~\circ~\circ~\circ~\circ~\circ~\circ~\circ\cr
}

\proclaim Table $-2$.~The norm $-2$  vectors of $II_{25,1}$. 

The following sets are in natural 1:1 correspondence:
\item{(1)} Orbits of norm $-2$ vectors in $II_{25,1}$ under $\Aut(II_{25,1})$.
\item{(2)} Orbits of norm $-2$ vectors $u$ of $D$ under $\Aut(D)$. 
\item{(3)} 25 dimensional even bimodular lattices $L$. 

$L$ is isomorphic to $u^\perp$. Table $-2$ lists the 121 elements of
any of these three sets.

The {\it height} $t$ is the height of the norm $-2$ vector $u$ of $D$,
in other words $-(u,w)$ where $w$ is the Weyl vector of $D$. The
letter after the height is just a name to distinguish vectors of the
same height, and is the letter referred to in the column headed ``Norm
$-2$'s'' of table $-4$. An asterisk after the letter means that the
vector $u$ is of type 1, in other words the lattice $L$ is the sum of
a Niemeier lattice and $a_1$.

The column ``Roots'' gives the Dynkin diagram of the norm 2 vectors of
$L$ arranged into orbits under $\Aut(L)$. ``Group'' is the order of
the subgroup of $\Aut(D)$ fixing $u$. The group $\Aut(L)$ is a split
extension $R.G$ where $R$ is the Weyl group of the Dynkin diagram and
$G$ is isomorphic to the subgroup of $\Aut(D)$ fixing $u$.

``$S$'' is the maximal number of pairwise orthogonal roots of $L$. 

The column headed ``Norm 0 vectors'' describes the norm 0 vectors $z$
corresponding to each orbit of roots of $u^\perp$ where $u$ is in $D$,
as in 3.5.2.  A capital letter indicates that the corresponding norm 0
vector is twice a primitive vector, otherwise the norm 0 vector is
primitive. $x$ stands for a norm 0 vector of type the Leech
lattice. Otherwise the letter $a$, $d$, or $e$ is the first letter of
the Dynkin diagram of the norm 0 vector, and its height is given by
${\rm height}(u)-h+1$ where $h$ is the Coxeter number of the component of
the Dynkin diagram of $u$. 

For example, the norm $-2$ vector of type $23a$ has 3 components in
its root system, of Coxeter numbers 12, 12, and 6, and the letters are
$e$, $a$, and $d$, so the corresponding norm 0 vectors have Coxeter
numbers 12, 12, and 18 and hence are norm 0 vectors with Dynkin
diagrams $E_6^4$, $A_{11}D_7E_6$, and $D_{10}E_7^2$.

See 4.3 for more information. 

\halign{
\hfill#&#\hfill&~\hfill$#$&~~\hfill#&~~~\hfill#&\hfill#&#\cr
Height&&\hbox{Roots}&Group&$S$&~Norm 0&~vectors\cr
\cr
1&a*&a_1 &8315553613086720000 &1&X\cr
\cr
2&a &a_2 &991533312000 &1&x\cr
\cr
3&a &a_1^9 &92897280 &9&a\cr
\cr
4&a &a_2a_1^{12} &190080 &13&aa\cr
\cr
5&a*&a_1^{24}a_1 &244823040 &25&aA\cr
5&b &a_2^4a_1^9 &3456 &13&aa\cr
5&c &a_3a_1^{15} &40320 &17&aa\cr
\cr
6&a &a_2^9 &3024 &9&a\cr
6&b &a_3a_2^5a_1^6 &240 &13&aaa\cr
\cr
7&a*&a_2^{12}a_1 &190080 &13&aA\cr
7&b &a_3^3a_2^4a_1^3 &48 &13&aaa\cr
7&c &a_3^4a_1^8a_1 &384 &17&aad\cr
7&d &a_4a_2^6a_1^5 &240 &13&aaa\cr
7&e &d_4a_1^{21} &120960 &25&ad\cr
\cr
8&a &a_3^6a_2 &240 &13&ad\cr
8&b &a_4a_3^3a_2^3a_1^2 &12 &13&aaaa\cr
8&c &d_4a_2^9 &864 &13&aa\cr
\cr
9&a*&a_3^8a_1 &2688 &17&aA\cr
9&b &a_4^2a_3^4a_1 &16 &13&aaa\cr
9&c &a_4^3a_3a_2^2a_1^3 &12 &13&aaaa\cr
9&d &d_4a_3^4a_3a_1^3 &48 &17&aada\cr
9&e &a_5a_3^3a_2^4 &24 &13&aaa\cr
9&f &a_5a_3^4a_1^6 &48 &17&aaa\cr
\cr
10&a &d_4a_4^3a_2^3 &12 &13&aaa\cr
10&b &a_5a_4^2a_3^2a_2a_1 &4 &13&aaaaa\cr
\cr
11&a*&a_4^6a_1 &240 &13&aA\cr
11&b &d_4^4a_1^9 &432 &25&dd\cr
11&c &a_5d_4^2a_3^3 &24 &17&daa\cr
11&d &a_5a_5a_4^2a_3a_1 &4 &13&aaaaa\cr
11&e &a_5^2d_4a_3^2a_1^2a_1 &8 &17&aaaad\cr
11&f &a_5^3a_2^4 &48 &13&aa\cr
11&g &d_5a_3^6a_1 &48 &17&aad\cr
11&h &a_6a_4^2a_3^2a_2a_1 &4 &13&aaaaa\cr
\cr
12&a &a_5^4a_2 &24 &13&ad\cr
12&b &d_5a_4^4a_2 &8 &13&aaa\cr
12&c &a_6d_4a_4^3 &6 &13&aaa\cr
12&d &a_6a_5^2a_3a_2^2 &4 &13&aaaa\cr
\cr
13&a*&a_5^4d_4a_1 &48 &17&aaA\cr
13&b &d_5a_5^2d_4a_3a_1 &4 &17&aaada\cr
13&c &d_5a_5^3a_1^3a_1 &12 &17&aaae\cr
13&d*&d_4^6a_1 &2160 &25&dD\cr
13&e &a_6^2a_5a_4a_1^2 &4 &13&aaaa\cr
13&f &a_7a_5a_4^2a_3 &4 &13&aaaa\cr
13&g &a_7a_5d_4a_3^2a_1^2 &4 &17&aaaaa\cr
\cr
14&a &a_6a_6d_5a_4a_2 &2 &13&aaaaa\cr
14&b &a_6^3d_4 &12 &13&aa\cr
14&c &a_7a_6a_5a_4a_1 &2 &13&aaaaa\cr
\cr
15&a*&a_6^4a_1 &24 &13&aA\cr
15&b &d_5^3a_5a_3 &12 &17&ade\cr
15&c &d_6d_4^4a_1^3 &24 &25&ddd\cr
15&d &d_6a_5^2a_5a_3 &4 &17&aada\cr
15&e &a_7d_5^2a_3^2a_1 &4 &17&aaad\cr
15&f &a_7^2d_4^2a_1 &8 &17&aad\cr
15&g &a_8a_5^3 &6 &13&aa\cr
15&h &a_8a_6a_5a_3a_2 &2 &13&aaaaa\cr
\cr
16&a &a_7^3a_2 &12 &13&ad\cr
16&b &a_8a_6d_5a_4 &2 &13&aaaa\cr
\cr
17&a*&a_7^2d_5^2a_1 &8 &17&aaA\cr
17&b &e_6a_5^3d_4 &12 &17&aae\cr
17&c &a_7d_6d_5a_5 &2 &17&daaa\cr
17&d &a_7^2d_6a_3a_1 &4 &17&aada\cr
17&e &a_8a_7^2a_1 &4 &13&aaa\cr
17&f &a_9d_5a_5d_4a_1 &2 &17&aaaaa\cr
17&g &a_9a_7a_4^2 &4 &13&aaa\cr
\cr
18&a &e_6a_6^3 &6 &13&aa\cr
18&b &a_9a_8a_5a_2 &2 &13&aaaa\cr
\cr
19&a*&a_8^3a_1 &12 &13&aA\cr
19&b &d_6^3d_4a_1^3 &6 &25&ddd\cr
19&c &a_7e_6d_5^2a_1 &4 &17&eaad\cr
19&d &d_7a_7d_5a_5 &2 &17&aaad\cr
19&e &d_7a_7^2a_3a_1 &4 &17&aaad\cr
19&f &a_9a_7d_6a_1a_1 &2 &17&aaaad\cr
19&g &a_{10}a_7a_6a_1 &2 &13&aaaa\cr
\cr
20&a &a_8^2e_6a_2 &4 &13&aaa\cr
20&b &a_{10}a_8d_5 &2 &13&aaa\cr
\cr
21&a*&a_9^2d_6a_1 &4 &17&aaA\cr
21&b &a_{11}d_6a_5a_3 &2 &17&aaaa\cr
21&c &a_{11}a_8a_5 &2 &13&aaa\cr
21&d*&d_6^4a_1 &24 &25&dD\cr
21&e &a_9e_6d_6a_3 &2 &17&aaad\cr
\cr
23&a &d_7e_6^2a_5 &4 &17&ead\cr
23&b &d_8d_6^2d_4a_1 &2 &25&dddd\cr
23&c &a_9d_7^2 &4 &17&da\cr
23&d &a_9d_8a_7 &2 &17&daa\cr
23&e &a_{11}d_7d_5a_1 &2 &17&aaad\cr
\cr
24&a &a_{11}^2a_2 &4 &13&ad\cr
24&b &a_{12}e_6a_6 &2 &13&aaa\cr
\cr
25&a*&a_{11}d_7e_6a_1 &2 &17&aaaA\cr
25&b &a_{13}d_6d_5 &2 &17&aaa\cr
25&e*&e_6^4a_1 &48 &17&eE\cr
\cr
26&a &a_{13}a_{10}a_1 &2 &13&aaa\cr
\cr
27&a*&a_{12}^2a_1 &4 &13&aA\cr
27&b &e_7d_6^3 &3 &25&dd\cr
27&c &a_9a_9e_7 &2 &17&ada\cr
27&d &d_9a_9e_6 &2 &17&ada\cr
27&e &a_{11}d_9a_5 &2 &17&aad\cr
27&f &a_{14}a_9a_2 &2 &13&aaa\cr
\cr
29&a &a_{11}e_7e_6 &2 &17&daa\cr
29&d*&d_8^3a_1 &6 &25&dD\cr
\cr
31&a &d_8^2e_7a_1a_1 &2 &25&ddde\cr
31&b &d_{10}d_8d_6a_1 &1 &25&dddd\cr
31&c &a_{15}d_8a_1 &2 &17&aad\cr
\cr
33&a*&a_{15}d_9a_1 &2 &17&aaA\cr
33&b &a_{15}e_7a_3 &2 &17&aad\cr
33&c &a_{17}a_8 &2 &13&aa\cr
\cr
35&a &e_7^3d_4 &6 &25&de\cr
35&b &a_{13}d_{11} &2 &17&da\cr
\cr
36&a &a_{18}e_6 &2 &13&aa\cr
\cr
37&a*&a_{17}e_7a_1 &2 &17&aaA\cr
37&d*&d_{10}e_7^2a_1 &2 &25&ddD\cr
\cr
39&a &d_{12}e_7d_6 &1 &25&ddd\cr
\cr
45&d*&d_{12}^2a_1 &2 &25&dD\cr
\cr
47&a &d_{10}e_8e_7 &1 &25&edd\cr
47&b &d_{14}d_{10}a_1 &1 &25&ddd\cr
47&c &a_{17}e_8 &2 &17&da\cr
\cr
48&a &a_{23}a_2 &2 &13&ad\cr
\cr
51&a*&a_{24}a_1 &2 &13&aA\cr
\cr
61&d*&d_{16}e_8a_1 &1 &25&ddD\cr
61&e*&e_8^3a_1 &6 &25&eE\cr
\cr
63&a &d_{18}e_7 &1 &25&dd\cr
\cr
93&d*&d_{24}a_1 &1 &25&dD\cr
}

\proclaim Table $-4$.~The norm $-4$ vectors of $II_{25,1}$. 

There is a natural 1:1 correspondence between the elements of the
following sets:
\item{(1)} Orbits of norm $-4$ vectors $u$ in $II_{25,1}$ under 
$\Aut(II_{25,1})$. 
\item{(2)}
Orbits of norm $-4$ vectors in the fundamental domain $D$ of
$II_{25,1}$ under $\Aut(D)$.
\item{(3)}
Orbits of norm $-1$ vectors $v$ of $I_{25,1}$ under $\Aut(I_{25,1})$. 
\item{(4)}
25 dimensional unimodular positive definite lattices $L$.
\item{(5)} 
Unimodular lattices $L_1$ of dimension  at most 25 with no vectors of norm 1. 
\item{(6)}
25 dimensional even lattices $L_2$ of determinant 4. 

$L_1$ is the orthogonal complement of the norm 1 vectors of $L$, $L_2$
is the lattice of elements of $L$ of even norm, $L_2$ is isomorphic to
$u^\perp$, and $L$ is isomorphic to $v^\perp$. Table $-4$ lists the
665 elements of any of these sets.

The height $t$ is the height of the norm $-4$ vector $u$ of $D$, in
other words $-(u,w)$ where $w$ is the Weyl vector of $D$. The things
in table $-4$ are listed in increasing order of their height.

Dim is the dimension of the lattice $L_1$. A capital $E$ after the
dimension means that $L_1$ is even.

The column ``roots'' gives the Dynkin diagram of the norm 2 vectors of
$L_2$ arranged into orbits under $\Aut(L_2)$.

``Group'' gives the order of the subgroup of $\Aut(D)$ fixing $u$. The
group $\Aut(L)\cong \Aut(L_2)$ is of the form $2\times R.G$ where $R$
is the group generated by the reflections of norm 2 vectors of $L$,
$G$ is the group described in the column ``group'', and 2 is the group
of order 2 generated by $-1$. If $\dim(L_1)\le 24$ then $\Aut(L_1)$ is
of the form $R.G$ where $R$ is the reflection group of $L_1$ and $G$
is as above.

For any root $r$ of $u^\perp$ the vector $u+r$ is a norm $-2$ vector
of $II_{25,1}$. This vector $u$ can  be found as follows. Let $X$ be the
component of the Dynkin diagram of $u^\perp$ to which $u$ belongs and
let $h$ be the Coxeter number of $X$. Then $r+u$ is conjugate to a
norm $-2$ vector of $II_{25,1}$ in $D$ of height $t-h+1$ (or $t-h$ if the
entry under ``Dim'' is $24E$) whose letter is the letter corresponding
to $X$ in the column headed ``norm $-2$'s''. For example let $u$ be
the vector of height 6 and root system $a_2^2a_1^{10}$. Then the norm
$-2$ vectors corresponding to roots from the components $a_2$ or $a_1$
have heights $6-3+1$ and $6-2+1$ and letters $a$ and $b$, so they are
the vectors $4a$ and $5b$ of table $-2$.

If $\dim(L_1)\le 24$ then the column ``neighbors'' gives the two even
neighbors of $L_1+I^{24-\dim(L_1)}$. If $\dim(L_1)\le 23$ then both
neighbors are isomorphic so only one is listed, and if $L_1$ is a
Niemeier lattice then the neighbor is preceded by 2 (to indicate that
the corresponding norm 0 vector is twice a primitive vector).  If the
two neighbors are isomorphic then there is an automorphism of $L$
exchanging them.
\bigskip

\halign{
\hfill$#$&~\hfill$#$&#&~\hfill$#$&~\hfill$#$&~\hfill#&~\hfill$#$&\hfill$#$\cr
\hbox{Height}&\hbox{Dim}&&\hbox{Roots}&\hbox{Group}&\hbox{norm $-2$'s}&
\hbox{Neighbors}&\cr
1 &24&E & &8315553613086720000 &&2\Lambda&      \cr
\cr
2 &23& &a_1^2 &84610842624000 &a&\Lambda &      \cr
2 &24& & &1002795171840 &&\Lambda&A_1^{24}      \cr
\cr
3 &25& &a_1^2 &88704000 &a&      &      \cr
\cr
4 &24& &a_1^8 &20643840 &a&A_1^{24}    &A_1^{24}    \cr
4 &25& &a_2^2 &26127360 &a&      &      \cr
4 &25& &a_1^6 &138240 &a&      &      \cr
\cr
5 &24& &a_1^{12} &190080 &a& A_1^{24}   & A_2^{12}   \cr
5 &25& &a_2a_1^7 &5040 &aa&      &      \cr
5 &25& &a_1^{10} &1920 &a&      &      \cr
\cr
6 &23& &a_1^{16}a_1^2 &645120 &ca&  A_1^{24}  &      \cr
6 &24& &a_2^2a_1^{10} &5760 &ab&    A_2^{12}&A_2^{12}    \cr
6 &24& &a_1^{16} &43008 &c&    A_1^{24}& A_3^8   \cr
6 &25& &a_3a_1^8 &21504 &ac&      &      \cr
6 &25& &a_2^2a_1^8 &128 &ab&      &      \cr
6 &25& &a_2a_1^{10}a_1 &120 &abc&      &      \cr
6 &25& &a_1^8a_1^6 &1152 &bc&      &      \cr
\cr
7 &24& &a_2^4a_1^8 &384 &bb&    A_2^{12}&A_3^8    \cr
7 &24&E &a_1^{24} &244823040 &a&   2A_1^{24}&      \cr
7 &25& &a_2^5a_1^3 &720 &bb&      &      \cr
7 &25& &a_3a_2a_1^9 &72 &acb&      &      \cr
7 &25& &a_2^4a_1^4a_1^2 &24 &bba&      &      \cr
7 &25& &a_3a_1^{12} &1440 &ab&      &      \cr
7 &25& &a_2^3a_1^6a_1^3 &12 &bbb&      &      \cr
7 &25& &a_2^2a_1^{12} &144 &cb&      &      \cr
\cr
8 &22& &a_3a_1^{22} &887040 &ae&    A_1^{24}&      \cr
8 &23& &a_2^6a_1^6a_1^2 &1440 &bda&   A_2^{12} &      \cr
8 &24& &a_3^2a_1^{12} &768 &cc&    A_3^8& A_3^8   \cr
8 &24& &a_3a_2^4a_1^6 &96 &bbb&    A_3^8&A_3^8    \cr
8 &24& &a_2^8 &672 &a&    A_3^8&A_3^8    \cr
8 &24& &a_2^6a_1^6 &240 &bd&   A_2^{12} &A_4^6    \cr
8 &24& &a_1^{24} &138240 &e&    A_1^{24}&D_4^6   \cr
8 &25& &a_4a_1^{12} &1440 &ad&      &      \cr
8 &25& &a_3a_2^4a_1^4 &16 &bbb&      &      \cr
8 &25& &a_3^2a_1^8a_1^2 &64 &cbc&      &      \cr
8 &25& &a_3a_2^3a_1^6a_1 &12 &bbbc&      &      \cr
8 &25& &a_3a_2^3a_1^3a_1^3a_1 &6 &bbbbd&      &      \cr
8 &25& &a_2^4a_2^2a_1^4 &8 &bab&      &      \cr
8 &25& &a_3a_2^2a_1^4a_1^4a_1^2 &16 &bbcbd&      &      \cr
8 &25& &a_2^4a_2a_1^4a_1^2a_1 &8 &bbbdb&      &      \cr
8 &25& &a_2^4a_1^8a_1^2 &48 &bdc&      &      \cr
8 &25& &a_3a_1^{15}a_1 &720 &cce&      &      \cr
\cr
9 &24& &a_3^2a_2^4a_1^4 &16 &bbb&  A_3^8  &A_4^6    \cr
9 &24& &a_2^8a_1^4 &384 &dc&    A_2^{12}&A_5^4D_4    \cr
9 &25& &a_4a_2^3a_1^6a_1 &12 &bdbb&      &      \cr
9 &25& &a_3^2a_2^4a_1^2 &16 &bba&      &      \cr
9 &25& &a_3^2a_2^2a_2a_1^4a_1 &4 &bbcba&      &      \cr
9 &25& &a_3a_3a_2^2a_2a_1^2a_1^2a_1 &2 &bbbbbbb&      &      \cr
9 &25& &a_3a_2^6a_1^2 &6 &abb&      &      \cr
9 &25& &a_3^2a_2^2a_1^4a_1^4 &8 &bcbb&      &      \cr
9 &25& &a_3a_2^4a_2a_1^4a_1 &8 &bbbbc&      &      \cr
9 &25& &a_3a_2^2a_2^2a_2a_1^2a_1^2a_1 &2 &bbbdbbb&      &      \cr
\cr
10 &22& &a_3a_2^{10} &2880 &ac&   A_2^{12} &      \cr
10 &23& &a_3^4a_1^8a_1^2 &384 &cfa&  A_3^8  &      \cr
10 &23& &a_3^3a_2^4a_1^2a_1^2 &48 &bbea& A_3^8   &      \cr
10 &24& &a_3^4a_2^2a_1^2 &32 &bab&    A_4^6&A_4^6    \cr
10 &24& &a_4a_3a_2^4a_1^4 &16 &bdbc&    A_4^6& A_4^6   \cr
10 &24& &a_3^2a_3a_2^4a_1^2 &16 &bbbe&    A_3^8&A_5^4D_4    \cr
10 &24& &a_3^4a_1^4a_1^4 &48 &cdf&    A_3^8& A_5^4D_4   \cr
10 &24& &a_3^4a_1^8 &384 &cd&    A_3^8&D_4^6   \cr
10 &24&E &a_2^{12} &190080 &a&   2A_2^{12}&      \cr
10 &25& &a_3^5 &1920 &c&      &      \cr
10 &25& &d_4a_2^4a_1^6 &144 &bcd&      &      \cr
10 &25& &d_4a_3a_1^{12} &576 &ced&      &      \cr
10 &25& &a_5a_1^{15} &720 &cf&      &      \cr
10 &25& &a_4a_3a_2^4a_1^2 &8 &bdbb&      &      \cr
10 &25& &a_3^3a_3a_2a_1^3 &6 &bcbb&      &      \cr
10 &25& &a_4a_3a_2^2a_2a_1^2a_1^2a_1 &2 &bdbbbce&      &      \cr
10 &25& &a_3^3a_3a_1^4a_1^2 &24 &bcce&      &      \cr
10 &25& &a_3^2a_3^2a_1^4a_1^2 &16 &cbbd&      &      \cr
10 &25& &a_4a_2^6a_1^2 &6 &abe&      &      \cr
10 &25& &a_4a_3a_2^2a_1^4a_1^2a_1^2 &4 &bdbccf&      &      \cr
10 &25& &a_3^2a_3a_2^2a_2a_1^2a_1 &2 &bbbbbc&      &      \cr
10 &25& &a_3^2a_3a_2^2a_2a_1a_1a_1 &2 &bbbadbc&      &      \cr
10 &25& &a_4a_2^5a_1^5 &10 &bbc&      &      \cr
10 &25& &a_3^2a_2^6 &48 &da&      &      \cr
10 &25& &a_3^2a_2^4a_2^2 &8 &bba&      &      \cr
10 &25& &a_3^3a_2^2a_1^3a_1^2a_1 &12 &bbdcf&      &      \cr
10 &25& &a_3a_3a_3a_2^2a_1^2a_1a_1a_1a_1 &2 &cbbbcebfd&      &      \cr
10 &25& &a_3a_3a_2^2a_2^2a_2a_1^2a_1 &2 &bbbbbee&      &      \cr
10 &25& &a_3^2a_2^4a_1^4a_1^2 &8 &dbcf&      &      \cr
10 &25& &a_3^3a_1^{12} &48 &cf&      &      \cr
10 &25& &a_3a_2^6a_2a_1^3 &12 &dbce&      &      \cr
\cr
11 &24& &a_4a_3^2a_3a_2^2a_1^2 &4 &bbbbb&    A_4^6&A_5^4D_4    \cr
11 &24& &a_4^2a_2^4a_1^4 &16 &dca&    A_4^6& A_5^4D_4   \cr
11 &24& &a_3^6 &240 &a&    A_4^6&D_4^6   \cr
11 &24& &a_3^4a_2^4 &24 &be&   A_3^8 &A_6^4    \cr
11 &25& &a_5a_2^4a_2a_1^4 &8 &befb&      &      \cr
11 &25& &d_4a_3a_2^4a_1^4 &8 &bcda&      &      \cr
11 &25& &a_4a_3^2a_3a_2a_1^2a_1 &2 &bbbcba&      &      \cr
11 &25& &a_4^2a_2^2a_2a_1^4a_1 &4 &dcbbb&      &      \cr
11 &25& &a_4a_3^2a_2^4 &4 &bbb&      &      \cr
11 &25& &a_4a_3^3a_1^6 &6 &cbb&      &      \cr
11 &25& &a_4a_3^2a_2^2a_2a_1^2a_1 &2 &bbbfab&      &      \cr
11 &25& &a_4a_3a_3a_2a_2a_2a_1a_1a_1 &1 &bbbbccbba&      &      \cr
11 &25& &a_4a_3a_3a_2a_2a_2a_1a_1a_1 &1 &bbbbecbbb&      &      \cr
11 &25& &a_3^4a_3a_1^4 &8 &bab&      &      \cr
11 &25& &a_3^2a_3^2a_2^2a_2a_1 &2 &bbbbb&      &      \cr
11 &25& &a_3^2a_3a_3a_2^2a_2a_1 &2 &bbabda&      &      \cr
11 &25& &a_4a_3a_2^2a_2^2a_2a_1^2a_1 &2 &bbeccba&      &      \cr
11 &25& &a_3^2a_3^2a_2^2a_1^4 &4 &bbdb&      &      \cr
\cr
12 &22& &a_3^6a_3a_1^2 &96 &dag&   A_3^8 &      \cr
12 &23& &a_4^2a_3^2a_2^2a_1^2a_1^2 &8 &bcbha& A_4^6   &      \cr
12 &23& &a_4a_3^5a_1^2 &40 &aba&    A_4^6&      \cr
12 &24& &d_4a_4a_2^6 &24 &dca&    A_5^4D_4& A_5^4D_4   \cr
12 &24& &d_4a_3^4a_1^4 &32 &cde&    A_5^4D_4& A_5^4D_4   \cr
12 &24& &a_5a_3^3a_1^6a_1 &24 &cfec&    A_5^4D_4&A_5^4D_4    \cr
12 &24& &a_4^2a_3^2a_3a_1^2 &8 &bbbd&    A_5^4D_4& A_5^4D_4   \cr
12 &24& &a_5a_3^2a_2^4a_1 &16 &bebf&    A_5^4D_4& A_5^4D_4   \cr
12 &24& &d_4a_3^4a_1^4 &48 &cdc&    D_4^6 &A_5^4D_4  \cr
12 &24& &d_4^2a_1^{16} &1152 &eb&    D_4^6&D_4^6  \cr
12 &24& &a_4^2a_3^2a_2^2a_1^2 &4 &bcbh&   A_4^6 &A_6^4    \cr
12 &24& &a_3^4a_3^2a_1^4 &32 &fdg&    A_3^8&A_7^2D_5^2    \cr
12 &25& &a_5a_3^2a_3a_1^4a_1 &8 &cefdc&      &      \cr
12 &25& &a_5a_3^2a_2^2a_2a_1a_1 &2 &bebbdh&      &      \cr
12 &25& &d_4a_3^2a_3a_2^2a_1^2 &4 &bddac&      &      \cr
12 &25& &d_4a_4a_2^4a_1^4 &8 &dcae&      &      \cr
12 &25& &a_5a_3^2a_2^2a_1^2a_1^2a_1 &4 &bfbdhf&      &      \cr
12 &25& &a_5a_3a_3a_2^2a_1^2a_1a_1a_1 &2 &bfebhdee&      &      \cr
12 &25& &a_4^2a_3a_3a_2a_1^2a_1 &2 &bbcade&      &      \cr
12 &25& &a_4a_4a_3a_3a_2a_1a_1a_1 &1 &bbccbddd&      &      \cr
12 &25& &a_4^2a_3a_2^4 &4 &bbb&      &      \cr
12 &25& &a_4^2a_3^2a_1^4a_1^2 &4 &bche&      &      \cr
12 &25& &d_4a_3^2a_2^4a_1a_1 &8 &bdaeg&      &      \cr
12 &25& &d_4a_3^3a_1^6a_1a_1 &12 &cdegb&      &      \cr
12 &25& &a_4a_3^2a_3^2a_2a_1 &2 &abcbc&      &      \cr
12 &25& &a_4^2a_3a_2^2a_2a_1a_1a_1 &2 &bcbbfdh&      &      \cr
12 &25& &a_4a_4a_3a_2a_2a_2a_1a_1a_1 &1 &bbcabbhhd&      &      \cr
12 &25& &a_4a_3^4a_1^4 &8 &ach&      &      \cr
12 &25& &a_4a_3^2a_3a_3a_1^2a_1^2 &2 &bbcfde&      &      \cr
12 &25& &a_4a_3^3a_2^3a_1 &6 &bbag&      &      \cr
12 &25& &a_4a_3a_3a_3a_2a_2a_2a_1 &1 &bbbebbah&      &      \cr
12 &25& &a_3^4a_3a_3a_1^2 &8 &bfdc&      &      \cr
12 &25& &a_4a_3^2a_3a_2^2a_1^2a_1a_1 &2 &bccbhge&      &      \cr
12 &25& &a_4a_3^2a_3a_2^2a_1^2a_1^2 &2 &bcfbhh&      &      \cr
12 &25& &a_3^2a_3^2a_3a_2^2a_1^2 &4 &edbbh&      &      \cr
\cr
13 &24& &a_4^4a_1^4 &24 &cc&    A_5^4D_4&A_6^4    \cr
13 &24& &a_5a_4a_3a_3a_2a_2a_1 &2 &bebbddd&  A_5^4D_4  & A_6^4   \cr
13 &24& &a_4^2a_3^4 &16 &bb&   A_4^6 &A_7^2D_5^2    \cr
13 &24& &a_4^2a_4a_3a_2^2a_1^2 &4 &ccahb&  A_4^6  & A_7^2D_5^2   \cr
13 &24&E &a_3^8 &2688 &a& 2A_3^8  &      \cr
13 &25& &a_6a_2^6a_1^3 &12 &dhd&      &      \cr
13 &25& &a_4^4a_1^2 &24 &ca&      &      \cr
13 &25& &a_5a_4a_3^2a_2a_1^2 &2 &bfbdc&      &      \cr
13 &25& &d_4a_4a_3^3a_1^2 &6 &bdac&      &      \cr
13 &25& &d_4^2a_2^6 &72 &cc&      &      \cr
13 &25& &a_5a_4a_3a_2a_2a_2a_1a_1 &1 &bebhdeda&      &      \cr
13 &25& &a_5a_4a_3a_2a_2a_2a_1a_1 &1 &bebhdddd&      &      \cr
13 &25& &d_4a_4a_3a_3a_2a_2a_1a_1 &1 &bdaaeccb&      &      \cr
13 &25& &a_4^2a_4a_3a_2^2 &2 &bcbd&      &      \cr
13 &25& &a_4a_4a_4a_3a_2a_1a_1a_1 &1 &bccbdbcd&      &      \cr
13 &25& &a_5a_3^3a_2^3 &6 &bbd&      &      \cr
13 &25& &a_5a_3^2a_3a_2^2a_2 &2 &abbhc&      &      \cr
13 &25& &a_5a_4a_2^2a_2^2a_2a_1^2 &2 &behdfd&      &      \cr
13 &25& &a_5a_3^2a_3a_2^2a_1^2a_1 &2 &bbbedd&      &      \cr
13 &25& &a_4a_4a_3a_3a_3a_2a_1 &1 &bcbbbdc&      &      \cr
13 &25& &a_4a_4a_3a_3a_3a_2a_1 &1 &bbbabdb&      &      \cr
13 &25& &a_4^3a_2^3a_1^3 &6 &cdb&      &      \cr
13 &25& &d_4a_3^3a_2^3a_2 &6 &bacg&      &      \cr
13 &25& &a_4^2a_3a_3a_2^2a_2a_1 &2 &ebaddd&      &      \cr
13 &25& &a_5a_2^9 &72 &cf&      &      \cr
\cr
14 &21& &d_4a_3^7 &336 &ag& A_3^8   &      \cr
14 &22& &a_4^4a_3a_2^2 &16 &aab& A_4^6   &      \cr
14 &23& &a_5a_4^2a_3a_3a_1^2a_1 &4 &bbddaf& A_5^4D_4   &      \cr
14 &23& &d_4a_4^3a_2^2a_1^2 &12 &caca&   A_5^4D_4 &      \cr
14 &23& &a_5d_4a_3^3a_1^3a_1^2 &12 &dfega& A_5^4D_4   &      \cr
14 &23& &a_5^2a_3a_2^4a_1^2 &16 &efda&   A_5^4D_4 &      \cr
14 &23& &d_4^2a_3^4a_1^2 &96 &dcd&   D_4^6&      \cr
14 &24& &a_5a_4^3a_1^3 &12 &cbe&   A_6^4 & A_6^4   \cr
14 &24& &a_5^2a_3^2a_2^2 &8 &eda&    A_6^4&A_6^4    \cr
14 &24& &a_6a_3^3a_2^3 &12 &bhd&    A_6^4& A_6^4   \cr
14 &24& &d_4a_4^2a_4a_2^2 &4 &caac&   A_5^4D_4 & A_7^2D_5^2   \cr
14 &24& &a_5d_4a_3^2a_3a_1^2a_1 &4 &dfecbg&  A_5^4D_4  & A_7^2D_5^2   \cr
14 &24& &a_5^2a_3^2a_1^2a_1^2a_1^2 &8 &febcg&    A_5^4D_4& A_7^2D_5^2   \cr
14 &24& &a_5a_4^2a_3a_3a_1 &4 &bbddf&    A_5^4D_4&A_7^2D_5^2    \cr
14 &24& &d_4^2a_3^4 &32 &dc&   D_4^6& A_7^2D_5^2   \cr
14 &24& &a_4^3a_3^3 &12 &bh&    A_4^6&A_8^3    \cr
14 &24& &a_3^8 &384 &g&    A_3^8& D_6^4  \cr
14 &25& &d_5a_3^2a_2^4a_1^2 &8 &bgbb&      &      \cr
14 &25& &d_5a_3^3a_1^8 &48 &cgc&      &      \cr
14 &25& &a_6a_3^2a_3a_2^2a_1 &2 &bhhcf&      &      \cr
14 &25& &a_5^2a_3a_3a_1^4 &8 &efee&      &      \cr
14 &25& &a_6a_3^2a_3a_2a_1^2a_1^2 &2 &bhhdeg&      &      \cr
14 &25& &d_4a_4^3a_1^3a_1 &6 &cabc&      &      \cr
14 &25& &a_5d_4a_3a_3a_2^2a_1 &2 &deeccb&      &      \cr
14 &25& &a_5a_5a_3a_2^2a_1^2a_1a_1 &2 &efddefg&      &      \cr
14 &25& &a_5a_4^2a_3a_2a_1^2 &2 &bbfdf&      &      \cr
14 &25& &a_5a_4a_4a_3a_2a_1a_1 &1 &cbbdabe&      &      \cr
14 &25& &d_4a_4^2a_3^2a_1^2 &4 &baeb&      &      \cr
14 &25& &d_4^2a_3^2a_3a_1^4 &8 &dcbb&      &      \cr
14 &25& &a_5a_4^2a_3a_1^2a_1^2a_1 &2 &cbhcee&      &      \cr
14 &25& &a_5d_4a_3^2a_1^4a_1^2a_1 &4 &dfegcg&      &      \cr
14 &25& &a_5a_4a_3a_3a_3a_2 &1 &bbhddc&      &      \cr
14 &25& &a_5a_4a_3a_3a_3a_2 &1 &bbeddc&      &      \cr
14 &25& &a_5a_4^2a_2^2a_2a_1^2 &2 &cbddf&      &      \cr
14 &25& &a_4^4a_2^2 &8 &ba&      &      \cr
14 &25& &d_4a_4a_4a_3a_2a_2a_1a_1 &1 &caaecbbg&      &      \cr
14 &25& &a_5a_4a_3a_3a_3a_1a_1a_1 &1 &bbhedbfg&      &      \cr
14 &25& &a_5a_3^2a_3^2a_3a_1 &4 &dddcg&      &      \cr
14 &25& &a_5a_3^2a_3a_3a_3a_1 &2 &bheddf&      &      \cr
14 &25& &a_5a_4a_3^2a_2^2a_1^2a_1 &2 &cbhdge&      &      \cr
14 &25& &a_4^2a_4a_3^2a_2a_1 &2 &bbdbf&      &      \cr
14 &25& &a_4a_4a_4a_3a_3a_2a_1 &1 &abbdhcf&      &      \cr
14 &25& &a_4^2a_4a_3a_2^2a_2a_1 &2 &bahbdf&      &      \cr
14 &25& &d_4a_3^4a_3a_1^4 &8 &fegg&      &      \cr
\cr
15 &24& &a_5^2a_4^2a_1^2 &4 &bdc&  A_6^4  & A_7^2D_5^2   \cr
15 &24& &a_6a_4a_4a_3a_2a_1a_1 &2 &chhceaa& A_6^4   & A_7^2D_5^2   \cr
15 &24& &a_5^2a_4a_3a_2^2 &4 &bfdf&   A_5^4D_4 & A_8^3   \cr
15 &24& &a_4^4a_3^2 &16 &hb&    A_4^6& A_9^2D_6   \cr
15 &25& &d_5a_3^5 &20 &ab&      &      \cr
15 &25& &d_5a_4a_3^2a_2^2a_1^2 &2 &bgbba&      &      \cr
15 &25& &a_6a_4^2a_3a_1^2a_1 &2 &chdcb&      &      \cr
15 &25& &a_6a_4a_4a_2a_2a_1a_1a_1 &1 &chhefacc&      &      \cr
15 &25& &a_5d_4a_4^2a_1^2a_1 &2 &abeab&      &      \cr
15 &25& &a_5a_5a_4a_3a_2a_1 &1 &bbddea&      &      \cr
15 &25& &a_6a_4a_3a_3a_2a_2a_1 &1 &bhdcfga&      &      \cr
15 &25& &a_6a_4a_3^2a_2^2a_1 &2 &bhdfc&      &      \cr
15 &25& &d_4^2a_4^2a_2^2 &4 &acb&      &      \cr
15 &25& &a_5^2a_4a_3a_1^2a_1^2 &2 &bedcc&      &      \cr
15 &25& &a_5d_4a_4a_3a_2a_2a_1 &1 &abecbba&      &      \cr
15 &25& &a_5^2a_3^2a_3a_1^2 &2 &bdac&      &      \cr
15 &25& &a_5a_4^2a_4a_2a_1^2 &2 &adfec&      &      \cr
15 &25& &a_5a_4^2a_4a_2a_1^2 &2 &addca&      &      \cr
15 &25& &a_5a_4a_4a_4a_2a_1a_1 &1 &bhddeca&      &      \cr
15 &25& &a_4^5 &5 &d&      &      \cr
15 &25& &a_5a_4a_4a_3a_3a_2 &1 &bhdcab&      &      \cr
15 &25& &d_4a_4^2a_3^2a_3 &2 &bccb&      &      \cr
\cr
16 &21& &d_4a_4^5 &40 &ab&   A_4^6 &      \cr
16 &22& &a_5^2d_4a_3^2a_3 &8 &ceba&    A_5^4D_4&      \cr
16 &22& &a_5^3a_3a_3a_1^3 &12 &ecad&   A_5^4D_4 &      \cr
16 &22& &d_4^4a_3a_1^6 &144 &bdc&   D_4^6&      \cr
16 &23& &a_6a_5a_4a_3a_2a_1^2a_1 &2 &bhdecah&  A_6^4  &      \cr
16 &23& &a_5^3a_4a_1^2a_1 &6 &daag&    A_6^4&      \cr
16 &24& &d_5d_4a_3^4 &16 &dgb&  A_7^2D_5^2  &   A_7^2D_5^2 \cr
16 &24& &d_5a_5a_3^2a_3a_1^2a_1 &4 &fgbceb&  A_7^2D_5^2  &  A_7^2D_5^2  \cr
16 &24& &a_5^2d_4^2a_1^2 &16 &cef&  A_7^2D_5^2  & A_7^2D_5^2   \cr
16 &24& &a_7a_3^4a_1^4 &16 &fge&  A_7^2D_5^2  &  A_7^2D_5^2  \cr
16 &24& &d_5a_4a_4^2a_2^2 &4 &cbba&  A_7^2D_5^2  & A_7^2D_5^2   \cr
16 &24& &a_6d_4a_4^2a_2 &4 &ahcb&  A_7^2D_5^2  &  A_7^2D_5^2  \cr
16 &24& &a_6a_5a_4a_3a_2a_1 &2 &bhdech&  A_6^4  &  A_8^3  \cr
16 &24& &a_5^2d_4a_3^2a_1^2 &4 &eegd&   A_5^4D_4 & A_9^2D_6  \cr
16 &24& &a_5a_5a_4^2a_3 &4 &dddf&    A_5^4D_4&  A_9^2D_6  \cr
16 &24& &a_5^2d_4a_3^2a_1^2 &8 &eebd&   A_5^4D_4 &  D_6^4 \cr
16 &24& &d_4^4a_1^8 &48 &bc&  D_4^6 & D_6^4  \cr
16 &24&E &a_4^6 &240 &a&  2A_4^6 &      \cr
16 &25& &a_7a_3^4a_1^2 &8 &fff&      &      \cr
16 &25& &a_7a_3a_3a_3a_2^2a_1^2 &2 &efggch&      &      \cr
16 &25& &a_5d_4^3a_1^3 &36 &bcb&      &      \cr
16 &25& &d_5a_4^2a_3a_3a_1a_1 &2 &bbcbdb&      &      \cr
16 &25& &d_5d_4a_3^3a_1^3a_1 &6 &dgbec&      &      \cr
16 &25& &d_5a_5a_3a_2^4a_1 &4 &egcad&      &      \cr
16 &25& &a_6a_5a_4a_2a_2a_1a_1 &1 &bhdcagh&      &      \cr
16 &25& &d_5a_4a_3^4 &4 &bbb&      &      \cr
16 &25& &d_5a_4^2a_3a_2^2a_1^2 &2 &cbcae&      &      \cr
16 &25& &a_6d_4a_4a_3a_2a_1a_1 &1 &ahcgafe&      &      \cr
16 &25& &a_5^2a_5a_3a_2 &2 &deeb&      &      \cr
16 &25& &a_6a_5a_3a_3a_2a_2 &1 &bhfebc&      &      \cr
16 &25& &a_6a_5a_3a_3a_2a_1a_1a_1 &1 &bhegchhh&      &      \cr
16 &25& &a_6a_4^2a_3^2a_1 &2 &bcee&      &      \cr
16 &25& &a_5a_5a_5a_2^2a_1^2a_1 &2 &defchd&      &      \cr
16 &25& &a_5^2d_4a_3a_2^2 &4 &cfcb&      &      \cr
16 &25& &a_5a_5d_4a_3a_2^2 &2 &ccdga&      &      \cr
16 &25& &a_6d_4a_3^2a_2^2a_2 &2 &ahgab&      &      \cr
16 &25& &a_5a_5a_4a_4a_2a_1 &1 &ddddcg&      &      \cr
16 &25& &a_5a_5a_4a_4a_2a_1 &1 &ddadad&      &      \cr
16 &25& &a_6a_4a_4a_3a_2a_2a_1 &1 &bddecah&      &      \cr
16 &25& &a_5d_4a_4^2a_3a_1 &2 &edccb&      &      \cr
16 &25& &a_5d_4^2a_3^2a_1a_1a_1 &2 &cebecf&      &      \cr
16 &25& &a_5^2a_4a_3^2a_1^2 &2 &haeh&      &      \cr
16 &25& &a_5^2a_4a_3a_3a_1^2 &2 &ddebd&      &      \cr
16 &25& &a_5a_5a_4a_3a_3a_1a_1 &1 &eddefdg&      &      \cr
16 &25& &a_5^2a_3^4 &8 &cf&      &      \cr
16 &25& &a_5^2a_4a_3a_2^2a_1^2 &2 &hdcch&      &      \cr
16 &25& &d_4a_4^2a_4a_3^2 &2 &hcbb&      &      \cr
\cr
17 &24& &a_7a_4^2a_3^2 &4 &bfc&   A_7^2D_5^2 &  A_8^3  \cr
17 &24& &a_6^2a_4a_3a_1^2 &4 &hebb& A_7^2D_5^2   &A_8^3    \cr
17 &24& &a_6a_5a_5a_3a_2 &2 &dddcg&   A_6^4 &  A_9^2D_6  \cr
17 &24& &a_6^2a_3^2a_2^2 &4 &hah&  A_6^4  &  A_9^2D_6  \cr
17 &24& &a_5^4 &24 &a& A_6^4   & D_6^4  \cr
17 &25& &a_7a_4a_4a_3a_2a_1 &1 &bfgcgb&      &      \cr
17 &25& &a_6^2a_3^2a_2a_1 &2 &hcfa&      &      \cr
17 &25& &a_7a_4^2a_2^2a_2a_1 &2 &bfhea&      &      \cr
17 &25& &d_5d_4a_4^2a_2^2 &2 &abbb&      &      \cr
17 &25& &d_5a_5a_4a_3a_3a_1 &1 &bbbaab&      &      \cr
17 &25& &d_5a_5a_4a_3a_2a_2a_1 &1 &bbbaedb&      &      \cr
17 &25& &a_6a_5d_4a_3a_2a_1 &1 &ecdafb&      &      \cr
17 &25& &a_6a_5a_4^2a_1^2 &2 &ddeb&      &      \cr
17 &25& &a_6a_5a_4a_3a_3 &1 &dcgbc&      &      \cr
17 &25& &a_6a_5a_4a_3a_3 &1 &dcfcc&      &      \cr
17 &25& &a_6a_5a_4a_3a_3 &1 &dceac&      &      \cr
17 &25& &a_5^2d_4a_4a_3 &2 &cdbb&      &      \cr
17 &25& &a_6a_4a_4a_4a_3a_1 &1 &deffab&      &      \cr
\cr
18 &21& &a_5^3d_4d_4a_1 &12 &bcab&  A_5^4D_4  &      \cr
18 &22& &a_6^2a_4^2a_3 &4 &caa&  A_6^4  &      \cr
18 &23& &a_6d_5a_4a_4a_2a_1^2 &2 &bhaaba&  A_7^2D_5^2  &      \cr
18 &23& &a_6^2d_4a_4a_1^2 &4 &ceba&  A_7^2D_5^2  &      \cr
18 &23& &d_5a_5^2d_4a_1^2a_1^2 &4 &ebcfa&  A_7^2D_5^2  &      \cr
18 &23& &a_7a_5a_4^2a_1^2a_1 &4 &dfcag&  A_7^2D_5^2  &      \cr
18 &23& &a_7d_4^2a_3^2a_1^2 &8 &cgfa&  A_7^2D_5^2  &      \cr
18 &23& &d_5^2a_3^4a_1^2 &16 &gea&   A_7^2D_5^2 &      \cr
18 &24& &a_7a_5^2a_2^2 &8 &ffa&   A_8^3 & A_8^3   \cr
18 &24& &a_7a_5a_4^2a_1 &4 &dfcg&   A_7^2D_5^2 & A_9^2D_6   \cr
18 &24& &a_7a_5d_4a_3a_1a_1a_1 &2 &eggfdfc&A_7^2D_5^2    &A_9^2D_6    \cr
18 &24& &a_6d_5a_4a_4a_2 &2 &bhaab&  A_7^2D_5^2  &  A_9^2D_6  \cr
18 &24& &d_5^2a_3^4 &16 &gb&  A_7^2D_5^2  & D_6^4  \cr
18 &24& &d_5a_5^2d_4a_1^2 &4 &ebbc& A_7^2D_5^2   &  D_6^4 \cr
18 &24& &a_5^2a_5d_4a_3a_1 &4 &gbcdb& A_5^4D_4   &  A_{11}D_7E_6 \cr
18 &24& &a_5^4a_1^4 &48 &cb&   A_5^4D_4 &E_6^4  \cr
18 &25& &d_6a_3^4a_3a_1^2 &8 &ddcd&      &      \cr
18 &25& &a_7a_5^2a_1^4 &8 &fge&      &      \cr
18 &25& &a_7a_5d_4a_2^2a_1 &2 &egfbc&      &      \cr
18 &25& &a_7d_4^2a_3a_1^4 &4 &cgef&      &      \cr
18 &25& &a_7a_5a_4a_3a_2 &1 &dfchb&      &      \cr
18 &25& &a_6^2a_5a_2a_1a_1 &2 &deaed&      &      \cr
18 &25& &a_7a_5a_4a_3a_1a_1a_1 &1 &dgchfeg&      &      \cr
18 &25& &a_6d_5a_4a_3a_2a_1a_1 &1 &bhaebfd&      &      \cr
18 &25& &d_5a_5a_5a_4a_1a_1 &1 &dbcacc&      &      \cr
18 &25& &a_7a_5a_3^2a_3a_1 &2 &dggfd&      &      \cr
18 &25& &a_7a_5a_3^2a_3a_1 &2 &dfhfg&      &      \cr
18 &25& &a_6a_6a_4a_3a_2a_1 &1 &cdchbg&      &      \cr
18 &25& &a_6a_5a_5a_4a_1 &1 &aeecc&      &      \cr
18 &25& &d_5a_5d_4a_3^2a_1^2a_1 &2 &ebcefb&      &      \cr
18 &25& &d_5d_4^2a_3^2a_3 &4 &cbbc&      &      \cr
18 &25& &d_5a_5a_4^2a_3a_1 &2 &dbadf&      &      \cr
18 &25& &d_5a_5a_4^2a_3a_1 &2 &dbabb&      &      \cr
18 &25& &a_6a_5d_4a_4a_2a_1 &1 &cgeabf&      &      \cr
18 &25& &a_6a_5a_4a_4a_3 &1 &cfacg&      &      \cr
18 &25& &a_5^2d_4^2a_3a_1^2 &2 &bgff&      &      \cr
18 &25& &a_6a_4^2a_4^2a_1 &2 &bcag&      &      \cr
\cr
19 &24& &a_8a_4^2a_3^2 &4 &hhb&   A_8^3 &   A_9^2D_6 \cr
19 &24& &a_7a_6a_5a_2a_1 &2 &dfceb& A_8^3   &   A_9^2D_6 \cr
19 &24& &a_6a_6a_5a_4a_1 &2 &eeaha&  A_6^4  &A_{11}D_7E_6   \cr
19 &24&E &d_4^6 &2160 &d&  2D_4^6&      \cr
19 &24&E &a_5^4d_4 &48 &aa& 2A_5^4D_4  &      \cr
19 &25& &d_6a_4^2a_4a_2^2 &2 &addc&      &      \cr
19 &25& &a_8a_4a_4a_3a_2a_1 &1 &hghbfb&      &      \cr
19 &25& &a_8a_4^2a_2^2a_2a_1 &2 &hhgeb&      &      \cr
19 &25& &a_7a_6a_4a_2a_2a_1 &1 &dfheeb&      &      \cr
19 &25& &d_5^2a_4a_4a_2^2 &2 &bbec&      &      \cr
19 &25& &a_6d_5d_4a_4a_2 &1 &bcaec&      &      \cr
19 &25& &a_6d_5a_5a_3a_2a_1 &1 &bdabca&      &      \cr
19 &25& &a_7a_5^2a_3a_1^2 &2 &dcab&      &      \cr
19 &25& &d_5a_5^3a_1 &3 &aaa&      &      \cr
19 &25& &a_7a_5^2a_2^2a_1^2 &2 &dcdb&      &      \cr
19 &25& &a_7a_5a_4a_4a_2 &1 &achgc&      &      \cr
19 &25& &a_7d_4a_4a_4a_3 &1 &ccefb&      &      \cr
19 &25& &a_6^2a_5a_3a_2 &2 &fbad&      &      \cr
19 &25& &a_6a_5^3 &6 &cb&      &      \cr
19 &25& &a_6a_5^2a_5 &2 &ecb&      &      \cr
19 &25& &a_6a_5a_5d_4a_2 &1 &babcc&      &      \cr
19 &25& &d_5a_5^2a_4a_2^2 &2 &dadc&      &      \cr
\cr
20 &20& &d_5a_5^4 &16 &ad&  A_5^4D_4  &      \cr
20 &20& &d_5d_4^5 &120 &dc& D_4^6  &      \cr
20 &21& &a_6^3d_4a_2 &6 &aaa&  A_6^4  &      \cr
20 &22& &a_7d_5a_5a_3a_3a_1 &2 &bgecad& A_7^2D_5^2   &      \cr
20 &22& &d_5^2a_5^2a_3 &8 &bba&  A_7^2D_5^2  &      \cr
20 &22& &a_7^2a_3^2a_3a_1^2 &8 &gdae& A_7^2D_5^2   &      \cr
20 &23& &a_7a_6^2a_1^2a_1^2 &4 &ecga& A_8^3   &      \cr
20 &23& &a_7^2a_4a_3a_1^2 &4 &faea&  A_8^3  &      \cr
20 &23& &a_8a_5^2a_2^2a_1^2 &4 &dhba&  A_8^3  &      \cr
20 &24& &a_7^2d_4a_1^4 &8 &gff&  A_9^2D_6  &   A_9^2D_6 \cr
20 &24& &d_6a_5^2a_3^2 &8 &edd&  A_9^2D_6  &   A_9^2D_6 \cr
20 &24& &a_8a_5^2a_3 &4 &dge&  A_9^2D_6  & A_9^2D_6   \cr
20 &24& &d_6a_5^2a_3^2 &4 &edc&  D_6^4 &A_9^2D_6    \cr
20 &24& &d_6d_4^3a_1^6 &12 &bcb&  D_6^4 &D_6^4   \cr
20 &24& &a_7d_5a_5a_3a_1a_1a_1 &2 &cgefecd& A_7^2D_5^2   &A_{11}D_7E_6   \cr
20 &24& &d_5d_5a_5^2a_1^2 &4 &bced& A_7^2D_5^2   &  A_{11}D_7E_6 \cr
20 &24& &a_7d_5d_4a_3a_3 &2 &bgecf& A_7^2D_5^2   &  A_{11}D_7E_6 \cr
20 &24& &a_6a_6d_5a_4 &2 &aaeb&   A_7^2D_5^2 & A_{11}D_7E_6  \cr
20 &24& &d_5^2a_5^2a_1^2 &8 &cbc&  A_7^2D_5^2  &E_6^4  \cr
20 &24& &a_6^2a_5^2 &4 &ch&  A_6^4  &  A_{12}^2 \cr
20 &24& &d_4^6 &48 &c&  D_4^6 & D_8^3  \cr
20 &25& &a_8a_5^2a_1^2a_1^2 &2 &dhfg&      &      \cr
20 &25& &d_6a_5^2a_3a_1^2a_1a_1 &2 &eddfeb&      &      \cr
20 &25& &d_6a_5d_4a_3^2a_1 &2 &ccdcd&      &      \cr
20 &25& &a_7a_7a_3a_2^2a_1^2 &2 &fggbg&      &      \cr
20 &25& &a_8d_4a_4a_3a_3 &1 &chbff&      &      \cr
20 &25& &d_5^2a_5d_4a_1^2a_1 &2 &bbecb&      &      \cr
20 &25& &a_7d_5a_4^2a_1a_1 &2 &cfbcf&      &      \cr
20 &25& &a_7a_6a_5a_3 &1 &ecfe&      &      \cr
20 &25& &a_6^2d_5a_3a_1^2 &2 &aedd&      &      \cr
20 &25& &a_7d_5a_4a_3a_3 &1 &bfbfc&      &      \cr
20 &25& &a_7a_6a_5a_2a_1a_1a_1 &1 &echbgfe&      &      \cr
20 &25& &a_7a_6a_4a_4a_1 &1 &ecbad&      &      \cr
20 &25& &a_6^2a_6a_3a_1 &2 &abee&      &      \cr
20 &25& &a_6a_6a_6a_3a_1 &1 &caceg&      &      \cr
20 &25& &a_6d_5a_5a_4a_2a_1 &1 &aeebad&      &      \cr
20 &25& &a_6a_6a_5d_4a_1 &1 &abfhd&      &      \cr
20 &25& &a_7a_5^2a_4a_1^2 &2 &ehag&      &      \cr
20 &25& &d_5a_5a_5a_4^2 &2 &fbdb&      &      \cr
\cr
21 &24& &a_8a_6a_5a_2a_1 &2 &ehbga&  A_8^3  & A_{11}D_7E_6  \cr
21 &24& &a_7^2a_4^2 &4 &cg& A_7^2D_5^2   & A_{12}^2  \cr
21 &25& &d_6a_6a_4a_4a_2 &1 &cdcdd&      &      \cr
21 &25& &a_8a_6a_4a_3a_1 &1 &eggbb&      &      \cr
21 &25& &a_7a_6d_5a_2a_1a_1 &1 &aecfba&      &      \cr
21 &25& &a_6d_5d_5a_4a_2 &1 &baacc&      &      \cr
21 &25& &a_7d_5a_5a_4a_1 &1 &acbdb&      &      \cr
21 &25& &a_7a_6a_5a_4 &1 &cgbe&      &      \cr
\cr
22 &21& &a_7d_5^2d_4a_3 &4 &beac& A_7^2D_5^2   &      \cr
22 &22& &a_8a_6^2a_3 &4 &bba&   A_8^3 &      \cr
22 &23& &a_8a_7a_5a_1^2a_1 &2 &cgeac&  A_9^2D_6  &      \cr
22 &23& &a_8a_6d_5a_2a_1^2 &2 &abhba&   A_9^2D_6 &      \cr
22 &23& &a_7d_6a_5a_3a_1^2a_1 &2 &dgdfab& A_9^2D_6   &      \cr
22 &23& &a_9a_5a_4^2a_1^2 &4 &fgba& A_9^2D_6   &      \cr
22 &23& &d_6d_5a_5^2a_1^2 &4 &bdcd&  D_6^4 &      \cr
22 &23& &d_5^4a_1^2 &48 &bd&  D_6^4 &      \cr
22 &24& &a_9a_5d_4a_3a_1a_1 &2 &gfffeb& A_9^2D_6   & A_{11}D_7E_6  \cr
22 &24& &a_7d_6a_5a_3a_1 &2 &dgcfe&  A_9^2D_6  & A_{11}D_7E_6  \cr
22 &24& &a_8a_6d_5a_2 &2 &abhb& A_9^2D_6   &A_{11}D_7E_6   \cr
22 &24& &d_6d_5a_5^2 &4 &bdc&  D_6^4 & A_{11}D_7E_6  \cr
22 &24& &d_5^4 &48 &b& D_6^4  & E_6^4 \cr
22 &24& &a_8a_7a_4a_3 &2 &chbg&  A_8^3  & A_{12}^2  \cr
22 &24& &a_7d_5^2a_3^2 &4 &eed& A_7^2D_5^2   & D_8^3  \cr
22 &24& &a_7^2d_4^2 &8 &fd&  A_7^2D_5^2  & D_8^3  \cr
22 &24&E &a_6^4 &24 &a& 2A_6^4  &      \cr
22 &25& &e_6d_4^2a_3^3 &12 &cbc&      &      \cr
22 &25& &e_6a_5a_4^2a_3a_1 &2 &dbace&      &      \cr
22 &25& &e_6a_5^2a_2^4 &8 &fba&      &      \cr
22 &25& &d_7a_3^6 &24 &ge&      &      \cr
22 &25& &a_9a_5d_4a_2^2 &2 &gfgb&      &      \cr
22 &25& &a_9a_5a_4a_3a_1a_1 &1 &ffbgbc&      &      \cr
22 &25& &a_9a_5a_3^2a_3 &2 &fggf&      &      \cr
22 &25& &a_8a_7a_4a_2a_1 &1 &chbbc&      &      \cr
22 &25& &d_6d_5a_5a_3a_3a_1 &1 &bdcdbe&      &      \cr
22 &25& &a_7d_6a_3^2a_3a_1^2 &2 &dgfeb&      &      \cr
22 &25& &a_8d_5a_5a_3a_1 &1 &agffe&      &      \cr
22 &25& &a_8d_5a_5a_2a_2a_1 &1 &ahfbab&      &      \cr
22 &25& &a_8a_6a_5a_3 &1 &bbgg&      &      \cr
22 &25& &a_8a_6a_5a_3 &1 &cbef&      &      \cr
22 &25& &d_6a_5^2d_4a_3 &2 &bcdb&      &      \cr
22 &25& &d_6a_5^3a_1^3 &6 &cdb&      &      \cr
22 &25& &a_8a_5^2d_4 &2 &bfe&      &      \cr
22 &25& &a_8a_6a_4a_4a_1 &1 &cbbbc&      &      \cr
22 &25& &a_7a_7a_5a_3a_1 &1 &ghefc&      &      \cr
22 &25& &d_5^3a_3^3 &6 &ec&      &      \cr
22 &25& &a_7d_5d_4^2a_3 &2 &bffc&      &      \cr
\cr
23 &24& &a_8^2a_3^2 &4 &hb&  A_9^2D_6  &  A_{12}^2 \cr
23 &24& &a_9a_6a_5a_2 &2 &cgbc&  A_9^2D_6  &A_{12}^2   \cr
23 &24& &a_7^3 &12 &a&   A_8^3 & D_8^3  \cr
23 &25& &d_7a_4^4 &4 &bd&      &      \cr
23 &25& &a_9a_6a_4a_3 &1 &cfgb&      &      \cr
23 &25& &a_8d_5^2a_2^2 &2 &ebe&      &      \cr
23 &25& &d_6a_6a_6a_4 &1 &accf&      &      \cr
23 &25& &d_6a_6^2a_4 &2 &bce&      &      \cr
23 &25& &a_8a_7d_4a_3 &1 &fbbb&      &      \cr
23 &25& &a_7^2d_5a_3 &2 &baa&      &      \cr
23 &25& &a_6^2d_5^2 &2 &cb&      &      \cr
\cr
24 &20& &a_7^2d_5d_5 &4 &cda&   A_7^2D_5^2 &      \cr
24 &21& &a_8^2d_4a_4 &4 &baa&   A_8^3 &      \cr
24 &22& &a_9a_7d_4a_3a_1 &2 &fffad&  A_9^2D_6  &      \cr
24 &22& &a_7^2d_6a_3 &4 &cfa&  A_9^2D_6  &      \cr
24 &22& &d_6^2d_4^2a_3a_1^2 &4 &cbdb&  D_6^4 &      \cr
24 &24& &d_7a_5^2a_5a_1 &4 &cdea&  A_{11}D_7E_6 & A_{11}D_7E_6  \cr
24 &24& &a_9d_5^2a_1^2a_1 &4 &efec&  A_{11}D_7E_6 &A_{11}D_7E_6   \cr
24 &24& &a_7e_6d_4a_3^2 &4 &bgce&   A_{11}D_7E_6& A_{11}D_7E_6  \cr
24 &24& &e_6a_6^2a_4 &4 &eaa&  A_{11}D_7E_6 & A_{11}D_7E_6  \cr
24 &24& &e_6d_5a_5^2a_1^2 &4 &cbca& E_6^4 & A_{11}D_7E_6  \cr
24 &24& &d_6^2d_4^2a_1^4 &4 &cbb& D_6^4  &  D_8^3 \cr
24 &24& &a_7^2d_6a_1^2 &4 &dfd&  A_9^2D_6  &  D_8^3 \cr
24 &24& &a_7^2d_5d_4 &4 &fde&    A_7^2D_5^2&  A_{15}D_9 \cr
24 &25& &d_7a_5^2d_4a_1a_1 &2 &bdeeb&      &      \cr
24 &25& &a_9a_7a_3^2a_1 &2 &fgcd&      &      \cr
24 &25& &d_6d_5^2a_5a_1 &2 &bcba&      &      \cr
24 &25& &a_7d_6d_5a_3a_1a_1 &1 &cedeeb&      &      \cr
24 &25& &a_8a_7a_6a_1 &1 &aebd&      &      \cr
24 &25& &a_7d_6a_5d_4a_1 &1 &cfdfb&      &      \cr
24 &25& &a_7^2a_7a_1^2 &2 &edd&      &      \cr
24 &25& &a_8a_7d_4a_4 &1 &bfgb&      &      \cr
24 &25& &a_7a_7d_5a_4 &1 &ccgb&      &      \cr
24 &25& &a_7^2d_5a_4 &2 &cea&      &      \cr
\cr
25 &24& &a_{10}a_6a_5a_1 &2 &hgbb&  A_{11}D_7E_6 &  A_{12}^2 \cr
25 &24& &a_8a_7^2 &4 &eb&    A_8^3&  A_{15}D_9 \cr
25 &24&E &a_7^2d_5^2 &8 &aa&  2A_7^2D_5^2 &      \cr
25 &25& &d_7a_6a_6a_2a_2 &1 &adedc&      &      \cr
25 &25& &a_{10}a_6a_4a_2a_1 &1 &hgcea&      &      \cr
25 &25& &e_6a_6d_5a_4a_2 &1 &acaea&      &      \cr
25 &25& &e_6a_6^2d_4 &2 &bca&      &      \cr
25 &25& &a_7e_6a_5a_4a_1 &1 &acaeb&      &      \cr
25 &25& &a_{10}a_5a_4a_4 &1 &gbcb&      &      \cr
25 &25& &a_8d_6a_6a_2 &1 &dbfc&      &      \cr
25 &25& &a_9a_6d_5a_2a_1 &1 &bfbeb&      &      \cr
25 &25& &a_9a_6^2a_2 &2 &agc&      &      \cr
25 &25& &a_9d_5a_5a_4 &1 &bbbe&      &      \cr
\cr
26 &19& &a_7^2d_6d_5 &4 &dae&  A_7^2D_5^2  &      \cr
26 &21& &a_9d_6a_5d_4 &2 &cfea& A_9^2D_6   &      \cr
26 &23& &a_9d_6d_5a_1^2a_1 &2 &cffab& A_{11}D_7E_6  &      \cr
26 &23& &a_{10}a_6d_5a_1^2 &2 &bbga&  A_{11}D_7E_6 &      \cr
26 &23& &a_8e_6a_6a_2a_1^2 &2 &ahaba& A_{11}D_7E_6  &      \cr
26 &23& &d_7a_7d_5a_3a_1^2 &2 &edeea&  A_{11}D_7E_6 &      \cr
26 &23& &e_6d_6a_5^2a_1^2 &4 &dbea&  A_{11}D_7E_6 &      \cr
26 &23& &e_6d_5^3a_1^2 &12 &bce&  E_6^4&      \cr
26 &24& &a_9^2a_2^2 &8 &ga&  A_{12}^2 & A_{12}^2  \cr
26 &24& &d_7a_7d_5a_3 &2 &eddc&  A_{11}D_7E_6 &D_8^3   \cr
26 &24& &a_9d_6a_5a_3 &2 &dfbd&  A_9^2D_6  & A_{15}D_9  \cr
26 &24& &a_9a_8a_5 &2 &ebc&  A_9^2D_6  &  A_{15}D_9 \cr
26 &24& &a_7^2d_5^2 &8 &ce&  A_7^2D_5^2  &  D_{10}E_7^2\cr
26 &25& &d_7d_5^2a_3a_3 &2 &bdba&      &      \cr
26 &25& &a_{10}d_5a_5a_2a_1 &1 &bgbbb&      &      \cr
26 &25& &a_7d_6^2a_3 &2 &bcc&      &      \cr
26 &25& &a_9d_6d_4a_3a_1 &1 &cfbeb&      &      \cr
26 &25& &a_9a_7a_6 &1 &egb&      &      \cr
26 &25& &a_9a_7a_6 &1 &efb&      &      \cr
26 &25& &a_7d_6d_5a_5a_1 &1 &dffeb&      &      \cr
\cr
27 &24& &a_8^2a_7 &4 &ga&  A_8^3  &  A_{17}E_7 \cr
27 &25& &a_8d_7a_4a_4 &1 &dbde&      &      \cr
27 &25& &a_9^2a_3a_1^2 &2 &baa&      &      \cr
27 &25& &a_{10}a_6^2a_1 &2 &eca&      &      \cr
\cr
28 &20& &a_9^2d_5a_1^2 &4 &fac&  A_9^2D_6  &      \cr
28 &20& &d_6^3d_5a_1^2 &6 &bdb&  D_6^4 &      \cr
28 &22& &a_9e_6d_5a_3a_1 &2 &cfead& A_{11}D_7E_6  &      \cr
28 &22& &a_9d_7a_5a_3 &2 &dfca&   A_{11}D_7E_6&      \cr
28 &22& &a_{11}d_5a_5a_3a_1 &2 &fbeae&  A_{11}D_7E_6 &      \cr
28 &22& &e_6^2a_5^2a_3 &8 &bae&  E_6^4&      \cr
28 &23& &a_{10}a_9a_2a_1^2a_1 &2 &bgaaf& A_{12}^2  &      \cr
28 &23& &a_{11}a_7a_4a_1^2 &2 &gcaa& A_{12}^2  &      \cr
28 &24& &d_8d_4^4 &8 &cb&  D_8^3 &D_8^3   \cr
28 &24& &a_9d_7a_5a_1a_1 &2 &efede& A_{11}D_7E_6  & A_{15}D_9  \cr
28 &24& &a_{11}d_5d_4a_3 &2 &fbeb& A_{11}D_7E_6  & A_{15}D_9  \cr
28 &24& &a_9a_7d_6a_1 &2 &fefc&  A_9^2D_6  &  A_{17}E_7 \cr
28 &24& &a_9a_7d_6a_1 &2 &fbfc&  A_9^2D_6  & D_{10}E_7^2 \cr
28 &24& &d_6^2d_6d_4a_1^2 &2 &bbbb& D_6^4  &D_{10}E_7^2  \cr
28 &24&E &a_8^3 &12 &a& 2A_8^3  &      \cr
28 &25& &d_8a_5^2a_5a_1 &2 &ddbe&      &      \cr
28 &25& &a_{11}a_7a_2^2a_1^2 &2 &gbaf&      &      \cr
28 &25& &a_{11}d_5a_4a_3 &1 &fcbb&      &      \cr
28 &25& &e_6d_6d_5a_5a_1 &1 &cceab&      &      \cr
28 &25& &d_7d_6a_5a_5 &1 &cdcb&      &      \cr
28 &25& &a_8a_7e_6a_1a_1 &1 &aeedc&      &      \cr
28 &25& &a_9e_6a_4^2a_1 &2 &cgbc&      &      \cr
28 &25& &a_{10}a_8a_4a_1 &1 &bbae&      &      \cr
28 &25& &a_9^2a_4a_1^2 &2 &gaf&      &      \cr
28 &25& &a_9a_8d_5a_1 &1 &fbcd&      &      \cr
\cr
29 &24& &a_{11}a_8a_3 &2 &bca&  A_{12}^2 &  A_{15}D_9 \cr
29 &25& &a_{10}d_6a_6 &1 &fbe&      &      \cr
\cr
30 &21& &d_7a_7e_6d_4 &2 &cada&  A_{11}D_7E_6 &      \cr
30 &22& &a_{10}^2a_3 &4 &ba& A_{12}^2  &      \cr
30 &23& &d_7^2a_7a_1^2 &4 &dcd& D_8^3  &      \cr
30 &23& &d_8a_7^2a_1^2 &4 &ddd&  D_8^3 &      \cr
30 &24& &d_8a_7^2 &4 &dd&  D_8^3 & A_{15}D_9  \cr
30 &24& &a_{11}d_6a_5a_1 &2 &fbba& A_{11}D_7E_6  &A_{17}E_7   \cr
30 &24& &a_{10}e_6a_6 &2 &agb& A_{11}D_7E_6  & A_{17}E_7  \cr
30 &24& &d_7a_7e_6a_3 &2 &caed& A_{11}D_7E_6  & D_{10}E_7^2 \cr
30 &24& &a_9e_6d_6a_1 &2 &efea& A_{11}D_7E_6  & D_{10}E_7^2 \cr
30 &24& &e_6^2d_5^2 &8 &ca&  E_6^4& D_{10}E_7^2 \cr
30 &25& &d_8a_7d_5a_3 &1 &cbdd&      &      \cr
30 &25& &a_{10}a_9a_4 &1 &bca&      &      \cr
30 &25& &d_7a_7d_5^2 &2 &cae&      &      \cr
\cr
31 &24& &a_{12}a_7a_4 &2 &gbf& A_{12}^2  & A_{17}E_7  \cr
31 &24&E &d_6^4 &24 &d& 2D_6^4 &      \cr
31 &24&E &a_9^2d_6 &4 &aa& 2A_9^2D_6  &      \cr
31 &25& &a_{12}a_6a_5 &1 &gba&      &      \cr
31 &25& &a_8e_6^2a_2^2 &2 &aaa&      &      \cr
31 &25& &a_{10}e_6d_5a_2 &1 &ebba&      &      \cr
31 &25& &a_{11}a_8d_4 &1 &bea&      &      \cr
31 &25& &a_{11}a_7d_5 &1 &bba&      &      \cr
31 &25& &a_8a_8d_7 &1 &dcb&      &      \cr
31 &25& &a_8^2d_7 &2 &ca&      &      \cr
\cr
32 &18& &a_9^2d_7 &4 &da&  A_9^2D_6  &      \cr
32 &18& &d_7d_6^3 &6 &db& D_6^4  &      \cr
32 &20& &a_{11}e_6d_5a_3 &2 &ebaa& A_{11}D_7E_6  &      \cr
32 &21& &a_{12}a_8d_4 &2 &bba&  A_{12}^2 &      \cr
32 &22& &d_8d_6^2a_3a_1^2 &2 &bbdb&D_8^3   &      \cr
32 &24& &d_8d_6^2a_1^2a_1^2 &2 &bbba&D_8^3   &D_{10}E_7^2  \cr
32 &24& &d_6^4 &8 &b&  D_6^4 & D_{12}^2 \cr
32 &25& &a_9d_8a_5a_1a_1 &1 &dfecb&      &      \cr
32 &25& &a_{12}a_7d_4 &1 &bbf&      &      \cr
32 &25& &a_9d_7d_6a_1 &1 &ceeb&      &      \cr
\cr
34 &19& &a_{11}d_7d_6 &2 &cea&  A_{11}D_7E_6 &      \cr
34 &19& &e_6^3d_6 &12 &ae& E_6^4 &      \cr
34 &23& &a_{11}d_8a_3a_1^2 &2 &dbca& A_{15}D_9  &      \cr
34 &23& &a_{13}a_8a_1^2a_1 &2 &caac& A_{15}D_9  &      \cr
34 &23& &d_9a_7^2a_1^2 &4 &eea&  A_{15}D_9 &      \cr
34 &24& &a_{13}d_6a_3a_1 &2 &bbcb&A_{15}D_9   &A_{17}E_7   \cr
34 &24& &d_9a_7^2 &4 &ed& A_{15}D_9  & D_{10}E_7^2 \cr
34 &24& &a_{11}d_7d_5 &2 &eee& A_{11}D_7E_6  & D_{12}^2 \cr
34 &25& &e_7a_7d_5a_5 &1 &cbca&      &      \cr
34 &25& &e_7a_7^2a_3a_1 &2 &dcab&      &      \cr
34 &25& &d_9a_7d_5a_3 &1 &ddeb&      &      \cr
34 &25& &a_{13}a_7a_3a_1 &1 &cfcc&      &      \cr
34 &25& &d_7^2e_6a_3 &2 &aca&      &      \cr
34 &25& &d_8a_7e_6a_3 &1 &edda&      &      \cr
34 &25& &a_{11}d_6^2 &2 &cb&      &      \cr
\cr
35 &24& &a_{11}^2 &4 &a& A_{12}^2  & D_{12}^2 \cr
35 &25& &a_{12}d_7a_4 &1 &ebc&      &      \cr
\cr
36 &20& &d_8^2d_5d_4 &2 &bda& D_8^3  &      \cr
36 &22& &a_{13}d_7a_3a_1 &2 &ebab& A_{15}D_9  &   \cr
36 &24& &a_9e_7a_7 &2 &cfa&  D_{10}E_7^2& A_{17}E_7  \cr
36 &24& &e_7d_6d_6d_4a_1 &2 &bbbaa&D_{10}E_7^2  & D_{10}E_7^2 \cr
36 &24& &d_8^2d_4^2 &2 &bb&  D_8^3 &  D_{12}^2\cr
36 &25& &a_{12}a_{10}a_1 &1 &aab&      &      \cr
\cr
37 &24&E &e_6^4 &48 &e& 2E_6^4&      \cr
37 &24&E &a_{11}d_7e_6 &2 &aaa& 2A_{11}D_7E_6 &      \cr
37 &25& &a_{13}e_6a_4a_1 &1 &baba&      &      \cr
\cr
38 &17& &a_{11}d_8e_6 &2 &dae&  A_{11}D_7E_6 &      \cr
38 &21& &a_{11}d_9d_4 &2 &dea&  A_{15}D_9 &      \cr
38 &23& &d_9e_6^2a_1^2 &4 &add& D_{10}E_7^2 &      \cr
38 &23& &a_9e_7e_6a_1^2 &2 &aecd& D_{10}E_7^2 &      \cr
38 &23& &a_{14}e_6a_2a_1^2 &2 &bfaa&  A_{17}E_7 &      \cr
38 &23& &a_{11}e_7a_5a_1^2 &2 &cbba& A_{17}E_7  &      \cr
38 &24& &a_{11}d_9a_3 &2 &eeb&  A_{15}D_9 &  D_{12}^2\cr
38 &24& &a_{12}a_{11} &2 &af&  A_{12}^2 & A_{24}  \cr
38 &25& &d_9d_7a_7 &1 &cdb&      &      \cr
38 &25& &a_{11}d_8d_5 &1 &dbc&      &      \cr
\cr
40 &20& &a_{15}d_5d_5 &2 &bba&  A_{15}D_9 &      \cr
40 &22& &d_8e_7d_6a_3a_1 &1 &bbada&D_{10}E_7^2  &      \cr
40 &22& &d_{10}d_6^2a_3 &2 &bbd& D_{10}E_7^2 &      \cr
40 &22& &a_{15}d_6a_3 &2 &bca& A_{17}E_7  &      \cr
40 &24& &d_{10}d_6^2a_1^2 &2 &bba&D_{10}E_7^2  &D_{12}^2  \cr
40 &24&E &a_{12}^2 &4 &a&  2A_{12}^2&      \cr
40 &25& &d_{10}a_9a_5 &1 &dbb&      &      \cr
40 &25& &a_{15}d_5a_4 &1 &bca&      &      \cr
40 &25& &e_7e_6^2a_5 &2 &aaa&      &      \cr
40 &25& &a_9e_7d_7 &1 &aca&      &      \cr
40 &25& &a_{11}e_7d_5a_1 &1 &aeba&      &      \cr
40 &25& &a_{14}a_9 &1 &ac&      &      \cr
\cr
41 &24& &a_{15}a_8 &2 &ac& A_{15}D_9  &A_{24}   \cr
\cr
42 &21& &a_{13}e_7d_4 &2 &aba& A_{17}E_7  &      \cr
\cr
43 &24& &a_{16}a_7 &2 &fa&  A_{17}E_7 & A_{24}  \cr
43 &24&E &d_8^3 &6 &d&  2D_8^3&      \cr
\cr
44 &16& &d_9d_8^2 &2 &db& D_8^3  &      \cr
44 &20& &e_7^2d_6d_5 &2 &bad& D_{10}E_7^2 &      \cr
44 &24& &d_8^2d_8 &2 &ba&  D_8^3 & D_{16}E_8 \cr
44 &24& &d_8^3 &6 &a&  D_8^3 & E_8^3 \cr
\cr
46 &23& &d_{11}a_{11}a_1^2 &2 &ebd& D_{12}^2 &      \cr
46 &24& &a_{15}d_8 &2 &cb& A_{15}D_9  &D_{16}E_8  \cr
46 &25& &d_{11}a_7e_6 &1 &dab&      &      \cr
\cr
47 &25& &a_{16}d_7 &1 &ca&      &      \cr
\cr
48 &18& &d_{10}e_7d_7a_1 &1 &abda&  D_{10}E_7^2&      \cr
48 &18& &a_{17}d_7a_1 &2 &cac&  A_{17}E_7 &      \cr
48 &22& &d_{10}^2a_3a_1^2 &2 &bdb& D_{12}^2 &      \cr
48 &24& &a_{15}e_7a_1 &2 &bcc&  A_{17}E_7 & D_{16}E_8 \cr
48 &24& &d_{10}e_7d_6a_1 &1 &abaa& D_{10}E_7^2 & D_{16}E_8 \cr
48 &24& &d_8e_7^2a_1^2 &2 &aaa&  D_{10}E_7^2& D_{16}E_8 \cr
48 &25& &a_{13}d_{10}a_1 &1 &bcb&      &      \cr
\cr
49 &24&E &a_{15}d_9 &2 &aa& 2A_{15}D_9 &      \cr
\cr
50 &15& &a_{15}d_{10} &2 &ba& A_{15}D_9  &      \cr
\cr
52 &20& &d_{12}d_8d_5 &1 &bad&  D_{12}^2&      \cr
52 &23& &a_{19}a_4a_1^2 &2 &caa& A_{24}  &      \cr
52 &24& &d_{12}d_8d_4 &1 &bab&  D_{12}^2& D_{16}E_8 \cr
\cr
55 &24&E &d_{10}e_7^2 &2 &dd& 2D_{10}E_7^2&      \cr
55 &24&E &a_{17}e_7 &2 &aa& 2A_{17}E_7 &      \cr
55 &25& &a_{19}d_5 &1 &aa&      &      \cr
\cr
56 &14& &d_{11}e_7^2 &2 &da&  D_{10}E_7^2&      \cr
56 &21& &a_{20}d_4 &2 &aa&  A_{24} &      \cr
\cr
58 &25& &e_8a_{11}e_6 &1 &aac&      &      \cr
58 &25& &a_{11}d_{13} &1 &bb&      &      \cr
\cr
60 &24& &e_8d_8^2 &2 &aa& E_8^3 &D_{16}E_8  \cr
\cr
62 &23& &a_{15}e_8a_1^2 &2 &cbd&D_{16}E_8  &      \cr
\cr
64 &22& &e_8e_7^2a_3 &2 &aae& E_8^3 &      \cr
64 &22& &d_{14}e_7a_3a_1 &1 &abda& D_{16}E_8 &      \cr
\cr
67 &24&E &d_{12}^2 &2 &d& 2D_{12}^2&      \cr
\cr
68 &12& &d_{13}d_{12} &1 &db& D_{12}^2 &      \cr
68 &20& &d_{12}e_8d_5 &1 &aad& D_{16}E_8 &      \cr
68 &24& &d_{12}^2 &2 &b& D_{12}^2 & D_{24} \cr
\cr
71 &24& &a_{23} &2 &a& A_{24}  &  D_{24}\cr
\cr
76 &16&E&d_{16}d_9 &1 &bd& D_{16}E_8 &      \cr
76 &16&E&d_9e_8^2 &2 &ea& E_8^3 &      \cr
76 &24& &d_{16}d_8 &1 &ba&D_{16}E_8  &D_{24}  \cr
76 &24&E &a_{24} &2 &a&2A_{24}  &      \cr
\cr
91 &24&E &e_8^3 &6 &e& 2E_8^3&      \cr
91 &24&E &d_{16}e_8 &1 &dd&2D_{16}E_8 &      \cr
\cr
92 & 8&E&d_{17}e_8 &1 &da&D_{16}E_8  &      \cr
\cr
100 &20& &d_{20}d_5 &1 &ad&D_{24}  &      \cr
\cr
139 &24&E &d_{24} &1 &d& 2D_{24}&      \cr
\cr
140 & 0&E&d_{25} &1 &d& D_{24} &      \cr
}

\proclaim Table 3 Orbits of $B/2B$, $B$ a Niemeier lattice. 

See section 3.6 and the last part of section 3.7. 

For each even unimodular lattice $B$ of dimension at most 24 we list
the orbits of $B/2B$ under $\Aut(B)$. If $\hat b$ is in $B/2B$ then
the orbit of $\hat b$ is denoted by the Dynkin diagram of $B_{\hat
b}$, where $B^{\hat b}$ is the set of vectors of $B$ that have even
inner product with $\hat b$. The orbits are written on rows according
to the minimal norm of a representative $b$ of $\hat b$ in $B$. The
integer by each orbit gives the ``non-reflection part'' of the
subgroup of $\Aut(B)$ fixing that orbit. Each orbit of norm 8 or norm
12 corresponds to an odd lattice $A$ with no vectors of norm 1 one of
whose even neighbors is $B$; we list the other even neighbor of $A$
next to the symbol for this orbit.

(This table has been omitted as is is long and boring, and all the
information in it can be reconstructed from tables $-2$ and $-4$. In
the unlikely event that you need the information in it, either ask me
for a copy or reconstruct it yourself from the list of type 2 vectors
of norms $-2$ and $-4$ given in tables $-2$ and $-4$.)
\bye